\documentclass[10pt]{article}

\usepackage{amssymb,amsmath,amsthm,amsfonts}
\usepackage[pdftex,pagebackref,colorlinks,linktocpage,citecolor=black,linkcolor=black]{hyperref}
\usepackage{graphicx,graphics}
\usepackage{epstopdf}
\usepackage{comment}
\usepackage{algorithm,algorithmic}
\usepackage{verbatim}
\usepackage{indentfirst}
\usepackage{mathtools}
\usepackage{cite}
\usepackage{url}
\usepackage[usenames,dvipsnames,svgnames,table]{xcolor}

\allowdisplaybreaks
\numberwithin{equation}{section}
\numberwithin{algorithm}{section}
\newtheorem{theorem}{Theorem}[section]
\newtheorem{lemma}{Lemma}[section]

\theoremstyle{definition}

\newtheorem{assumption}{Assumption}[section]
\theoremstyle{remark}
\newtheorem{remark}{Remark}[section]

\def\cT {\mathcal T}
\def\bold {\boldsymbol}
\def\Om {\Omega}
\def\J {\mathcal{J}}

\def\to{\rightarrow}
\def\eps {\varepsilon}
\def\A {{\mathcal A}}

\def\p {\partial}

\newcommand{\dx}{\,{\rm d}x}
\newcommand{\ds}{\,{\rm d}s}

\allowdisplaybreaks

\begin{document}
\title{Adaptive Approximations of Inclusions in a Semilinear Elliptic Problem Related to Cardiac Electrophysiology\thanks{The work of B. Jin is supported by Hong Kong RGC General Research Fund (Project 14306824) and a start-up fund from The Chinese University of Hong Kong. The work of Y. Xu is partially supported by the National Natural Science Foundation of China (Projects 12250013, 12261160361 and 12271367), the Science and Technology Commission of Shanghai Municipality (Projects 20JC1413800 and 22ZR1445400) and General Research Fund (Projects KF202318 and KF202468) from Shanghai Normal University.}}
\author{
Bangti Jin\thanks{Department of Mathematics, The Chinese University of Hong Kong, Shatin, N.T., Hong Kong, P.R. China. (b.jin@cuhk.edu.hk, bangti.jin@gmail.com; fengruwang@cuhk.edu.hk) }\and
Fengru Wang\footnotemark[2] \and Yifeng Xu\thanks{Department of Mathematics and Scientific Computing Key Laboratory of Shanghai Universities, Shanghai Normal University, Shanghai 200234, P.R. China. (yfxu@shnu.edu.cn, mayfxu@gmail.com)}
}

\date{}

\maketitle

\begin{abstract}
{In this work, we investigate the numerical reconstruction of inclusions in a semilinear elliptic
equation arising in the mathematical modeling of cardiac ischemia. We propose an adaptive finite
element method for the resulting constrained minimization problem that is relaxed by a phase-field
approach. The \textit{a posteriori} error estimators of the adaptive algorithm consist of three
components, i.e., the state variable, the adjoint variable and the complementary relation. Moreover,
using tools from adaptive finite element analysis and nonlinear optimization, we establish the strong
convergence for a subsequence of adaptively generated discrete solutions to a solution of the
continuous optimality system. Several numerical examples are presented to illustrate the convergence
and efficiency of the adaptive algorithm.}
{semilinear elliptic equation, inclusion recovery, \textit{a posteriori} error estimator, adaptive finite element, convergence}
\end{abstract}

\section{Introduction}

Ischemic / coronary heart disease that results from a restriction in blood supply to the heart represents the most common form of heart disease around the world (see the website \url{https://www.cdc.gov/heartdisease/coronary_ad.htm} of Centers for Disease Control and Prevention for the US and the reports \cite{Daponte:2022} and \cite{Wang:2023} for Europe and China, respectively), and also it is the leading cause of heart attack. Thus the early detection of the disease is of primary importance. It is often inferred from the electrical activity of the heart by means of body surface or intracardiac measurements, and the bidomain / monodomain model represents the predominant mathematical model \cite{ColliPavarino:2014}. In the past two decades, several numerical investigations dealing with ischemia identification from measurements of surface potentials have been performed by casting the problem into an optimization framework: a stationary model involving the heart-torso coupling \cite{NielsenLysaker:2007}, a nonstationary monodomain model from an isolated heart \cite{LysakerNielsen:2006}, and ischemia identification from intracardiac electrograms \cite{AlvarezAlonso:2012}.

Following \cite{beretta2016}, we consider the insulated monodomain model in the steady state by ignoring
the coupling with the ionic model but instead adopting a phenomenological description of the ion current,
which leads to the study of a Neumann boundary value problem for a semilinear elliptic equation.
Over an open bounded polygonal domain $\Omega\subset\mathbb{R}^2 $, the governing equation reads
\begin{equation}\label{diffequ}
\left\{\begin{aligned}
        -{\nabla}\cdot(\widetilde{\sigma}{\nabla}y)+ \chi_{\Omega\setminus\overline{\omega}} y^3 &= f,& & \mbox{in}~\Omega, \\
        \dfrac{\partial y}{\partial {n}} &= 0,& & \mbox{on}~\partial\Omega,
    \end{aligned}
\right.
\end{equation}
where $\omega\subsetneq\Omega$ is an open subdomain (representing the region occupied by ischemia in cardiac electrophysiology), $\chi_{\Omega\setminus\overline{\omega}}$ is the characteristic function of the set $\Omega\setminus\overline{\omega}$, $f\in L^2(\Omega)$ is a fixed source, ${n}$ is the unit outward normal vector to the boundary $\partial\Omega$ and the conductivity $\widetilde{\sigma}({x})$ is piecewise constant:
\[
    \widetilde{\sigma}({x}) = \left\{
        \begin{array}{ll}
        \sigma, & {x}\in \omega,\\
        1, & {x}\in \Omega\setminus\overline{\omega},
        \end{array}
    \right.
\]
with $0<\sigma \ll 1$. In the context of cardiac electrophysiology \cite{ColliPavarino:2014}, the state
variable $y$ represents the distribution of the electric transmembrane potential, $\widetilde{\sigma}$
is the tissue conductivity. According to the
experimental observation, in an ischemic or infarcted region, cells are not excitable and thus the conductivity
$\sigma$ is substantially different from that of healthy tissues \cite{Strinstra:2005}. Note that ion transport circumvents
ischemic areas, so that also the ionic membrane current is multiplied by $\chi_{\Omega\setminus \omega}$.
The cubic polynomial term $y^3$, denoting the ionic current across the cell membrane, represents one significant
empirical choice, and describes the macroscopic behavior of excitable cells. Finally, $f$ represents a
current stimulus applied to the tissue, usually in a confined region (and for a short time interval),
expressing the initial electrical stimulus, related to the so-called pacemaker potential. In practice, the ischemia region
$\omega$ is also assumed to be well separated from $\p\Omega$.
\begin{assumption}\label{ass01}
    Given $d_0>0$, $\chi_\omega = 0~\text{a.e. in}~\Omega^{d_0}$ with  $\Omega^{d_0} = \{{x}\in \Omega:~\text{dist}({x},\partial \Omega)\leq d_0\}$.
\end{assumption}
Throughout, we make the following assumption on the source $f$. This condition was also employed in \cite{beretta2018a}. It guarantees the Fr\'{e}chet differentiability of the functional $\J_\eps$ in the minimization problem of interest (see \eqref{min_G-L}), then the associated optimality system (see \eqref{optsys}), and the convergence of finite element approximations with uniform refinements.
\begin{assumption}\label{ass02}
    The source $f$ is greater than or equal to a constant $m>0$, i.e.,
    $f \geq m$ a.e. in $\Omega$.
\end{assumption}

Let $u=\chi_\omega$, $a(u)=1- (1-\sigma)u$ and $b(u)= 1- u$. Then the variational
formulation for problem \eqref{diffequ} is to find $y:=y(u)\in H^1(\Omega)$ such that
\begin{equation}\label{vp}
    (a(u) {\nabla} y , {\nabla} \phi )
    + (b(u) y^3, \phi ) = (f,\phi), \quad \forall \phi \in H^1(\Omega),
\end{equation}
where the notation $(\cdot,\cdot)$ denotes the $L^2(\Omega)$ inner product. We use the notation $y(u)$
to indicate its dependence on $u$. The unique
solvability of problem \eqref{vp} has been established in \cite{beretta2018a} using
Minty-Browder theorem \cite[Theorem 26.A]{zeidler1990}, and the monotonicity of the term $(b(u) y^3, \phi )$.
In this work, we are concerned with the inverse problem of recovering the inclusion
$\omega$ from the over-posed noisy boundary observation $y^\delta$ on $\partial\Omega$ with a noise level
$\delta = \|y^\delta-y\|_{L^2(\p\Omega)}$,
where $\delta>0$ measures the accuracy of $y^\delta$. Like most practical inverse
problems, it is  ill-posed in the sense that the solution may be very sensitive with respect to the perturbation of the data $y^\delta$ \cite{beretta2018a}. One first investigation of the inverse problem was provided in \cite{beretta2016}, where an anasymptotic expansion of the perturbed
electrical potential was derived and a reconstruction procedure was given. Beretta et al \cite{beretta2017c} proposed a reconstruction algorithm
based on topological derivatives to detect the
position of small inclusions. To remedy the inherent
ill-posedness, we employ the classic Tikhonov regularization \cite{EnglHankeNeubauer2000,ItoJin:2015}, which leads to
the following PDE constrained minimization problem:
\begin{equation}\label{min_bv}
    \min_{u\in X} \big\{\J(u) = \tfrac{1}{2} \|y(u) - y^\delta\|^2_{L^2(\p\Omega)} + \alpha |u|_{\mathrm{TV}(\Omega)}\big\},
\end{equation}
where $X = \{ v\in \mathrm{BV}(\Omega):v\in\{0,1\}~\text{a.e. in}~\Omega,~v = 0~\text{a.e. in}~\Omega^{d_0}\}$
and $|\cdot|_{\rm TV}$ denotes the total variation seminorm, which promotes piecewise constancy of the reconstructed
conductivity \cite{ChambolleCaselles:2010}, and then taking the minimizer (if it exists) as an approximation. The scalar $\alpha>0$ is the regularization parameter.
Under Assumption \ref{ass02}, the existence of a minimizer to \eqref{min_bv}, the stability of a minimizer with
respect to the data $y^\delta$ and the convergence of the sequence of minimizers as $\alpha\to 0^+$ have been established
in \cite{beretta2018a}, using a compactness argument.

Numerically, the non-differentiablity of the functional $\J$ and the discrete nature of the admissible set $X$ pose
significant challenges. One powerful idea is to
replace the total variation penalty $|\cdot|_{\rm TV}$ in $\J$ by a smooth approximation. \cite{beretta2018a} proposed regularizing $\frac{1}{2} \|y(u) -y^\delta\|^2_{L^2(\p\Omega)}$ with
a Ginzburg-Landau energy. That is, with the relaxed admissible set $\mathcal{A}=\{v\in
H^1(\Omega):~0\leq v\leq 1~\text{a.e. in}~\Omega,~v=0~\text{a.e. in}~\Omega^{d_0}\}$ and $\eps>0$, one employs the following constrained minimization problem:
\begin{equation}\label{min_G-L}
    \min_{u\in \A} \left\{\J_\eps (u) = \frac{1}{2} \|y(u) - y^\delta\|^2_{L^2(\p\Omega)} + \alpha \left( \eps \int_\Omega |{\nabla} u|^2 \dx + \frac{1}{\eps} \int_\Omega u(1-u)  \dx \right) \right\}.
\end{equation}
An analysis of problem \eqref{min_G-L} is provided in \cite{beretta2018a}: (i) For a fixed $\varepsilon$, \eqref{min_G-L} has a minimizer which is stable with respect to the measurement $y^{\delta}$. (ii) \eqref{min_bv} and \eqref{min_G-L} are connected via $\Gamma$-convergence,
i.e., $\J_{\eps_k}$ $\Gamma$-converges to $\J$ in the $L^1(\Omega)$-norm as $\varepsilon_k\to0$
and the sequence of minimizers $\{u_k^\ast\}_{k\geq0}$ to $\{\J_{\eps_k}\}_{k\geq0}$ has a
subsequence converging to a minimizer $u^\ast$ to $\J$ in the $L^1(\Omega)$-norm. Berettta et al \cite{beretta2018a} analyzed the convergence of the finite element method (for uniform refinements), and designed an
algorithm based on a parabolic obstacle problem to find the discrete minimizer.

This work builds on the work \cite{beretta2018a}. Due to the presence of the jumps of $\tilde{\sigma}$ and the degeneracy
of the nonlinear term $y^3$ across the interface $\partial\omega$, the solution $y(u)$ of problem \eqref{diffequ} naturally possesses local singularities. This results in poor
numerical performance of FEM over uniform mesh refinements and then affects the reconstruction
accuracy. So there is still room for improvement in the numerical efficiency for resolving
\eqref{min_G-L}, which we aim to explore in this work. One natural option is to incorporate adaptive techniques into the process of computing.
Adaptive finite element methods (AFEMs) typically consist of the following successive loops:
\begin{equation}\label{afem_loop}
  \mbox{SOLVE}\rightarrow\mbox{ESTIMATE}\rightarrow\mbox{MARK}\rightarrow\mbox{REFINE}.
\end{equation}
The module SOLVE outputs the numerical solution of \eqref{min_G-L} by a finite element method over a given mesh. The \textit{a posteriori} error estimate in the module ESTIMATE identifies the region with local singularities. Elements over there are selected in the module MARK for refinement, which is performed in the module REFINE. Over the past four decades, AFEM has been established as a very popular technique in scientific and engineering computing since it is able to achieve the desired accuracy at a reduced computational cost. For forward problems, the theory of AFEM is well developed by now \cite{ainsworth2000,carstensen2014,nochetto2009,verfurth2013}. In recent years, several researchers have contributed to the applications and analysis of AFEM in the context of several linear and nonlinear inverse problems, e.g., inverse scattering \cite{BeilinaJohnson:2005,BeilinaClason:2006}, flux reconstruction \cite{xu2015a}, inverse Robin problem \cite{xu2015} and electrical impedance tomography  \cite{JinXu:2020,JinXuZou:2016}. The major difference between the AFEM for inverse problems and that for direct problems lies in the fact that, due to the ill-posedness, an estimate of
the difference between computed and exact coefficients is replaced by \textit{a posteriori} estimates of the accuracy of either the Lagrangian \cite{BeilinaJohnson:2005,BeilinaClason:2006} or the Tikhonov functional \cite{BeilinaKlibanov:2010,KaltenbacherKirchner:2014,ClasonKaltenbacherWachsmuth:2016} or the optimality system \cite{xu2015,xu2015a,JinXuZou:2016,JinXu:2020}. This departure brings significant challenges in the convergence analysis. It is also worth mentioning that to treat the nonsmoothness of the total variation penalty, Bartels and his collaborators \cite{Bartels:2015,BartelsToveyWassmer:2022} proposed a duality approach for AFEM to minimize a total variation regularized formulation arising from the Rudin–Osher–Fatemi model (for image denoising), for which the dual problem is smooth and strongly convex subject to the constraint. However, the duality technique does not lead to a convex dual problem for the functional $\J$ in \eqref{min_bv}, due to the nonlinear forward map given by \eqref{vp}, and thus the approach in \cite{Bartels:2015,BartelsToveyWassmer:2022} does not apply directly to the present setting.

In this work, we employ AFEM to recover the inclusion $u$ through problem \eqref{min_G-L} subject to the PDE constraint \eqref{vp}. An adaptive algorithm of the form \eqref{afem_loop}, based on the first-order necessary optimality system, is proposed in Algorithm \ref{afem_CE}. The algorithm involves \textit{a posteriori} error estimators for the state $y$, the adjoint $p$ and the control $u$ (via the variational inequality). Due to the nonconvexity of problem \eqref{min_G-L}, to prove the reliability of the error estimators without any further assumptions along the traditional line \cite{ainsworth2000, verfurth2013} is still open. Nevertheless, we adopt related arguments for forward problems \cite{GanterPraetorius:2022,morin2008,siebert2011} to establish its convergence: the sequence of discrete solutions $(u^\ast_k,y^\ast_k,p^\ast_k)$ by Algorithm \ref{afem_CE} contains a subsequence strongly convergent in the $H^1(\Omega)$-topology to a triplet $(u^\ast,y^\ast,p^\ast)$ solving the optimality system associated with problem \eqref{min_G-L} subject to \eqref{vp}; see Theorem \ref{thm:conv} for the precise statement. Moreover, some computable quantities, serving as the error estimators in Algorithm \ref{afem_CE}, arise naturally from the analysis. To the best of our knowledge, this work represents the first work on the adaptive FEM for solving an inverse problem associated with a semilinear elliptic problem.

The proof of Theorem \ref{thm:conv} is lengthy and technical, and involves three crucial ingredients. The overall proof strategy proceeds as follows. First, by utilizing techniques in the nonlinear PDE-constrained optimization, we establish an auxiliary convergence: the sequence of discrete triplets $\{(u^\ast_k,y^\ast_k,p^\ast_k)\}_{k\geq0}$ generated by Algorithm \ref{afem_CE} contains a subsequence $\{(u^\ast_{k_j},y^\ast_{k_j},p^\ast_{k_j})\}_{j\geq0}$ converging to a $H^1(\Omega)$-limiting triplet $(u_\infty^\ast,y_\infty^\ast,p_\infty^\ast)$ (see Theorems \ref{thm:conv_medmin} and \ref{thm:conv_med_costate}). The convergence of this kind is first observed for two-point boundary value problems in the seminal work \cite{babuska:1984}. Second, inspired by the approach in the work \cite{GanterPraetorius:2022} (see also
\cite{morin2008, siebert2011}), we show that three subsequences of associated estimators $\{\eta_{k_j,1}\}_{j\geq0}$, $\{\eta_{k_j,2}\}_{j\geq0}$ and $\{\eta_{k_j,3}\}_{j\geq0}$ in the module ESTIMATE all tend to zero in Theorem \ref{thm:conv_est}. This is achieved by establishing a vanishing limit of the maximal error indicators and two structural properties for the estimators $\{\eta_{k_j,i}\}_{j\geq0}$ for $i=1,2,3$: stability on non-refined elements and reduction on refined elements. These important auxiliary results are shown in Lemmas \ref{lem:conv_aux5_maxest}, \ref{lem:stab_est_unref} and \ref{lem:reduction_est_ref}. Third and last, we prove the desired subsequence convergence by verifying that $ (u_\infty^\ast,y_\infty^\ast,p_\infty^\ast)$ solves the first-order necessary optimality system \eqref{optsys} (see Lemma \ref{lem:conv_aux4_vp}). Note that the presence of a semilinear term requires a fairly delicate treatment in the convergence analysis when compared with that for linear elliptic problems. In addition, we have adapted recent advances on the AFEM convergence analysis \cite{GanterPraetorius:2022} for direct problems to the context of inverse problems. Thus, the analysis strategy differ markedly from existing works \cite{JinXu:2020,JinXuZou:2016,xu2015,xu2015a} on the convergence of AFEM for linear and nonlinear PDE inverse problems. These also represent the main technical novelty of the work.

The rest of this paper is organized as follows. In Section \ref{sec:alg}, we propose an adaptive algorithm for problem \eqref{min_G-L} subject to \eqref{vp}. In Section \ref{sec:conv}, we give the technical proof of the convergence of the adaptive algorithm. We present in Section \ref{sec:numer} some numerical experiments to illustrate the algorithm. Throughout, the notation $c$ denotes a generic constant, which
may differ at each occurrence, but it is always independent of the mesh size and other quantities of interest unless specified otherwise.

\section{Adaptive algorithm}\label{sec:alg}

In this section, we develop an adaptive algorithm for recovering the inclusion $u$. Let $\cT_0$ be a shape regular conforming triangulation of the domain
$\overline{\Omega}$ into closed triangles and $\mathbb{T}$ be the set of all possible conforming triangulations of $\overline{\Omega}$
obtained from $\cT_0$ by the successive use of the newest vertex bisection algorithm \cite{nochetto2009}. Then the set $\mathbb{T}$ is uniformly shape regular, i.e., the
shape-regularity of any $\mathcal{T}\in\mathbb{T}$ is bounded by a constant depending only on $\cT_0$ \cite{nochetto2009,traxler:1997}. Over any triangulation $\cT\in\mathbb{T}$, we define a continuous piecewise linear finite element space
\begin{equation*}
   V_\cT = \left\{v\in C(\overline{\Omega}): v|_T\in P_1(T)~ \forall T\in\cT\right\},
\end{equation*}
where $P_1(T)$ consists of all linear functions on the element $T$. The space $V_\cT$ is used for approximating the potential $y$ and the discrete admissible set $\mathcal{A}_\mathcal{T}$ for the inclusion $u$ is taken to be
$\mathcal{A}_\cT:= V_\cT \cap \mathcal{A}.$

Given a discrete approximation $u_\cT \in \mathcal{A}$, the finite element approximation of problem \eqref{vp} is to find $y_\cT=y_\cT(u_\cT) \in V_\cT$ such that
\begin{equation}\label{disvp}
    (a(u_\cT) {\nabla} y_\cT, {\nabla} \phi_\cT) + (b(u_\cT) y_\cT^3, \phi_\cT) = (f,\phi_\cT), \quad \forall\phi_\cT \in V_\cT.
\end{equation}
Now the numerical approximation of problem \eqref{min_G-L} is given by
\begin{equation}\label{min_G-L_dis}
    \min_{u_\cT\in \A_\cT} \left\{\J_{\eps,\cT} (u_\cT) = \frac{1}{2} \|y_\cT(u_\cT) - y^\delta\|^2_{L^2(\p\Omega)} + \alpha \left( \eps \int_\Omega |{\nabla} u_\cT|^2 \dx + \frac{1}{\eps} \int_\Omega\!\! u_\cT(1\!-\!u_\cT)  \dx \right) \right\}.
\end{equation}
The minimizer $u^\ast_\cT \in \A_\cT$ to problem \eqref{min_G-L_dis}, the related discrete state $y_\cT^\ast\in V_\cT$ and the discrete adjoint $p^\ast_\cT \in V_\cT$ satisfy the following discrete optimality system \cite{beretta2018a}
\begin{subequations}\label{optsys_dis}
    \begin{align}
            &(a(u^\ast_\cT){\nabla} y^\ast_\cT , {\nabla} \phi_\cT)
    + (b(u^\ast_\cT) (y^\ast_\cT)^3,\phi_\cT ) = (f, \phi_\cT), \quad \forall\phi_\cT \in V_\cT,\label{optsys_dis-1}\\
        &(a(u_\cT^\ast){\nabla} p_\cT^\ast,{\nabla} \psi_\cT ) + 3 (b(u_\cT^\ast) (y^\ast_\cT)^2 p^\ast_\cT, \psi_\cT) =  (y_\cT^\ast - y^\delta, \psi_\cT)_{L^2({\p\Omega})}, \quad \forall \psi_\cT \in V_\cT , \label{optsys_dis-2}\\
        &( (1-\sigma) (v_\cT-u^\ast_\cT){\nabla} y^\ast_\cT ,{\nabla} p^\ast_\cT)+ (v_\cT-u_\cT^\ast, (y^\ast_\cT)^3 p^\ast_\cT )\nonumber\\
         &\quad + 2 \alpha \varepsilon ({\nabla} u^\ast_\cT ,{\nabla} (v_\cT - u_\cT^\ast)) +\alpha\varepsilon^{-1} (1-2u^\ast_\cT,v_\cT - u_\cT^\ast) \geq 0, \quad \forall v_\cT \in \A_\cT.\label{optsys_dis-3}
        \end{align}
\end{subequations}
This system forms the basis for deriving the \textit{a posteriori} error estimators.

To give an adaptive algorithm for problem \eqref{min_G-L} subject to \eqref{vp}, we introduce some notation. The collection of all edges (respectively all interior edges) in the triangulation $\cT\in\mathbb{T}$ is denoted by $\mathcal{F}_{\cT}$ (respectively
$\mathcal{F}_{\cT}(\Omega)$) and its restriction to the boundary $\p\Omega$ by $\mathcal{F}_{
\cT}(\p\Omega)$, respectively. An edge $F$ has a fixed normal unit vector ${n}_{F}$ in $\overline{\Omega}$ with ${n}_{F}={n}$ on $\partial\Omega$. Over any triangulation $\cT\in\mathbb{T}$, we define a piecewise constant meshsize function $h_\cT:\overline{\Omega}\rightarrow \mathbb{R}_{+}$ by
\begin{equation}\label{meshsize_def}
    h_\cT |_T := h_T = |T|^{1/2},\quad \forall T\in\cT.
\end{equation}
Since $\cT$ is shape regular, $h_T$ is equivalent to the diameter of any $T\in\cT$.

For the solution triplet $(u_\cT^\ast, y_\cT^\ast, p_\cT^\ast)$ of the system \eqref{optsys_dis}, we define three element residuals on each element $T\in\cT$ and three jump residuals on each edge $F\in\mathcal{F}_\cT$ respectively by
\begin{align*}
    R_{T,1}(u_\cT^\ast, y_\cT^\ast) &: = {\nabla} \cdot (a(u_\cT^\ast) {\nabla} y_\cT^\ast) - b(u_{\cT}^\ast)(y_\cT^\ast)^3 + f,\\
    R_{T,2}(u_\cT^\ast, y_\cT^\ast, p_\cT^\ast) &: = {\nabla} \cdot (a(u_\cT^\ast){\nabla}p_\cT^\ast) - 3b(u_{\cT}^\ast)(y_\cT^\ast)^2 p_\cT^\ast,\\
    R_{T,3}(u_\cT^\ast, y_\cT^\ast, p_\cT^\ast) &: = (1-\sigma){\nabla} y_\cT^\ast \cdot{\nabla} p_\cT^\ast + (y_\cT^\ast)^3 p_\cT^\ast + {\alpha}\varepsilon^{-1} (1-2u_\cT^\ast),\\
    J_{F,1}(u_{\cT}^{\ast},y^{\ast}_{\cT}) &:= [a(u_\cT^\ast) \nabla y_\cT^\ast \cdot{n}_F],\\
    J_{F,2}(u_{\cT}^{\ast},y^{\ast}_{\cT},p^{\ast}_{\cT}) &:=
    \left\{\begin{array}{ll}
                        [a(u_{\cT}^{\ast})\nabla p_{\cT}^{\ast}\cdot{n}_{F}],\quad&
                        \mbox{for} ~~F\in\mathcal{F}_{\cT}(\Omega),\\ [1ex]
                        a(u_{\cT}^{\ast})\nabla p_{\cT}^{\ast}\cdot{n}-(y_{\cT}^{\ast}-y^{\delta}),\quad&
                        \mbox{for} ~~ F\in\mathcal{F}_{\cT}(\p\Omega),\end{array}\right.\\
    J_{F,3}(u^{\ast}_{\cT}) &:=2\alpha\varepsilon[\nabla u_{\cT}^{\ast}\cdot{n}_{F}],
\end{align*}
where $[\cdot]$ denotes the jump across an interior edge $F\in\mathcal{F}(\Omega)$ while it is equal to the value of the involved function on a boundary edge $F\in\mathcal{F}_\cT(\p\Omega)$. Then for any
$T\in \cT$, we define three error indicators
\begin{align*}
    \eta_{\cT,1}^{2}(u_{\cT}^{\ast},y_{\cT}^{\ast},T)
    & :=h_{T}^{2}\|R_{T,1}(u^{\ast}_{\cT},y^{\ast}_{\cT})\|_{L^{2}(T)}^{2}
    +\sum_{F\subset\partial T}h_{T}\|J_{F,1}(u^{\ast}_{\cT},y^{\ast}_{\cT})\|_{L^{2}(F)}^{2},\\
    \eta_{\cT,2}^{2}(u_{\cT}^{\ast},y_{\cT}^{\ast},p_{\cT}^{\ast},T)
    &:=h_{T}^{2}\|R_{T,2}(u^{\ast}_{\cT},y_{\cT}^{\ast},p^{\ast}_{\cT})\|_{L^{2}(T)}^{2}
    +\sum_{F\subset\partial T}h_{T}\|J_{F,2}(u^{\ast}_{\cT},y_{\cT}^{\ast},p^{\ast}_{\cT})\|_{L^{2}(F)}^{2},\\
    \eta_{\cT,3}^{2}(u_{\cT}^{\ast},y_{\cT}^{\ast},p_{\cT}^{\ast},T)
    &:=h_{T}^{2}\|R_{T,3}(u_{\cT}^\ast,y_{\cT}^{\ast},p_{\cT}^{\ast})\|^{2}_{L^{2}(T)}+\sum_{F\subset\partial T}h_{T}\|J_{F,3}(u_{\cT}^{\ast})\|^{2}_{L^{2}(F)}.
\end{align*}
Then for any collection of elements $\mathcal{M}\subseteq \cT$, we define
\begin{align*}
    \eta_{\cT,1}^{2}(\mathcal{M})&: = \sum_{T\in\mathcal{M}}\eta_{\cT,1}^{2}(u_{\cT}^{\ast},y_{\cT}^{\ast},T),\quad
    \eta_{\cT,2}^{2}(\mathcal{M}): = \sum_{T\in\mathcal{M}}\eta_{\cT,2}^{2}(u_{\cT}^{\ast},y_{\cT}^{\ast},p_{\cT}^{\ast},T),\\
    \eta_{\cT,3}^{2}(\mathcal{M})&: = \sum_{T\in\mathcal{M}}\eta_{\cT,3}^{2}(u_{\cT}^{\ast},y_{\cT}^{\ast},p_{\cT}^{\ast},T).
\end{align*}
When $\mathcal{M}=\cT$, we drop $\cT$ from these estimators.

\begin{remark}
The obtained \textit{a posteriori} estimators contain the small parameter $\eps$ as $\eps^{-1}$, in analogous to the sharp interface limit of Allen-Chan and Cahn-Hilliard type models. Thus, it is of interest to analyze the robustness of the estimators with respect to $\eps$. In the context of \textit{a posteriori} error estimators for the time-dependent Allen-Cahn model, this issue has been investigated extensively; see, e.g.,  \cite{Kessler:2004,BartelsMuller:2011,BartelsMuller:2011b}, and various strategies to derive reliability estimates that are explicit and polynomially dependent on $\eps^{-1}$ have been proposed and analyzed. Note that there is still no reliability result for the AFEM algorithm proposed in the present work.
\end{remark}

\begin{algorithm}
\caption{AFEM for problem \eqref{min_G-L} subject to the constraint \eqref{vp}}\label{afem_CE}
    \begin{algorithmic}[1]
  \STATE {(INITIALIZE)} Specify an initial mesh $\cT_{0}$, choose the tolerance \texttt{TOL} and $K$
  \FOR {$k=0:K-1$}
  \STATE {(SOLVE)}
        Solve problem \eqref{disvp}--\eqref{min_G-L_dis} over $\cT_k$ for $(u_k^\ast, y_k^\ast) \in \mathcal{A}_k \times V_k$ and \eqref{optsys_dis} for $p_k^\ast \in V_k$.
  \STATE {(ESTIMATE)}
        Compute the error estimators $\eta_{k,1}^2$, $\eta_{k,2}^2$ and $\eta_{k,3}^2$.

  \STATE Check the stopping criterion: if $\eta_{k,1}^2+\eta_{k,2}^2+\eta_{k,3}^2\leq \texttt{TOL}$, then break.

  \STATE {(MARK)} Let $\widetilde{\eta}_k = \max_{1\leq i\leq 3}(\eta_{k,i})$. Mark a subset $\mathcal{M}_k\subseteq\cT_k$ such that $\mathcal{M}_k$ contains at least one element $\widetilde{T}_{k}\in \cT_k$ with the largest error indicator:
  \begin{equation}\label{marking}
    \widetilde{\eta}_{k}(\widetilde{T}_k) = \max_{T\in\cT_k} \widetilde{\eta}_{k}(T).
  \end{equation}

  \STATE {(REFINE)} Refine each element $T$ in $\mathcal{M}_{k}$ by the newest vertex bisection algorithm to get $\cT_{k+1}$.

  \ENDFOR

  \STATE Output $(u_k^*,y_k^*,p_k^*)$.

\end{algorithmic}
\end{algorithm}

Now we can present an adaptive algorithm of the form \eqref{afem_loop} for problem \eqref{min_G-L}
subject to \eqref{vp}; see Algorithm \ref{afem_CE} for details. The dependence on the triangulation $\mathcal{T}_k$ is indicated by the iteration number $k$ in the subscript, e.g., $V_{k}$ for $V_{\cT_{k}}$. In the algorithm, we prescribe a maximum iteration number at the beginning and include a stopping test after the ESTIMATE module. In the
module MARK, we adopt a separate marking by comparing the three error
estimators and utilizing the dominant one in the marking strategy. Most widely
used strategies, including the maximum and practical D\"{o}rfler's strategies \cite{siebert2011},
fulfill condition \eqref{marking}. A typical example in the module MARK reads
$$
    \eta_{k,1}^2(\mathcal{M}_k)\geq \theta \eta_{k,1}^2 \quad \mbox{and}\quad \max_{T\in\cT_k\setminus\mathcal{M}_k}\eta_{k,1}(u_k^\ast,y_k^\ast,T)\leq\min_{T\in\mathcal{M}_k}\eta_{k,1}(u_k^\ast,y_k^\ast,T),
$$
with the parameter $\theta\in(0,1]$ if $\eta_{k,1}^2 = \max_{1\leq i\leq 3}\eta_{k,i}^2$ at the $k$-th iteration.

It is worth noting that the error estimators in Algorithm \ref{afem_CE} are not  provably reliable, i.e., they do not provide upper bounds on the error. This is mainly attributed to the nonconvexity of the functional $\J_{\eps,k}$. Nonetheless, we can still discuss the convergence of the adaptive algorithm. Let $p^*\in H^1(\Omega)$ be the related adjoint state of $(u^*,y^*)$. The triplet $(u^*,y^*,p^*)$ satisfies the following optimality system associated with problem \eqref{min_G-L} subject to \eqref{vp}:
\begin{equation}\label{optsys}
    \left\{\begin{aligned}
            &(a(u^\ast) {\nabla} y^\ast, {\nabla} \phi)
    + (b(u^\ast) (y^\ast)^3,\phi) =(f, \phi), \quad \forall\phi \in H^1(\Omega),\\
        &(a(u^\ast) {\nabla} p^\ast, {\nabla} \psi) + 3 (b(u^\ast) (y^\ast)^2 p^\ast, \psi) = (y^\ast - y^\delta, \psi)_{L^2(\p\Omega)}, \quad \forall \psi \in H^1(\Omega) , \\
        &((1-\sigma) (v-u^\ast) {\nabla} y^\ast, {\nabla} p^\ast) + (v-u^\ast, (y^\ast)^3 p^\ast) \\
         &\quad + 2 \alpha \varepsilon ({\nabla} u^\ast, {\nabla} (v - u^\ast)) + \alpha\varepsilon^{-1} (1-2u^\ast,v - u^\ast) \geq 0, \quad \forall v \in \A.
        \end{aligned}
        \right.
\end{equation}
Then we have the following convergence result.
\begin{theorem}\label{thm:conv}
The sequence of solution triplets $\{(u_k^\ast,y_k^\ast,p_k^\ast)\}_{k\geq0}$ generated by Algorithm
\ref{afem_CE} contains a subsequence $\{(u_{k_j}^\ast,y_{k_j}^\ast,p_{k_j}^\ast)\}_{j\geq0}$ that converges to a
solution $(u^\ast,y^\ast,p^\ast)$ of the system \eqref{optsys} in the sense
    \begin{equation}\label{eq:conv}
        \|u^\ast_{k_j}-u^\ast\|_{H^1(\Omega)},\, \|y^\ast_{k_j}-y^\ast\|_{H^1(\Omega)},\, \|p^\ast_{k_j}-p^\ast\|_{H^1(\Omega)}\textcolor{blue}{ \to 0}\quad \text{as}~j\to \infty.
    \end{equation}
\end{theorem}

The proof of Theorem \ref{thm:conv} is lengthy and fairly technical, and thus it is deferred to Section \ref{sec:conv} below. In the analysis, a vanishing limit of the sequences $\eta_{k,i}$ ($1\leq i\leq 3$) of the estimators up to a subsequence is crucial. This is also precisely the motivation for the error estimators in Algorithm \ref{afem_CE}.

\begin{remark}
Theorem \ref{thm:conv} gives only a subsequence convergence. This is due to a lack of uniqueness of the solutions to the optimal system, or equivalently strong convexity of the objective $\mathcal{J}_\eps$, which in itself is attributed to the nonlinearity of the forward map $u\mapsto y(u)|_{\partial\Omega}$. Note that the strong convexity of the objective plays a crucial role in the a posteriori error analysis of PDE optimal control problems in establishing reliability; see e.g., \cite{HintHoppeIliashKieweg:2008,LiuYan:2001,KohlsRoschSiebert:2014}. The absence of strong convexity brings enormous technical challenges in the convergence analysis of the proposed adaptive algorithm, including establishing the reliability of the algorithm.
\end{remark}

\section{Proof of Theorem \ref{thm:conv}}\label{sec:conv}
The starting point of the analysis
is an auxiliary problem posed over a limiting set as first observed in \cite{babuska:1984}. With the sequences $\{V_k\}_{k\geq0}$ and $\{\A_k\}_{k\geq0}$ generated by Algorithm
\ref{afem_CE}, we define
\[
    V_\infty: = \overline{\bigcup_{k\geq 0}V_k}~ \text{in}~H^1(\Omega)\text{-norm}\quad\mbox{and}\quad
    \A_\infty: = \overline{\bigcup_{k\geq 0}\A_k}~\text{in}~H^1(\Omega)\text{-norm}.
\]
Then $\A_\infty$ is a closed convex subset of $\A$ \cite[Lemma 4.1]{JinXuZou:2016}, and $V_\infty$ is a closed
subspace of $H^1(\Omega)$. Hence $V_\infty$ is a Hilbert space and separable. Over
the set $\A_\infty$, we consider a limiting minimization problem:
\begin{equation}\label{medmin}
    \underset{u_\infty\in \A_\infty}\min \left\{\J_{\eps,\infty} (u_\infty) = \frac{1}{2} \|y_\infty - y^\delta\|^2_{L^2(\p\Omega)} + \alpha \left( \eps \int_\Omega |{\nabla} u_\infty|^2 \dx + \frac{1}{\eps} \int_\Omega u_\infty(1-u_\infty)  \dx \right) \right\},
\end{equation}
where $y_\infty: = y_\infty(u_\infty)\in V_\infty$ solves
\begin{equation}\label{medvp_state}
     (a(u_\infty){\nabla} y_\infty,{\nabla} \phi_\infty)
    + (b(u_\infty) y_\infty^3, \phi_\infty) = (f, \phi_\infty), \quad \forall \phi_\infty \in V_\infty.
\end{equation}

We use the following elementary result in the convergence analysis.
\begin{lemma}\label{lem:sol_medvp}
Under Assumptions \ref{ass01}--\ref{ass02}, for any $u_\infty\in \A_\infty$, there exists a
nonempty set $\Omega_\infty \subseteq \Omega$ and a constant $Q_\infty$ such that
$b(u_\infty) y_\infty^2(u_\infty) \geq Q_\infty$ almost everywhere in $ \Omega_\infty$.
\end{lemma}

\begin{proof}
The proof is inspired by \cite[Proposition 6.6]{beretta2017b}. By the definition of
$\A_\infty$, there exists a sequence $\{u_k\in\mathcal{A}_k\}_{k\geq 0}$ such that $u_k\to u_\infty$ in $H^1(\Omega)$. We may extract a subsequence, still denoted by $\{u_k\}_{k\geq0}$, that converges to $u_\infty$ almost everywhere (a.e.) in $\Omega$. Since $u_k=0$ in $\Omega^{d_0}$ for each $k\geq0$ (cf. the definition $\A_k = V_k \cap \A$), $u_\infty=0$ a.e. in $\Omega^{d_0}$. Now assume that $b(u_\infty) y_\infty^2(u_\infty)=0$ a.e. in $\Omega$. This and \eqref{medvp_state} imply
\[
   (a(u_\infty){\nabla}y_\infty,{\nabla} \phi_\infty)= (f, \phi_\infty), \quad \forall\phi_\infty \in V_\infty.
\]
Setting $\phi_\infty = c$ in the identity yields $\int_{\Omega}f \dx=0$, which contradicts Assumption \ref{ass02}.
\end{proof}

In Section \ref{ssec:limiting}, we prove that the sequence of discrete
solutions generated by Algorithm \ref{afem_CE} contains a subsequence convergent to a minimizer, the associated state of problem
\eqref{medmin}-\eqref{medvp_state} and the solution of the adjoint problem (cf. \eqref{medvp_costate} below). This also
implies the existence of a minimizer to problem \eqref{medmin}-\eqref{medvp_state}. In Section \ref{ssec:estimator}, we
establish the convergence to zero for the associated subsequence of estimators. Finally, the main result follows once the
auxiliary minimizer, the related state and the solution of  \eqref{medvp_costate} are shown to solve the optimality system \eqref{optsys}
associated with problem \eqref{min_G-L} in Section \ref{ssec:conv-overall}.

\subsection{The properties of the limit problem}\label{ssec:limiting}
Over each sequence of meshes $\{\cT_k\}_{k\geq0}$, discrete sets $\{\A_k\}_{k\geq0}$ and discrete spaces $\{V_k\}_{k\geq0}$ generated by Algorithm \ref{afem_CE}, the constraint \eqref{disvp} yields a sequence of discrete forward maps $\{u_k \mapsto y_k(u_k)\}_{k\geq0}$. This naturally motivates the study on the continuity of $y_k(u_k)$ as $k\to\infty$.
\begin{lemma}\label{lem:cont_control-state}
    If the sequence $\{u_k\in\mathcal{A}_k\}_{k\geq0}$ converges to some $u_\infty\in \A_\infty$ in $L^1(\Omega)$, then
\begin{equation}\label{cont_control-state}
    \lim_{k\to\infty}\| y_{\infty}(u_\infty) - y_k(u_k)\|_{H^1(\Omega)}=0.
\end{equation}
\end{lemma}
\begin{proof}
For the given sequence $\{u_k\}_{k\geq0}$ and $u_\infty$, we define a sequence of discrete operators  $A_k:V_k\to V_k^\ast$ and $A_\infty : V_\infty \to V_\infty^\ast$ respectively by
\begin{align*}
  \langle A_k v_k, \phi_k \rangle &= (a(u_k){\nabla} v_k, {\nabla} \phi_k) + (b(u_k) y_k^3, \phi_k),  \quad \forall \phi_k \in V_k,\\
    \langle A_\infty v_\infty,
    \phi_\infty \rangle &= (a(u_\infty){\nabla} v_\infty, {\nabla} \phi_\infty)
    + (b(u_\infty) v_\infty^3, \phi_\infty),  \quad \forall \phi_\infty \in V_\infty.
\end{align*}
The function $t\mapsto \langle A_k (v_k+t w_k), \phi_k \rangle$ (respectively $t\mapsto \langle A_\infty (v_\infty+t w_\infty), \phi_\infty \rangle$) is continuous
from $\mathbb{R}$ into $\mathbb{R}$ for all $v_k, w_k, \phi_k \in V_k$ (respectively all $v_\infty,w_\infty, \phi_\infty\in V_\infty$). By the proof of \cite[Proposition 2.1]{beretta2018a} and Assumption \ref{ass01}, the operator $A_k$ (respectively $A_\infty$) is strictly monotone and coercive. Since $V_k$ and $V_\infty$ both are separable Hilbert spaces, Browder-Minty theorem \cite[Theorem 26.A]{zeidler1990} implies that problem \eqref{disvp} over each $\cT_k$ ($k\geq0$) (respectively \eqref{medvp_state}) has a unique solution in $V_k$ (respectively $V_\infty$). Furthermore, we have the following growth condition from below (cf. the proof of \cite[Proposition 2.1]{beretta2018a})
\[
    \langle A_k v_k , v_k \rangle - (f,v_k) \geq \tfrac{\sigma}{c} \|v_k\|^2_{H^1(\Omega)} - \|f\|_{L^2(\Omega)} \|v_k\|_{H^1(\Omega)} - \tfrac{\sigma^2 |\Omega|}{4},\quad \forall k\geq 0,
\]
where the constant $c$ depends only on $\Omega$ and $d_0$. Thus there exists a constant $R>0$, independent of $k$, such that
each solution $y_k(u_k)$ of \eqref{disvp} satisfies $\|y_k(u_k)\|_{H^1(\Omega)} \leq R$, i.e., $\{\|y_k(u_k)\|_{H^1(\Omega)}\}_{k \geq 0}$ is uniformly bounded. Since $V_\infty$ is closed, this assertion, the $L^1(\Omega)$ convergence of $\{u_k\}_{k\geq0}$ to $u_\infty$ and Sobolev compact embedding  \cite[Theorem 6.3]{AdamsFournier:2003} imply that there exist some $\overline{y}\in V_\infty$ and two subsequences $\{u_{k_j}\}_{j\geq 0}$ and $\{y_{k_j}\equiv y_{k_j}(u_{k_j})\}_{j\geq 0}$ such that
\begin{equation}\label{pf:cont_control-state02}
    u_{k_j} \to u_\infty \quad \text{a.e. in}~\Omega, \quad y_{k_j}(u_{k_j}) \rightharpoonup \overline{y}\quad \text{in}~H^1(\Omega), \quad y_{k_j} \to \overline{y} \quad \text{in}~L^4(\Omega), \quad \mbox{as}~j\to \infty.
\end{equation}
It remains to show that $\overline{y}$ solves \eqref{medvp_state} and the $H^1(\Omega)$ weak convergence in \eqref{pf:cont_control-state02} can be lifted to a strong one. This is proved in two steps.

\medskip

\noindent \textit{Step 1.} For any $\phi_\infty \in V_\infty$, the density of $\bigcup_{k\geq0}V_k$ in $V_\infty$ ensures the existence of an $H^1(\Omega)$ convergent sequence $\{\phi_k\}_{k\geq0}$.
In the identity $a(u_{k_j}) {\nabla} \phi_{k_j} - a(u_{\infty}) {\nabla} \phi_{\infty} = a(u_{k_j}) ({\nabla} \phi_{k_j} - {\nabla} \phi_{\infty} ) + ( a(u_{k_j}) -a(u_\infty) ) {\nabla} \phi_{\infty}$,
the convergence of $\{\phi_{k_j}\}_{j\geq 0}$ in $H^1(\Omega)$ gives
\[
    \lim_{j\to\infty}\| a(u_{k_j}) ({\nabla} \phi_{k_j} - {\nabla} \phi_{\infty} )\|_{L^2(\Omega)}  \leq \lim_{j\to\infty}\|  {\nabla} \phi_{k_j} - {\nabla} \phi_{\infty} \|_{L^2(\Omega)} = 0,
\]
and the pointwise convergence in \eqref{pf:cont_control-state02}, $|u_{k_j} - u_\infty| \leq 1$ and Lebesgue's dominated convergence theorem \cite[Theorem 1.19, p. 28]{Evans:2015} yield
\[
    \lim_{j\to\infty}\| ( a(u_{k_j}) -a(u_\infty) ) {\nabla} \phi_{\infty} \|_{L^2(\Omega)} \leq \lim_{j\to\infty}\| ( u_{k_j} - u_\infty ) {\nabla} \phi_{\infty} \|_{L^2(\Omega)} = 0.
\]
Thus, we have
$\|a(u_{k_j}) {\nabla} \phi_{k_j} - a(u_{\infty}) {\nabla} \phi_{\infty} \|_{L^2(\Omega)}
    \to 0$ for any $\phi_\infty\in V_\infty$,
which, together with the $H^1(\Omega)$ weak convergence in \eqref{pf:cont_control-state02}, yields
\begin{equation}\label{pf:cont_control-state03}
  \lim_{j\to\infty}  ( a(u_{k_j}) {\nabla} y_{k_j} , {\nabla} \phi_{k_j}  ) = ( a(u_{\infty}) {\nabla} \overline{y}, {\nabla} \phi_{\infty}) ,\quad \forall \phi_\infty\in V_\infty.
\end{equation}
The generalized H\"{o}lder inequality, the continuous Sobolev embedding $H^1(\Omega)\hookrightarrow L^p(\Omega)$ (for any $ p  \geq2$ in 2D) and Lebesgue's dominated convergence theorem imply as $j\to\infty$,
\begin{align*}
    \left| ( b(u_{k_j}) \phi_{k_j} - b(u_\infty) \phi_\infty ) , y_{k_j}^3)\right| &\leq \| b(u_k) \phi_{k_j} - b(u_\infty) \phi_\infty\|_{L^4(\Omega)} \| y_{k_j} \|_{L^4(\Omega)}^3 \\
    &\leq c (\| b(u_{k_j}) ( \phi_{k_j} - \phi_\infty )\|_{L^4(\Omega)} + \| ( b(u_{k_j}) - b(u_\infty) )\phi_\infty \|_{L^4(\Omega)}) \\
    &\leq c ( \| \phi_{k_j} - \phi_\infty \|_{H^1(\Omega)}  +  \| ( u_{k_j} - u_\infty ) \phi_\infty\|_{L^4(\Omega)} ) \to 0.
\end{align*}
Similarly,  we have as $j\to\infty$,
\begin{align*}
\left| ( b(u_\infty) \phi_\infty , y_{k_j}^3 - \overline{y}^3 )\right|
     \leq  \|y_{k_j} - \overline{y} \|_{L^4(\Omega)}  \| y_{k_j}^2 + y_{k_j}\overline{y} + \overline{y}^2\|_{L^2(\Omega)} \| \phi_\infty \|_{L^4(\Omega)} \to 0.
\end{align*}
Then collecting these two limits leads to
\begin{equation}\label{pf:cont_control-state04}
  \lim_{j\to\infty}  (b(u_{k_j}) y_{k_j}^3 , \phi_{k_j}) = ( b(u_{\infty}) y_{\infty}^3 , \phi_{\infty}), \quad \forall \phi_\infty \in V_\infty.
\end{equation}
The limits \eqref{pf:cont_control-state03}-\eqref{pf:cont_control-state04} and the $H^1(\Omega)$ convergence of $\phi_{k_j}$ to $\phi_{\infty}$ yield that $\overline{y}$ solves problem \eqref{medvp_state}.

\medskip
\noindent \textit{Step 2.} Prove $\|\nabla ( y_{k_j} - \overline{y}) \|_{L^2(\Omega)} \to 0$. Setting $\phi_\cT = y_{k_j}$ in \eqref{disvp} over $\cT_{k_j}$ and then using \eqref{pf:cont_control-state02} give
\begin{equation}\label{pf:cont_control-state06}
    \begin{aligned}
    & \quad  \|a(u_{k_j})^{\frac{1}{2}} {\nabla} y_{k_j}\|^2_{L^2(\Omega)} + \|{b(u_{k_j})}^{\frac{1}{4}} y_{k_j}\|^4_{L^4(\Omega)} = (f,  y_{k_j}) \\
    & \to (f , \overline{y}) = \|a(u_{\infty})^{\frac{1}{2}} {\nabla} \overline{y}\|^2_{L^2(\Omega)} + \|{b(u_\infty)}^{\frac{1}{4}} \overline{y}\|^4_{L^4(\Omega)}, \quad \mbox{as}~j\to \infty.
    \end{aligned}
\end{equation}
Next, by the generalized H\"{o}lder inequality and the $L^4(\Omega)$ convergence in \eqref{pf:cont_control-state02},  we have
\begin{align*}
    &\quad\lim_{j\to\infty}\left| \|{b(u_{k_j})}^{\frac{1}{4}} y_{k_j}\|^4_{L^4(\Omega)} - \|{b(u_{k_j})}^{\frac{1}{4}} \overline{y}\|^4_{L^4(\Omega)} \right| \\
    &\leq \lim_{j\to\infty}(\| y_{k_j}\|_{L^4(\Omega)} + \| \overline{y}\|_{L^4(\Omega)} ) ( \|y_{k_j}\|_{L^4(\Omega)}^2 + \|\overline{y}\|_{L^4(\Omega)}^2 )\| y_{k_j}- \overline{y}\|_{L^4(\Omega)} = 0.
    \end{align*}
Using Lebesgue's dominated convergence theorem, we find
\[
   \lim_{j\to\infty} \left| \|{b(u_{k_j})}^{\frac{1}{4}} \overline{y}\|^4_{L^4(\Omega)} - \|{b(u_\infty)}^{\frac{1}{4}} \overline{y}\|^4_{L^4(\Omega)} \right| \leq \lim_{j\to\infty}\| (u_{k_j} - u_\infty) \overline{y}^4\|_{L^1(\Omega)} = 0.
\]
These two limits yield $\|{b(u_{k_j})}^{\frac{1}{4}} y_{k_j}\|^4_{L^4(\Omega)} \to \|{b(u_\infty)}^{\frac{1}{4}} \overline{y}\|^4_{L^4(\Omega)} $. This and \eqref{pf:cont_control-state06} give
\begin{equation*}
   \lim_{j\to\infty} \|a(u_{k_j})^{\frac{1}{2}} {\nabla} y_{k_j}\|^2_{L^2(\Omega)}  = \|a(u_{\infty})^{\frac{1}{2}} {\nabla} \overline{y}\|^2_{L^2(\Omega)}.
\end{equation*}
Moreover, Lebesgue's dominated convergence theorem implies \begin{align*}
    \lim_{j\to\infty}\|(a(u_{k_j}) - a(u_\infty)){\nabla}\overline{y}\|_{L^2(\Omega)} = 0.
\end{align*}
By the $H^1(\Omega)$ weak convergence in \eqref{pf:cont_control-state02} again, there holds
$$\lim_{j\to\infty} (a(u_{k_j}) \nabla y_{k_j},\nabla \overline{y} ) = (a(u_{\infty}) \nabla \overline{y},\nabla \overline{y} ).$$
Similarly,
$$\lim_{j\to\infty}\|a(u_{k_j})^{\frac{1}{2}}{\nabla}\overline{y}\|_{L^2(\Omega)}^2 = \|a(u_{\infty})^{\frac{1}{2}}{\nabla}\overline{y}\|_{L^2(\Omega)}^2.$$
The last three limits and the identity
\[
    \|a(u_{k_j})^{\frac{1}{2}}\nabla( y_{k_j} - \overline{y} )\|_{L^2(\Omega)}^2 =  \|a(u_{k_j})^{\frac{1}{2}} {\nabla} y_{k_j}\|^2_{L^2(\Omega)} - 2 (a(u_{k_j}) \nabla y_{k_j} ,\nabla \overline{y} ) + \|a(u_{k_j})^{\frac{1}{2}}{\nabla}\overline{y}\|_{L^2(\Omega)}^2
\]
imply the $H^1(\Omega)$ convergence of the sequence $\{y_{k_j}\}_{j\geq 0}$. Since problem \eqref{medvp_state} has a unique solution, the $H^1(\Omega)$ convergence for the whole sequence with  $\overline{y} = y_\infty(u_\infty)$ follows by a standard subsequence contradiction argument.
\end{proof}

\begin{theorem}\label{thm:conv_medmin}
Let $\{\A_k\times V_k\}_{k\geq 0}$ be a sequence of discrete admissible sets and FEM spaces generated by Algorithm \ref{afem_CE}. Then the sequence of adaptively-generated minimizers $\{(u_k^\ast, y_k^\ast)\}_{k\geq 0}$ of problem \eqref{disvp}--\eqref{min_G-L_dis} contains a subsequence $\{(u_{k_j}^\ast, y_{k_j}^\ast)\}_{j\geq 0}$ convergent to a minimizer $(u^\ast_\infty, y^\ast_\infty)$ of problem \eqref{medmin}--\eqref{medvp_state} in $H^1(\Omega)\times H^1(\Omega)$:
    \begin{equation}\label{conv_medmin}
        u_{k_j}^\ast \to u_\infty^\ast\quad\text{a.e. in}~\Omega, \quad u_{k_j}^\ast \to u_\infty^\ast,\quad y_{k_j}^\ast \to y_\infty^\ast\quad \text{in}~H^1(\Omega), \quad \mbox{as}~j\to \infty.
    \end{equation}
\end{theorem}
\begin{proof}
It follows from Lemma \ref{lem:cont_control-state} and the triangle inequality that
$\|y_k(0) - y_\infty(0)\|_{H^1(\Omega)}\to 0$. Then $\{\|y_k(0)\|_{H^1(\Omega)}\}_{k\geq 0}$
is uniformly bounded. Since $0 \in \A_k$ for all $k\geq 0$ and the functional $\J_{\eps,k}$ attains its minimum value at $u_k^\ast$ over $\A_k$, $0\leq \J_{\eps,k}(u_k^\ast) \leq  \J_{\eps,k}(0) \leq c$. This, the uniform bound $\|u_k^\ast\|_{L^2(\Omega)}\leq c$
(due to the box constraint in $\A_k$) and the boundedness of $\Omega$ imply that $\{\|u_k^\ast\|_{H^1(\Omega)}\}_{k\geq0}$ is uniformly bounded. Now we may extract a subsequence $\{u_{k_j}^\ast\}_{j\geq 0}$ and obtain some $u^\ast_\infty \in \A_\infty$ such that as $j\to\infty$,
\begin{equation}\label{pf:conv_med01}
        u^\ast_{k_j}\rightharpoonup u^\ast_{\infty}~\text{in}~ H^1(\Omega),
        \quad u^\ast_{k_j}\to u^\ast_{\infty}~\text{in}~L^2(\Omega),\quad
        u^\ast_{k_j}\to u^\ast_{\infty}~\text{a.e. in}~\Omega.
\end{equation}
By Lemma \ref{lem:cont_control-state} again and the trace theorem, we deduce that the subsequence $\{y_{k_j}^\ast\}_{j\geq 0}$ converges to $y^\ast_\infty$, the unique solution of problem \eqref{medvp_state} with $u_\infty = u^\ast_{\infty}$, in the $H^1(\Omega)$ norm, i.e.,
\begin{equation}\label{pf:conv_med02}
  \lim_{j\to\infty} \|y_{k_j}^\ast-y_\infty^\ast \|_{L^2(\partial\Omega)} \leq c\lim_{j\to\infty} \|y_{k_j}^\ast-y_\infty^\ast \|_{H^1(\Omega)}= 0.
\end{equation}
Moreover, the pointwise convergence in \eqref{pf:conv_med01} and Lebesgue's dominated convergence theorem imply
\begin{equation}\label{pf:conv_med03}
   \lim_{j\to\infty}( u^\ast_{k_j},1-u^\ast_{k_j}) = (u^\ast_{\infty},1-u^\ast_\infty).
\end{equation}
Meanwhile, for any fixed $u_\infty \in \A_\infty$, there exists a sequence $\{u_k\in\A_k\}_{k\geq0}$ such that $u_k\to u_\infty$ in $H^1(\Omega)$. The preceding argument leads to
$\lim_{k\to\infty} \|y_k(u_k)-y_\infty(u_\infty)\|_{L^2(\p\Omega)} = 0$ and $
     \lim_{k\to\infty}  (u_k,1-u_k) = (u_\infty,1-u_\infty)$.
By the relations \eqref{pf:conv_med02}--\eqref{pf:conv_med03}, the weak lower semicontinuity of the $H^1(\Omega)$-seminorm and the minimizing property  of $u_k^*$ to $\J_{\eps,k}(u_k^\ast)$ over $\A_k$, we obtain for any $u_\infty\in \A_\infty$
\begin{align}
&\J_{\eps,\infty}(u^\ast_\infty) \leq \liminf_{j\to\infty} \J_{\eps,k_j}(u^\ast_{k_j}) \leq \limsup_{j\to\infty}\J_{\eps,k_j}(u^\ast_{k_j}) \nonumber\\
\leq&\limsup_{k\to\infty}\J_{\eps,k}(u^\ast_{k}) \leq
\limsup_{k\to\infty}\J_{\eps,k}(u_{k}) = \J_{\eps,\infty}(u_\infty).\label{pf:conv_med04}
\end{align}
Thus, $u_\infty^\ast \in \A_\infty$ is a minimizer of $\J_{\eps,\infty}$ over $\A_\infty$. The first and third convergences in \eqref{conv_medmin} follow directly from \eqref{pf:conv_med01} and \eqref{pf:conv_med02}, respectively. By taking $u_\infty=u_\infty^\ast$ in \eqref{pf:conv_med04} and noting the relations \eqref{pf:conv_med02}--\eqref{pf:conv_med03}, we get
$\|{\nabla}u_{k_j}^\ast\|_{L^2(\Omega)}^2 \to \|{\nabla}u_{\infty}^\ast\|_{L^2(\Omega)}^2$ as $j\to\infty$. This and the weak convergence in \eqref{pf:conv_med01} yield the second assertion in \eqref{conv_medmin}.
\end{proof}

With a minimizing pair $(u_\infty^\ast, y_\infty^\ast)$ of the limiting problem \eqref{medmin}--\eqref{medvp_state}, we next study a limiting adjoint problem:  find $p_\infty^\ast\in V_\infty$ such that
\begin{equation}\label{medvp_costate}
    (a(u^\ast_\infty){\nabla} p_\infty^\ast,{\nabla} \psi_\infty) + 3(b(u_\infty^\ast) (y_\infty^\ast)^2 p_\infty^\ast, \psi_\infty) = (y_\infty^\ast-y^\delta, \psi_\infty)_{L^2(\p\Omega)}, \quad \forall \psi_\infty \in V_\infty.
\end{equation}
By Lemma \ref{lem:sol_medvp}, the following version of Poincar\'{e} inequality \cite[Appendix]{beretta2016}
\begin{equation}\label{eqn:Poincare}
 \|u\|_{H^1(\Omega)} \leq c_{\rm p}(\|\nabla u\|_{L^2(\Omega)} + \|u\|_{L^2(\Omega_\infty)}),\quad \forall u\in H^1(\Omega)
\end{equation}
with $c_{\rm p}$ depending on $\Omega$, and Lax-Milgram lemma,
problem \eqref{medvp_costate} has a unique solution  $p_\infty^*\in V_\infty$. Further, we define an auxiliary discrete
problem associated with \eqref{medvp_costate}: find $\widetilde{p}_k \in V_k$ such that
\begin{equation}\label{disvp_aux_costate}
    (a(u^\ast_\infty){\nabla} \widetilde{p}^\ast_k, {\nabla} \psi_k) + 3(b(u_\infty^\ast) {(y_\infty^\ast)^2} \widetilde{p}^\ast_k, \psi_k) = (y_\infty^\ast-y^\delta, \psi_k)_{L^2(\partial\Omega)}, \quad \forall \psi_k \in V_k.
\end{equation}
It is the Galerkin approximation of \eqref{medvp_costate} over $V_k$. By Cea's lemma and the density of $\bigcup_{k\geq0}V_k$ in $V_\infty$ with respect to the $H^1(\Omega)$ norm, we deduce
\begin{equation}\label{conv_aux1_costate}
   \lim_{k\to\infty} \|p_\infty^\ast-\widetilde{p}^\ast_k\|_{H^1(\Omega)}\leq \lim_{k\to\infty}c\inf_{\psi_k\in V_k} \|p_\infty^\ast - \psi_k\|_{H^1(\Omega)}= 0.
\end{equation}
\begin{theorem}\label{thm:conv_med_costate}
Under the condition of Theorem \ref{thm:conv_medmin}, the subsequence of discrete adjoint solutions $\{p^\ast_{k_j}\}_{j\geq0}$ generated by Algorithm \ref{afem_CE} converges to the solution $p_\infty^\ast$ of problem \eqref{medvp_costate}:
    \begin{equation}\label{conv_med_costate}
        \|p^\ast_{k_j} - p_\infty^\ast\|_{H^1(\Omega)} \to 0, \quad \mbox{as}~j\to \infty.
    \end{equation}
\end{theorem}
\begin{proof}
Let $\overline{u}_{k_j}^* = u_{k_j} - u_\infty^*$, $\overline{y}_{k_j}^* = {y}_{k_j}^*(u_{k_j}^*)-y_\infty^*(u_\infty^*)$ and $\widehat{p}_{k_j}^*=p_{k_j}^\ast -
\widetilde{p}_{k_j}^*$. By the triangle inequality and the estimate \eqref{conv_aux1_costate}, it suffices to show $\|\widehat{p}_{k_j}^\ast\|_{H^1(\Omega)}\to 0$. Note that problem \eqref{disvp_aux_costate} is a perturbed version of \eqref{optsys_dis-2} over $V_k$ with $u^\ast_\infty$ in place of $u_k^\ast$. Subtracting these two equations and then setting $\psi_{k_j}=\widehat{p}_{k_j}^*$ lead to
\begin{align*}
      (a(u^\ast_{k_j}){\nabla}\widehat{p}^\ast_{k_j},&{\nabla}\widehat{p}^\ast_{k_j})+ 3(b(u^\ast_{k_j}) (y_{k_j}^\ast)^2,(\widehat{p}^\ast_{k_j})^2) =
            (y_{k_j}^\ast - y^\ast_\infty,\widehat{p}^\ast_{k_j})_{L^2(\p\Omega)} \\
             & \quad + ((a(u_\infty^\ast)-a(u_{k_j}^\ast)) {\nabla} \widetilde{p}^\ast_{k_j} , {\nabla}\widehat{p}^\ast_{k_j}) + 3 ((b(u_\infty^\ast) (y^\ast_\infty)^2 - b(u_{k_j}^\ast) (y_{k_j}^\ast)^2) \widetilde{p}_{k_j}^\ast,\widehat{p}^\ast_{k_j}).
\end{align*}
By Lemma \ref{lem:sol_medvp}, since $a(u_{k_j})\geq \sigma$, we have
\begin{align*}
      \sigma \|{\nabla}\widehat{p}^\ast_{k_j}\|_{L^2(\Omega)}^2
        +3Q_\infty\|\widehat{p}^\ast_{k_j}\|_{L^2(\Omega_\infty)}^2 &\leq  \sigma\|{\nabla}\widehat{p}^\ast_{k_j}\|_{L^2(\Omega)}^2
        +3 (b(u^\ast_\infty)(y^\ast_\infty)^2,(\widehat{p}^\ast_{k_j})^2)\\
        & \leq (\overline{y}_{k_j}^\ast,\widehat{p}^\ast_{k_j})_{L^2(\p\Omega)}
              + ((1-\sigma)\overline{u}_{k_j}^\ast{\nabla} \widetilde{p}^\ast_{k_j},{\nabla}\widehat{p}^\ast_{k_j})\\
              & \quad + 3((b(u_\infty^\ast) (y^\ast_\infty)^2 - b(u_{k_j}^\ast) (y_{k_j}^\ast)^2) \widetilde{p}_{k_j}^\ast,\widehat{p}^\ast_{k_j}) \\
              &\quad+ 3(b(u_\infty^\ast) (y_\infty^\ast)^2- b(u_{k_j}^\ast)(y_{k_j}^\ast)^2,(\widehat{p}^\ast_{k_j})^2):=\sum_{i=1}^4{\rm I}_i.
\end{align*}
Now it follows from the trace theorem and the triangle inequality that
\begin{align*}
  {\rm I}_1 \leq c\|\overline{y}_{k_j}^\ast\|_{H^1(\Omega)}\|\widehat{p}_{k_j}^\ast \|_{H^1(\Omega)}\quad \mbox{and}\quad
  {\rm I}_2  \leq (\|{\nabla} (\widetilde{p}_{k_j}^\ast - p_\infty^\ast) \|_{L^2(\Omega)}  + \|\overline{u}_{k_j}^\ast{\nabla}  p^\ast_\infty\|_{L^2(\Omega)})\|{\nabla}\widehat{p}^\ast_{k_j}\|_{L^2(\Omega)}.
\end{align*}
By the generalized H\"{o}lder's inequality and Sobolev embedding $H^1(\Omega) \hookrightarrow L^p(\Omega)$ (for any $p \geq 2$ in 2D),
\begin{align*}
  {\rm I}_3&= 3(b(u_{k_j}^\ast)((y_\infty^\ast)^2-(y_{k_j}^\ast)^2) \widetilde{p}^\ast_{k_j},\widehat{p}^\ast_{k_j})
             + 3((b(u^\ast_\infty) - b(u^\ast_{k_j}))(y_\infty^\ast)^2  \widetilde{p}_{k_j}^\ast, \widehat{p}^\ast_{k_j})\\
       &\leq 3 \|\overline{y}_{k_j}^\ast\|_{L^4(\Omega)}\|y_{k_j}^\ast+y_\infty^\ast\|_{L^4(\Omega)}
              \|\widetilde{p}_{k_j}^\ast\|_{L^4(\Omega)}\|\widehat{p}^\ast_{k_j} \|_{L^4(\Omega)} + 3 \|\overline{u}_{k_j}^\ast(y_\infty^\ast)^2\|_{L^2(\Omega)}\|\widetilde{p}_{k_j}^\ast\|_{L^4(\Omega)}\|\widehat{p}^\ast_{k_j}\|_{L^4(\Omega)}\\
       &\leq  c(\|\overline{y}_{k_j}^\ast\|_{H^1(\Omega)}\|y_{k_j}^\ast+y_\infty^\ast\|_{H^1(\Omega)}\|\widetilde{p}_{k_j}^\ast\|_{H^1(\Omega)} + \|\overline{u}_{k_j}^\ast(y_\infty^\ast)^2\|_{L^2(\Omega)}\|\widetilde{p}_{k_j}^\ast\|_{H^1(\Omega)})
              \|\widehat{p}^\ast_{k_j}\|_{H^1(\Omega)},\\
       {\rm I}_4
      &\leq c(\|\overline{y}_{k_j}^\ast\|_{H^1(\Omega)}\|y_{k_j}^\ast+y_\infty^\ast\|_{H^1(\Omega)} +  \|\overline{u}_{k_j}^\ast(y_\infty^\ast)^2\|_{L^2(\Omega)})\|\widehat{p}^\ast_{k_j}\|_{H^1(\Omega)}^2.
\end{align*}
By Theorem \ref{thm:conv_medmin} and \eqref{conv_aux1_costate}, $\|\overline{y}_{k_j}^\ast\|_{H^1(\Omega)}$, $\|\overline{u}_{k_j}^\ast(y_\infty^\ast)^2\|_{L^2(\Omega)}\to0$ as $j\to \infty$ and $\{y_{k_j}^\ast\}_{j\geq0}$, $\{\widetilde{p}_{k_j}^\ast\}_{j\geq0}$ are uniformly bounded in $H^1(\Omega)$. Hence, for sufficiently large $j$, Poincar\'{e} inequality (cf. \eqref{eqn:Poincare}) implies
\[
        \|\widehat{p}^\ast_{k_j}\|_{H^1(\Omega)}^2
        \leq \frac{c_6\|\overline{y}_{k_j}^\ast\|_{H^1(\Omega)}+\|{\nabla} ( \widetilde{p}_{k_j}^* - p_\infty^\ast ) \|_{L^2(\Omega)}  + \|\overline{u}_{k_j}^\ast{\nabla}  p^\ast_\infty\|_{L^2(\Omega)} + c_7 \|\overline{u}_{k_j}^\ast(y_\infty^\ast)^2\|_{L^2(\Omega)}}{c_4 - c_5(\|\overline{y}_{k_j}^\ast\|_{H^1(\Omega)}+ \|\overline{u}_{k_j}^\ast(y_\infty^\ast)^2\|_{L^2(\Omega)})}\|\widehat{p}^\ast_{k_j}\|_{H^1(\Omega)}.
\]
The desired result follows from Theorem \ref{thm:conv_medmin}, \eqref{conv_aux1_costate}, $\|\overline{u}_{k_j}^\ast(y_\infty^\ast)^2\|_{L^2(\Omega)}\to0$ and
$\|\overline{u}_{k_j}^\ast{\nabla}  p^\ast_\infty\|_{L^2(\Omega)} \to 0$ as $j\to \infty$
(due to dominated convergence theorem and pointwise convergence of $u_{k_j}^*$ to $u_\infty^*$ in \eqref{conv_medmin}).
\end{proof}

\subsection{Zero limit of error estimators}\label{ssec:estimator}
In view of Theorems \ref{thm:conv_medmin} and \ref{thm:conv_med_costate}, to prove Theorem \ref{thm:conv},
it suffices to show that the limiting triplet $(u^\ast_\infty,y^\ast_\infty,p^\ast_\infty)$ solves the system \eqref{optsys}.
This depends on the properties of modules ESTIMATE, MARK and REFINE in Algorithm \ref{afem_CE}. Now we analyze the
behaviour of the error estimators in the adaptive loop \eqref{afem_loop}. We denote by $\omega(T)$ the union of the element $T$
and its neighbours in some $\cT\in\mathbb{T}$ and by $\omega(F)$ the union of elements in $\cT$ sharing
a common edge with $F\in\mathcal{F}_\cT$. Let $\cT_{k}^{+}:=\bigcap_{l\geq k}\cT_{l}$, i.e., the set of all elements
not refined after the $k$-th iteration, $\Omega_k^+:=\bigcup_{T\in\cT^+_k}T$ and $\Omega_{k}^{0}:=\bigcup_{T\in\cT_k\setminus\cT^{+}_{k}}T$. By the definition of $\mathcal{T}_k^+$, each element in the set $\cT_k\setminus\cT^{+}_{k}$ is refined (at least once) after the $k$-th iteration, and thus over the set $\Omega_k^0$, the meshsize function $h_k$ (cf. \eqref{meshsize_def}) has the following
property \cite[Corollary 3.3]{siebert2011}:
\begin{equation}\label{eqn:conv_zero_mesh}
    \lim_{k\rightarrow\infty}\|h_{k}\chi^0_{k}\|_{L^\infty(\Omega)}=0,
\end{equation}
where $\chi_k^0$ is the characteristic function of $\Omega_k^0$. Below
we adopt the shorthand notation $\overline{y}_{k_j,k_l}^* = y_{k_j}^* - y_{k_l}^*$, $\overline{y}_{k_j}^* = y_{k_j}^*-y_\infty$ and similarly $\overline{u}_{k_j,k_l}^*$, $\overline{u}_{k_j}^*$, $\overline{p}_{k_j,k_l}^*$ and $\overline{p}_{k_j}^*$.  Moreover, we have the following inverse estimates \cite[Lemmas 12.1 and 12.8]{ErnGuermond:2021}.

\begin{lemma}\label{lem:inverse}
There exists $c$ such that for all $p, r \in [1,\infty]$, $\ell,m\in\{0,1\}$ with $m\geq \ell$,
$$| v|_{W^{m,p}(T)}
\leq  ch^{\ell-m+d(\frac1p-\frac1r)}_T
|v|_{W^{\ell,r}(T)},\\.
\quad \forall v\in P_1(T),\ T \in \mathcal{T}.$$
Moreover, there is $c$ such that for all $p, r \in [1,\infty]$,
$$\|\nabla v\|_{L^p(F)}\leq ch^{-\frac1p +d(\frac1p-\frac1r)}_T
\|\nabla v\|_{L^r(T)}, \quad v \in  P_1(T), \ T\in \mathcal{T} ,\ F \subset\partial T.$$
\end{lemma}

The first auxiliary result is the following local stability for the three error indicators.

\begin{lemma}\label{lem:stab_est}
Let $\{(u_k^\ast,y_k^\ast, p_k^\ast)\}_{k\geq0}$ be the sequence generated by Algorithm \ref{afem_CE}. Then for any $T\in\cT_k$, there hold
\begin{align}\label{stab_est1}
        \eta_{k,1}(u^\ast_k,y_k^\ast,T)&\leq c\big(\| \nabla y_k^\ast\|_{L^2(\omega(T))} +h_T\|y_k^\ast\|_{L^6(T)}^3 + h_T\|f\|_{L^2(T)}\big),\\
    \label{stab_est2}
      \eta_{k,2}(u^\ast_k,y_k^\ast,p^\ast_k,T) &\leq c \big( \| \nabla {p}^\ast_k\|_{L^2(\omega(T))} + h_T\|y^\ast_k\|^2_{L^6(T)} \|p^\ast_k\|_{L^6(T)} + h_{T}^{1/2}\|y^\ast_k-y^\delta\|_{L^2(\p T\cap\partial\Omega)} \big),\\
      \eta_{k,3}(u^\ast_k,y_k^\ast,p^\ast_k,T) &\leq c\big( \| \nabla u_{k}^\ast\|_{L^2(\omega(T))}+\| \nabla y_k^\ast\|_{L^2(T)}\| \nabla p_k^\ast\|_{L^2(T)}\nonumber\\
      & \qquad +h_T^{2/3}\|y_k^\ast\|^3_{L^6(T)}\|p_k^\ast\|_{L^6(T)}+
        h_T \|2u^\ast_k - 1\|_{L^2(T)} \big).    \label{stab_est3}
    \end{align}
\end{lemma}
\begin{proof}
The identity $R_{T,1}(u^\ast_k,y_k^\ast) = (1-(1-\sigma)) \nabla u_k^\ast\cdot \nabla y_k^\ast -(1-u^\ast_k)(y_k^\ast)^3 + f$,  H\"{o}lder's inequality and the inverse estimate in Lemma \ref{lem:inverse} imply
        \begin{align*}
        h_T\|R_{T,1}(u^\ast_k,y_k^\ast)\|_{L^2(T)} &\leq h_T\| \nabla u_k^\ast\|_{L^4(T)}\| \nabla y_k^\ast \|_{L^4(T)} + h_T\| y_k^\ast \|_{L^6(T)}^3 + h_T\|f\|_{L^2(T)} \\
         & \leq c h_T h_T^{-1} \| \nabla u_k^\ast\|_{L^2(T)}\| \nabla y_k^\ast \|_{L^2(T)} + h_T \| y_k^\ast \|_{L^6(T)}^3 + h_T\|f\|_{L^2(T)}\\
         &\leq c \|u_{k}^\ast\|_{L^\infty(T)}\| \nabla y_k^\ast \|_{L^2(T)} + h_T \| y_k^\ast \|_{L^6(T)}^3 + h_T\|f\|_{L^2(T)}\\
         &\leq c\| \nabla y_k^\ast \|_{L^2(T)} + h_T \| y_k^\ast \|_{L^6(T)}^3 + h_T\|f\|_{L^2(T)}  .
        \end{align*}
    Further, by the inverse estimate in Lemma \ref{lem:inverse} again, we get
    \[
        h_T^{1/2}\|J_{F,1}(u_k^\ast,y^\ast_k)\|_{L^2(F)} \leq c \| \nabla y^\ast_k\|_{L^2(\omega(F))},\quad \forall F \subset \partial T.
    \]
     The inequality \eqref{stab_est1} follows from these two estimates. A similar argument yields
     \begin{align*}
        h_T\|R_{T,2}(u^\ast_k,y^\ast_k,p^\ast_k)\|_{L^2(T)} &\leq c \| \nabla p_k^\ast\|_{L^2(T)} + 3 h_T\|y^\ast_k\|^2_{L^6(T)} \|p^\ast_k\|_{L^6(T)},\\
        h_T^{1/2}\|J_{F,2}(u^\ast_k,y^\ast_k,p^\ast_k)\|_{L^2(F)}
        &\leq c\| \nabla p^\ast_k\|_{L^2(\omega(F))}+h_{T}^{1/2}\|y^\ast_k-y^\delta\|_{L^2(F\cap\partial\Omega)},\quad \forall F\subset \partial T,
     \end{align*}
which give the estimate \eqref{stab_est2}. Finally, by H\"{o}lder's inequality and  Lemma \ref{lem:inverse}, we get
\begin{align}
    & \quad h_T\|R_{T,3}(u_k^\ast,y_k^\ast,p_k^\ast)\|_{L^2(T)} \nonumber\\
    &\leq h_T \| \nabla y_k^\ast\cdot \nabla p_k^\ast\|_{L^2(T)}  + h_T \|y_k^\ast\|^3_{L^8(T)}\|p_k^\ast\|_{L^8(T)} + h_T\alpha \varepsilon^{-1} \|2u^\ast_k - 1\|_{L^2(T)} \label{eqn:RT3-inverse0} \\
    &\leq c \left(h_T h_T^{-1}\| \nabla y_k^\ast\cdot \nabla p_k^\ast\|_{L^1(T)}+h_T h_T^{-1/3}\|y_k^\ast\|^3_{L^6(T)}\|p_k^\ast\|_{L^6(T)}+ h_T \|2u^\ast_k - 1\|_{L^2(T)} \right)\label{eqn:RT3-inverse1}\\
    &\leq c\left(\| \nabla y_k^\ast\|_{L^2(T)}\| \nabla p_k^\ast\|_{L^2(T)}+ h_T^{2/3}\|y_k^\ast\|^3_{L^6(T)}\|p_k^\ast\|_{L^6(T)}+
        h_T \|2u^\ast_k - 1\|_{L^2(T)} \right),\nonumber\\
        &\quad h_T^{1/2} \| J_{F,3},
         (u_k^\ast)\|_{L^2(F)} \leq c \| \nabla u_{k}^\ast\|_{L^2(\omega(F))}.\nonumber
\end{align}
Thus, the desired estimate \eqref{stab_est3} holds. This completes the proof of the lemma.
\end{proof}

Lemma \ref{lem:stab_est} enables studying the limit of the error indicator $\eta_{k,i}(T)$ over the marked elements.

\begin{lemma}\label{lem:conv_aux5_maxest}
Let $\{(u_{k_j}^\ast,y_{k_j}^\ast,p_{k_j}^\ast)\}_{j\geq 0}$ be the convergent subsequence in Theorems \ref{thm:conv_medmin} and \ref{thm:conv_med_costate} of the sequence $\{(u_k^\ast,y_k^\ast,p_k^\ast)\}_{k\geq0}$ generated by Algorithm \ref{afem_CE}. Then there holds
    \begin{equation}\label{conv_aux5_maxest}
        \lim_{j\to\infty}\Big(\max_{T\in\mathcal{M}_{k_j}}\eta_{k_j,1}(u^\ast_{k_j},y_{k_j}^\ast,T)+\max_{T\in\mathcal{M}_{k_j}}\eta_{k_j,2}(u^\ast_{k_j},y_{k_j}^\ast,p^\ast_{k_j},T)
        +\max_{T\in\mathcal{M}_{k_j}}\eta_{k_j,3}(u^\ast_{k_j},y_{k_j}^\ast,p^\ast_{k_j},T)\Big)=0.
    \end{equation}
\end{lemma}
\vspace{-.3cm}
\begin{proof}
Let $T_j^{i}$ ($i=1,2,3$) be the element with the largest error indicator $\eta_{k_j,i}(T)$ among $\mathcal{M}_{k_j}$. Since $\mathcal{M}_{k_j}\subset\cT_{k_j}^0$, by the identity  \eqref{eqn:conv_zero_mesh}, we have
\begin{equation}\label{pf:aux_conv_maxest}
    |T_{j}^i| \leq |\omega(T_{j}^i)| \leq c h_{T_{j}^i}^2 \leq c \|h_{k_j}\chi_{k_j}^0\|^2_{L^\infty(\Omega)} \to 0 \quad \mbox{as}~j \to \infty.
\end{equation}
It follows from Lemma \ref{lem:stab_est}, Sobolev embedding theorem and the trace theorem that
\begin{align*}
        \eta_{k_j,1}(u_{k_j}^\ast,y_{k_j}^\ast,T_{j}^1) &\leq c \big(\| \nabla \overline{y}_{k_j}^\ast\|_{H^1(\Omega)}+\|y_\infty\|_{H^1(\omega(T_j^1))}+h_{T^1_j}\|y^\ast_{k_j}\|^3_{H^1(\Omega)}+h_{T^1_j}\|f\|_{L^2(\Omega)}\big),\\
        \eta_{k_j,2}(u_{k_j}^\ast,y_{k_j}^\ast,p_{k_j}^\ast,T_{j}^2) &\leq c \big( \| \nabla\overline{p}_{k_j}^\ast\|_{H^1(\Omega)}+\|p_\infty\|_{H^1(\omega(T_j^2))} + h_{T^2_j}\|y^\ast_{k_j}\|^2_{H^1(\Omega)}\|p^\ast_{k_j}\|_{H^1(\Omega)}\\
        &\qquad+h_{T^2_j}^{1/2}\big(\|y_{k_j}^\ast\|_{H^1(\Omega)}+\|y^\delta\|_{L^2(\p\Omega)}\big)
        \big).
        \end{align*}
By the estimate \eqref{pf:aux_conv_maxest}, the absolute continuity of $\|\cdot\|_{H^1(\Omega)}$ with respect to the Lebesgue measure implies that $\|y_\infty\|_{H^1(\omega(T_j^1))}$ and $\|p_\infty\|_{H^1(\omega(T_j^2))} $ both tend to zero.
In view of Theorems \ref{thm:conv_medmin} and \ref{thm:conv_med_costate}, $\| \nabla\overline{y}_{k_j}^\ast\|_{H^1(\Omega)}$ and $\| \nabla \overline{p}_{k_j}^\ast\|_{H^1(\Omega)}$ tend to zero as $j\to \infty$, from which we also deduce that both $\{\|y^{\ast}_{k_j}\|_{H^1(\Omega)}\}_{j\geq 0}$ and $\{\|p^{\ast}_{k_j}\|_{H^1(\Omega)}\}_{j\geq 0}$ are uniformly bounded. This and the estimate \eqref{pf:aux_conv_maxest} imply $ h_{T^1_j}\|y^\ast_{k_j}\|^3_{H^1(\Omega)}$, $h_{T^1_j}\|f\|_{L^2(\Omega)}$, $h_{T^2_j}\|y^\ast_{k_j}\|^2_{H^1(\Omega)}\|p^\ast_{k_j}\|_{H^1(\Omega)}$ and $h_{T^2_j}^{1/2}(\|y_{k_j}^\ast\|_{H^1(\Omega)}+\|y^\delta\|_{L^2(\p\Omega)})$ all tend to zero as $j\to \infty$. Hence, $$
\lim_{j\to\infty}\max_{T\in\mathcal{M}_{k_j}}\eta_{k_j,1}(u^\ast_{k_j},y_{k_j}^\ast,T) =\lim_{j\to\infty}\max_{T\in\mathcal{M}_{k_j}}\eta_{k_j,2}(u^\ast_{k_j},y_{k_j}^\ast,p^\ast_{k_j},T)=0.$$ For the third limit, we make use of \eqref{stab_est3} and get
        \begin{align*}
              \eta_{k_j,3}(u^\ast_{k_j},y_{k_j}^\ast,p^\ast_{k_j},T) &\leq c\big( \|\nabla \overline{u}_{k_j}^\ast\|_{L^2(\Omega)} +  \|\nabla u_{\infty}\|_{L^2(\omega(T_j^3))} +
        h_{T_j^3} (2\|u^\ast_{k_j}\|_{L^2(\Omega)}+|\Omega|^{1/2}) \\ &\quad + h_{T_j^3}^{2/3}\|y_{k_j}^\ast\|^3_{L^6(\Omega)}\|p_{k_j}^\ast\|_{L^6(\Omega)} +(\| \nabla\overline{y}_{k_j}^\ast\|_{L^2(\Omega)}+\|\nabla y_{\infty}\|_{L^2(\omega(T_j^3))})\| \nabla p_{k_j}^\ast\|_{L^2(\Omega)} \big).
        \end{align*}
The proof is completed by appealing to \eqref{pf:aux_conv_maxest} again and repeating the preceding argument.
\end{proof}

Below we establish the stability on non-refined elements (Lemma \ref{lem:stab_est_unref}) and reduction on refined elements (Lemma \ref{lem:reduction_est_ref}) for the estimators. These two properties and the behaviour of the maximal error indicators in Lemma \ref{lem:conv_aux5_maxest} ensure the existence of a null sequence of estimators associated with  the sequence $\{(u^\ast_{k_j},y^\ast_{k_j},p^\ast_{k_j})\}_{j\geq 0}$. The analysis strategy is inspired by the works \cite{GanterPraetorius:2022,morin2008,siebert2011} for direct problems. The proof of the following stability results is lengthy and deferred to the appendix.

\begin{lemma}\label{lem:stab_est_unref}
    Let $\{(u_{k_j}^\ast,y_{k_j}^\ast,p_{k_j}^\ast)\}_{j\geq 0}$ be the convergent subsequence in Theorems \ref{thm:conv_medmin} and \ref{thm:conv_med_costate} of the sequence $\{(u_k^\ast,y_k^\ast,p_k^\ast)\}_{k\geq0}$ generated by Algorithm \ref{afem_CE}. Then for $j<l$, there hold
    \begin{align}
        \eta_{k_l,1}(\cT_{k_j}\cap \cT_{k_l})  &\leq \eta_{k_j,1}(\cT_{k_j}\cap \cT_{k_l}) + c (\|\overline{u}_{{k_l},k_j}^\ast\|_{H^1(\Omega)} + \|\overline{y}_{{k_l},k_j}^\ast \|_{H^1(\Omega)} \nonumber\\
         & \quad + \|{\nabla}\overline{y}_{k_j}^\ast\|_{L^2(\Omega)}  + \|\overline{u}_{{k_l},k_j}^\ast {\nabla} y_{\infty}^\ast\|_{L^2(\Omega)}),\label{stab_est_unref1}\\
    \eta_{k_l,2}(\cT_{k_l}\cap \cT_{k_j}) &\leq \eta_{k_j,2}(\cT_{k_l}\cap \cT_{k_j}) + c(\|\overline{u}_{k_l,k_j}^\ast\|_{H^1(\Omega)} + \|\overline{y}_{k_l,k_j}^\ast\|_{H^1(\Omega)} + \|\overline{p}_{k_l,k_j}^\ast \|_{H^1(\Omega)} \nonumber \\
       &\quad +   \|{\nabla}\overline{p}_{k_j}^\ast\|^2_{L^2(\Omega)}
       + \|\overline{u}_{k_l,k_j}^\ast {\nabla} p_{\infty}^\ast\|^2_{L^2(\Omega)} ),\label{stab_est_unref2}\\
        \eta_{k_l,3}(\cT_{k_j}\cap \cT_{k_l}) &\leq \eta_{k_j,3}(\cT_{k_j}\cap \cT_{k_l}) + c (\|\overline{u}_{k_l,k_j}^\ast \|_{H^1(\Omega)} + \|\overline{y}_{k_l,k_j}^\ast \|_{H^1(\Omega)} + \|\overline{p}_{k_l,k_j}^\ast \|_{H^1(\Omega)}).\label{stab_est_unref3}
    \end{align}
\end{lemma}

\begin{remark}
The convergence analysis of the algorithm is restricted to the case $d=2$. This restriction is mainly due to the use of the inverse inequality in the proof of Lemmas \ref{lem:stab_est} and \ref{lem:stab_est_unref}. For example, in the proof of Lemma \ref{lem:stab_est} for $d=3$, the inverse estimate implies that in the step from \eqref{eqn:RT3-inverse0} to \eqref{eqn:RT3-inverse1}, the term for bounding $h_T\|R_{T,3}(u_k^\ast,y_k^\ast,p_k^\ast)\|_{L^2(T)}$ becomes $h_T h_T^{-3/2}\| \nabla y_k^\ast\cdot \nabla p_k^\ast\|_{L^1(T)}$, which involves a negative power of the local meshsize and thus the convergence analysis does not work any more. This flaw also occurs to \eqref{pf:stab_est_unref1_09} in the proof of Lemma \ref{lem:stab_est_unref} for $d=3$ (see the appendix). This might be due to a limitation of the proof technique, and there might be alternative proof techniques that can overcome the restriction.  Alternatively one may modify the estimators to weaker norms for $d=3$.
\end{remark}

The next lemma gives an important contraction property of the estimators.
\begin{lemma}\label{lem:reduction_est_ref}
    Let $\{(u_{k_j}^\ast,y_{k_j}^\ast,p_{k_j}^\ast)\}_{j\geq 0}$ be the convergent subsequence given in Theorems \ref{thm:conv_medmin} and \ref{thm:conv_med_costate} of the sequence $\{(u_k^\ast,y_k^\ast,p_k^\ast)\}_{k\geq0}$ generated by Algorithm \ref{afem_CE}. Then for $j<l$, there exist constants $q_i\in (0,1)$, $i=1,2,3$, such that
    \begin{align}
        \eta^2_{k_{l},1}(\cT_{k_l}\setminus \cT_{k_j}) &\leq q_{1}\eta^2_{k_j,1}(\cT_{k_j}\setminus \cT_{k_l}) + c(\|\overline{u}_{k_l,k_j}^\ast \|^2_{H^1(\Omega)} + \|\overline{y}_{k_l,k_j}^\ast\|^2_{H^1(\Omega)}\nonumber \\
         &\quad + \|{\nabla}\overline{y}_{k_j}^\ast\|^2_{L^2(\Omega)}
         +\|\overline{u}_{k_l,k_j}^\ast {\nabla} y_{\infty}^\ast\|_{L^2(\Omega)}^2),\label{reduction_est_ref1}\\
        \eta^2_{k_l,2}(\cT_{k_l}\setminus \cT_{k_j}) &\leq q_{2}\eta_{k_j,2}^2(\cT_{k_j}\setminus \cT_{k_l}) + c(\|\overline{u}_{k_l,k_j}^\ast \|^2_{H^1(\Omega)} + \|\overline{y}_{k_l,k_j}^\ast\|^2_{H^1(\Omega)} + \|\overline{p}_{k_l,k_j}^\ast\|^2_{H^1(\Omega)}        \nonumber\\
        &\quad +   \|{\nabla}\overline{p}_{k_j}^\ast \|^2_{L^2(\Omega)}
       + \|\overline{u}_{k_l,k_j}^\ast {\nabla} p_{\infty}^\ast\|^2_{L^2(\Omega)}),\label{reduction_est_ref2}\\
        \eta_{k_l,3}^2(\cT_{k_l}\setminus \cT_{k_j}) &\leq q_{3}\eta_{k_j,3}^2(\cT_{k_j}\setminus \cT_{k_l}) + c (\|\overline{u}_{k_l,k_j}^\ast\|^2_{H^1(\Omega)} + \|\overline{y}_{k_l,k_j}^\ast\|^2_{H^1(\Omega)} + \|\overline{p}_{k_l,k_j}^\ast \|^2_{H^1(\Omega)}).
        \label{reduction_est_ref3}
    \end{align}
\end{lemma}
\begin{proof}
By applying the triangle inequality on each element $T \in \cT_{k_l}\setminus\cT_{k_j}$, we derive
\[
        \eta_{k_l,1}(u_{k_l}^\ast,y^\ast_{k_l},T) \leq \eta_{k_l,1}(u_{k_j}^\ast,y^\ast_{k_j},T) + \mathrm{I}
\]
with the term ${\rm I}$ given by
\begin{align*}
        \mathrm{I}:= &(h_T^2\|{\nabla} \cdot(a(u_{k_l}^\ast){\nabla} y_{k_l}^\ast - a(u_{k_j}^\ast){\nabla} y_{k_j}^\ast) + (b(u_{k_j}^\ast)(y_{k_j}^\ast)^3-b(u_{k_l}^\ast)(y_{k_l}^\ast)^3)\|_{L^2(T)}^2  \\
        &+ \sum_{F\subset\partial T}h_T\|[(a(u_{k_l}^\ast){\nabla} y_{k_l}^\ast -a(u_{k_j}^\ast){\nabla} y_{k_j}^\ast) \cdot{n}_F]\|^2_{L^2(F)})^{1/2}.
\end{align*}
The argument in the proof of Lemma \ref{lem:stab_est_unref} in the appendix yields
\[
    \mathrm{I} \leq c(\|\overline{u}_{k_l,k_j}^\ast \|_{H^1(\omega(T))} + \|\overline{y}_{k_l,k_j}^\ast \|_{H^1(\omega(T))} + \|{\nabla}\overline{y}_{k_j}^\ast\|_{L^2(\omega(T))} +\|\overline{u}_{k_l,k_j}^\ast {\nabla} y_{\infty}^\ast\|_{L^2(\omega(T))}).
\]
Using Young's inequality with $\delta>0$, summing over all elements $T\in\cT_{k_l}\setminus\cT_{k_j}$ and the finite overlapping property of the neighborhood $\omega_T$, we further get
\begin{align*}
   \eta_{k_l,1}^2(u_{k_l}^\ast,y^\ast_{k_l},\cT_{k_l}\setminus\cT_{k_j}) \leq &(1+\delta)\eta_{k_j,1}^2(u_{k_j}^\ast,y^\ast_{k_j},\cT_{k_l}\setminus\cT_{k_j})
    +c(1+\delta^{-1}) \big( \|\overline{u}_{k_l,k_j}^\ast \|_{H^1(\Omega)}^2 + \|\overline{y}_{k_l,k_j}^\ast \|_{H^1(\Omega)}^2 \\
    &\quad + \|{\nabla}\overline{y}_{k_j}^\ast\|^2_{L^2(\Omega)}  +\|\overline{u}_{k_l,k_j}^\ast {\nabla} y_{\infty}^\ast\|_{L^2(\Omega)}^2\big).
\end{align*}
Since $\cT_{k_l}$ is a refinement of $\cT_{k_j}$, each element $T\in \cT_{k_l}\setminus\cT_{k_j}$ is generated from some $T'\in\cT_{k_j}\setminus\cT_{k_l}$ by at least one bisection. Moreover, since each component of $a(u^\ast_{k_j}){\nabla}y^{\ast}_{k_j}$ is linear on $T'$, $J_{F,1}(u^\ast_{k_j},y^\ast_{k_j})$ is zero across $F\in\mathcal{F}_{k_l}$ also pertaining to the interior of $T'$ and $h^2_T=|T|\leq \frac{1}{2} |T'| = \frac{1}{2}h^2_{T'}$. Then we arrive at
\[
\eta_{k_l,1}^2(u^\ast_{k_j},y^\ast_{k_j},\cT_{k_l}\setminus\cT_{k_j}) \leq 2^{-1/2} {u^\ast_{k_j},y^\ast_{k_j},\eta_{k_j,1}^2(\cT_{k_j}\setminus\cT_{k_l})}.
\]
Now for small enough $\delta>0$, we get the first assertion with $q_{1} = (1+\delta)
/\sqrt{2} \in (0,1)$. The other two estimates follow similarly, and hence the details are omitted.
\end{proof}

Lastly, we establish the vanishing property of the subsequences of the estimators. We borrow the techniques in \cite{GanterPraetorius:2022,morin2008,siebert2011} for direct problems.
\begin{theorem}\label{thm:conv_est}
    For the convergent subsequence $\{(u^\ast_{k_j}, y^\ast_{k_j}, p^\ast_{k_j})\}_{j\geq0}$ in Theorems \ref{thm:conv_medmin} and \ref{thm:conv_med_costate}, the three sequences of associated estimators converge to zero as $j\to \infty$, i.e.,
    \begin{equation}\label{conv_est}
      \lim_{j\to\infty}\eta_{k_j,1} =\lim_{j\to\infty}\eta_{k_j,2}
      =\lim_{j\to\infty}\eta_{k_j,3}=0.
    \end{equation}
\end{theorem}
\begin{proof}
The key ingredient of the proof is an enlarged set defined by $$\cT^+_{k_j\to k_l}:=\{T\in \cT_{k_l}\cap\cT_{k_j}~|~T\cap \Omega_{k_j}^+\neq \emptyset\}\supseteq\cT_{k_j}^+$$
for $j<l$ and an elementary yet critical disjoint splitting of $\cT_l$ as
\[
\begin{aligned}
    \cT_{k_l}&=\big(\cT_{k_l}\setminus\cT_{k_j}\big)\cup \big(\cT_{k_l}\cap\cT_{k_j}\big)
				= \big(\cT_{k_l}\setminus\cT_{k_j}\big)\cup \big(\cT_{k_l}\cap((\cT_{k_j}\setminus\cT^+_{k_j\to k_l})\cup \cT^+_{k_l \to k_l})\big)\\
				& = \big(\cT_{k_l}\setminus\cT_{k_j}\big) \cup \big(\cT_{k_l}\cap (\cT_{k_j}\setminus\cT^+_{k_j\to k_l}) \big)  \cup \big(\cT_{k_l}\cap\cT^+_{k_j\to k_l}\big) \\
       &=
\big(\cT_{k_l}\setminus\cT_{k_j}\big)\cup \big(\cT_{k_l}\cap (\cT_{k_j}\setminus\cT^+_{k_j\to k_l})\big)  \cup \cT^+_{k_j\to k_l}.
			\end{aligned}
    \]
Note that the set $\cT^+_{k_j\to k_l}$ is well defined since $\cT_{k_l}$ is a refinement of $\cT_{k_j}$ and $\cT_{k_j}^+:=\bigcap_{m\geq k_{j}}\cT_{m}\subset\cT_{k_j}\cap\cT_{k_l}$. By the definition, $\cT^+_{k_j\to k_l}$ includes $\cT_{k_j}^{+}$ and its neighboring elements in $\cT_{k_j}$ that are not refined at least until the $k_l$-th iteration, while the set $\cT_{k_l}\setminus\cT_{k_j}$ contains the elements in $\cT_{k_l}$ generated during the adaptive refinement process from the mesh level $\cT_{k_j}$ to $\cT_{k_l}$. All elements in the remaining set $\cT_{k_l}\cap (\cT_{k_j}\setminus\cT^+_{k_j\to k_l})$ will eventually be refined after the $k_l$-th iteration (see Step 2 below for more detail). Moreover, the uniform shape-regularity of the sequence $\{\cT_k\}_{k\geq0}$ given by Algorithm \ref{afem_CE} implies
    \begin{equation}\label{thm:conv_est_pf01}
        \# \cT^+_{k_j\to k_l}\leq c \# \cT_{k_j}^+,
    \end{equation}
    with the constant $c$ depending only on $\cT_0$. Thus $\eta_{k_l,i}$ ($1\leq i \leq 3$) can be rewritten as
    \begin{equation}\label{thm:conv_est_pf02}
        \eta_{k_l,i}^2 = \eta_{k_l,i}^2(\cT_{k_l}\setminus\cT_{k_j}) + \eta_{k_l,i}^2(\cT_{k_l}\cap (\cT_{k_j}\setminus\cT^+_{k_j \to k_l})) + \eta_{k_l,i}^2(\cT^+_{k_j \to k_l}),\quad i=1,2,3.
    \end{equation}
Next we treat the three terms in the right hand side of \eqref{thm:conv_est_pf02} in four steps, which comprise the rest of the proof. The overall proof strategy is as follows. (i) The reduction property in Lemma \ref{lem:reduction_est_ref} and the auxiliary convergence in Theorems \ref{thm:conv_medmin} and \ref{thm:conv_med_costate} imply that $\eta_{k_l,i}^2(\cT_{k_l}\setminus\cT_{k_j})$ decreases compared to $\eta_{k_j,i}^2(\cT_{k_j}\setminus\cT_{k_l})$ up to a Cauchy sequence in Step 1. (ii) Inspired by the proof of \cite[Proposition 4.2]{morin2008}, we resort to the local stability in Lemma \ref{lem:stab_est} and a geometric observation on $\cT_{k_l}^0\cap\cT_{k_j}^0$ (see \eqref{thm:conv_est_pf06} below) to show that for each fixed $j\in\mathbb{N}_0$, $\eta_{k_l,i}^2(\cT_{k_l}\cap (\cT_{k_j}\setminus\cT^+_{k_j \to k_l}))\to 0$ as $l\to \infty$ in Steps 2-3. (iii) Similar to \cite[Proposition 3.7]{siebert2011}, using assumption \eqref{marking} in the module MARK of Algorithm \ref{afem_CE}, Lemma \ref{lem:conv_aux5_maxest} and \eqref{thm:conv_est_pf01}, we deduce $\eta_{k_l,i}^2(\cT^+_{k_j \to k_l})\to 0$ as $l\to\infty$ with $j$ fixed in Step 4. Finally, Lemmas \ref{lem:stab_est_unref} and \ref{lem:reduction_est_ref}, and the argument in the proof of \cite[Theorem 3.1]{GanterPraetorius:2022} yield the desired zero limits. Due to the result in Step 1, this analysis differs markedly from that for direct problems \cite{morin2008}, where three contributions in an alternative decomposition of the error estimator are separately shown to vanish in the limit.

\smallskip\noindent \textit{Step 1.} Lemma \ref{lem:reduction_est_ref} implies
    \begin{align}
        \eta^2_{k_{l},1}(\cT_{k_l}\setminus \cT_{k_j}) &\leq q_{1}\eta^2_{k_j,1}(\cT_{k_j}\setminus \cT_{k_l}) + c(\|\overline{u}_{k_l,k_j}^\ast \|^2_{H^1(\Omega)} + \|\overline{y}_{k_l,k_j}^\ast\|^2_{H^1(\Omega)}+ \|{\nabla}\overline{y}_{k_j}^\ast\|^2_{L^2(\Omega)}\nonumber
        \\
        &\quad+\|\overline{u}_{k_l,k_j}^\ast {\nabla} y_{\infty}^\ast\|_{L^2(\Omega)}^2)  := q_{1}\eta^2_{k_j,1}(\cT_{k_j}\setminus \cT_{k_l}) + \epsilon_1(k_l,k_j), \label{thm:conv_est_pf03} \\
        \eta^2_{k_l,2}(\cT_{k_l}\setminus \cT_{k_j}) &\leq q_{2}\eta_{k_j,2}^2(\cT_{k_j}\setminus \cT_{k_l}) + c(\|\overline{u}_{k_l,k_j}^\ast \|^2_{H^1(\Omega)} + \|\overline{y}_{k_l,k_j}^\ast\|^2_{H^1(\Omega)} + \|\overline{p}_{k_l,k_j}^\ast\|^2_{H^1(\Omega)}  +   \|{\nabla}\overline{p}_{k_j}^\ast \|^2_{L^2(\Omega)} \nonumber\\
        &\quad + \|\overline{u}_{k_l,k_j}^\ast {\nabla} p_{\infty}^\ast\|^2_{L^2(\Omega)}) := q_{2}\eta^2_{k_j,2}(\cT_{k_j}\setminus \cT_{k_l}) + \epsilon_2(k_l,k_j), \label{thm:conv_est_pf04} \\
        \eta_{k_l,3}^2(\cT_{k_l}\setminus \cT_{k_j}) &\leq q_{3}\eta_{k_j,3}^2(\cT_{k_j}\setminus \cT_{k_l}) + c (\|\overline{u}_{k_l,k_j}^\ast\|^2_{H^1(\Omega)} + \|\overline{y}_{k_l,k_j}^\ast\|^2_{H^1(\Omega)} + \|\overline{p}_{k_l,k_j}^\ast \|^2_{H^1(\Omega)}) \nonumber \\
        &:= q_{3}\eta^2_{k_j,3}(\cT_{k_j}\setminus \cT_{k_l}) + \epsilon_3(k_l,k_j). \label{thm:conv_est_pf05}
    \end{align}
    Theorems \ref{thm:conv_medmin} and \ref{thm:conv_med_costate} and Lebesgue's dominated convergence theorem imply that $\{u^\ast_{k_j}\}_{j\geq0}$, $\{y^\ast_{k_j}\}_{j\geq0}$ and $\{p^\ast_{k_j}\}_{j\geq0}$ are all convergent in $H^1(\Omega)$ while $\{u^\ast_{k_j}\nabla y_\infty^\ast\}_{j\geq0}$ and $\{u^\ast_{k_j}\nabla p_\infty^\ast\}_{j\geq0}$ are both convergent in $L^2(\Omega)$. Hence the terms $\epsilon_i(k_l,k_j)$ ($i=1,2,3$) converge to zero respectively as $j$ and $l$ tend to $\infty$.

\medskip
\noindent\textit{Step 2.} Now we deal with $\eta_{k_l,i}(\cT_{k_l}\cap(\cT_{k_j}\setminus \cT_{k_j\to k_l}^+))$, $i=1,2,3$.  Since $\cT_{k_j}^+\subseteq\cT_{k_j\to k_l}^+$ for $j<l$, we have $\cT_{k_j}\setminus \cT_{k_j\to k_l}^+ \subseteq \cT_{k_j}\setminus \cT_{k_j}^+ = \cT_{k_j}^0$ and $\cT_{k_l}\cap (\cT_{k_j}\setminus\cT^+_{k_j \to k_l})\subseteq \cT_{k_l}\cap\cT_{k_j}^0$. Since the set $\cT_{k_l}\cap\cT_{k_j}^0$ consists of elements that will eventually be refined after the $k_j$-th iteration and also belong to $\cT_{k_l}$, we deduce that any element $T\in \cT_{k_l}\cap\cT_{k_j}^0$ is not refined until after the $k_l$-th iteration. Thus the definition of the set $\cT_{k_l}^0$ implies  \begin{equation*} \cT_{k_l}\cap\cT_{k_j}^0=\cT_{k_l}^0\cap\cT_{k_j}^0 \supseteq \cT_{k_l}\cap (\cT_{k_j}\setminus\cT^+_{k_j \to k_l}).
\end{equation*}
When $l$ is sufficiently large, due to the property $\eqref{eqn:conv_zero_mesh}$, the meshsize of each element in $\cT_{k_l}^0$ is smaller than that of each one in $\cT_{k_j}^0$ for each fixed $j\in\mathbb{N}_0$. Therefore, by letting $\Omega(\cT^0_{k_j}\cap\cT^0_{k_l}):=\bigcup_{T\in\cT^0_{k_j}\cap\cT^0_{k_l}}T$, we conclude that $\cT_{k_j}^0\cap\cT_{k_l}^0=\emptyset$ for sufficiently large $l$ and for each fixed $j\in\mathbb{N}_0$,
\begin{equation}\label{thm:conv_est_pf06}
\lim_{l\to\infty}\left|\Omega(\cT^0_{k_j}\cap\cT^0_{k_l})\right|=0.
\end{equation}
Squaring both sides of \eqref{stab_est1} in Lemma \ref{lem:stab_est}, summing the resulting estimate over the set $\cT_{k_l}\cap (\cT_{k_j}\setminus\cT^+_{k_j \to k_l})$ and  using the relation $\cT_{k_l}\cap (\cT_{k_j}\setminus\cT^+_{k_j \to k_l}) \subseteq \cT^0_{k_j}\cap\cT^0_{k_l}$ and the triangle inequality lead to
\[
    \begin{aligned}
        &\eta_{k_l,1}^2(\cT_{k_l}\cap (\cT_{k_j}\setminus\cT^+_{k_j \to k_l})) \leq c \sum_{T\in\cT_{k_l}\cap (\cT_{k_j}\setminus\cT^+_{k_j \to k_l})} \big(  \|{\nabla}y_{k_l}^\ast\|_{L^2(\omega(T))}^2 +h_T^2\|y_{k_l}^\ast\|_{L^6(T)}^6 + h_T^2\|f\|_{L^2(T)}^2 \big) \\
         \leq & c ( \|{\nabla}y_{k_l}^\ast\|_{L^2(\Omega(\cT_{k_l}^0\cap\cT_{k_j}^0))}^2 + \|y_{k_l}^\ast\|_{L^6(\Omega(\cT_{k_l}^0\cap\cT_{k_j}^0))}^6 + \|f\|^2_{L^2(\Omega(\cT_{k_l}^0\cap\cT_{k_j}^0))}) \\
        \leq & c \big( \|{\nabla}\overline{y}_{k_l}^\ast\|_{L^2(\Omega)}^2 + \|\overline{y}_{k_l}^\ast\|_{L^6(\Omega)}^6
        + \|{\nabla}y_\infty^\ast\|_{L^2(\Omega(\cT_{k_l}^0\cap\cT_{k_j}^0))}^2 + \|y_\infty^\ast\|^6_{L^6(\Omega(\cT_{k_l}^0\cap\cT_{k_j}^0))} +  \|f\|^2_{L^2(\Omega(\cT_{k_l}^0\cap\cT_{k_j}^0))} \big).
    \end{aligned}
\]
The first two terms tend to zero as $l\to\infty$ by Theorem \ref{thm:conv_medmin} and Sobolev embedding theorem. By the identity \eqref{thm:conv_est_pf06} and the absolute continuity of $\|\cdot\|_{L^2(\Omega)}$ and $\|\cdot\|_{L^6(\Omega)}$ with respect to the Lebesgue measure, the remaining three terms also vanish in the limit. Hence, for each fixed $j\in\mathbb{N}_0$, there holds
\begin{equation}\label{thm:conv_est_pf07}
   \lim_{l\to\infty} \eta_{k_l,1}^2(\cT_{k_l}\cap (\cT_{k_j}\setminus\cT^+_{k_j \to k_l}))\to 0.
\end{equation}

\noindent\textit{Step 3}. Vanishing limits of $\eta_{k_l,2}^2(\cT_{k_l}\cap (\cT_{k_j}\setminus\cT^+_{k_j \to k_l}))$ and $\eta_{k_l,3}^2(\cT_{k_l}\cap (\cT_{k_j}\setminus\cT^+_{k_j \to k_l}))$. Since $\cT_{k_l}\cap (\cT_{k_j}\setminus\cT^+_{k_j \to k_l})\subseteq \cT_{k_l}^0\cap \cT_{k_j}^0$, we have $h_T \leq \|h_{k_l}\chi^0_{k_l}\|_{L^\infty(\Omega)}$ for each $T\in \cT_{k_l}\cap (\cT_{k_j}\setminus\cT^+_{k_j \to k_l})$, over which, the inverse estimates in Lemma \ref{lem:inverse} imply
$$h_T\|y^\ast_{k_l}\|^2_{L^6(T)} \|p^\ast_{k_l}\|_{L^6(T)}\leq c \|y^\ast_{k_l}\|^2_{L^4(T)} \|p^\ast_{k_l}\|_{L^2(T)}.$$
This, \eqref{stab_est2} and the trace theorem yield
\begin{align*}
&\quad\eta_{k_l,2}^2(\cT_{k_l}\cap (\cT_{k_j}\setminus\cT^+_{k_j \to k_l})) \\
&\leq c\sum_{T\in\cT_{k_l}\cap (\cT_{k_j}\setminus\cT^+_{k_j \to k_l})}
( \| \nabla {p}^\ast_{k_l}\|^2_{L^2(\omega(T))} + \|y^\ast_{k_l}\|^4_{L^4(T)} \|p^\ast_{k_l}\|_{L^2(T)}^2 + \|h_{k_l}\chi^0_{k_l}\|_{L^\infty(\Omega)}\|y^\ast_{k_l}-y^\delta\|^2_{L^2(\p T\cap\partial\Omega)}) \\
& \leq c ( \|{\nabla} p_{k_l}^\ast\|_{L^2(\Omega(\cT_{k_l}^0\cap\cT_{k_j}^0))}^2 + \|y^\ast_{k_l}\|^4_{L^4(\Omega(\cT_{k_l}^0\cap\cT_{k_j}^0))} \|p^\ast_{k_l}\|_{L^2(\Omega(\cT_{k_l}^0\cap\cT_{k_j}^0))}^2\\
&\qquad + \|h_{k_l}\chi^0_{k_l}\|_{L^\infty(\Omega)} (\|y_{k_l}^\ast\|^2_{L^2(\partial\Omega)} + \|y^\delta\|_{L^2(\partial\Omega)}^2)) \\
& \leq c( \|{\nabla}\overline{p}_{k_l}^\ast\|_{L^2(\Omega)}^2 + \|{\nabla}p_\infty^\ast\|_{L^2(\Omega(\cT_{k_l}^0\cap\cT_{k_j}^0))}^2 + \|h_{k_l}\chi^0_{k_l}\|_{L^\infty(\Omega)}(\|y_{k_l}^\ast\|^2_{H^1(\Omega)} + \|y^\delta\|_{L^2(\partial\Omega)}^2) \\
&\qquad + ( \|\overline{y}_{k_l}^\ast \|_{L^4(\Omega)}^4 + \|y_\infty^\ast\|^4_{L^4(\Omega(\cT_{k_l}^0\cap\cT_{k_j}^0))}) ( \|\overline{p}_{k_l}^\ast\|_{L^2(\Omega)}^2 + \|p_\infty^\ast\|^2_{L^2(\Omega(\cT_{k_l}^0\cap\cT_{k_j}^0))})).
    \end{align*}
By Theorem \ref{thm:conv_medmin}, $\{\|y_{k_j}^\ast\|_{H^1(\Omega)}\}_{j\geq0}$ is uniformly bounded. Then by Theorems \ref{thm:conv_medmin} and \ref{thm:conv_med_costate}, Sobolev embedding theorem, the identity \eqref{thm:conv_est_pf06}, the absolute continuity of $\|\cdot\|_{L^2(\Omega)}$ and $\|\cdot\|_{L^4(\Omega)}$ with respect to the Lebesgue measure and \eqref{eqn:conv_zero_mesh}, all the terms tend to zero as $l\to\infty$ for each fixed $j$, which further implies
\begin{equation}\label{thm:conv_est_pf08}
    \lim_{l\to\infty}\eta_{k_l,2}^2(\cT_{k_l}\cap (\cT_{k_j}\setminus\cT^+_{k_j \to k_l}))= 0.
\end{equation}
By using the estimate \eqref{stab_est3} and the inverse estimate, we obtain
\[
    \begin{aligned}
        \eta_{k_l,3}^2(\cT_{k_l}\cap (\cT_{k_j}\setminus\cT^+_{k_j \to k_l})) & \leq c \Big( \|{\nabla} u_{k_l}^\ast\|_{L^2(\Omega(\cT_{k_l}^0\cap\cT_{k_j}^0))}^2 + \|{\nabla} y_{k_l}^\ast\|_{L^2(\Omega(\cT_{k_l}^0\cap\cT_{k_j}^0))}^2\|{\nabla} p_{k_l}^\ast\|_{L^2(\Omega(\cT_{k_l}^0\cap\cT_{k_j}^0))}^2 \\
        & \qquad + \|y^\ast_{k_l}\|^6_{L^6(\Omega(\cT_{k_l}^0\cap\cT_{k_j}^0))} \|p^\ast_{k_l}\|_{L^2(\Omega(\cT_{k_l}^0\cap\cT_{k_j}^0))}^2 + \|2 u^\ast_{k_l}-1\|_{L^2(\Omega(\cT_{k_l}^0\cap\cT_{k_j}^0))}^2\Big).
    \end{aligned}
\]
Hence for each fixed $j\in\mathbb{N}_0$, a similar argument leads to
\begin{equation}\label{thm:conv_est_pf09}
       \lim_{l\to\infty} \eta_{k_l,3}^2(\cT_{k_l}\cap (\cT_{k_j}\setminus\cT^+_{k_j \to k_l}))  =0.
    \end{equation}

\medskip
\noindent\textit{Step 4.} Vanishing limit of $\widetilde{\eta}_{k_l}^2(\cT_{k_l\to k_j}^+)$. We prove the assertion using assumption \eqref{marking} in the module MARK of Algorithm \ref{afem_CE}, which has not been invoked in Steps 1-3. Since the marking strategy in Algorithm \ref{afem_CE} is applied to $\widetilde\eta_{k_l} = \max_{1\leq i\leq 3}\eta_{k_l,i}$, then by \eqref{thm:conv_est_pf01}, \eqref{marking} and \eqref{conv_aux5_maxest}, there holds for each fixed $j\in\mathbb{N}_0$, as $l\to \infty$,
\begin{align}
     \widetilde{\eta}_{k_l}^2(\cT_{k_l\to k_j}^+) & \leq \#\cT_{k_j \to k_l}^+ \max_{ T \in \cT_{k_j \to k_l}^+ } \widetilde{\eta}_{k_l}^2(T) \leq c \#\cT_{k_j}^+ \max_{ T \in \mathcal{M}_{k_l} } \widetilde{\eta}_{k_l}^2(T) \nonumber\\
        &\leq c \#\cT_{k_j}^+ \Big(\max_{ T \in \mathcal{M}_{k_l} }  \eta_{k_l,1}^2(T) + \max_{ T \in \mathcal{M}_{k_l} } \eta_{k_l,2}^2(T) + \max_{ T \in \mathcal{M}_{k_l} } \eta_{k_l,3}^2(T)\Big)  \to 0, \label{thm:conv_est_pf10}
\end{align}
where we have omitted $y_{k_l}^\ast$, $p_{k_l}^\ast$ and $u_{k_l}^\ast$ in $\eta_{k_l,i}(T)$ ($1\leq i\leq 3$).
Since $\widetilde{\eta}_{k_j}: = \max_{1\leq i\leq 3} \eta_{k_j,i}$, the estimates \eqref{thm:conv_est_pf07}--\eqref{thm:conv_est_pf10} imply that for each $j\in \mathbb{N}_0$, there exists some $N(j)>j$ such that for each $l>N(j)$,
\begin{equation}\label{thm:conv_est_pf11}
\widetilde{\eta}_{k_{l}}^2(\cT_{k_l}\cap\cT_{k_j}) \leq 	(q' - \max_{1\leq i\leq 3}q_i) \widetilde{\eta}_{k_j}^2,
\end{equation}
with some $q' \in (\max_{1\leq i\leq 3}q_i,1)$. This and  \eqref{thm:conv_est_pf03}--\eqref{thm:conv_est_pf05} allow extracting a further subsequence $\{k_{j_n}\}_{n\geq0}$ satisfying $\widetilde{\eta}_{k_{j_{n+1}}}^2 \leq q' \widetilde{\eta}_{k_{j_{n}}}^2 + \max_{1\leq i\leq 3}\epsilon_{k_{j_{n}},i}$,
where $\epsilon_{k_{j_n},i} = \epsilon_i(k_{j_{n+1}},k_{j_{n}})\to 0$ as $n\to\infty$, cf.  Step 1. Then by \cite[Lemma 3.2]{GanterPraetorius:2022}, $\widetilde{\eta}_{k_{j_{n}}}\to0$ as $n\to\infty$, and thus
\begin{equation}\label{thm:conv_est_pf12}
                \lim_{n\to\infty} \eta_{k_{j_{n},i}} = 0,\quad 1 \leq i \leq 3.
            \end{equation}
Last, Lemmas
\ref{lem:stab_est_unref} and \ref{lem:reduction_est_ref} and the argument of Step 1 give $\eta_{k_j,i}^2\leq 2\eta_{k_{j_n},i}^2 + \xi_i(k_j,k_{j_n})$ ($1\leq i\leq 3$), where $\xi_i(k_j,k_{j_n})$ tends to zero as $k_j>k_{j_n} \to \infty$ and $n\to\infty$. This and \eqref{thm:conv_est_pf12} imply the desired assertion.
\end{proof}

\begin{remark}
In the proofs of Lemmas \ref{lem:stab_est}, \ref{lem:conv_aux5_maxest}, \ref{lem:stab_est_unref} and \ref{lem:reduction_est_ref} and Theorem \ref{thm:conv_est}, we resort to the approach developed for direct problems \cite{GanterPraetorius:2022,morin2008,nochetto2009,siebert2011}. The main changes lie in the decomposition \eqref{thm:conv_est_pf02} and the observation \eqref{thm:conv_est_pf06}. This idea has been adopted to analyze the convergence of AFEM for elliptic eigenvalue optimization \cite{LiXuZhu:2023}. Theorem \ref{thm:conv_est} motivates \textit{a posteriori} error estimators $\eta_{k,i}$ $(i=1,2,3)$ in the module ESTIMATE of Algorithm \ref{afem_CE} even if it holds only for the subsequences  $\{\eta_{k_j,i}\}_{j\geq0}$ ($i=1,2,3)$).
\end{remark}

\subsection{The proof of Theorem \ref{thm:conv}}\label{ssec:conv-overall}

Since $\eta_{k_j,i}$ ($1\leq i\leq 3$) does not provide reliable bounds for $\|u^\ast_{k_j}-u^\ast\|_{H^1(\Omega)}$, $\|y^\ast_{k_j}-y^\ast\|_{H^1(\Omega)}$ and $\|p^\ast_{k_j}-p^\ast\|_{H^1(\Omega)}$, Theorem \ref{thm:conv} is not a direct consequence of Theorem \ref{thm:conv_est} as in \cite{GanterPraetorius:2022, morin2008}. To prove the theorem, we introduce three residual functionals associated with the triplet $(u^\ast_k,
y^\ast_k, p^\ast_k)$ over $H^1(\Omega)$ and $\mathcal{A}$, defined by
\begin{align*}
    \langle \mathcal{R}_{1}(u^\ast_k,y^\ast_k) , \phi \rangle &: = (a(u^\ast_{k}){\nabla}y^\ast_{k},{\nabla}\phi) + (b(u^\ast_{k}) (y^\ast_{k})^3,\phi) - (f, \phi), \quad \forall  \phi \in H^1(\Omega),\\
    \langle \mathcal{R}_{2}(u^\ast_k,y^\ast_k, p^\ast_k) , \psi \rangle & : = (a(u^\ast_{k_j}){\nabla}p^\ast_{k},{\nabla}\psi) + 3(b(u^\ast_{k}) (y^\ast_{k})^2 p_{k}^\ast,\psi) - (y_{k}^\ast-y^\delta,\psi)_{L^2(\p\Omega)}, \quad \forall  \psi \in H^1(\Omega),\\
    \langle \mathcal{R}_{3}(u^\ast_k,y^\ast_k, p^\ast_k), v \rangle &: = ((1-\sigma) v {\nabla}y_{k}^\ast,{\nabla}p_{k}^\ast) + (v,(y_{k}^\ast)^3 p_{k}^\ast)  + 2 \alpha \varepsilon ({\nabla} u^\ast_k ,{\nabla} v) + \alpha\varepsilon^{-1} (1-2u^\ast_k, v),  \quad \forall  v \in \mathcal{A}.
\end{align*}

By bounding the residuals $\mathcal{R}_{i}$ on the convergent subsequence
$\{(u_{k_j}^\ast,y_{k_j}^\ast,p_{k_j}^\ast)\}_{j\geq 0}$, we prove the following
weak-type convergence using Theorem \ref{thm:conv_est}. The proof employs a Clement-type interpolation operator $\Pi_k: L^1(\Omega)\to V_k$ \cite{JinXu:2020}. Let $ \mathcal{N}_k$ be the set of all nodes of $\cT_k$ and $\{\phi_{x}\}_{\bold{x}\in\mathcal{N}_k} $ be the nodal basis functions in $V_k$. For each ${x}\in\mathcal{N}_k$, the support of $\phi_x$ is denoted by $\omega_x$, i.e., the union of all elements in $\cT_k$ with ${x}$ being a vertex. Then we define an interpolation operator $\Pi_k:L^1(\Om)\to V_k$ by
$\Pi_k v := \sum_{\bold{x}\in\mathcal{N}_k} \frac{1}{|\omega_x|}\int_{\omega_x}v \dx \phi_x$.
Then $\Pi_k v \in \mathcal{A}_k$ for $v \in\mathcal{A}$. Further, the following error estimates hold for any $v\in H^1(\Omega)$ \cite[Lemma 5.3]{JinXu:2020}
\begin{equation}\label{eqn:cn_int_est}
    \|v-\Pi_k v\|_{L^2(T)} \leq c h_T \|{\nabla}v\|_{L^2(\omega(T))},\quad \|v-\Pi_k v\|_{L^2(\p T)} \leq c h^{1/2}_T \|{\nabla}v\|_{L^2(\omega(T))},\quad \forall T\in \cT_k.
\end{equation}
\begin{lemma}\label{lem:conv_res}
For the convergent subsequence $\{(u_{k_j}^\ast,y_{k_j}^\ast,p_{k_j}^\ast)\}_{j\geq 0}$ given in Theorems \ref{thm:conv_medmin} and \ref{thm:conv_med_costate} of the sequence $\{(u_k^\ast,y_k^\ast,p_k^\ast)\}_{k\geq0}$ generated by Algorithm \ref{afem_CE}, there hold
\begin{align}
   \lim_{j\to\infty}\langle \mathcal{R}_{1}(u^\ast_{k_j},y^\ast_{k_j}) , \phi \rangle &= 0,\quad \forall \phi \in H^1(\Omega),\label{conv_res1}\\
   \lim_{j\to\infty}\langle \mathcal{R}_{2}(u^\ast_{k_j},y^\ast_{k_j}, p^\ast_{k_j}) , \psi \rangle& = 0,\quad\forall  \psi \in H^1(\Omega),\label{conv_res2}\\
   \liminf_{j\to\infty}\langle \mathcal{R}_{3}(u^\ast_{k_j},y^\ast_{k_j}, p^\ast_{k_j}), v -  u^\ast_{k_j} \rangle &\geq 0,\quad \forall  v \in \mathcal{A}.\label{conv_res3}
\end{align}
\end{lemma}
\begin{proof}
Since the triplet $(u^\ast_{k_j},y^\ast_{k_j},p^\ast_{k_j})$ solves the system \eqref{optsys_dis-1}--\eqref{optsys_dis-2}, the following two identities hold
    \begin{align*}
        \langle \mathcal{R}_1(u^\ast_{k_j}, y^\ast_{k_j}) , \phi \rangle &= \langle \mathcal{R}_1(u^\ast_{k_j}, y^\ast_{k_j}) , \phi - \phi_{k_j}\rangle,\quad  \forall\phi_{k_j}\in V_{k_j},\\
        \langle \mathcal{R}_2(u^\ast_{k_j}, y^\ast_{k_j},p^\ast_{k_j}) , \phi \rangle &= \langle \mathcal{R}_2(u^\ast_{k_j}, y^\ast_{k_j},p^\ast_{k_j}) , \psi - \psi_{k_j} \rangle,\quad \forall \psi_{k_j}\in V_{k_j}.
    \end{align*}
Upon using the operator $\Pi_{k_j}$ and elementwise integration by parts, we obtain for any $\phi,\psi\in H^1(\Omega)$,
\begin{align*}
        &|\langle \mathcal{R}_1(u^\ast_{k_j}, y^\ast_{k_j}) , \phi - \Pi_{k_j} \phi \rangle|\\
         =& \Big| \sum_{T\in\cT_{k_j}}\int_{T} - R_{T,1}(u^\ast_{k_j},y^\ast_{k_j}) (\phi - \Pi_{k_j}\phi) \dx + \sum_{F\in\mathcal{F}_{k_j}}\int_F J_{F,1}(u^\ast_{k_j},y^\ast_{k_j}) (\phi - \Pi_{k_j}\phi)  \ds    \Big| \\
        \leq& c \sum_{T\in\cT_{k_j}}\eta_{{k_j},1}(u_{k_j}^\ast,y_{k_j}^\ast,T) \|{\nabla}\phi\|_{L^2(\omega(T))}\leq c\eta_{{k_j},1}(u_{k_j}^*,y_{k_j}^*) \|\phi\|_{H^1(\Omega)},\\
        &|\langle \mathcal{R}_2(u^\ast_{k_j}, y^\ast_{k_j},p^\ast_{k_j}) , \psi - \Pi_{k_j} \psi \rangle|\\
         =& \Big| \sum_{T\in\cT_{k_j}}\int_{T} - R_{T,2}(u^\ast_{k_j},y^\ast_{k_j},p^\ast_{k_j}) (\psi - I_{k_j}\psi) \dx + \sum_{F\in\mathcal{F}_{k_j}}\int_F J_{F,2}(u^\ast_{k_j},y^\ast_{k_j},y^\ast_{k_j}) (\psi - \Pi_{k_j}\psi)  \ds    \Big| \\
        \leq &c \sum_{T\in\cT_{k_j}}\eta_{{k_j},2}(u_{k_j}^\ast,y_{k_j}^\ast,p_{k_j}^\ast,T) \|{\nabla}\psi\|_{L^2(\omega(T))}\leq c \eta_{{k_j},2}(u_{k_j}^*,y_{k_j}^*,p_{k_j}^*)\|\psi\|_{H^1(\Omega)}.
        \end{align*}
Therefore, Theorem \ref{thm:conv_est} implies the identities \eqref{conv_res1} and \eqref{conv_res2}. To prove the inequality \eqref{conv_res3}, from  \eqref{optsys_dis-3}, we deduce that for any $ v \in \mathcal{A}$,
\begin{align*}
        \langle \mathcal{R}_{3}(u^\ast_{k_j},y^\ast_{k_j}, p^\ast_{k_j}), v -  u^\ast_{k_j} \rangle  &= \langle \mathcal{R}_{3}(u^\ast_{k_j},y^\ast_{k_j}, p^\ast_{k_j}), v - \Pi_{k_j} v \rangle + \langle \mathcal{R}_{3}(u^\ast_{k_j},y^\ast_{k_j}, p^\ast_{k_j}), \Pi_{k_j}v -  u^\ast_{k_j}\rangle\\
        &\geq
        \langle \mathcal{R}_{3}(u^\ast_{k_j},y^\ast_{k_j}, p^\ast_{k_j}), v - \Pi_{k_j}v  \rangle.
    \end{align*}
Elementwise integration by parts and  \eqref{eqn:cn_int_est} give that for any $v\in\mathcal{A}$,
\begin{align*}
        | \langle \mathcal{R}_{3}(u^\ast_{k_j},y^\ast_{k_j}, p^\ast_{k_j}), v - \Pi_{k_j}v  \rangle | & = \Big| \sum_{T\in\cT_{k_j}}\int_T R_{T,3}(u_{k_j}^\ast,y^\ast_{k_j},p^\ast_{k_j}) ( v - \Pi_{k_j} v ) \dx + \sum_{F\in\mathcal{F}_{k_j}} \int_F J_{F,3}(u^\ast_{k_j})\ds \Big| \\
        & \leq c \sum_{T\in\cT_{k_j}}\eta_{{k_j},3}(u_{k_j}^\ast,y^\ast_{k_j},p^\ast_{k_j},T) \|{\nabla}v\|_{L^2(\omega(T))}\\ &\leq c \eta_{{k_j},3}(u_{k_j}^\ast,y^\ast_{k_j},p^\ast_{k_j})\|v\|_{H^1(\Omega)}.
    \end{align*}
Passing to the limit $j\to\infty$ and appealing to Theorem \ref{thm:conv_est} again conclude the proof.
\end{proof}

The next result gives the residual convergence for the convergent subsequence.
\begin{lemma}\label{lem:conv_aux3_vp}
For the convergent subsequence $\{(u_{k_j}^\ast,y_{k_j}^\ast,p_{k_j}^\ast)\}_{j\geq 0}$ in Theorems \ref{thm:conv_medmin} and \ref{thm:conv_med_costate} and the solution triplet $(u_\infty^\ast,y_\infty^\ast,p_\infty^\ast)$ of
problem \eqref{medmin}--\eqref{medvp_state} and \eqref{medvp_costate}, there hold as $j\to\infty$,
\begin{align}
   & (a(u^\ast_{k_j}){\nabla}y^\ast_{k_j},{\nabla}\phi) + (b(u^\ast_{k_j}) (y^\ast_{k_j})^3,\phi)
         \to( a(u^\ast_{\infty}){\nabla}y^\ast_{\infty},{\nabla}\phi) + (b(u^\ast_{\infty}) (y^\ast_{\infty})^3,\phi), \quad \forall\phi \in H^1(\Omega),\label{conv_aux3_vp1}\\
   &     (a(u^\ast_{k_j}){\nabla}p^\ast_{k_j},{\nabla}\psi) + 3(b(u^\ast_{k_j}) (y^\ast_{k_j})^2 p_{k_j}^\ast,\psi)\nonumber\\
    &\qquad\qquad         \to(a(u^\ast_{\infty}){\nabla}y^\ast_{\infty},{\nabla}\psi) + 3(b(u^\ast_{\infty}) (y^\ast_{\infty})^2 p_\infty^\ast,\psi), \quad \forall\psi \in H^1(\Omega),\label{conv_aux3_vp2}\\
   &     ((1-\sigma) (v - u_{k_j}^\ast) {\nabla}y_{k_j}^\ast,{\nabla}p_{k_j}^\ast) + (v - u_{k_j}^\ast, (y_{k_j}^\ast)^3 p_{k_j}^\ast)\nonumber\\
     &\qquad\qquad    \to ((1-\sigma) (v - u_{\infty}^\ast){\nabla}y_{\infty}^\ast,{\nabla}p_{\infty}^\ast )+ (v - u_{\infty}^\ast, (y_{\infty}^\ast)^3 p_{\infty}^\ast), \quad \forall v \in \A.\label{conv_aux3_vp3}
\end{align}
\end{lemma}
\begin{proof}
For any $\phi\in H^1(\Omega)$, direct computation leads to
    \begin{align*}
        |(a(u^\ast_{k_j}){\nabla}y^\ast_{k_j}-a(u^\ast_\infty){\nabla}y^\ast_\infty,{\nabla} \phi) |
         & = |((a(u^\ast_{k_j})-a(u_\infty^\ast)){\nabla}y^\ast_{k_j}, {\nabla} \phi) +( a(u^\ast_\infty){\nabla}\overline{y}_{k_j}^\ast,{\nabla} \phi)| \\
        & \leq \|\overline{u}_{k_j}^\ast{\nabla}\phi\|_{L^2(\Omega)}\|{\nabla}y_{k_j}^\ast\|_{L^2(\Omega)}
        +\|{\nabla}\overline{y}_{k_j}^\ast\|_{L^2(\Omega)}\|\nabla\phi\|_{L^2(\Omega)}.
    \end{align*}
By Theorem \ref{thm:conv_medmin}, the sequence $\{y_{k_j}^\ast\}_{j\geq0}$ is uniformly bounded in $H^1(\Omega)$.
Then by Lebesgue's dominated convergence theorem \cite[Theorem 1.19, p. 28]{Evans:2015} and Theorem \ref{thm:conv_medmin}, we obtain
\begin{equation}\label{pf:aux_conv_vp01}
       \lim_{j\to\infty}(a(u^\ast_{k_j}){\nabla}y^\ast_{k_j}-a(u^\ast_\infty){\nabla}y^\ast_\infty,{\nabla} \phi) = 0, \quad \forall\phi\in H^1(\Omega).
\end{equation}
Likewise, for any $\phi\in H^1(\Omega)$, we have the splitting
\[
(b(u_{k_j}^\ast)(y_{k_j}^\ast)^3 - b(u_{\infty}^\ast)(y_{\infty}^\ast)^3, \phi) =  ((b(u_{k_j}^\ast) - b(u_{\infty}^\ast))(y_{k_j}^\ast)^3 ,\phi) + (b(u_{\infty}^\ast)((y_{k_j}^\ast)^3-(y_\infty^\ast)^3), \phi),
    \]
    which upon repeating the preceding argument yields
    \[
      \lim_{j\to\infty}  |( (b(u_{k_j}^\ast) - b(u_{\infty}^\ast))(y_{k_j}^\ast)^3, \phi )|
        \leq \lim_{j\to\infty}\|\overline{u}_{k_j}^\ast\phi\|_{L^2(\Omega)}\|y_{k_j}^\ast\|^3_{L^6(\Omega)} = 0, \quad \forall\phi\in H^1(\Omega).
    \]
    Meanwhile, we deduce
\begin{align*}
    |( b(u_{\infty}^\ast)((y_{k_j}^\ast)^3-(y_\infty^\ast)^3), \phi)| &
    = |(\overline{y}_{k_j}^\ast((y^\ast_{k_j})^2+y_{k_j}^\ast y_{k_l}^\ast+(y^\ast_{\infty})^2),\phi)|
    \leq c\|\overline{y}_{k_j}^\ast\|_{L^4(\Omega)}\|\phi\|_{L^4(\Omega)}.
\end{align*}
By Sobolev embedding theorem and the argument for \eqref{pf:aux_conv_vp01}, we get
$$\lim_{j\to\infty}|(b(u_{\infty}^\ast)((y_{k_j}^\ast)^3-(y_\infty^\ast)^3), \phi)| = 0,\quad \forall\phi\in H^1(\Omega).$$
The preceding two results together yield
\begin{equation}\label{pf:aux_conv_vp02}
    \lim_{j\to\infty}(b(u_{k_j}^\ast)(y_{k_j}^\ast)^3 - b(u_{\infty}^\ast)(y_{\infty}^\ast)^3,\phi)= 0,\quad \forall\phi \in H^1(\Omega).
\end{equation}
The estimates \eqref{pf:aux_conv_vp01} and \eqref{pf:aux_conv_vp02} together imply \eqref{conv_aux3_vp1}.
Next we prove the estimate \eqref{conv_aux3_vp2}. The argument for  \eqref{pf:aux_conv_vp01} and Theorem \ref{thm:conv_med_costate} lead to
\begin{equation}\label{pf:aux_conv_vp03}
        \lim_{j\to\infty}(a(u^\ast_{k_j}){\nabla}p^\ast_{k_j},{\nabla}\psi) = (a(u^\ast_{\infty}){\nabla}p^\ast_{\infty},{\nabla}\psi), \quad \forall\psi\in H^1(\Omega).
    \end{equation}
Then we split the nonlinear term into
\begin{align*}
         (b(u^\ast_{k_j}) &(y^\ast_{k_j})^2 p_{k_j}^\ast  - b(u^\ast_{\infty}) (y^\ast_{\infty})^2 p_{\infty}^\ast,\psi) = (b(u^\ast_{k_j}) ((y_{k_j}^\ast)^2 - (y_{\infty}^\ast)^2)  p_{k_j}^\ast,\psi) \\
        & + ((b(u_{k_j}^\ast)-b(u^\ast_\infty)) (y_\infty^\ast)^2 p_{k_j}^\ast,\psi)
         + (b(u^\ast_\infty)(p_{k_j}^\ast-p_{\infty}^\ast) (y_\infty^\ast)^2,  \psi), \quad \forall\psi\in H^1(\Omega).
        \end{align*}
By Theorems \ref{thm:conv_medmin} and \ref{thm:conv_med_costate} and Sobolev embedding theorem, we deduce that the subsequences $\{y_{k_j}^\ast\}_{j\geq 0}$ and $\{p_{k_j}^\ast\}_{j\geq 0}$  are both uniformly bounded in $H^1(\Omega)$ and as $j \to \infty$,
\begin{align*}
        |(b(u^\ast_{k_j})((y_{k_j}^\ast)^2 - (y_{\infty}^\ast)^2)p_{k_j}^\ast, \psi)| & \leq c \|\overline{y}_{k_j}^\ast \|_{H^1(\Omega)} \|\psi\|_{H^1(\Omega)} \to 0,\\
        | ((b(u_{k_j}^\ast)-b(u^\ast_\infty)) (y_\infty^\ast)^2 p_{k_j}^\ast, \psi )| & \leq
        c\|\overline{u}_{k_j}^\ast(y_\infty^\ast)^2\|_{L^2(\Omega)} \|\psi\|_{H^1(\Omega)}\to 0,\\
        |(b(u^\ast_\infty)(p_{k_j}^\ast-p_{\infty}^\ast) (y_\infty^\ast)^2,\psi)| &\leq c \|\overline{p}_{k_j}^\ast\|_{H^1(\Omega)} \|y_\infty^\ast\|^2_{H^1(\Omega)} \|\psi\|_{H^1(\Omega)}\to0.
\end{align*}
Consequently,
\begin{equation}\label{pf:aux_conv_vp04}
    \lim_{j\to\infty}(b(u^\ast_{k_j}) (y^\ast_{k_j})^2 p_{k_j}^\ast,\psi) = (b(u^\ast_{\infty}) (y^\ast_{\infty})^2 p_{\infty}^\ast, \psi), \quad \forall\psi\in H^1(\Omega).
\end{equation}
Now the identity \eqref{conv_aux3_vp2} is direct from \eqref{pf:aux_conv_vp03} and \eqref{pf:aux_conv_vp04}. Last, for the assertion \eqref{conv_aux3_vp3}, we only show the convergence for the second term, since the first one can be proved similarly (also cf. \cite{JinXuZou:2016}). By the pointwise  convergence of $\{u_{k_j}^\ast\}_{j\geq 0}$ in Theorem \ref{thm:conv_medmin}, Lebesgue's dominated convergence theorem, Sobolev embedding theorem and the $H^1(\Omega)$ convergence of $\{p_{k_j}^\ast\}_{j\geq0}$ in Theorem \ref{thm:conv_med_costate}, as $j \to \infty$, we have
\begin{align*}
    |(u_{k_j}^\ast p_{k_j}^\ast,(y^\ast_{k_j})^3-(y_\infty^\ast)^3) |
    & \leq c\|p^\ast_{k_j}\|_{L^4(\Omega)}\|\overline{y}^\ast_{k_j}\|_{L^4(\Omega)}(\|y^\ast_{k_j}\|^2_{L^4(\Omega)} + \|y^\ast_{\infty}\|^2_{L^4(\Omega)})
     \leq c \|\overline{y}^\ast_{k_j}\|_{H^1(\Omega)}\to0,\\
    |((u_{k_j}^\ast-u_\infty^\ast) (y^\ast_{\infty})^3, p_{k_j}^\ast )| &\leq \|\overline{u}^\ast_{k_j}(y^\ast_\infty)^3\|_{L^2(\Omega)}\|p_{k_j}^\ast\|_{H^1 (\Omega)}\to0,\\
        |(u_{\infty}^\ast(y_{\infty}^\ast)^3,p_{k_j}^\ast - p_\infty^\ast)| &\leq \|y_\infty^\ast\|_{L^4(\Omega)}^3 \|\overline{p}_{k_j}^\ast\|_{L^4(\Omega)} \to 0.
\end{align*}
These three vanishing limits yield
$$(u_{k_j}^\ast (y^\ast_{k_j})^3, p^\ast_{k_j}) \to (u_{\infty}^\ast (y^\ast_{\infty})^3, p^\ast_{\infty})\quad \mbox{as } j\to\infty.$$
A similar argument gives
|$$  \lim_{j\to\infty}(v, (y_{k_j}^\ast)^3 p_{k_j}^\ast) = (v, (y_{\infty}^\ast)^3 p_{\infty}^\ast),\quad \forall v \in \A.$$
Hence, the last assertion follows.
\end{proof}

The next result shows that triplet $(u^\ast_\infty,y^\ast_\infty,p_\infty^*)$ verifies the
necessary optimality system \eqref{optsys}. Finally, the desired convergence in Theorem \ref{thm:conv} is direct from Theorems \ref{thm:conv_medmin}, \ref{thm:conv_med_costate} and Lemma \ref{lem:conv_aux4_vp}.
\begin{lemma}\label{lem:conv_aux4_vp}
The minimizing pair $(u^\ast_\infty,y^\ast_\infty)$ to problem \eqref{medmin}--\eqref{medvp_state} and
the associated adjoint $p_\infty^\ast$ of problem \eqref{medvp_costate} also solve
\begin{align}
         &(a(u_\infty^\ast){\nabla} y^\ast_\infty ,{\nabla} \phi)
    + (b(u_\infty^\ast) (y^\ast_\infty)^3, \phi) = (f, \phi), \quad \forall \phi \in H^1(\Omega),\label{conv_aux4_vp_state}\\
      &(a(u^\ast_\infty){\nabla} p_\infty^\ast,{\nabla} \psi) + 3(b(u_\infty^\ast) (y_\infty^\ast)^2 p_\infty^\ast, \psi) =     (y_\infty^\ast-y^\delta,\psi)_{L^2(\partial\Omega)},\quad \forall \psi \in H^1(\Omega),\label{conv_aux4_vp_costate}\\
   & ((1-\sigma) (v-u_\infty^\ast) {\nabla} y_\infty^\ast  ,{\nabla} p_\infty^\ast) + (v-u_\infty^\ast,(y_\infty^\ast)^3 p_\infty^\ast) \nonumber\\
        & \quad + 2 \alpha \varepsilon ({\nabla} u_\infty^\ast,{\nabla} (v - u_\infty^\ast)) + \alpha\varepsilon^{-1} (1-2u^\ast_\infty,v - u_\infty^\ast) \geq 0, \quad \forall v \in \A.\label{conv_aux4_vp_control}
\end{align}
\end{lemma}
\begin{proof}
The assertion \eqref{conv_aux4_vp_state} is direct from \eqref{conv_res1} in Lemma \ref{lem:conv_res}
and \eqref{conv_aux3_vp1} in Lemma \ref{lem:conv_aux3_vp}. Next, the trace theorem and the $H^1(\Omega)$
convergence of the subsequence $\{y^\ast_{k_j}\}_{j\geq0}$ (cf. Theorem \ref{thm:conv_medmin}) imply
\[
\lim_{j\to\infty}  (y_{k_j}^\ast-y^\delta, \psi)_{L^2(\partial\Omega)} = (y_\infty^\ast-y^\delta, \psi)_{L^2(\partial\Omega)},\quad \forall \psi\in H^1(\Omega),
\]
which, together with \eqref{conv_res1} in Lemma \ref{lem:conv_res} and \eqref{conv_aux3_vp2} in Lemma
\ref{lem:conv_aux3_vp}, yields the assertion \eqref{conv_aux4_vp_costate}. Last, by
the $H^1(\Omega)$ convergence of the subsequence $\{u^\ast_{k_j}\}_{j\geq0}$ in  \eqref{conv_medmin}, we have
for any $v\in\A$,
\begin{align*}
  \lim_{j\to\infty}  ({\nabla} u_{k_j}^\ast ,{\nabla} (v - u_{k_j}^\ast))  =  ({\nabla} u_\infty^\ast, {\nabla} (v - u_\infty^\ast)) \quad \mbox{and}\quad  \lim_{j\to \infty}(1-2u^\ast_{k_j},v - u_{k_j}^\ast)
        =   (1-2u^\ast_\infty,v - u_\infty^\ast).
\end{align*}
This, \eqref{conv_res3} in Lemma \ref{lem:conv_res} and \eqref{conv_aux3_vp3} in Lemma \ref{lem:conv_aux3_vp} conclude the inequality \eqref{conv_aux4_vp_control}.
\end{proof}

\section{Numerical experiments and discussions}\label{sec:numer}
Now we provide numerical experiments to show the convergence of Algorithm \ref{afem_CE}. Throughout we take a square domain $\Omega=(-1,1)^2$, and the initial mesh $\mathcal{T}_0$ to be a uniform triangulation with a $21 \times 21$ grid. We obtain the exact data $y^*$ by solving problem \eqref{diffequ} on a fine mesh with a $401 \times 401$ grid. We generate the noisy data $y^\delta$ by
\begin{equation*}
y^\delta (x) = y^*(x) + 1\% \|y^{*}\|_{L^{\infty}(\partial \Omega)} \xi(x),
\end{equation*}
where $\xi(x)$ is a random variable drawn from the standard normal distribution.
Numerically we investigate several scenarios involving inclusions of different shapes, including circles, ellipses, and multiple circles, and employ the experimental setup in \cite{beretta2018a}. For the case of one inclusion, we take two fixed sources for the numerical recovery, i.e., $f_1(x_1,x_2) = x_1$ and $f_2(x_1,x_2) = x_2$, and for the more complex cases involving multiple circular inclusions, we also use a third source $f_3(x_1,x_2) = \frac{x_1+x_2}{2}$.
(Note however that the choice does not satisfy the positivity condition in Assumption \ref{ass02}.)
Throughout, the conductivity $\sigma$ of the inclusion $\omega$ is fixed at $\sigma = \text{1e-4}$. The other parameters are shown in Table \ref{tab:para}, where $\varepsilon$ is the phase field parameter and $\alpha$ is the regularization parameter.
To solve the constrained optimization problem \eqref{min_G-L_dis} on the initial mesh $\mathcal{T}_0$, we employ the algorithm in \cite{beretta2018a} (i.e., parabolic obstacle method) with the initial guess $u=0$ and the BFGS \cite{ByrdLu:1995} (implemented by the IPOPT project) for the subsequent meshes $\{\mathcal{T}_{k}\}_{k>0}$ (with the $\mathcal{T}_k$-interpolant of results over the previous mesh as the initial guess), in order to facilitate fast convergence of the optimization process. In each MARK step, based on the computed error indicator $\widetilde{\eta}_k$, we fix $\theta = 0.65$ in the practical D\"{o}rfler's strategy.

\begin{table}[hbt!]
\centering
\caption{Experimental parameters for the examples.}
\label{tab:para}
\begin{tabular}{c|ccc}
\hline
Configuration&$\alpha$&$\varepsilon$\\
\hline
Circle&$\text{1.5e-3}$&$1/16\pi$\\
Ellipse&$\text{1.5e-3}$&$1/16\pi$\\
Two Circles&$\text{1e-3}$&$1/16\pi$\\
Four Circles&$\text{2e-3}$&$1/12\pi$\\
\hline
\end{tabular}
\end{table}

\newcommand{\fht}{.16\textwidth}
\begin{figure}[hbt!]
\centering
\setlength{\tabcolsep}{0pt}
\begin{tabular}{cccccc}
\includegraphics[height = \fht, trim = {3.5cm 0cm 3.5cm 0cm}, clip]{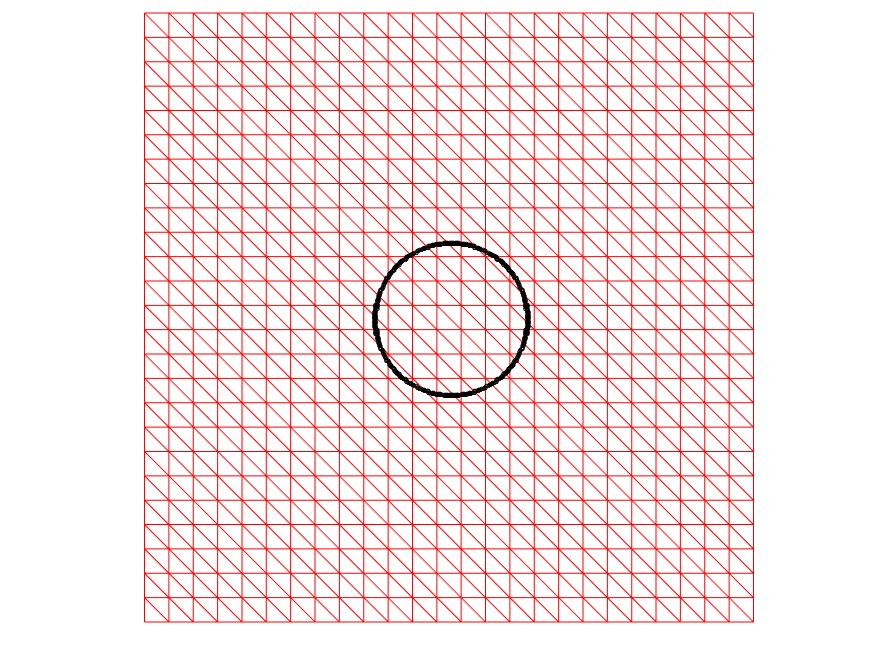}
&\includegraphics[height = \fht, trim = {3.5cm 0cm 3.5cm 0cm}, clip]{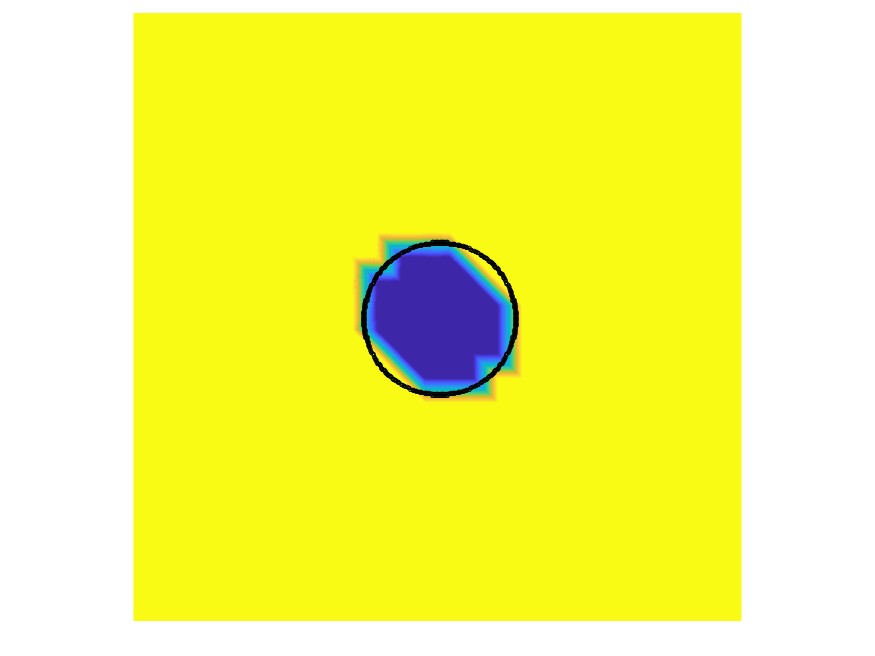}
&\includegraphics[height = \fht, trim = {3.5cm 0cm 3.5cm 0cm}, clip]{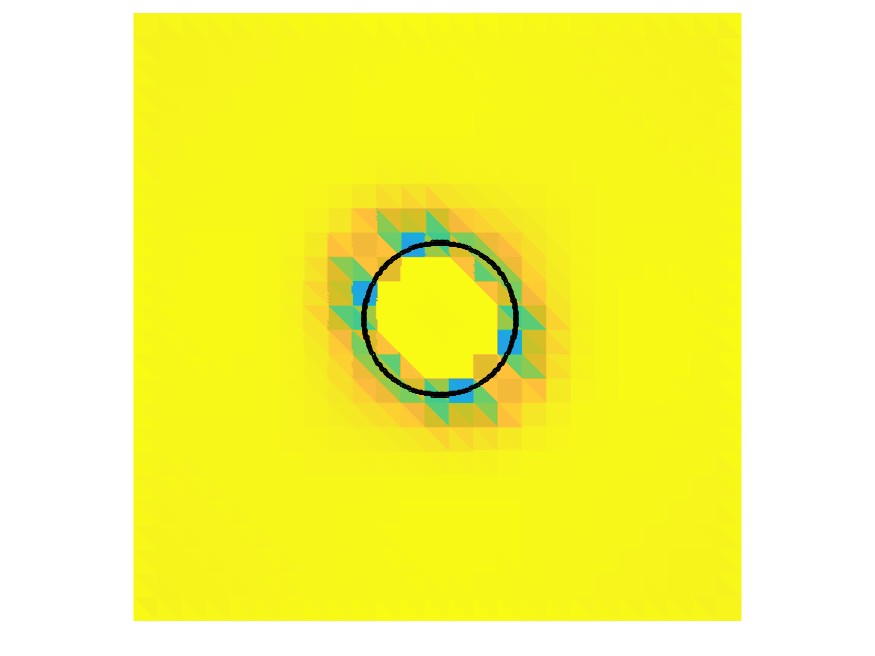}
&\includegraphics[height = \fht, trim = {3.5cm 0cm 3.5cm 0cm}, clip]{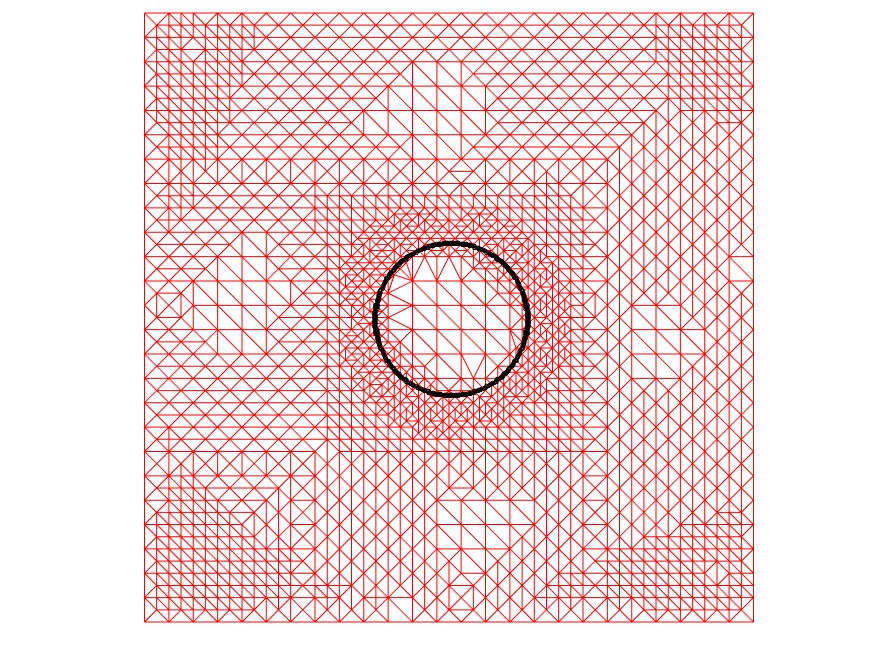}
&\includegraphics[height = \fht, trim = {3.5cm 0cm 3.5cm 0cm}, clip]{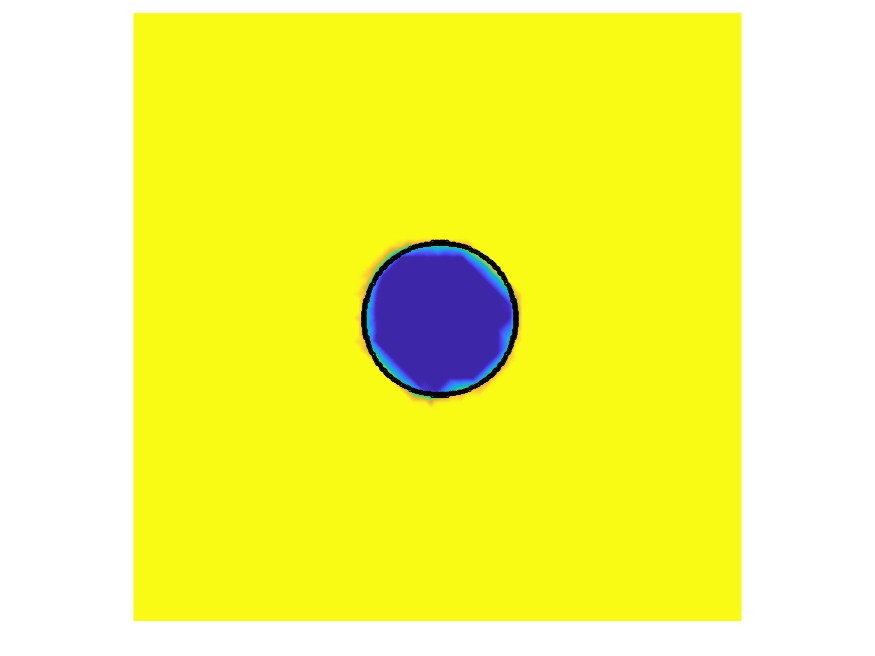}
&\includegraphics[height = \fht, trim = {3cm 0cm 3cm 0cm}, clip]{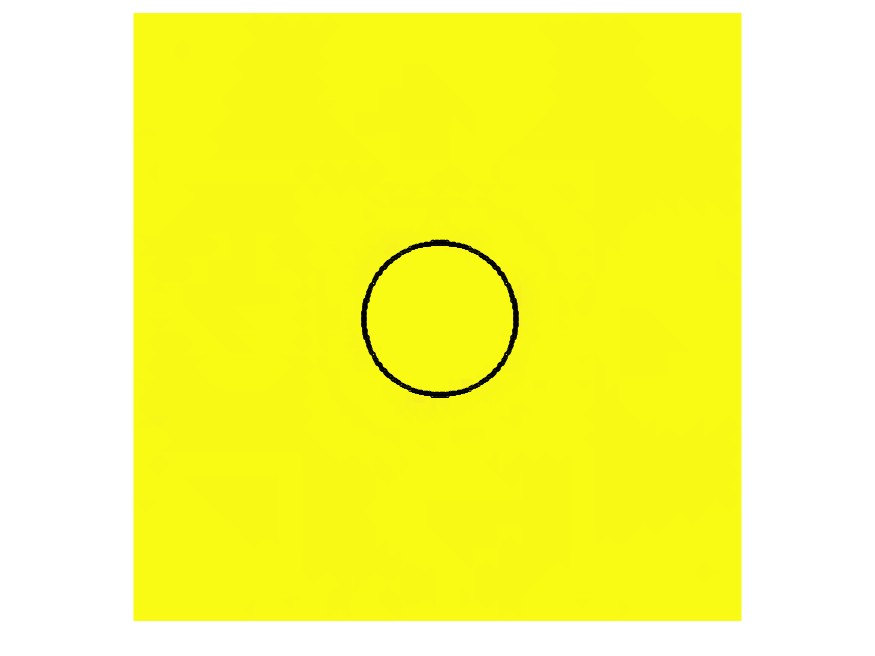}\\
\includegraphics[height = \fht, trim = {3.5cm 0cm 3.5cm 0cm}, clip]{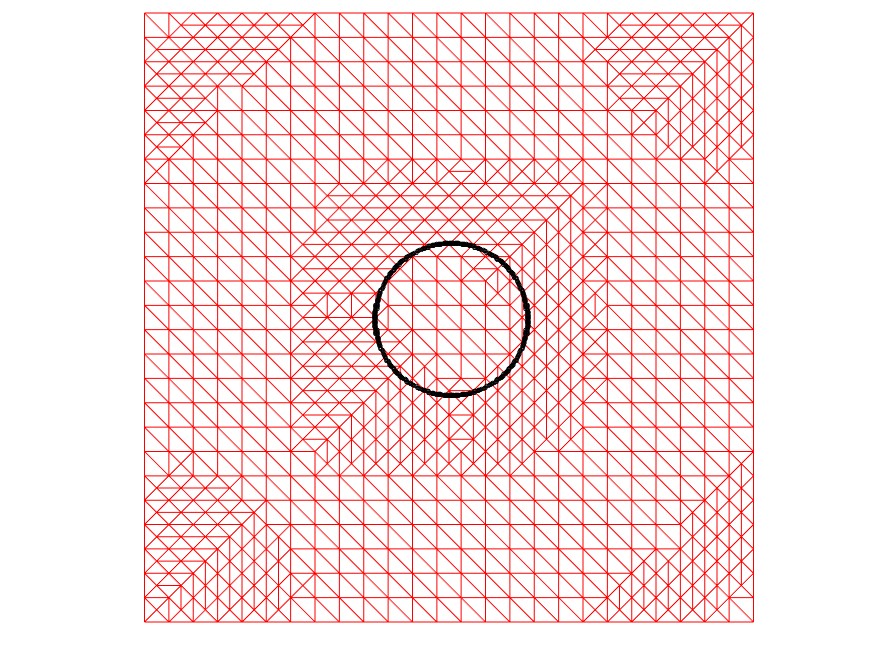}
&\includegraphics[height = \fht, trim = {3.5cm 0cm 3.5cm 0cm}, clip]{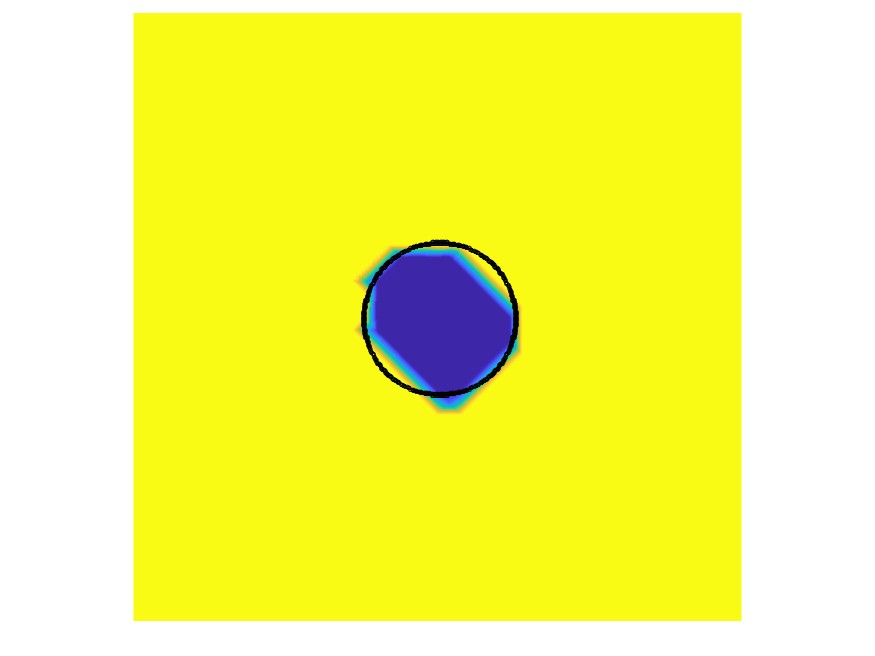}
&\includegraphics[height = \fht, trim = {3.5cm 0cm 3.5cm 0cm}, clip]{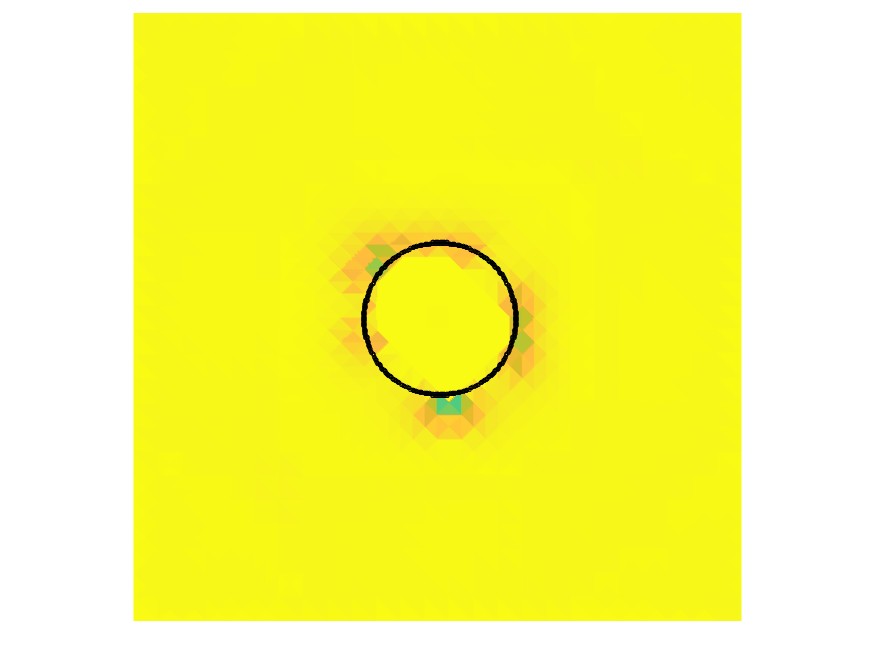}
&\includegraphics[height = \fht, trim = {3.5cm 0cm 3.5cm 0cm}, clip]{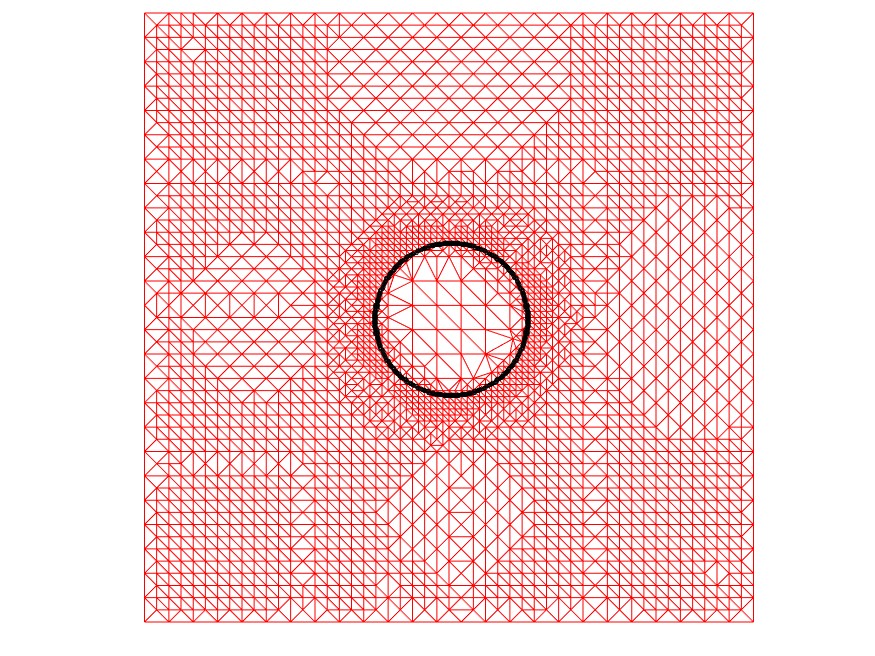}
&\includegraphics[height = \fht, trim = {3.5cm 0cm 3.5cm 0cm}, clip]{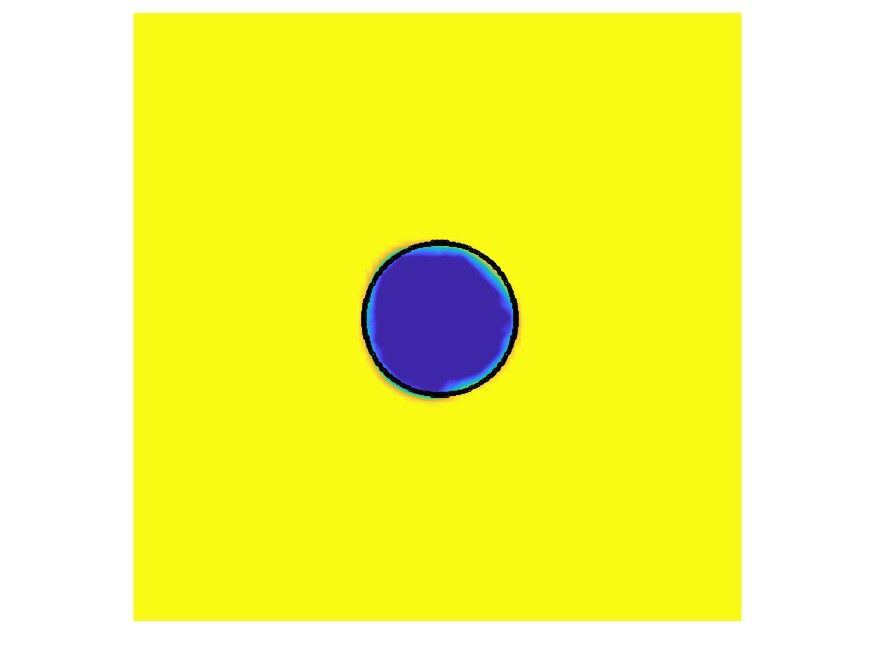}
&\includegraphics[height = \fht, trim = {3.5cm 0cm 3.5cm 0cm}, clip]{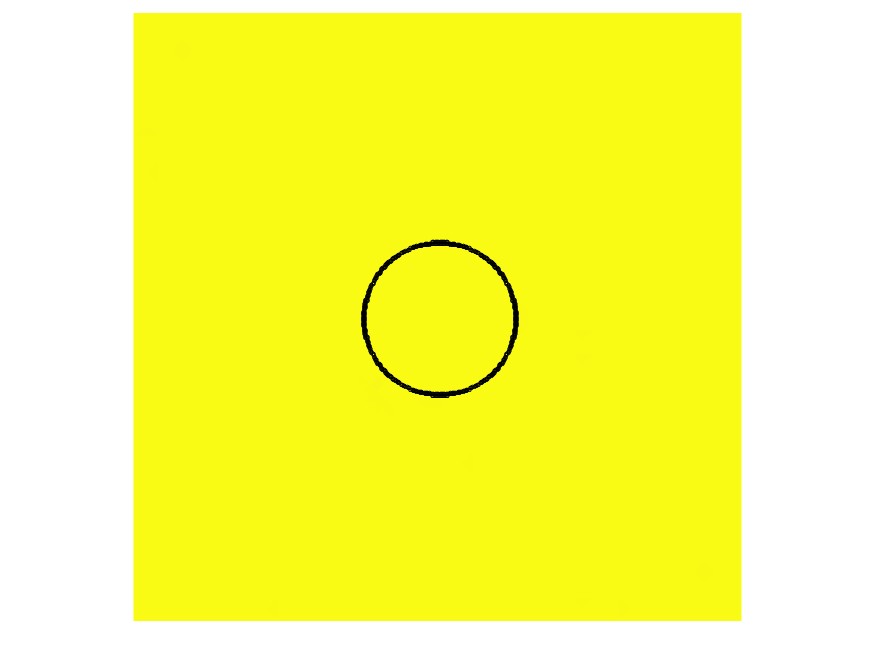}\\
\includegraphics[height = \fht, trim = {3.5cm 0cm 3.5cm 0cm}, clip]{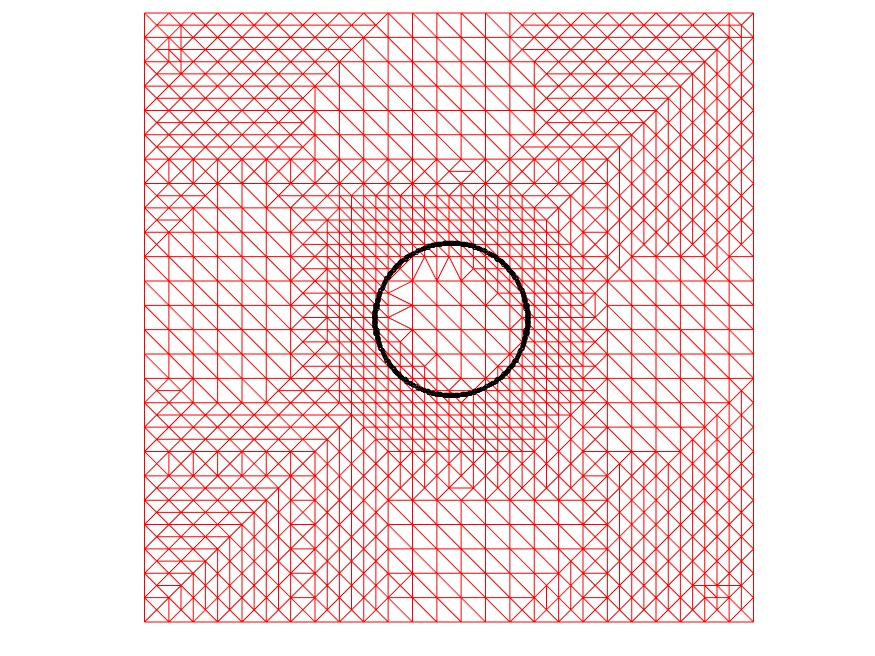}
&\includegraphics[height = \fht, trim = {3.5cm 0cm 3.5cm 0cm}, clip]{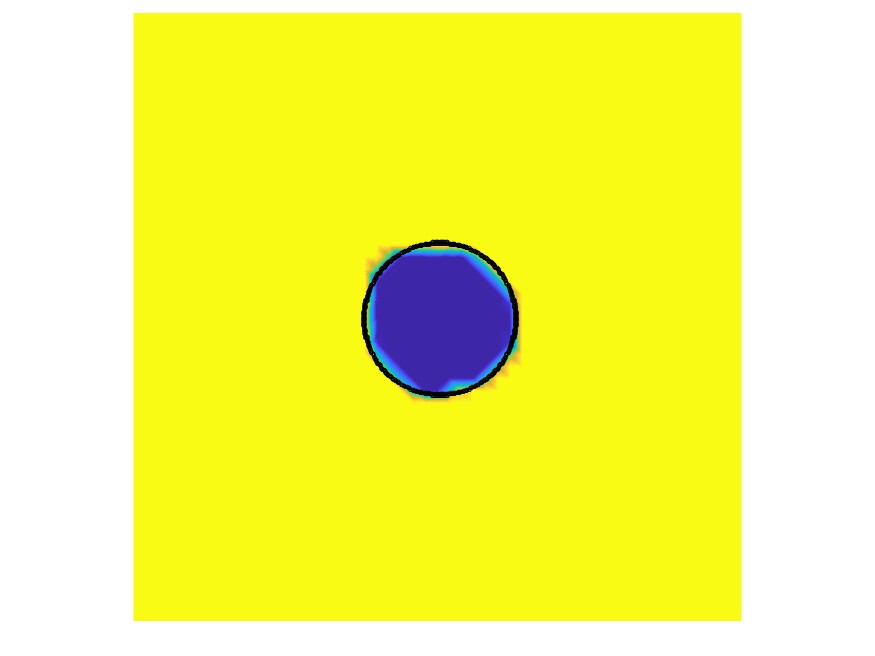}
&\includegraphics[height = \fht, trim = {3.5cm 0cm 3.5cm 0cm}, clip]{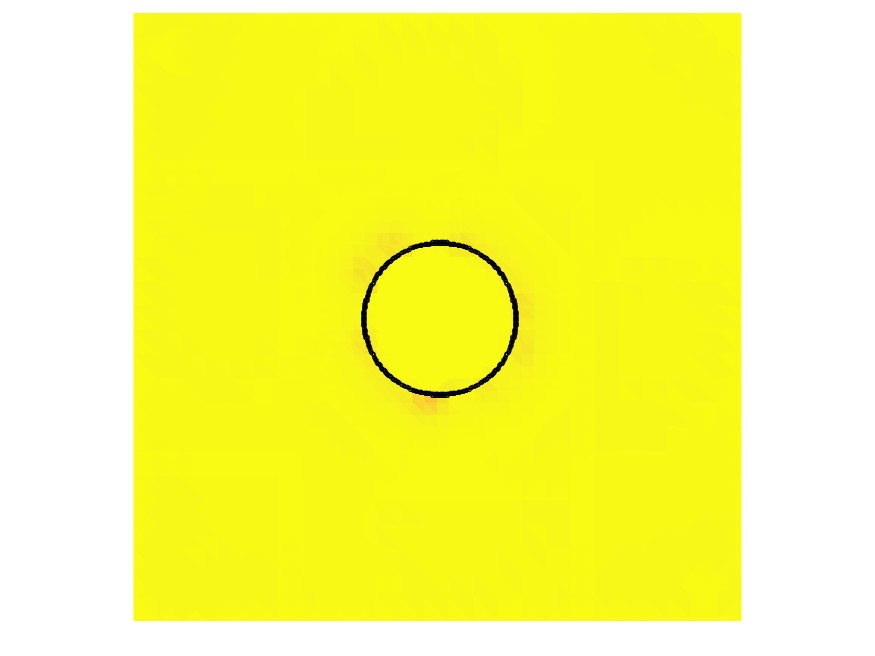}
&\includegraphics[height = \fht, trim = {3.5cm 0cm 3.5cm 0cm}, clip]{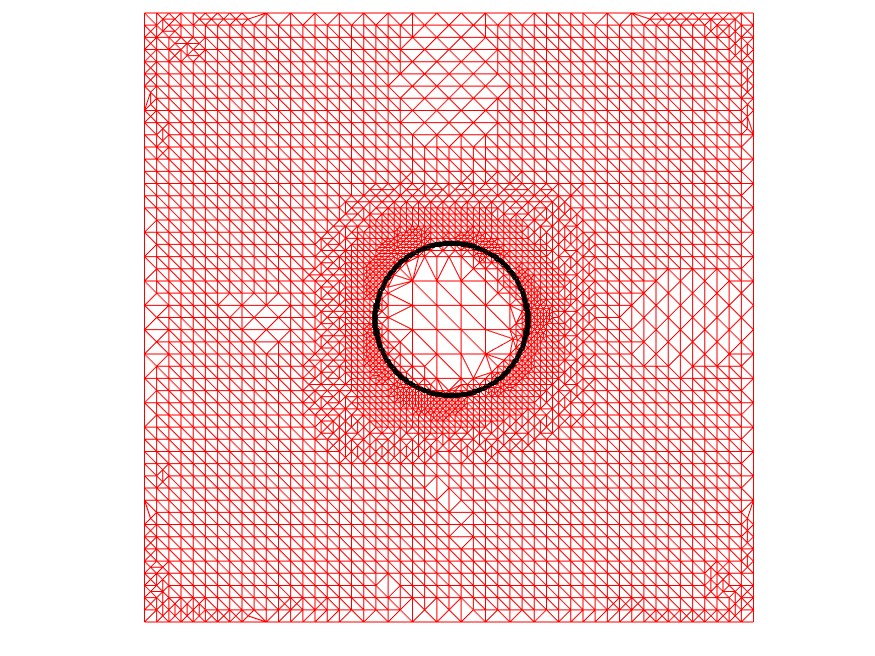}
&\includegraphics[height = \fht, trim = {3.5cm 0cm 3.5cm 0cm}, clip]{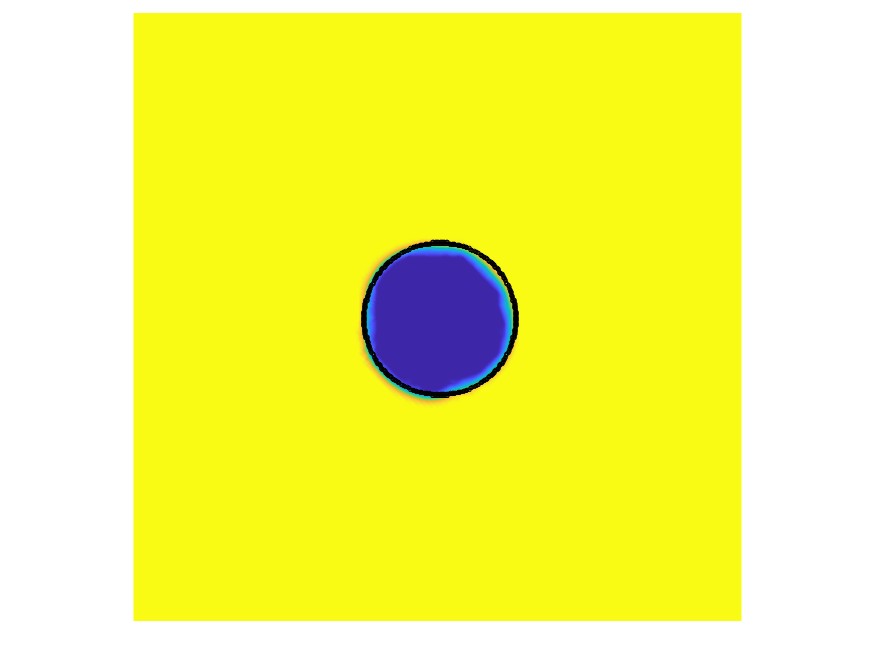}
&\includegraphics[height = \fht, trim = {3.5cm 0cm 3.5cm 0cm}, clip]{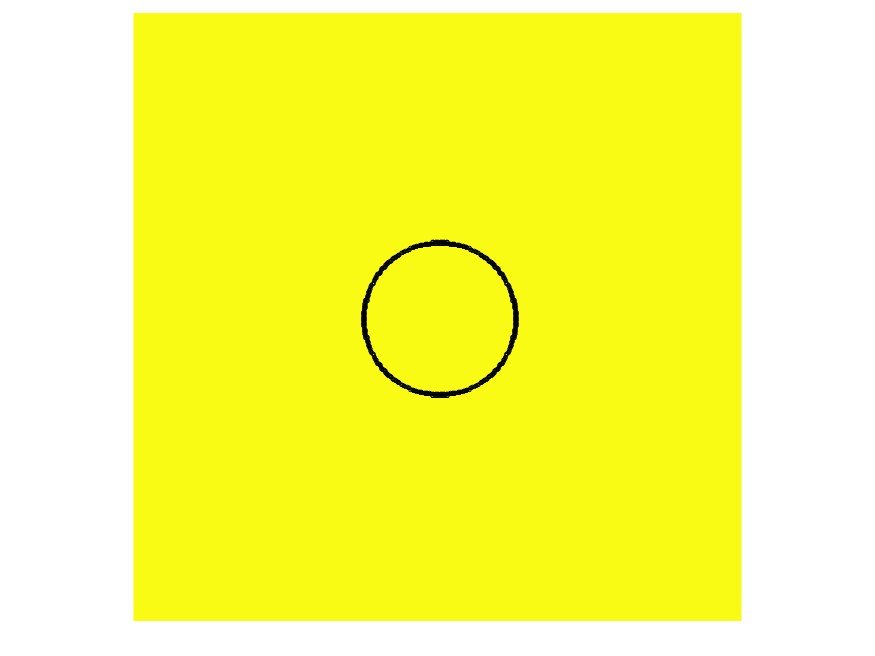}
\end{tabular}
\caption{The results by the adaptive method for the exact data $y^*$.
From left to the right are the mesh, recovered inclusion and error indicator function. The number of nodes for each step is
676, 938, 1304, 1822, 2550 and 3582.
}
\label{fig:circlefreeadaptive}
\end{figure}

\begin{figure}[hbt!]
\centering
\setlength{\tabcolsep}{0pt}
\begin{tabular}{cccccc}
\includegraphics[height = \fht, trim = {3.5cm 0cm 3.5cm 0cm}, clip]{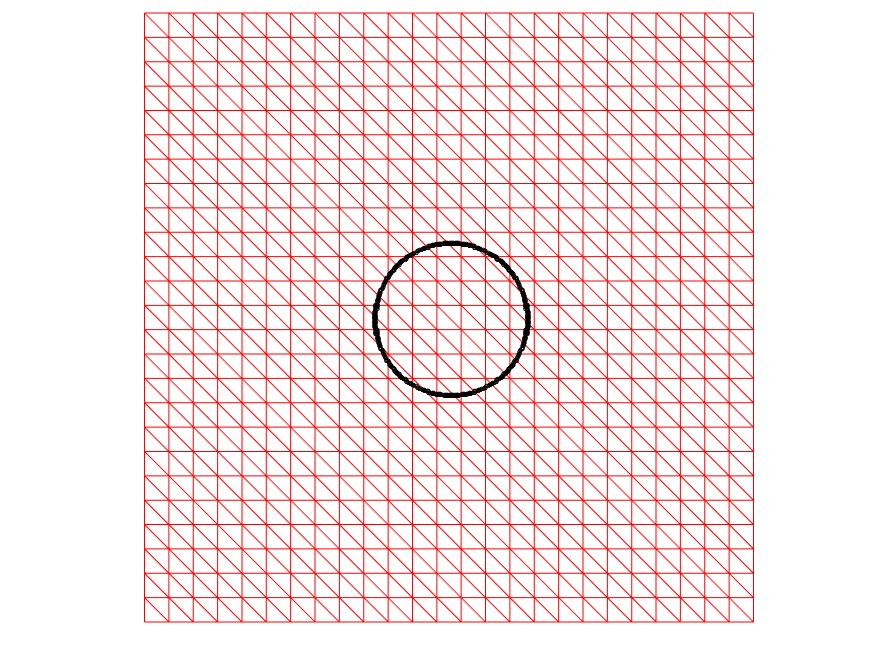}
&\includegraphics[height = \fht, trim = {3.5cm 0cm 3.5cm 0cm}, clip]{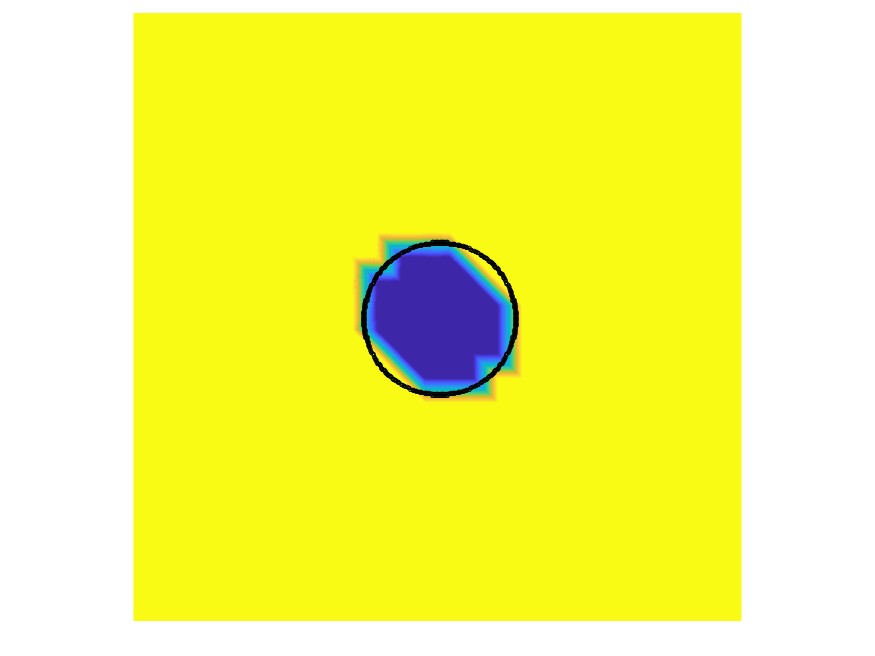}
&\includegraphics[height = \fht, trim = {3.5cm 0cm 3.5cm 0cm}, clip]{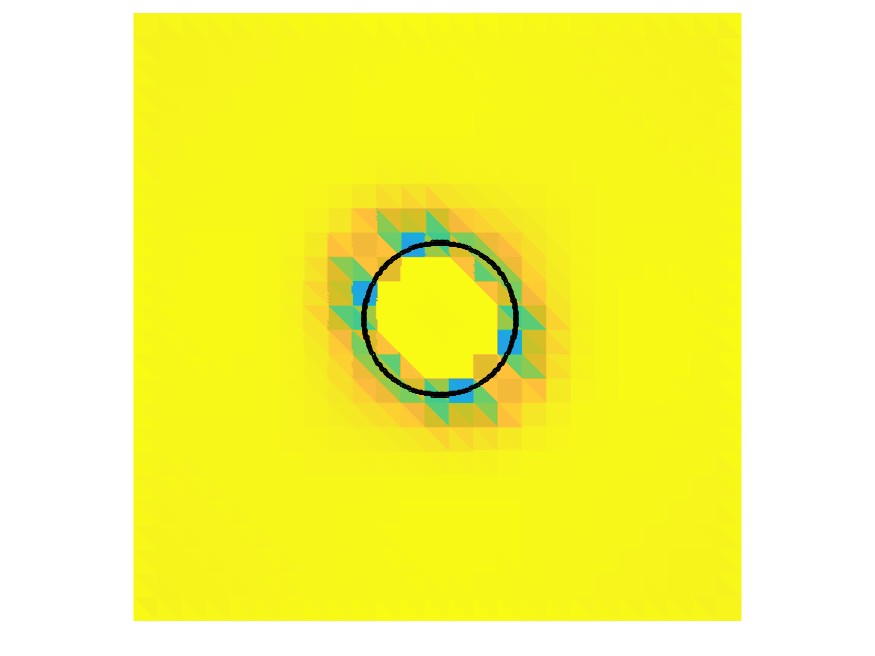}
&\includegraphics[height = \fht, trim = {3.5cm 0cm 3.5cm 0cm}, clip]{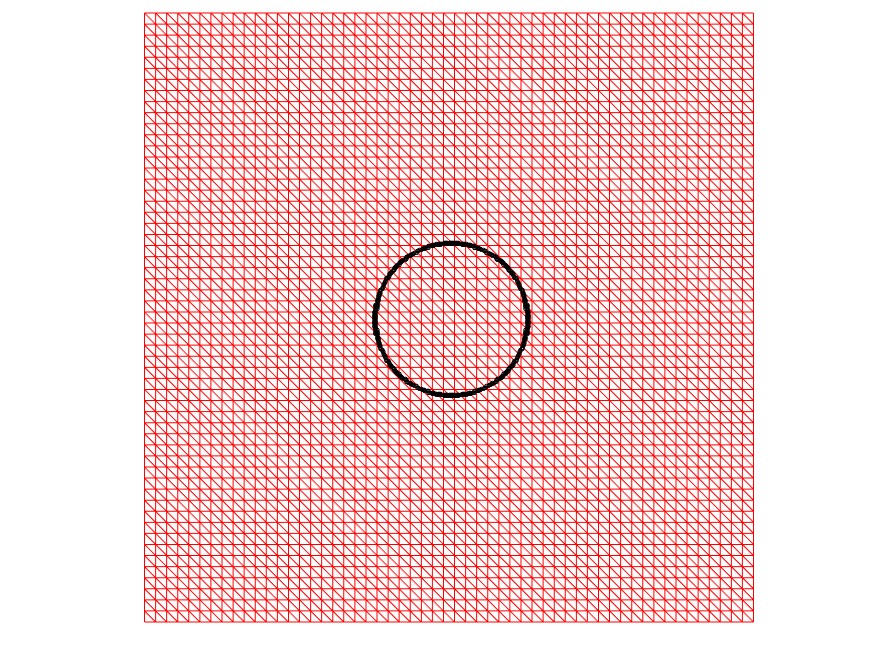}
&\includegraphics[height = \fht, trim = {3.5cm 0cm 3.5cm 0cm}, clip]{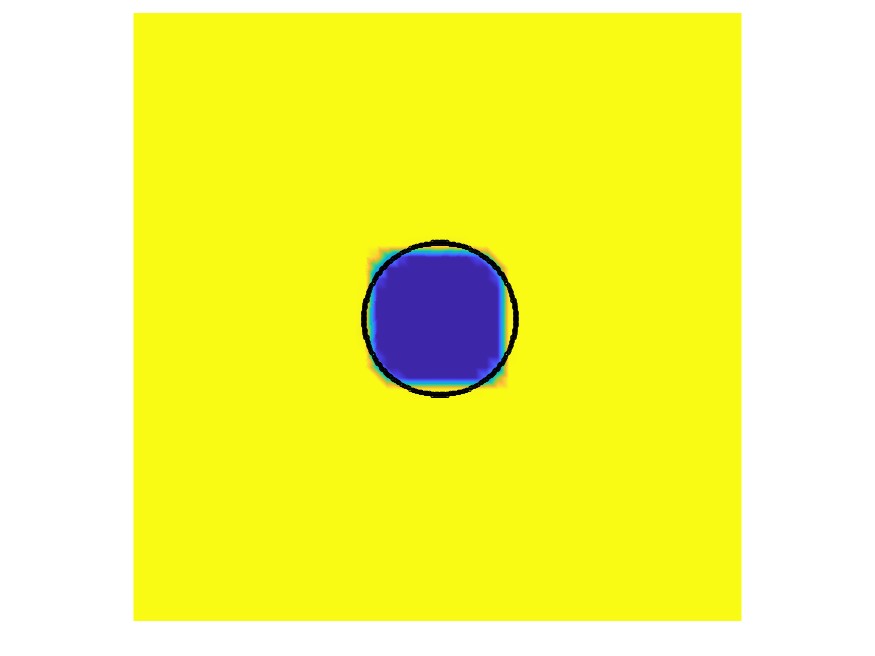}
&\includegraphics[height = \fht, trim = {3.5cm 0cm 3.5cm 0cm}, clip]{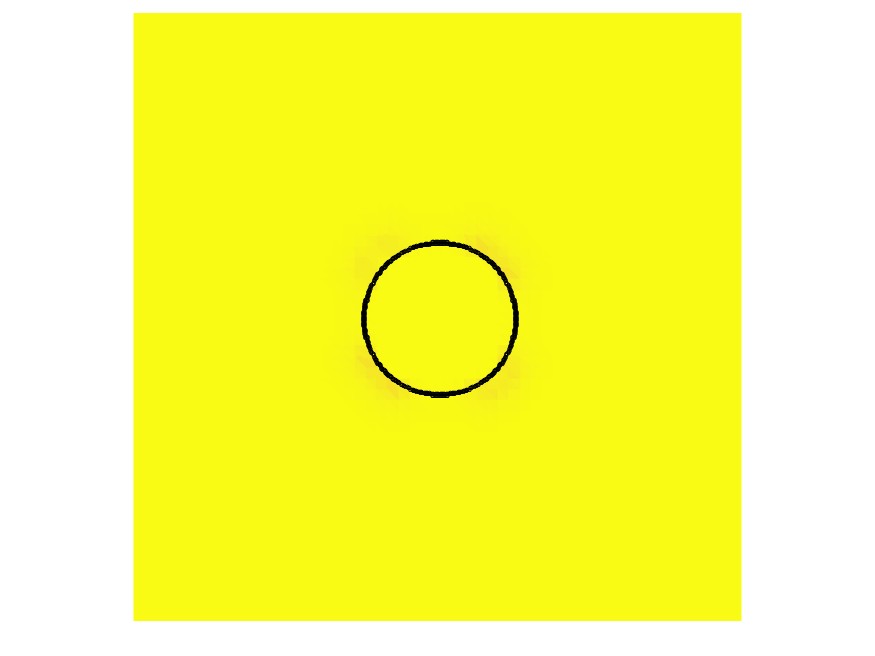}\\
\includegraphics[height = \fht, trim = {3.5cm 0cm 3.5cm 0cm}, clip]{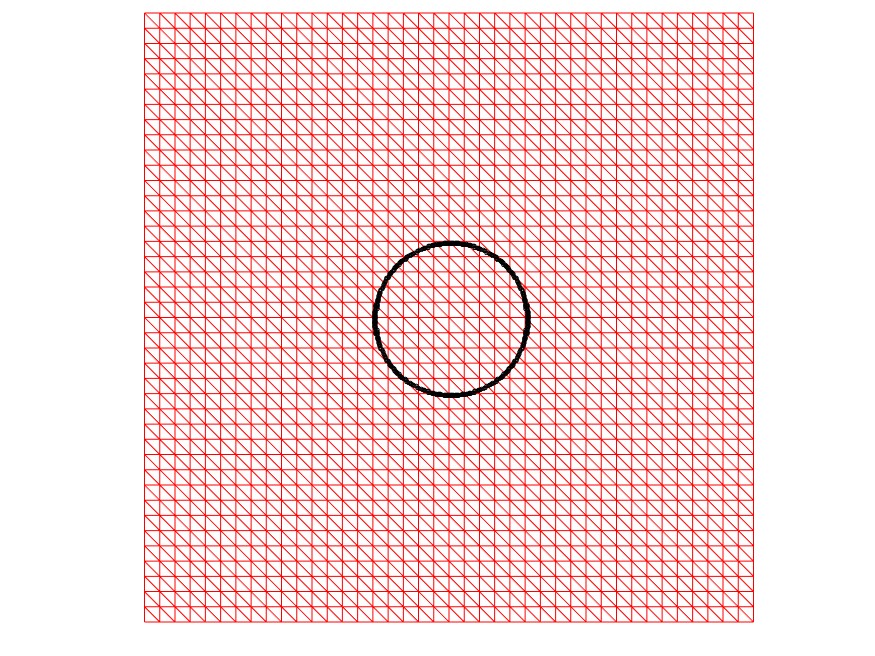}
&\includegraphics[height = \fht, trim = {3.5cm 0cm 3.5cm 0cm}, clip]{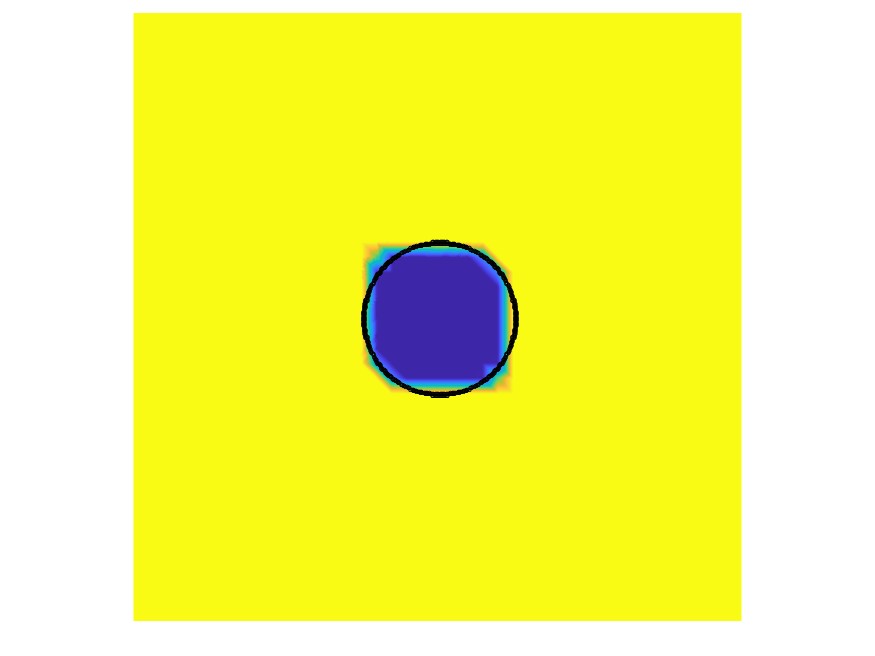}
&\includegraphics[height = \fht, trim = {3.5cm 0cm 3.5cm 0cm}, clip]{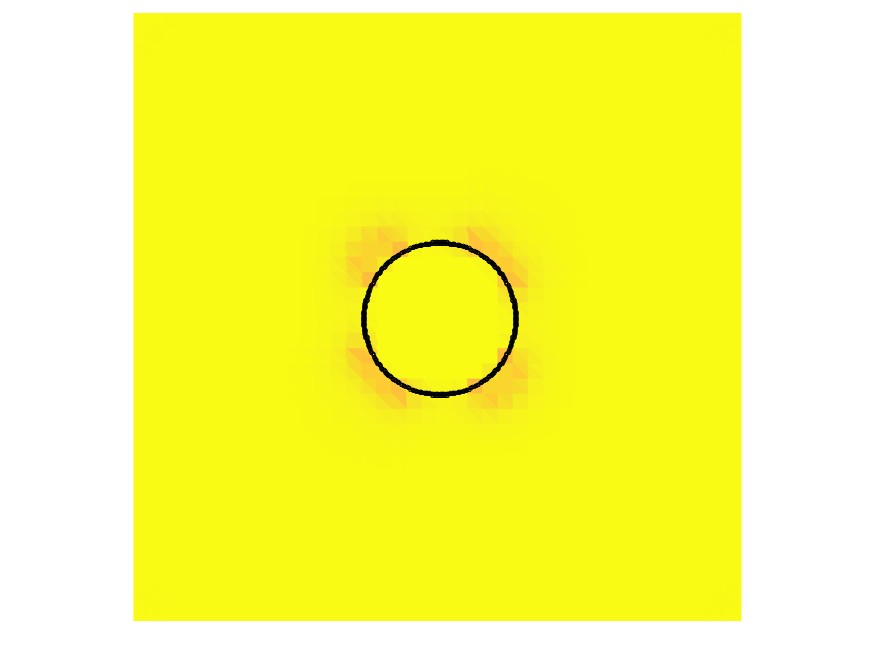}
&\includegraphics[height = \fht, trim = {3.5cm 0cm 3.5cm 0cm}, clip]{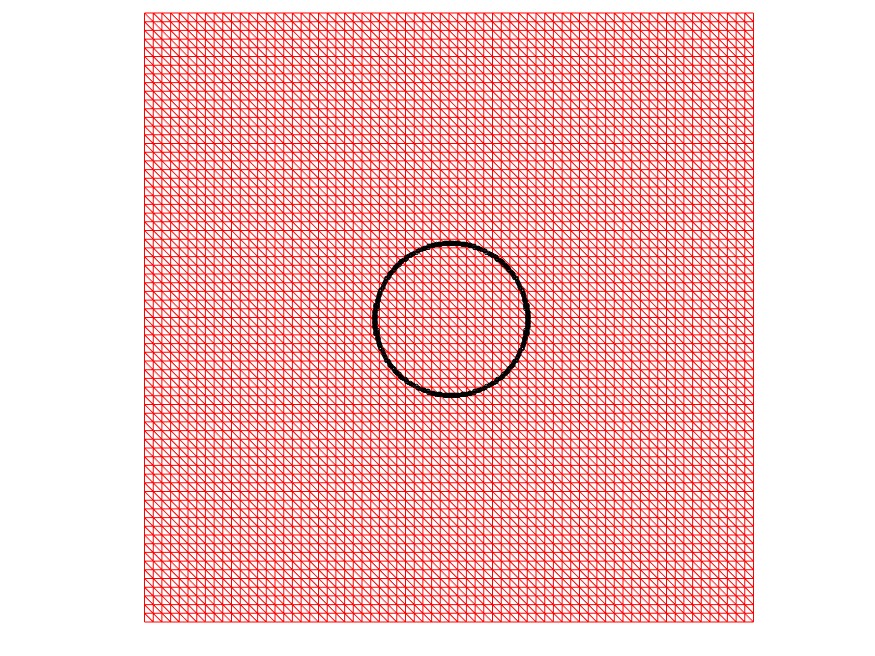}
&\includegraphics[height = \fht, trim = {3.5cm 0cm 3.5cm 0cm}, clip]{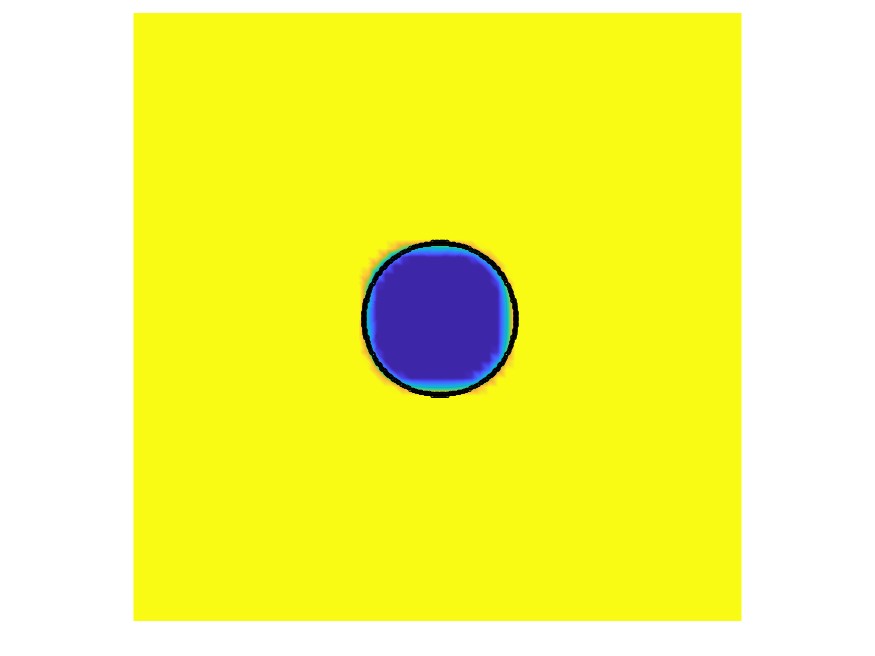}
&\includegraphics[height = \fht, trim = {3.5cm 0cm 3.5cm 0cm}, clip]{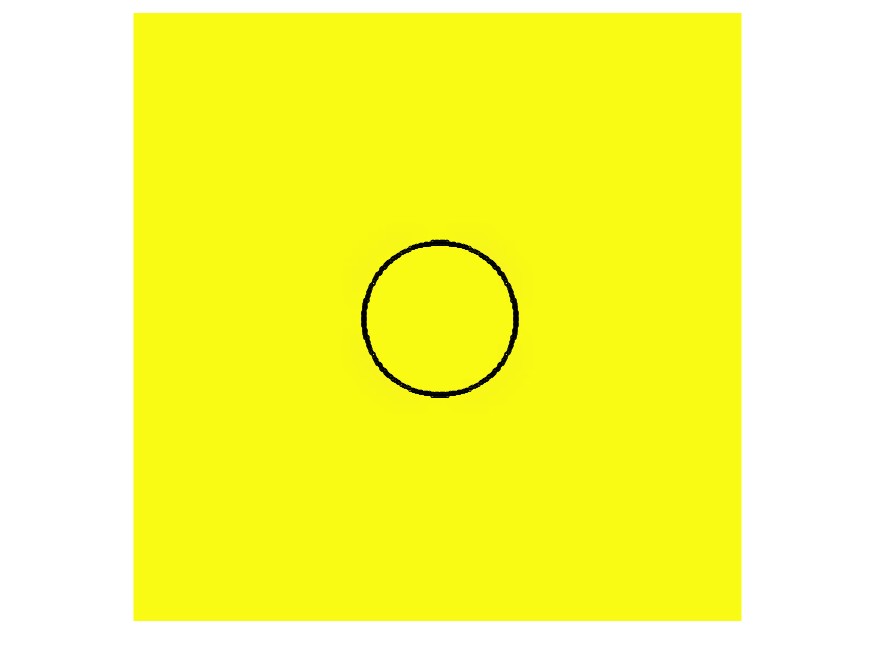}
\end{tabular}
\caption{The results by uniform mesh refinements for the exact data $y^*$.
From left to the right are the mesh, recovered inclusion and error indicator function.
The number of nodes at each step is 676, 1681, 3136 and 5041.
}
\label{fig:circlefreeuniform}
\end{figure}

\begin{figure}[hbt!]
\centering
\setlength{\tabcolsep}{0pt}
\begin{tabular}{cccccc}
\includegraphics[height = \fht, trim = {3.5cm 0cm 3.5cm 0cm}, clip]{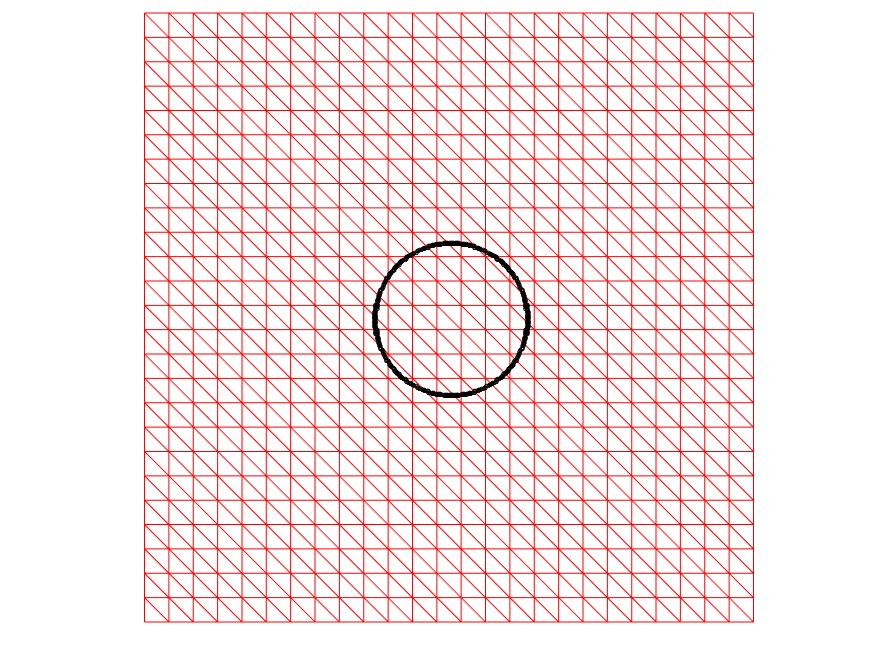}
&\includegraphics[height = \fht, trim = {3.5cm 0cm 3.5cm 0cm}, clip]{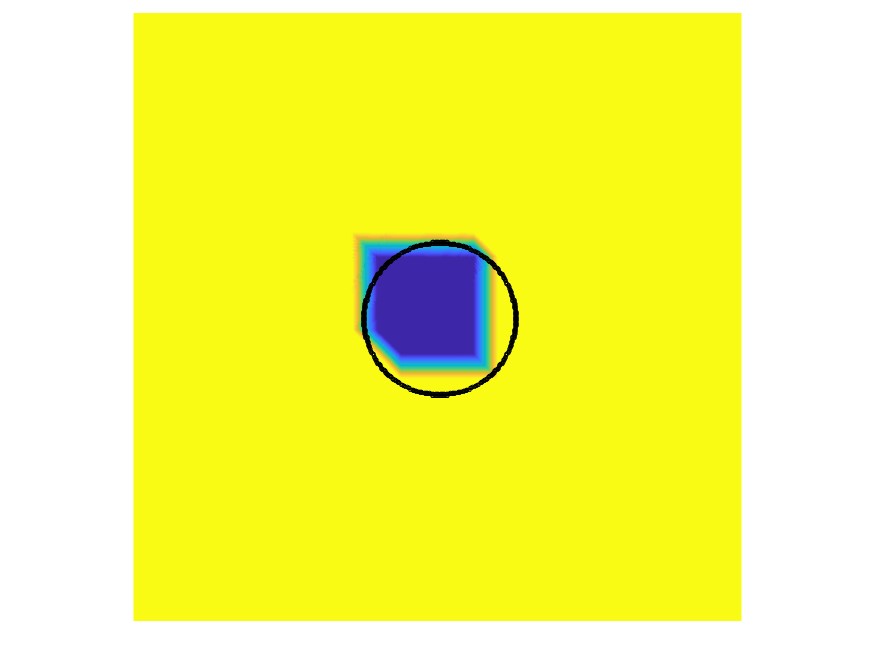}
&\includegraphics[height = \fht, trim = {3.5cm 0cm 3.5cm 0cm}, clip]{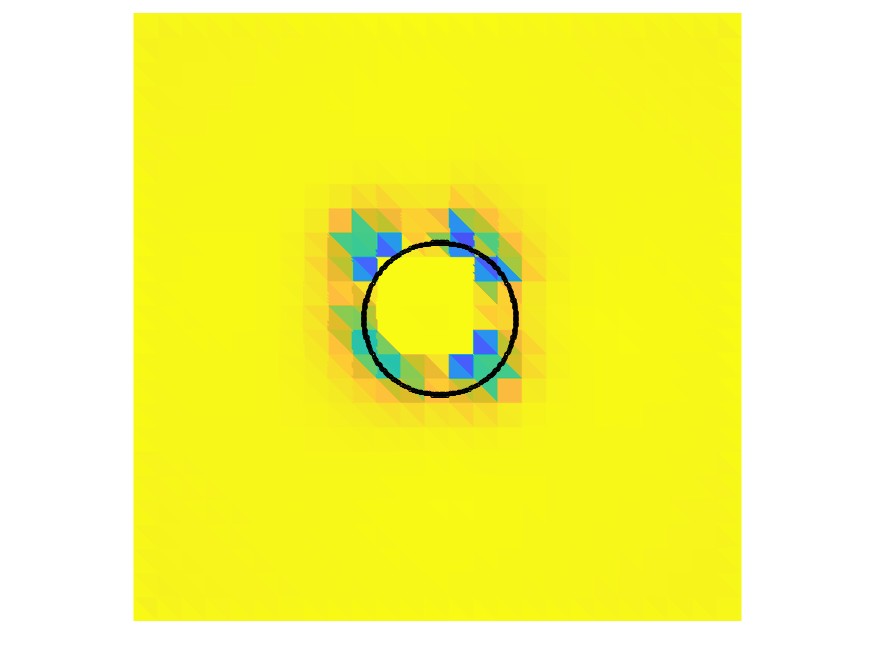}
&\includegraphics[height = \fht, trim = {3.5cm 0cm 3.5cm 0cm}, clip]{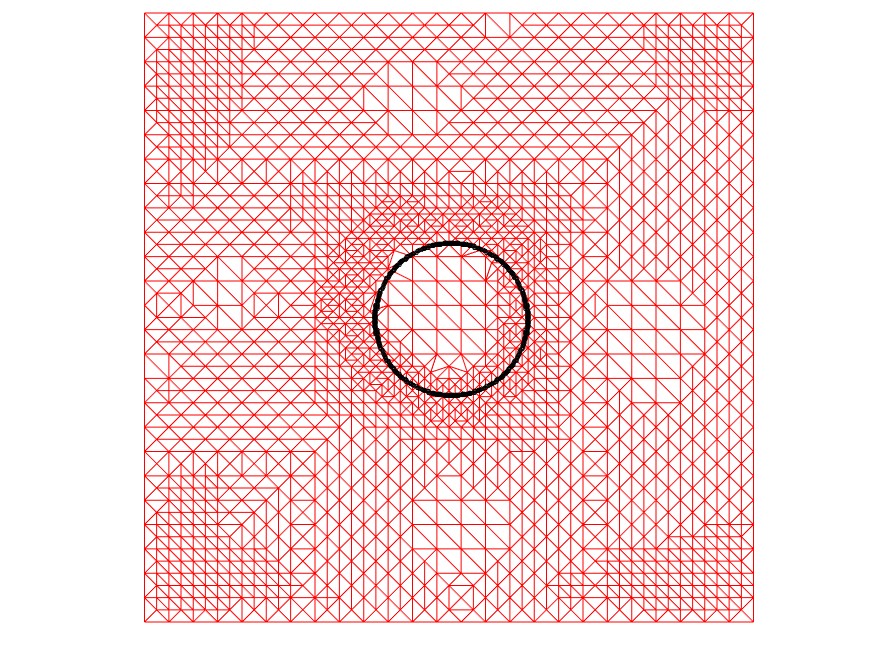}
&\includegraphics[height = \fht, trim = {3.5cm 0cm 3.5cm 0cm}, clip]{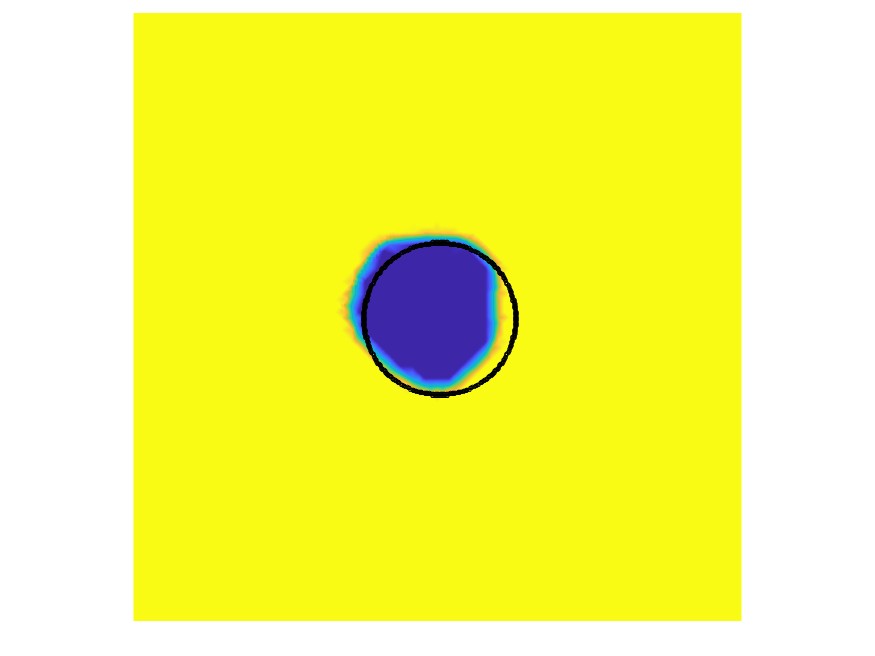}
&\includegraphics[height = \fht, trim = {3.5cm 0cm 3.5cm 0cm}, clip]{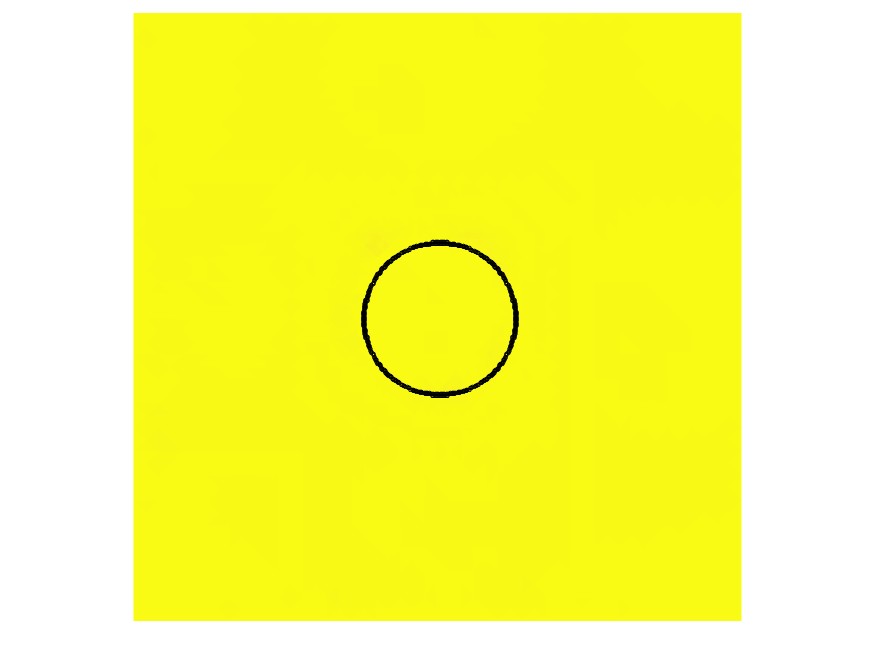}\\
\includegraphics[height = \fht, trim = {3.5cm 0cm 3.5cm 0cm}, clip]{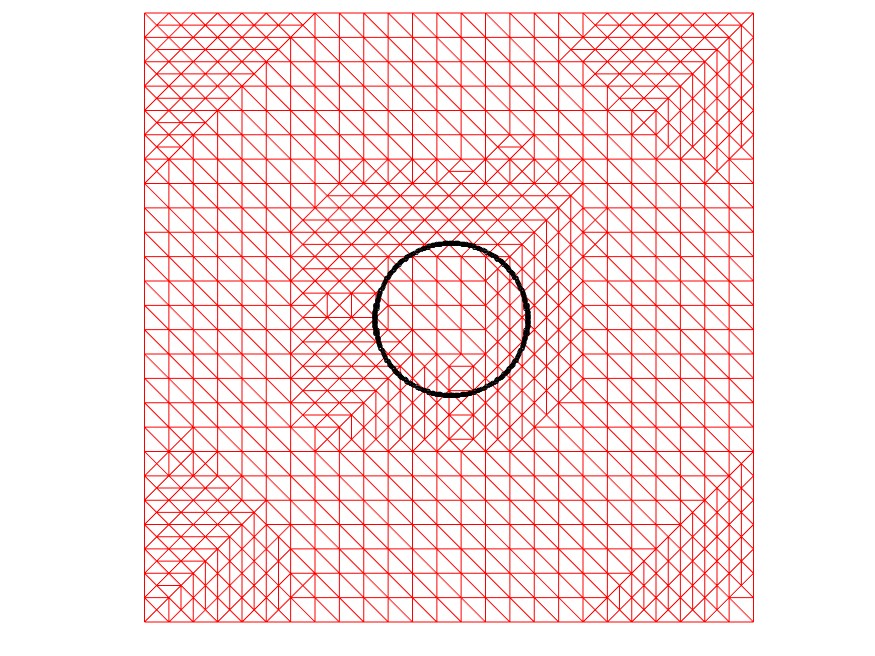}
&\includegraphics[height = \fht, trim = {3.5cm 0cm 3.5cm 0cm}, clip]{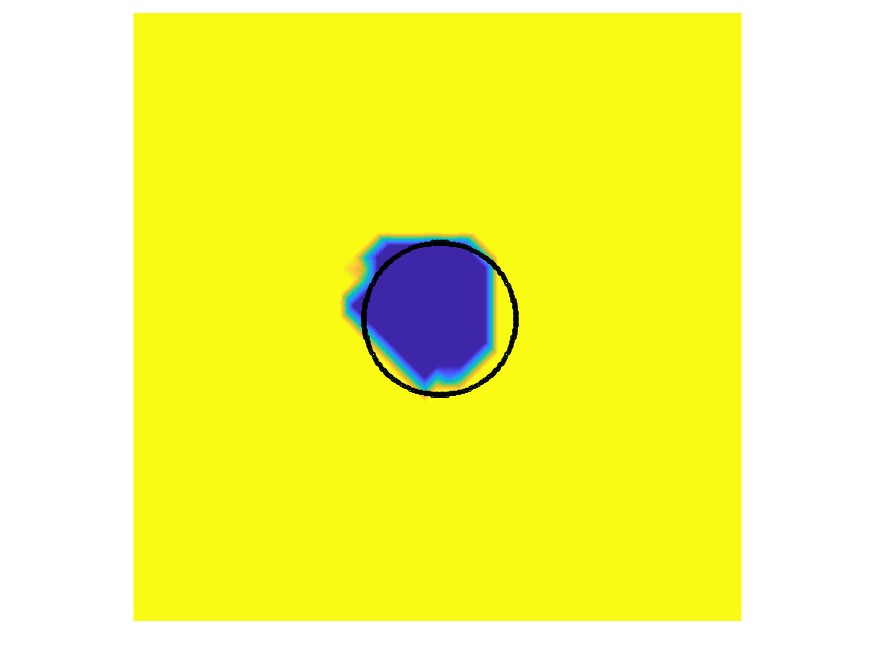}
&\includegraphics[height = \fht, trim = {3.5cm 0cm 3.5cm 0cm}, clip]{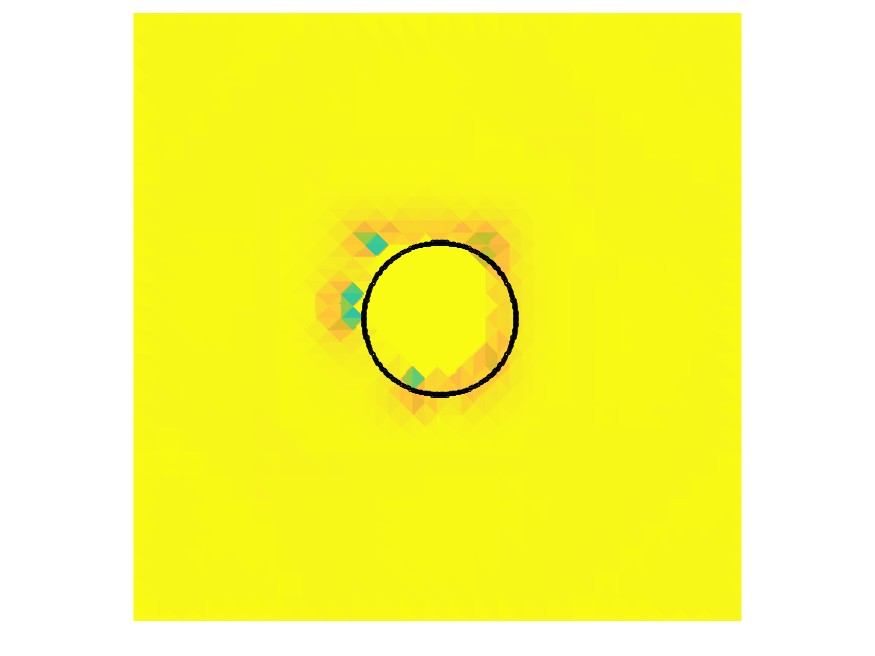}
&\includegraphics[height = \fht, trim = {3.5cm 0cm 3.5cm 0cm}, clip]{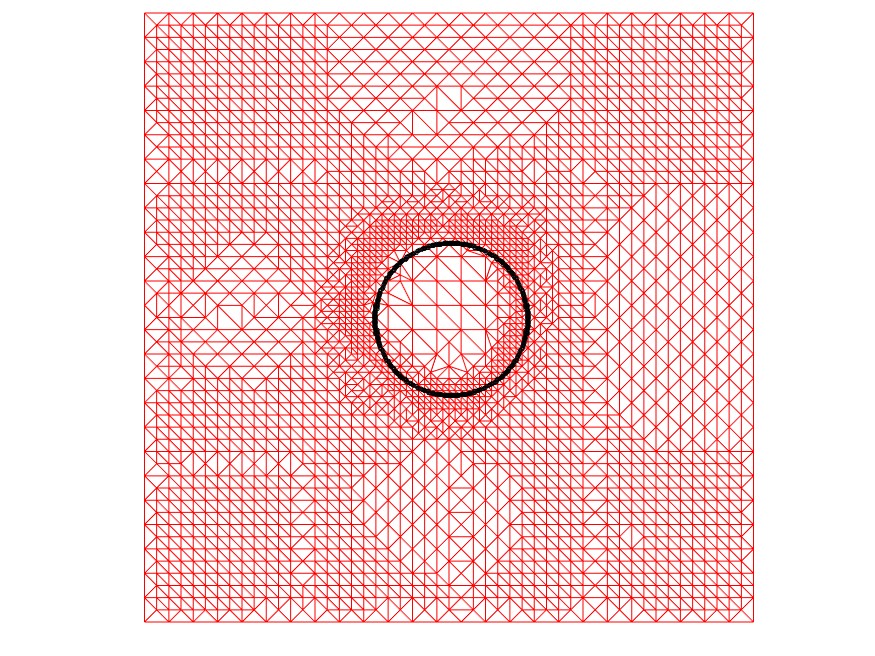}
&\includegraphics[height = \fht, trim = {3.5cm 0cm 3.5cm 0cm}, clip]{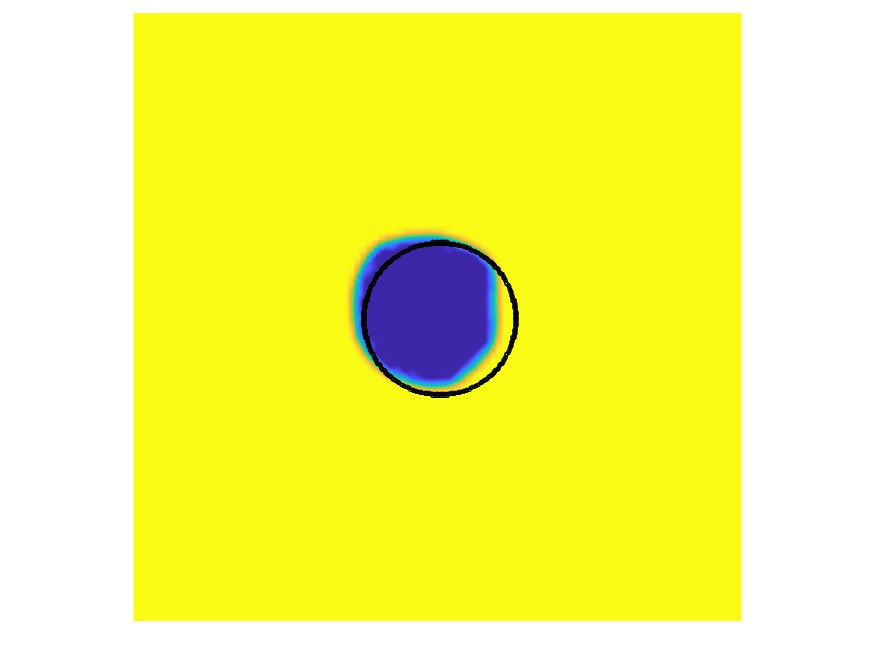}
&\includegraphics[height = \fht, trim = {3.5cm 0cm 3.5cm 0cm}, clip]{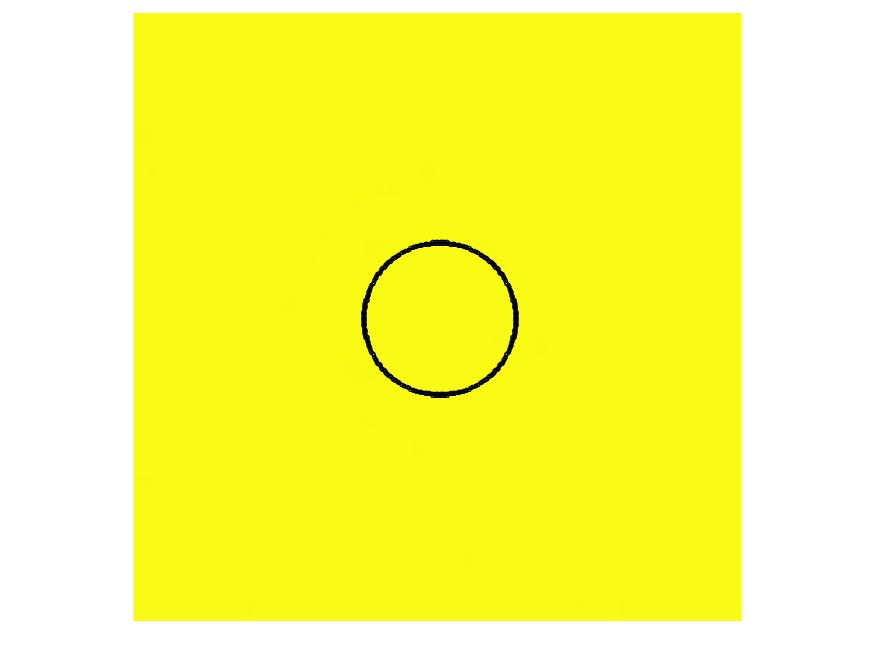}\\
\includegraphics[height = \fht, trim = {3.5cm 0cm 3.5cm 0cm}, clip]{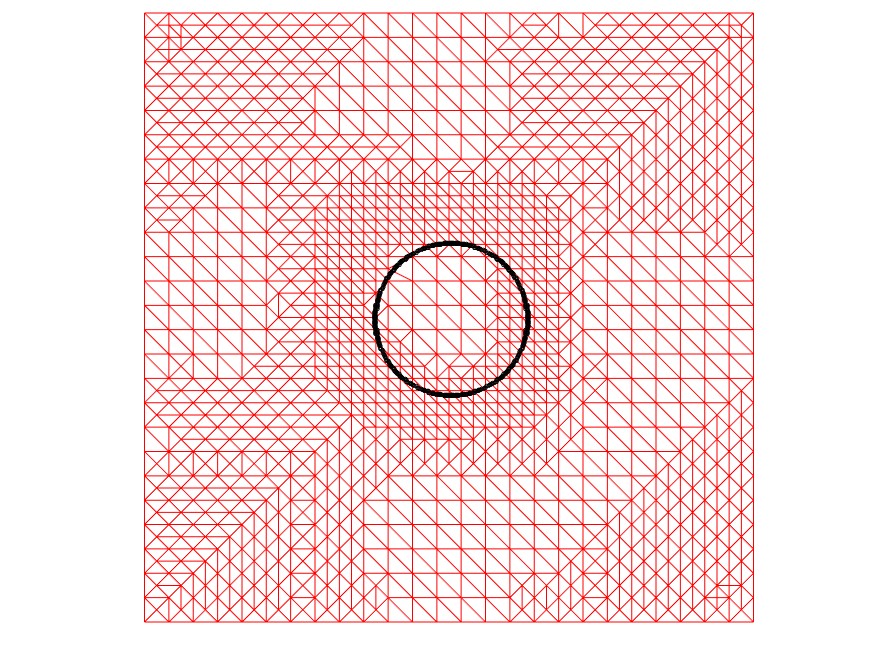}
&\includegraphics[height = \fht, trim = {3.5cm 0cm 3.5cm 0cm}, clip]{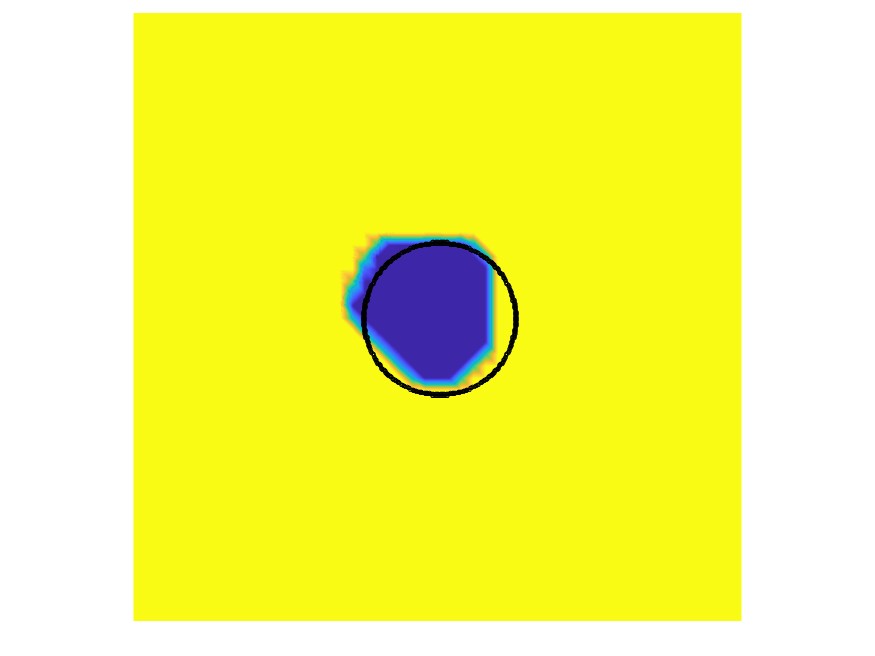}
&\includegraphics[height = \fht, trim = {3.5cm 0cm 3.5cm 0cm}, clip]{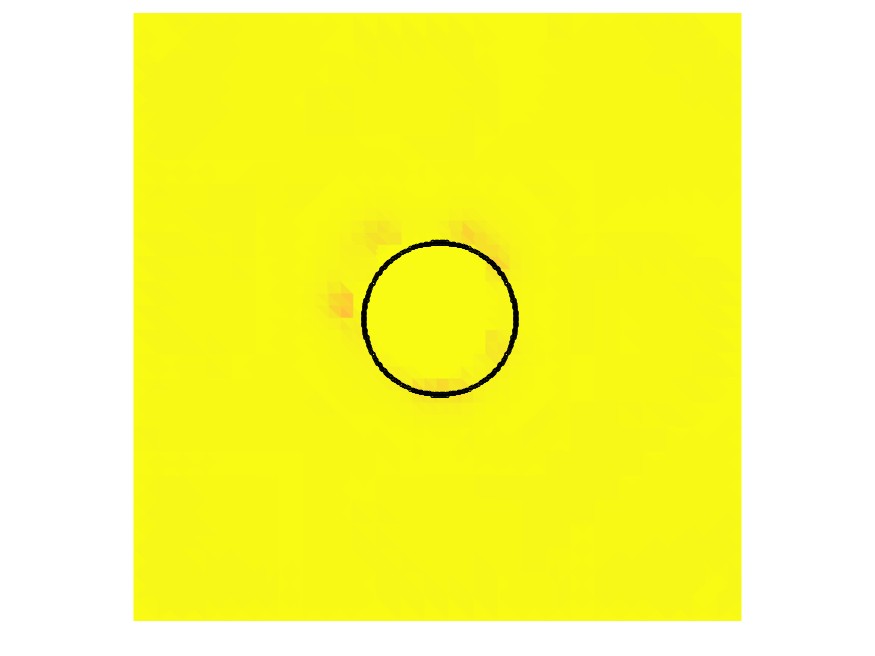}
&\includegraphics[height = \fht, trim = {3.5cm 0cm 3.5cm 0cm}, clip]{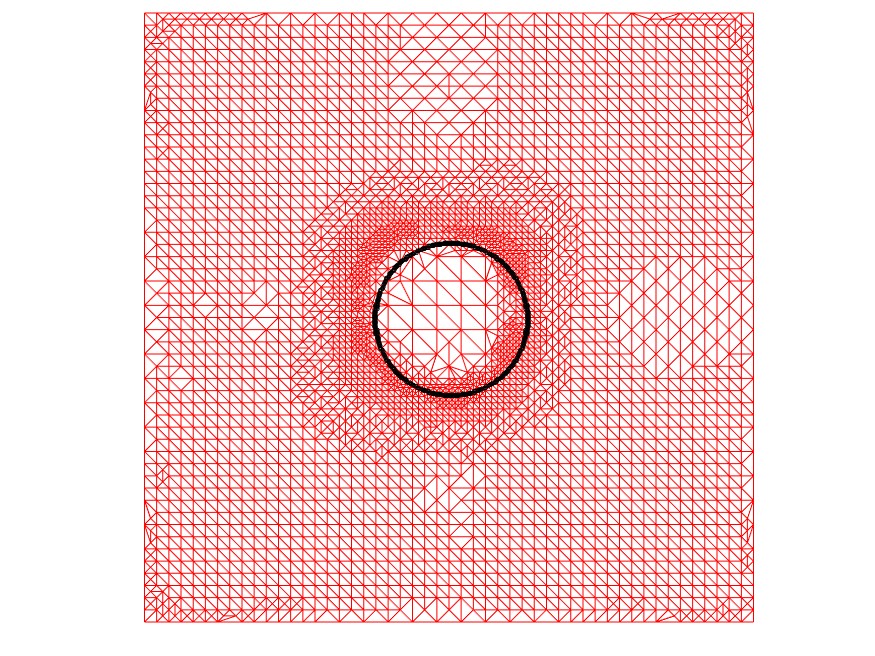}
&\includegraphics[height = \fht, trim = {3.5cm 0cm 3.5cm 0cm}, clip]{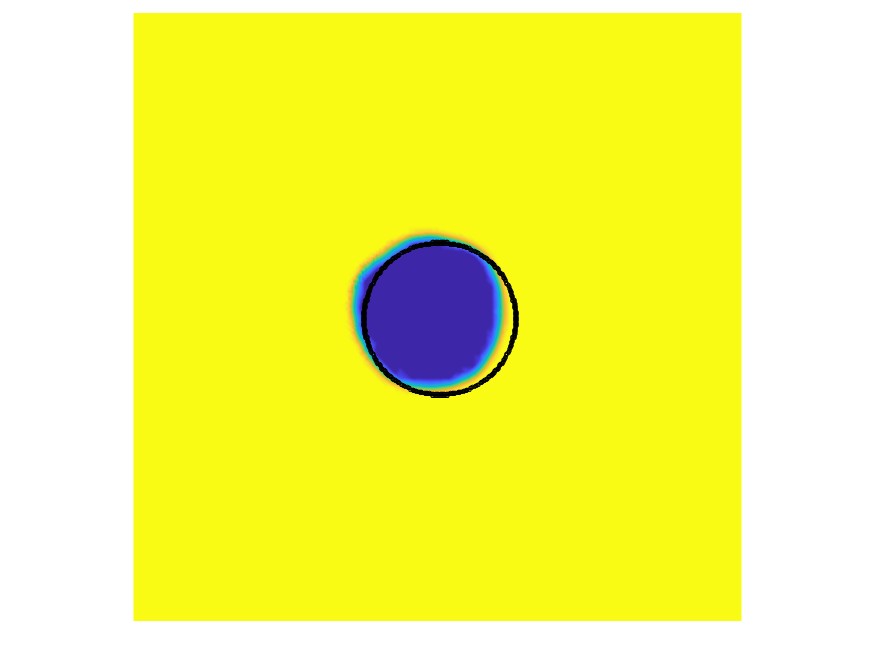}
&\includegraphics[height = \fht, trim = {3.5cm 0cm 3.5cm 0cm}, clip]{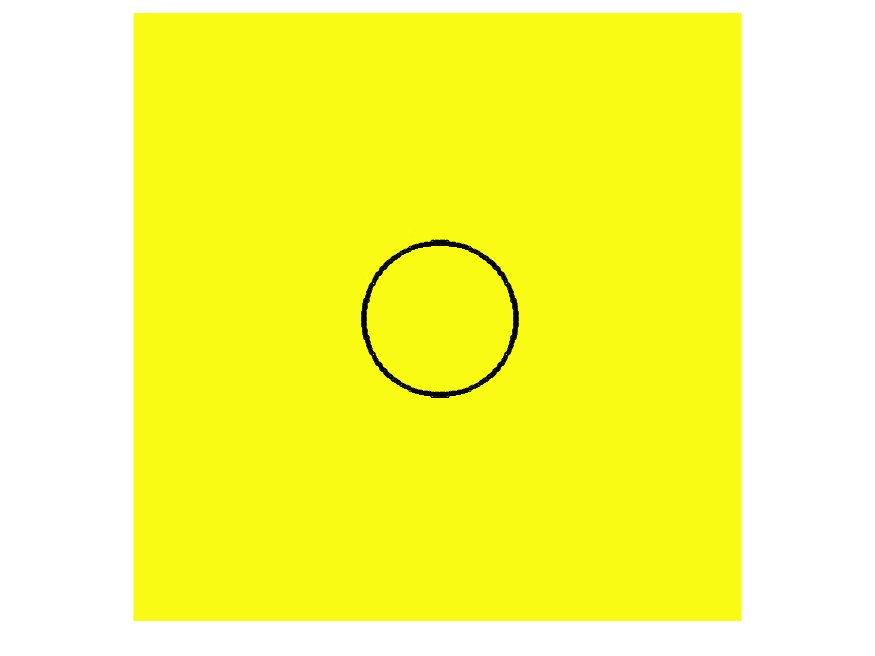}
\end{tabular}
\caption{
The results by the adaptive method for the noisy data $y^\delta$. From left to the right are the mesh, recovered inclusion and error indicator function.
The number of nodes for each step is
676, 933, 1289, 1799, 2515 and 3523.
}
\label{fig:circlenoiseadaptive}
\end{figure}

\begin{figure}[hbt!]
\centering
\setlength{\tabcolsep}{0pt}
\begin{tabular}{cccccc}
\includegraphics[height = \fht, trim = {3.5cm 0cm 3.5cm 0cm}, clip]{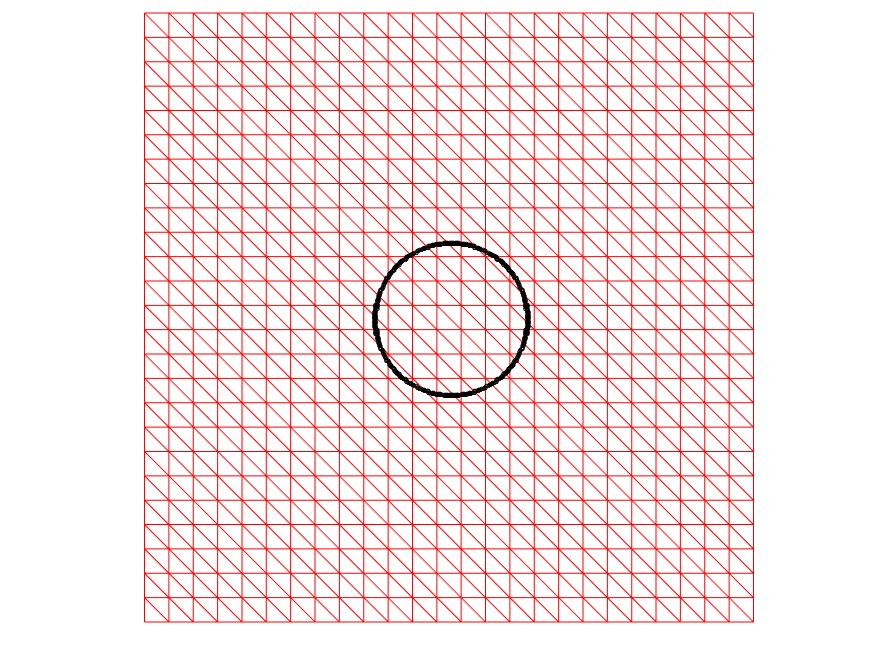}
&\includegraphics[height = \fht, trim = {3.5cm 0cm 3.5cm 0cm}, clip]{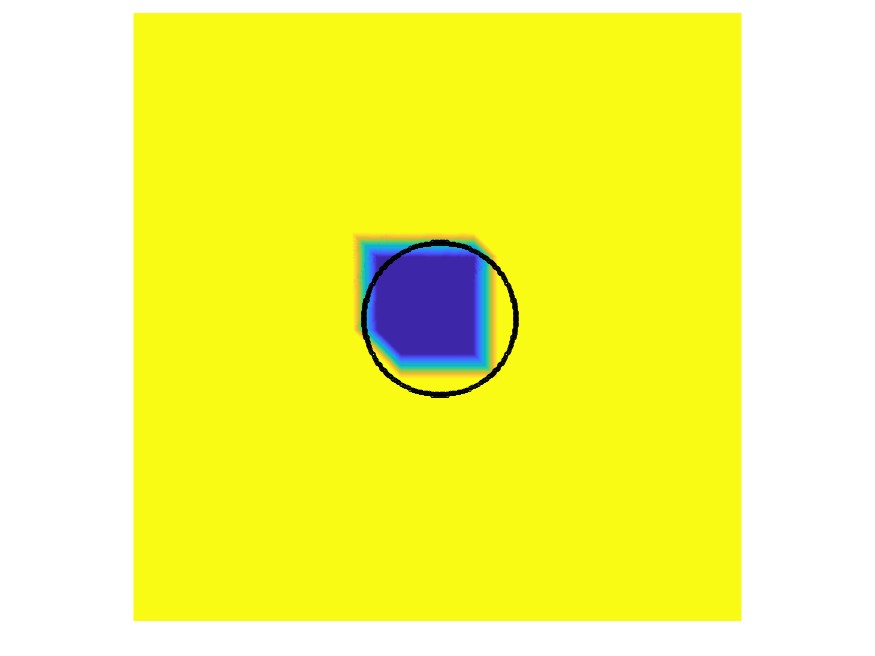}
&\includegraphics[height = \fht, trim = {3.5cm 0cm 3.5cm 0cm}, clip]{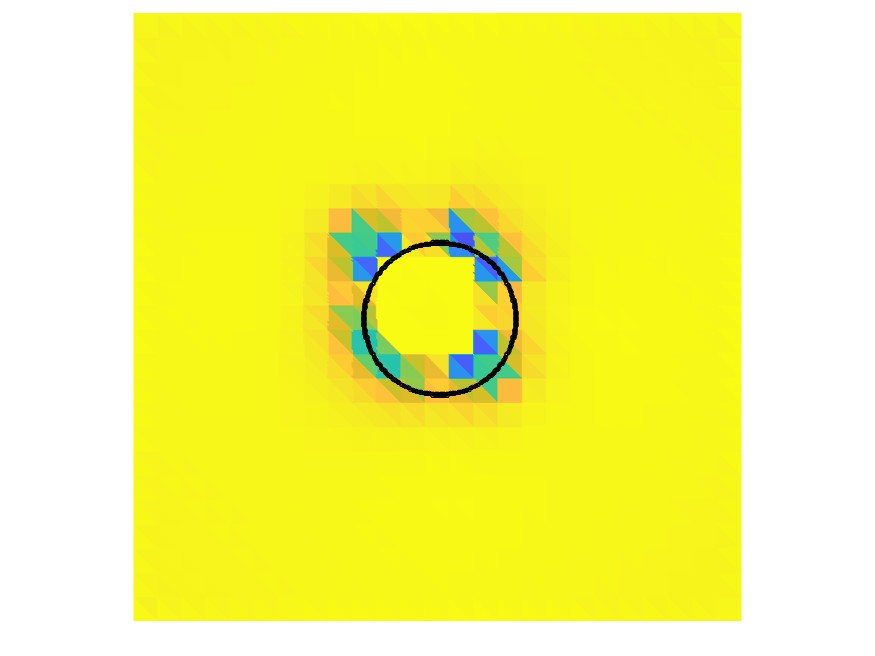}
&\includegraphics[height = \fht, trim = {3.5cm 0cm 3.5cm 0cm}, clip]{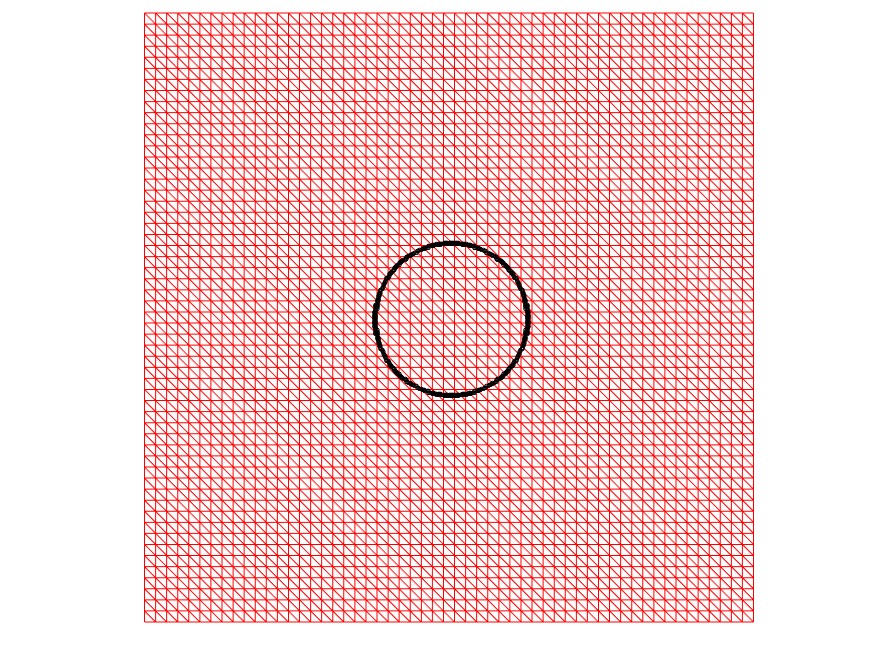}
&\includegraphics[height = \fht, trim = {3.5cm 0cm 3.5cm 0cm}, clip]{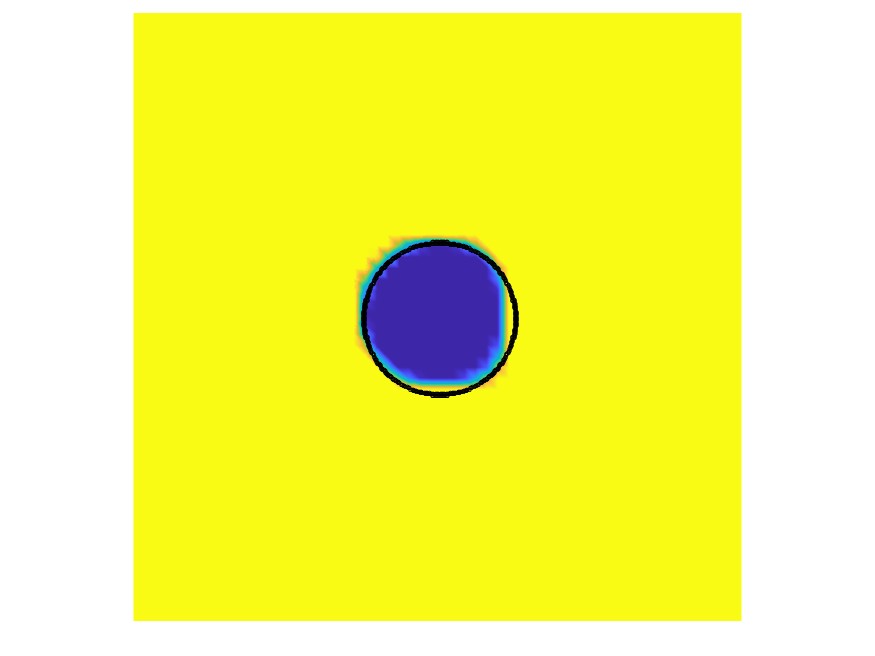}
&\includegraphics[height = \fht, trim = {3.5cm 0cm 3.5cm 0cm}, clip]{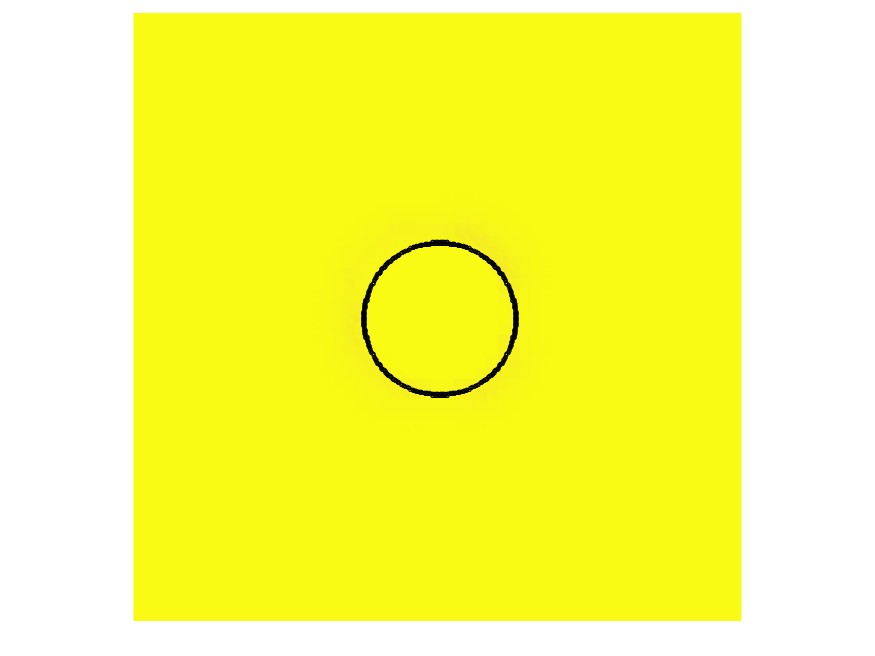}\\
\includegraphics[height = \fht, trim = {3.5cm 0cm 3.5cm 0cm}, clip]{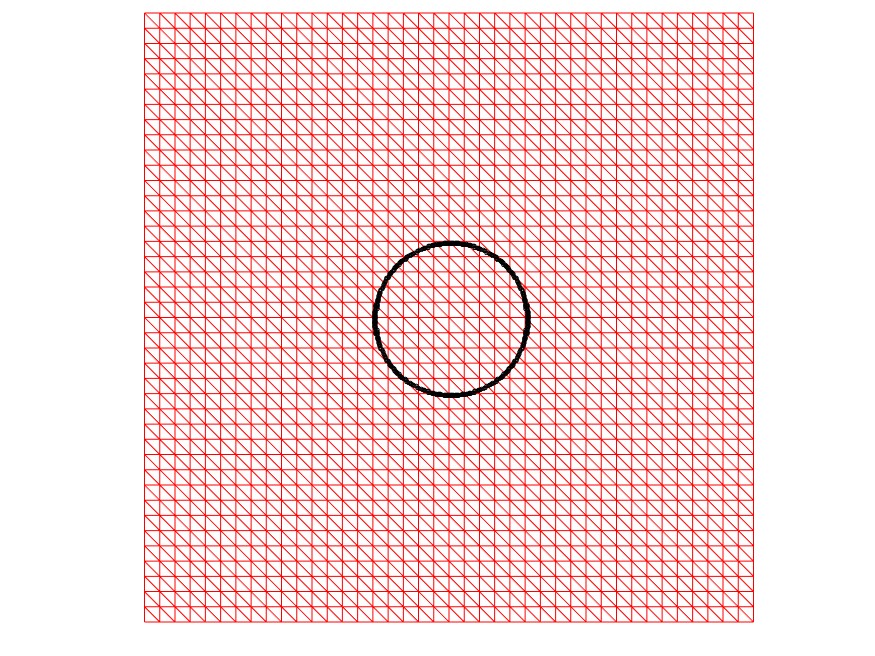}
&\includegraphics[height = \fht, trim = {3.5cm 0cm 3.5cm 0cm}, clip]{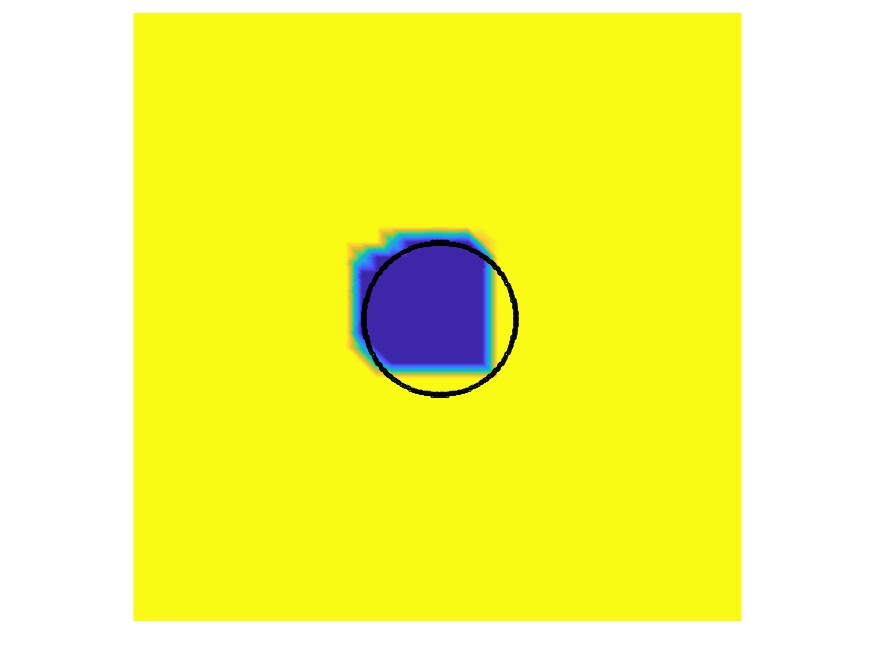}
&\includegraphics[height = \fht, trim = {3.5cm 0cm 3.5cm 0cm}, clip]{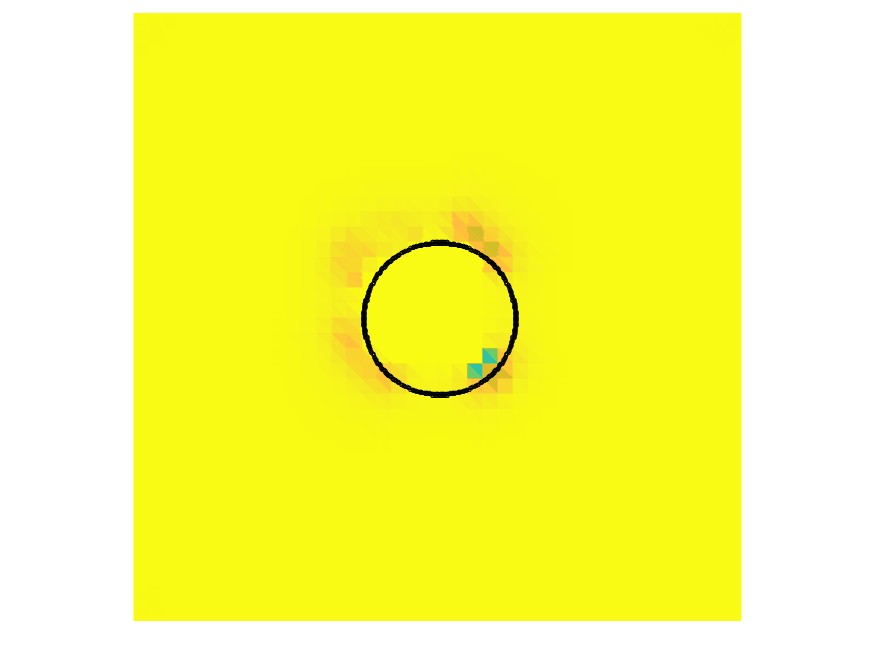}
&\includegraphics[height = \fht, trim = {3.5cm 0cm 3.5cm 0cm}, clip]{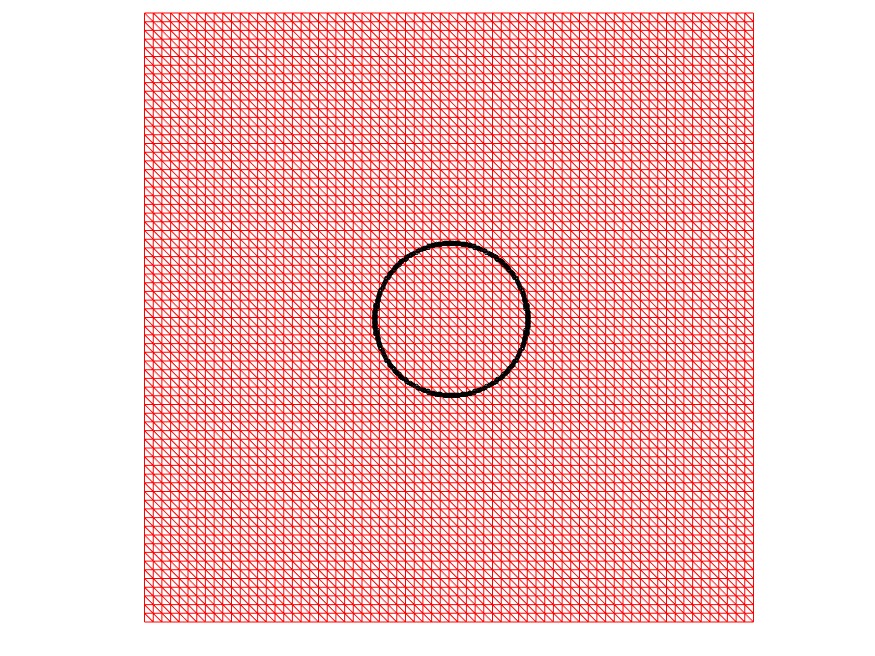}
&\includegraphics[height = \fht, trim = {3.5cm 0cm 3.5cm 0cm}, clip]{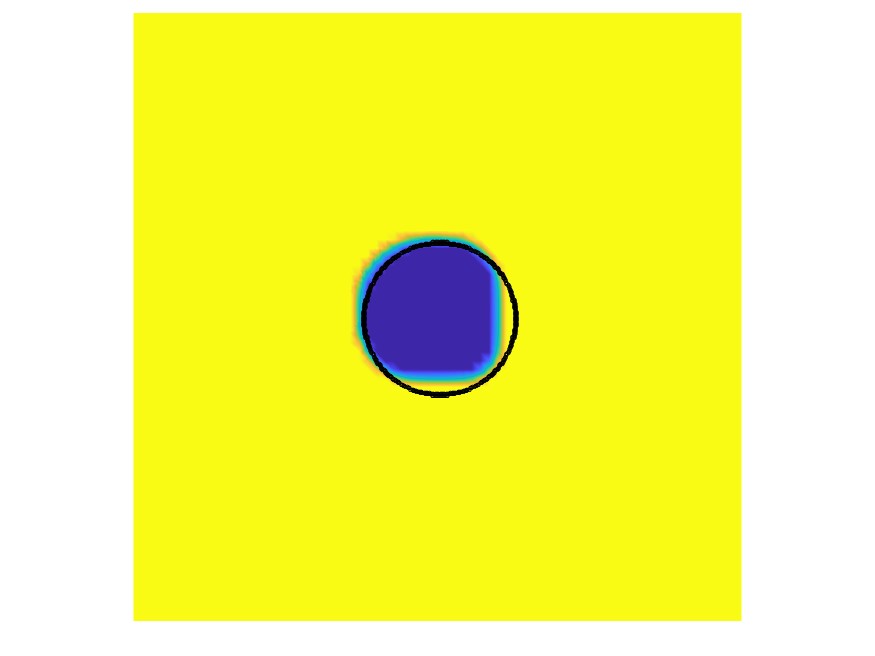}
&\includegraphics[height = \fht, trim = {3.5cm 0cm 3.5cm 0cm}, clip]{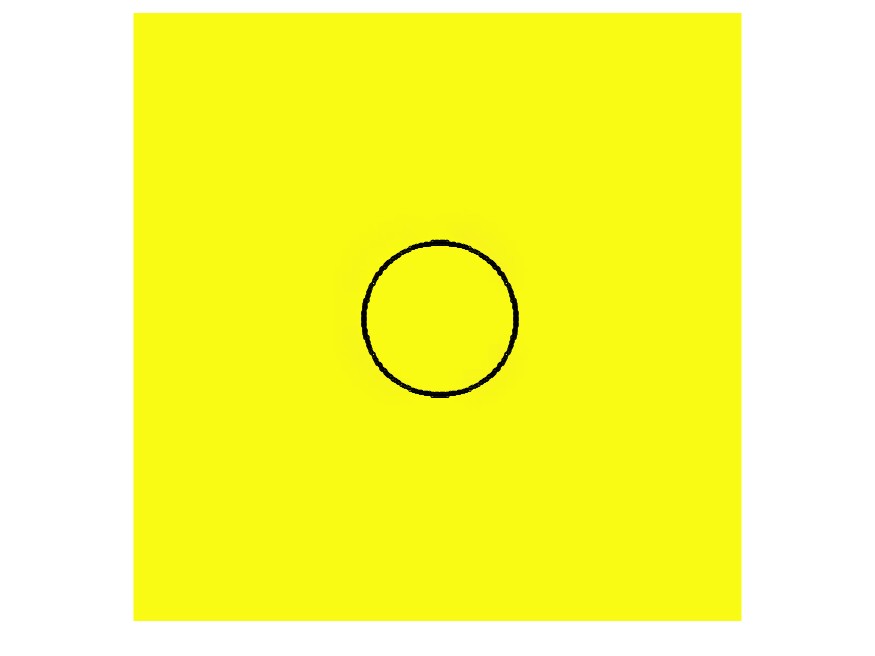}
\end{tabular}
\caption{The results by uniform mesh refinements for the noisy data $y^\delta$.
From left to the right are the mesh, recovered inclusion and error indicator function.
The number of nodes at each step is 676, 1681, 3136 and 5041.
}
\label{fig:circlenoiseuniform}
\end{figure}

The numerical results for the case of one circular inclusion are shown in Figs. \ref{fig:circlefreeadaptive} and \ref{fig:circlefreeuniform} for the exact data $y^*$, and Figs. \ref{fig:circlenoiseadaptive} and \ref{fig:circlenoiseuniform} for the noisy data $y^\delta$, respectively.
In each figure, from left to right, the plots show the resulting mesh, the reconstructed inclusion, and the error indicator $\widetilde{\eta}_{k}(T):=\max_{j=1,3}\eta_{k}(T)$ (abbreviating the arguments $y^*_k,u^*_k$ and $p_k^*$).
To facilitate the comparison of the indicator $\widetilde{\eta}_{k}$ on different meshes $\mathcal{T}_k$, the color bars have been normalized to the range $[0, \text{2e-3}]$.
On the initial mesh $\mathcal{T}_0$, the reconstructions can generally capture the location of the circular inclusion $\omega$, but the recovered shape suffer from pronounced deformations near the inclusion interface $\partial\omega$.
The \textit{a posteriori} error estimator $\widetilde\eta_k$ also takes relatively large values near the interface $\partial\omega$ (about the order 2e-3), and can correctly identify the region with large errors associated with the state (due to the discontinuous interface) and the conductivity, which are then used for the marking and the subsequent adaptive mesh refinement.
Indeed, based on the estimator $\widetilde\eta_k$, the adaptive algorithm can effectively perform the mesh refinement in the concerned region, and after several adaptive iterations, the algorithm can accurately resolve the circular inclusion.
Notably, the inclusion location remains fairly stable after multiple refinement loops, and the final mesh refinements take place mostly around the interface $\partial\omega$ when compared with the initial mesh $\mathcal{T}_0$, and partly extends also towards the external boundary, where the boundary observation $y^\delta$ is taken.
This is expected since the discretization error around the interface $\partial\omega$ becomes very small after several refinement steps and {the value of the estimator $\widetilde{\eta}_k$ is comparable everywhere so the mesh starts to be refined everywhere.}
It is observed that the adaptive algorithm not only enhances the accuracy of the shape recovery of the inclusions but also adapts to the accuracy of the data $y^\delta$ on the boundary $\partial\Omega$.
When the algorithm terminates, the estimator $\widetilde{\eta}_k$ also takes small value, indicating the convergence of the algorithm to a solution of the first-order optimality system \eqref{optsys_dis}, cf. Theorem \ref{thm:conv}. The convergence is also consistently observed for the uniform refinement strategy, but it requires more degree of freedoms in order to achieve comparable accuracy. The presence of noise does influence the reconstruction slightly. Nonetheless, the reconstruction remains reasonably accurate, and the preceding observations are still valid. For both exact data $y^*$ and noisy data $y^\delta$, the final reconstructions given by the adaptive and uniform refinements are similar to each other. Moreover, the comparison between the reconstructions for exact and noisy data in Figs. \ref{fig:circlefreeadaptive}--\ref{fig:circlenoiseuniform} indicates that the algorithm is robust with respect to data noise, which is a highly desirable feature for practical applications.

\begin{figure}[hbt!]
\centering
\setlength{\tabcolsep}{0pt}
\begin{tabular}{cccccc}
\includegraphics[height = \fht, trim = {3.5cm 0cm 3.5cm 0cm}, clip]{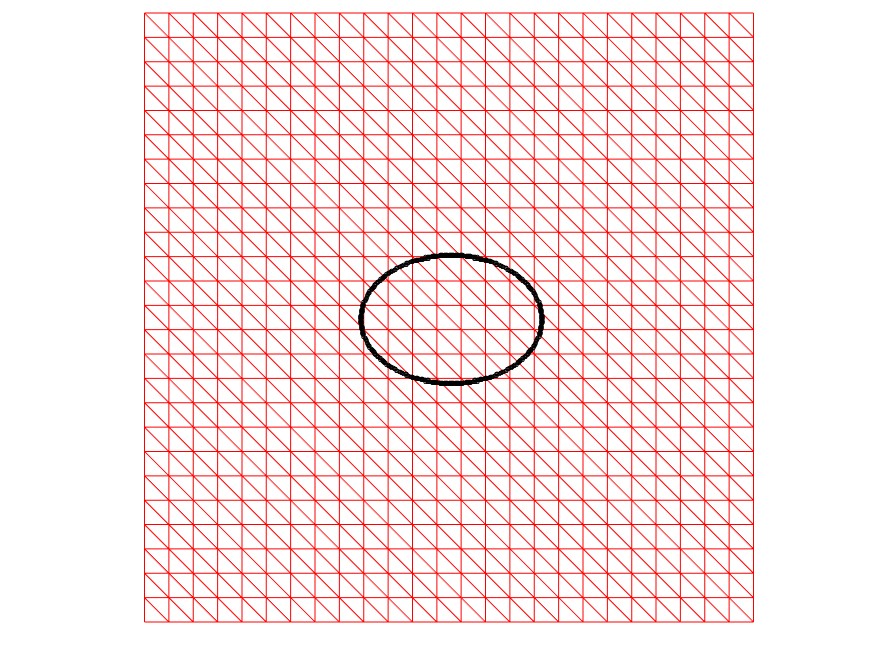}
&\includegraphics[height = \fht, trim = {3.5cm 0cm 3.5cm 0cm}, clip]{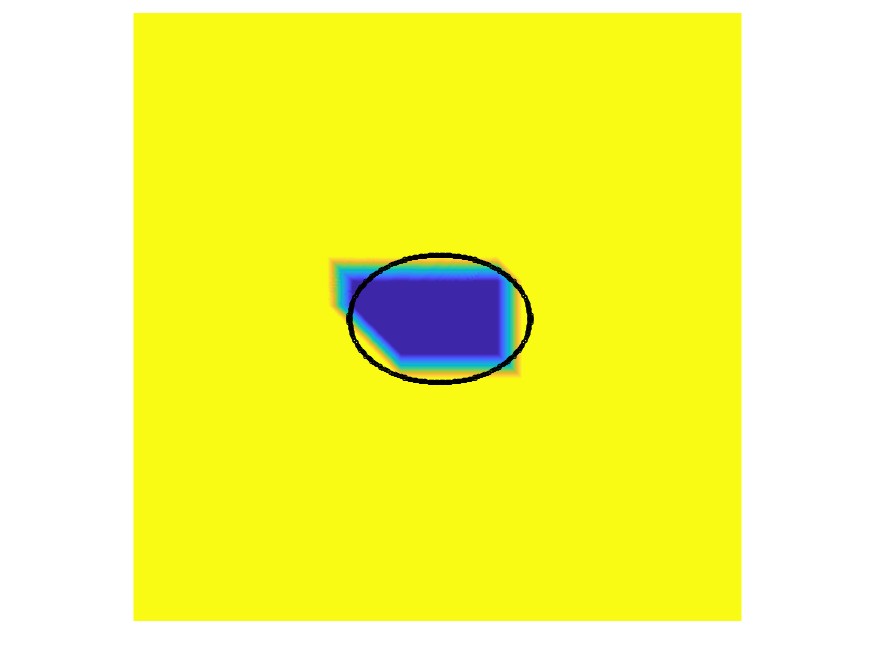}
&\includegraphics[height = \fht, trim = {3.5cm 0cm 3.5cm 0cm}, clip]{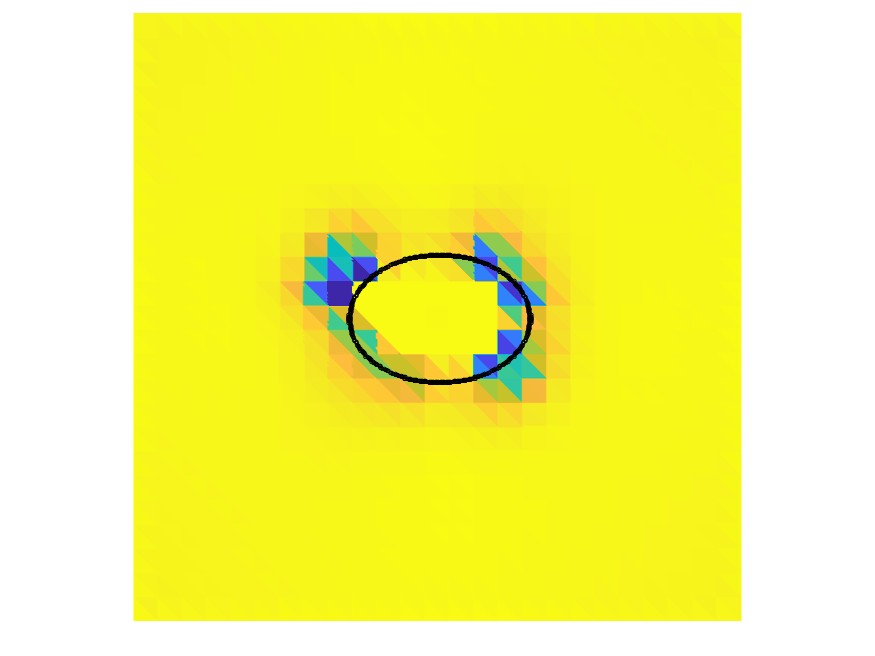}
&\includegraphics[height = \fht, trim = {3.5cm 0cm 3.5cm 0cm}, clip]{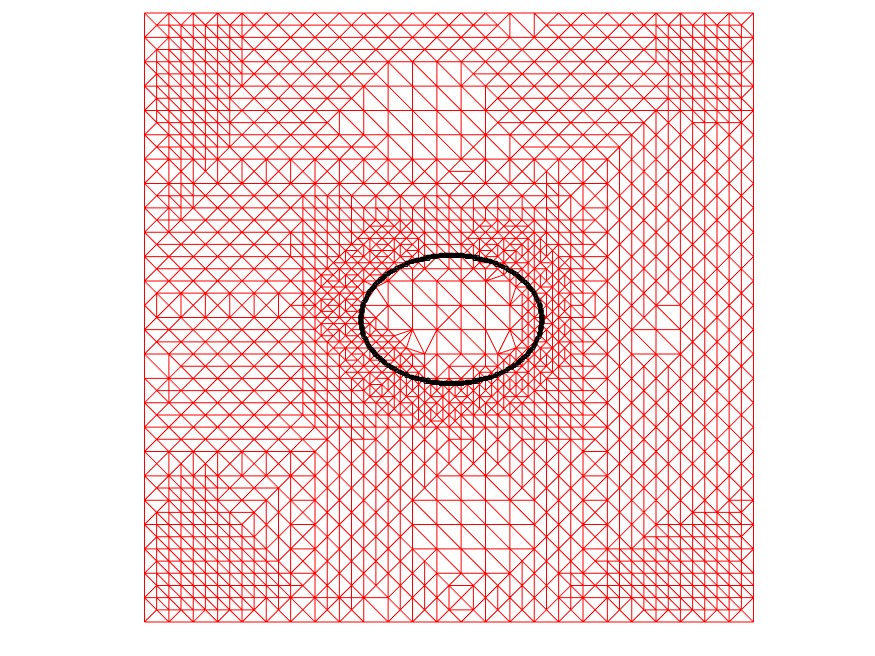}
&\includegraphics[height = \fht, trim = {3.5cm 0cm 3.5cm 0cm}, clip]{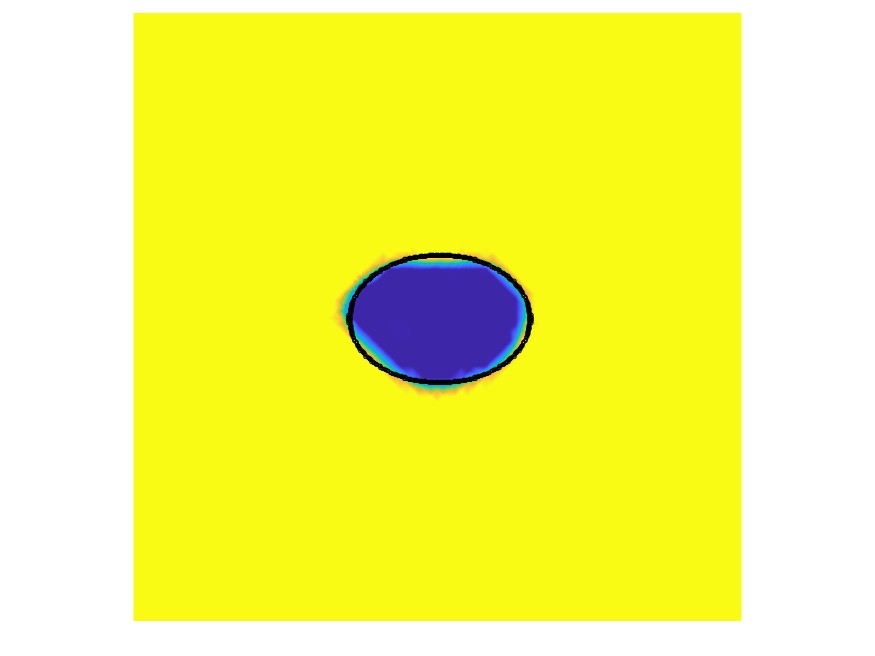}
&\includegraphics[height = \fht, trim = {3.5cm 0cm 3.5cm 0cm}, clip]{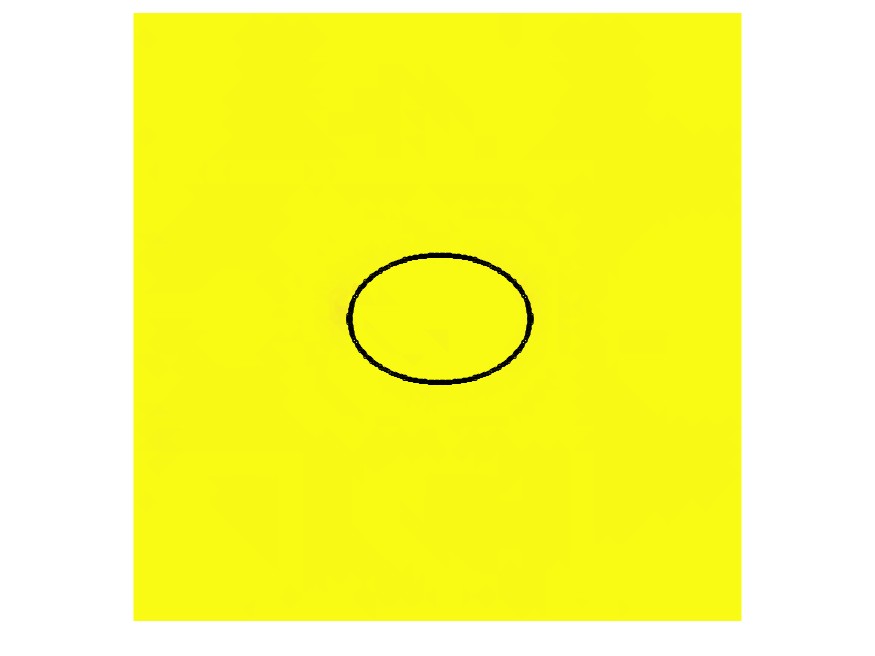}\\
\includegraphics[height = \fht, trim = {3.5cm 0cm 3.5cm 0cm}, clip]{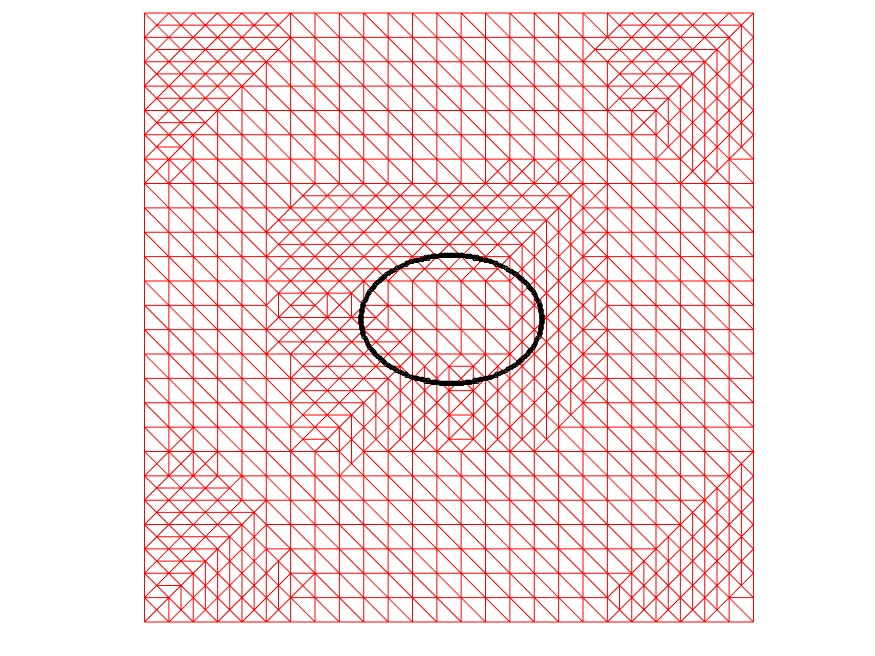}
&\includegraphics[height = \fht, trim = {3.5cm 0cm 3.5cm 0cm}, clip]{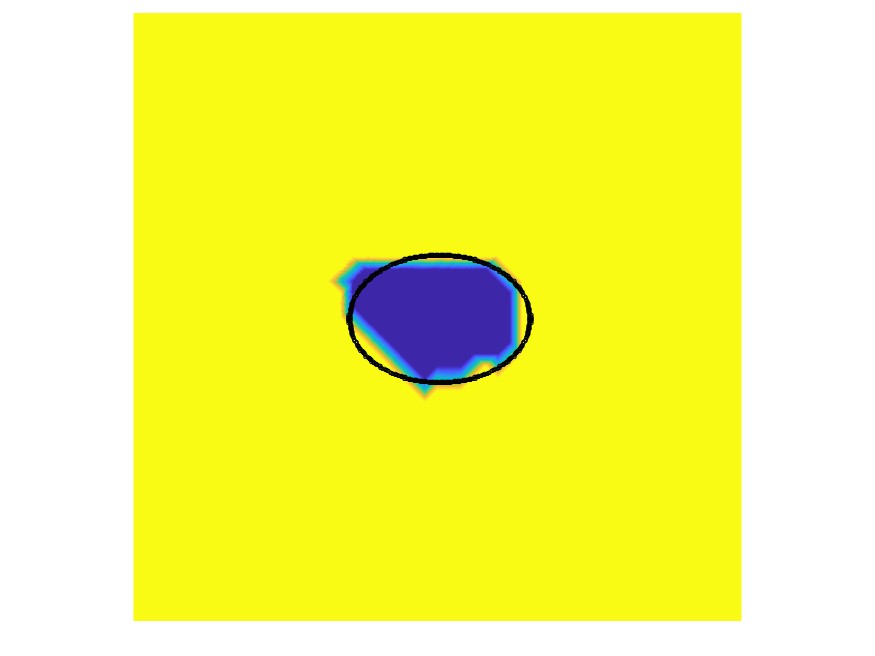}
&\includegraphics[height = \fht, trim = {3.5cm 0cm 3.5cm 0cm}, clip]{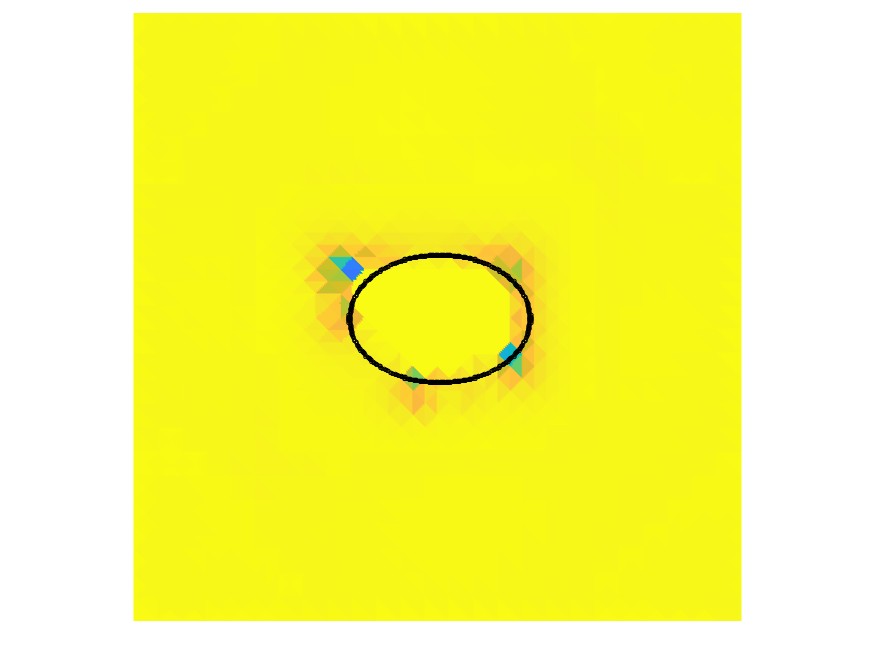}
&\includegraphics[height = \fht, trim = {3.5cm 0cm 3.5cm 0cm}, clip]{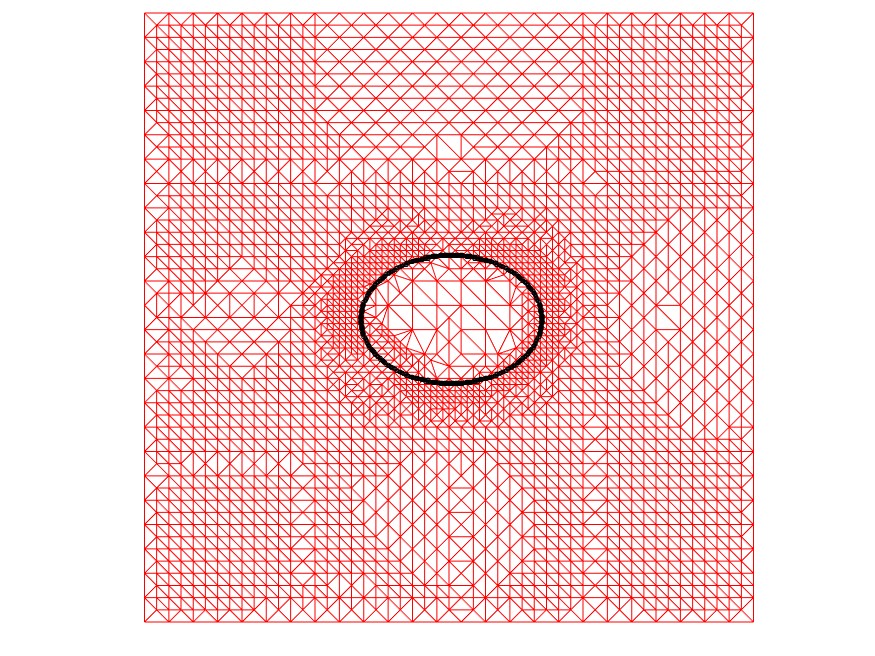}
&\includegraphics[height = \fht, trim = {3.5cm 0cm 3.5cm 0cm}, clip]{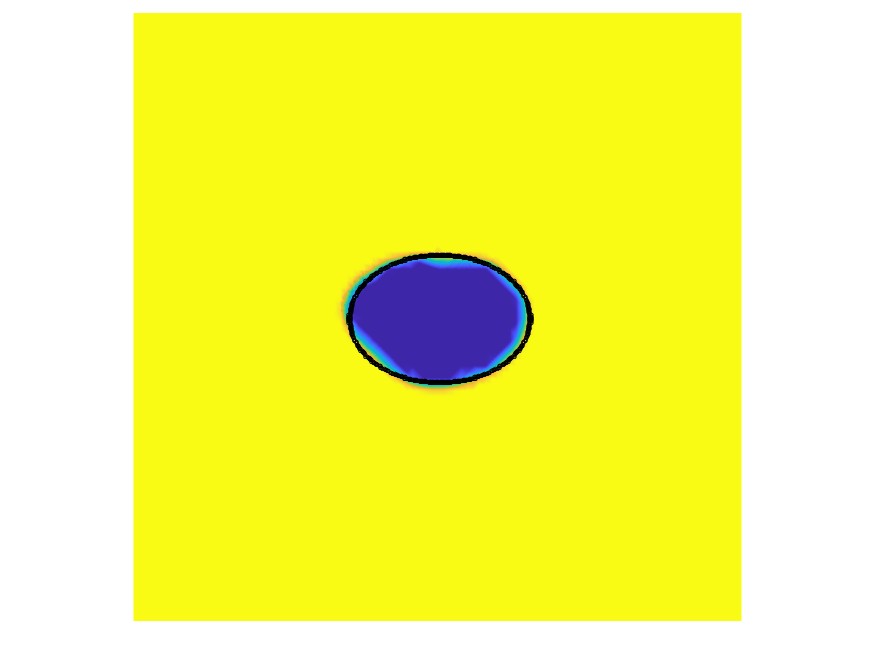}
&\includegraphics[height = \fht, trim = {3.5cm 0cm 3.5cm 0cm}, clip]{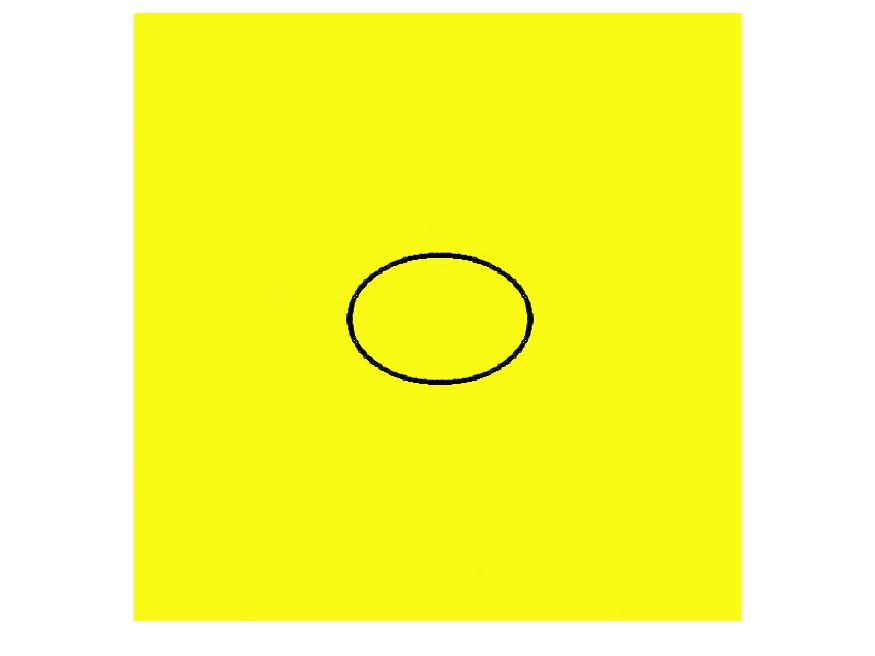}\\
\includegraphics[height = \fht, trim = {3.5cm 0cm 3.5cm 0cm}, clip]{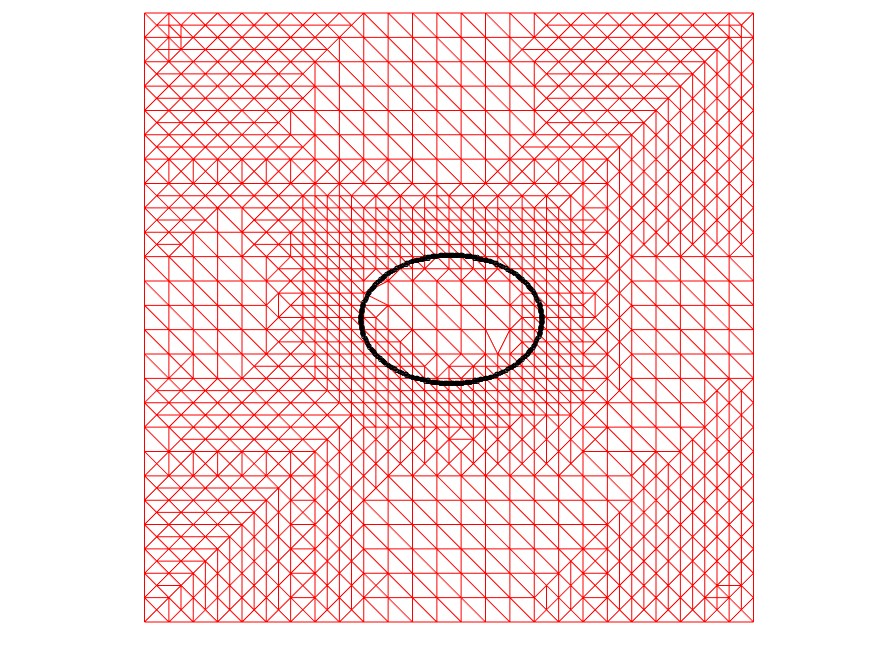}
&\includegraphics[height = \fht, trim = {3.5cm 0cm 3.5cm 0cm}, clip]{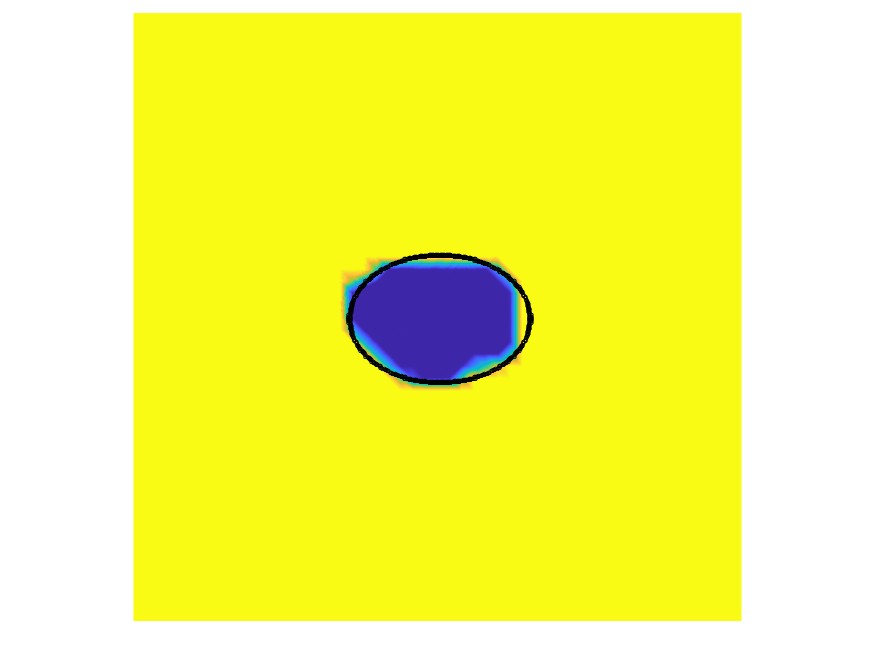}
&\includegraphics[height = \fht, trim = {3.5cm 0cm 3.5cm 0cm}, clip]{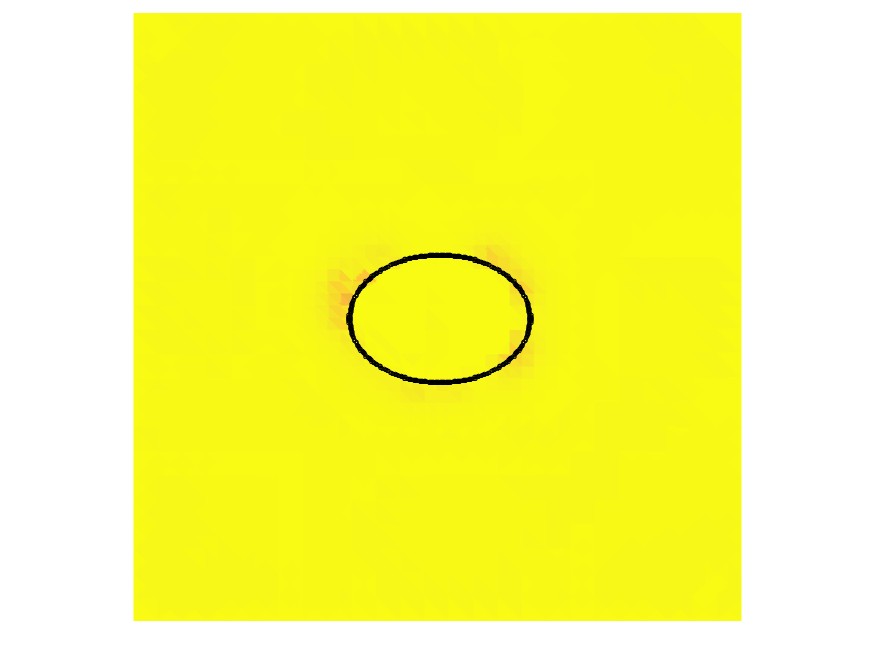}
&\includegraphics[height = \fht, trim = {3.5cm 0cm 3.5cm 0cm}, clip]{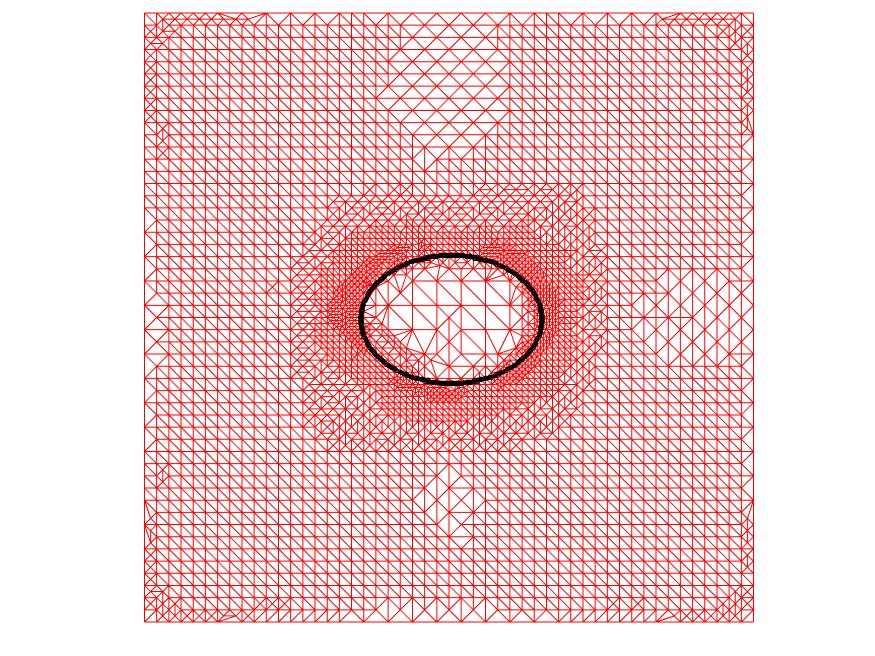}
&\includegraphics[height = \fht, trim = {3.5cm 0cm 3.5cm 0cm}, clip]{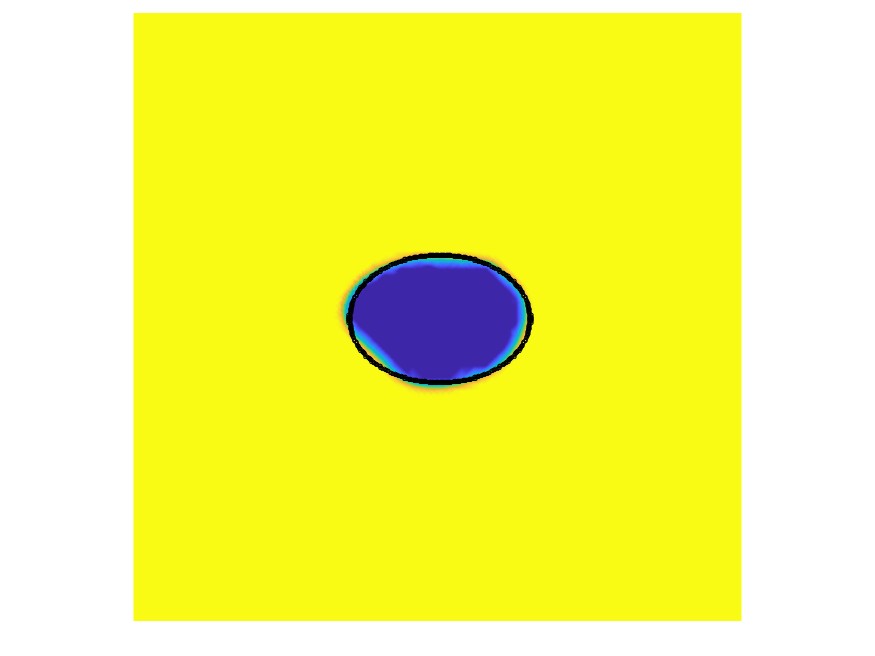}
&\includegraphics[height = \fht, trim = {3.5cm 0cm 3.5cm 0cm}, clip]{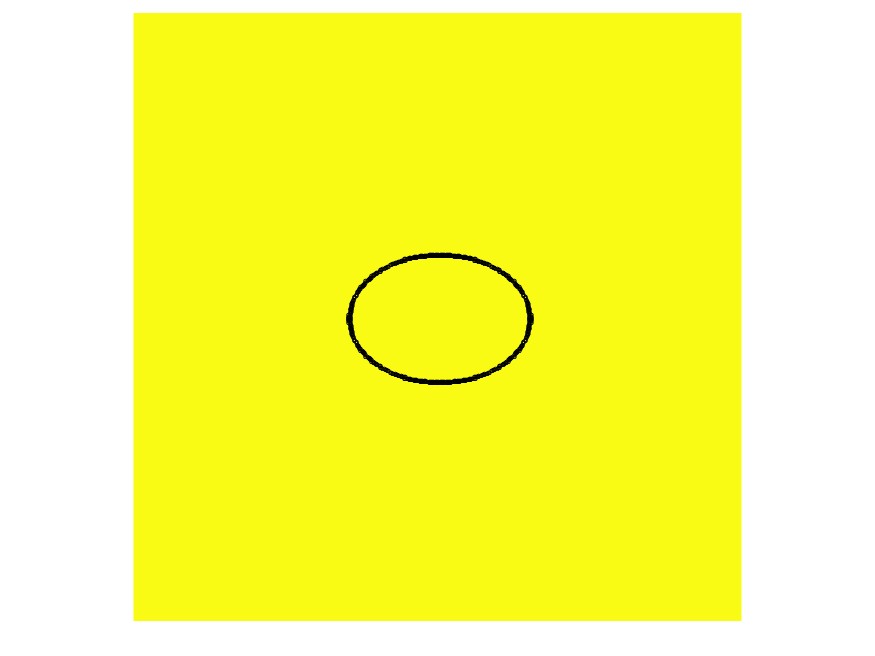}
\end{tabular}
\caption{
The results by the adaptive method for the noisy data $6^\delta$.
From left to the right are the mesh, recovered inclusion and error indicator function.
The number of nodes for each step is
676, 939, 1303, 1823, 2557 and 3586.
}
\label{fig:ellipnoiseadaptive}
\end{figure}

\begin{figure}[hbt!]
\centering
\setlength{\tabcolsep}{0pt}
\begin{tabular}{cccccc}
\includegraphics[height = \fht, trim = {3.5cm 0cm 3.5cm 0cm}, clip]{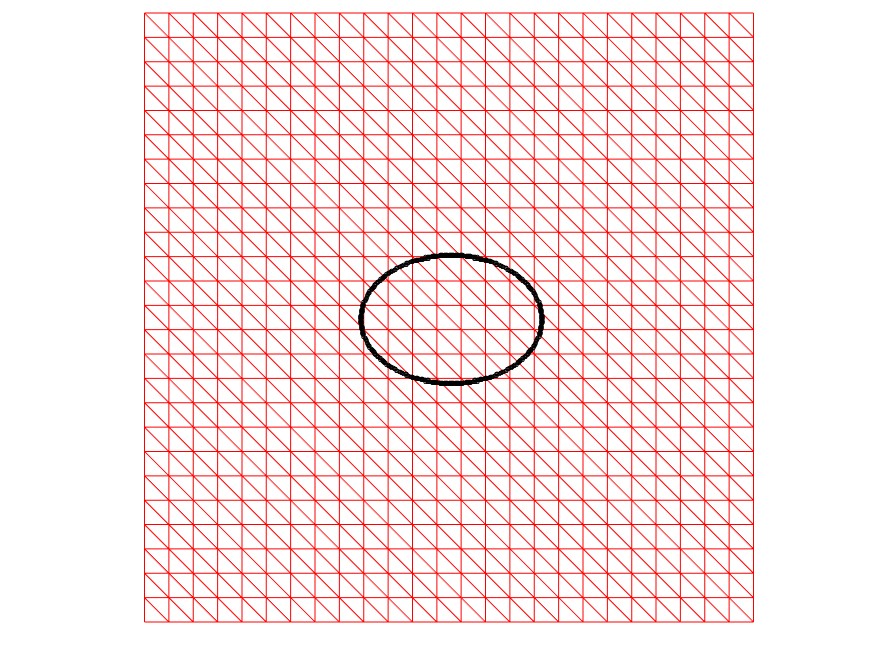}
&\includegraphics[height = \fht, trim = {3.5cm 0cm 3.5cm 0cm}, clip]{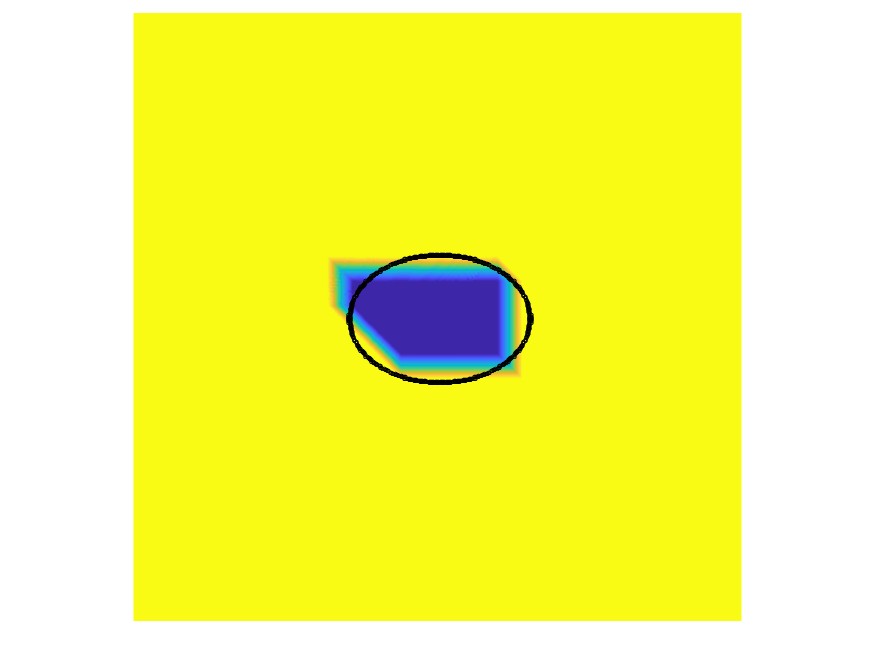}
&\includegraphics[height = \fht, trim = {3.5cm 0cm 3.5cm 0cm}, clip]{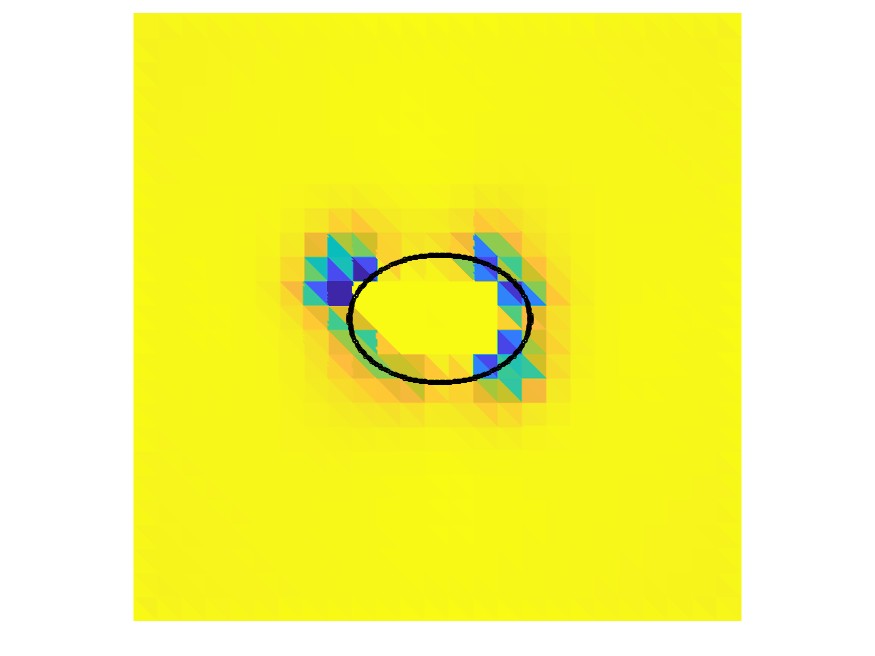}
&\includegraphics[height = \fht, trim = {3.5cm 0cm 3.5cm 0cm}, clip]{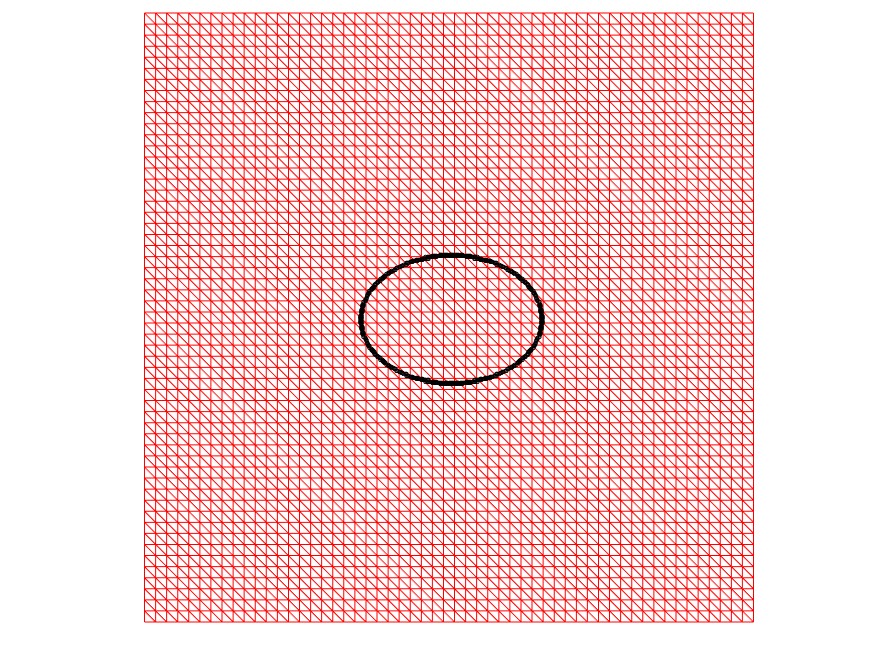}
&\includegraphics[height = \fht, trim = {3.5cm 0cm 3.5cm 0cm}, clip]{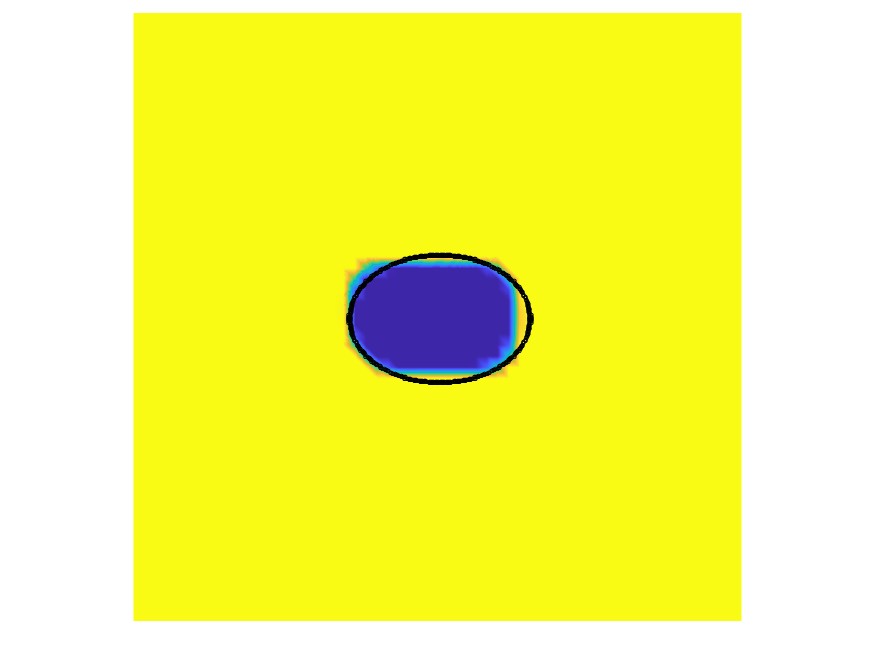}
&\includegraphics[height = \fht, trim = {3.5cm 0cm 3.5cm 0cm}, clip]{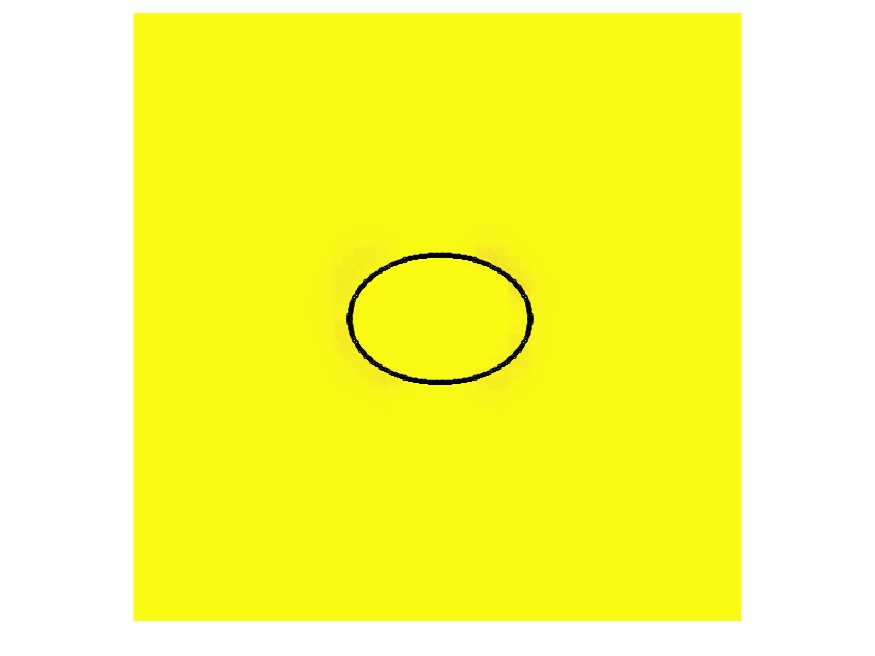}\\
\includegraphics[height = \fht, trim = {3.5cm 0cm 3.5cm 0cm}, clip]{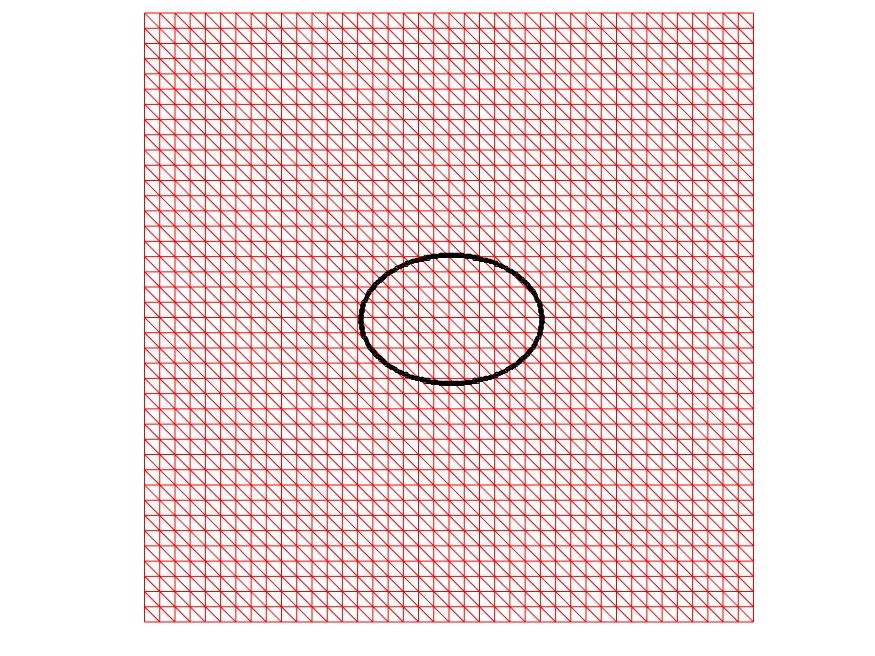}
&\includegraphics[height = \fht, trim = {3.5cm 0cm 3.5cm 0cm}, clip]{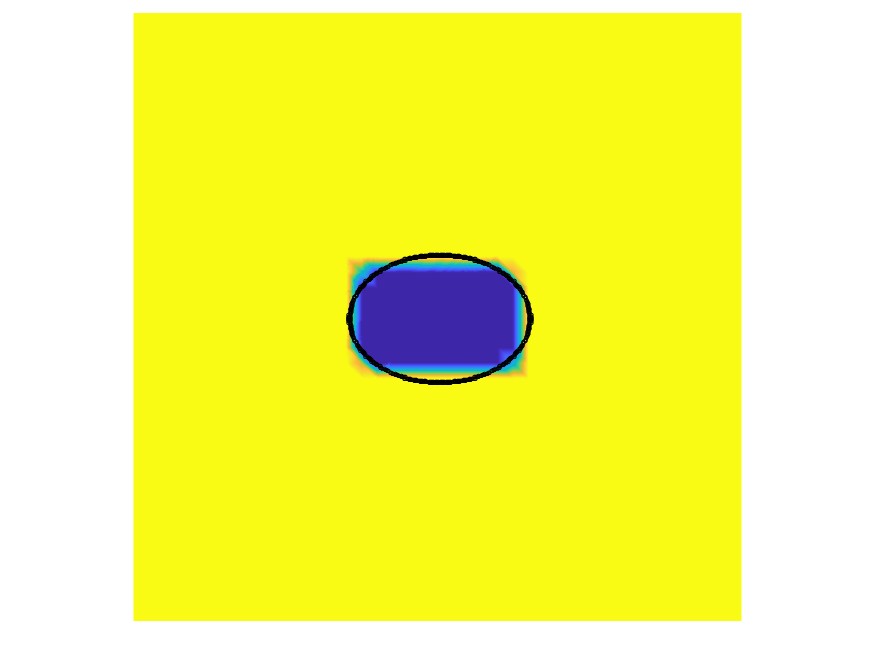}
&\includegraphics[height = \fht, trim = {3.5cm 0cm 3.5cm 0cm}, clip]{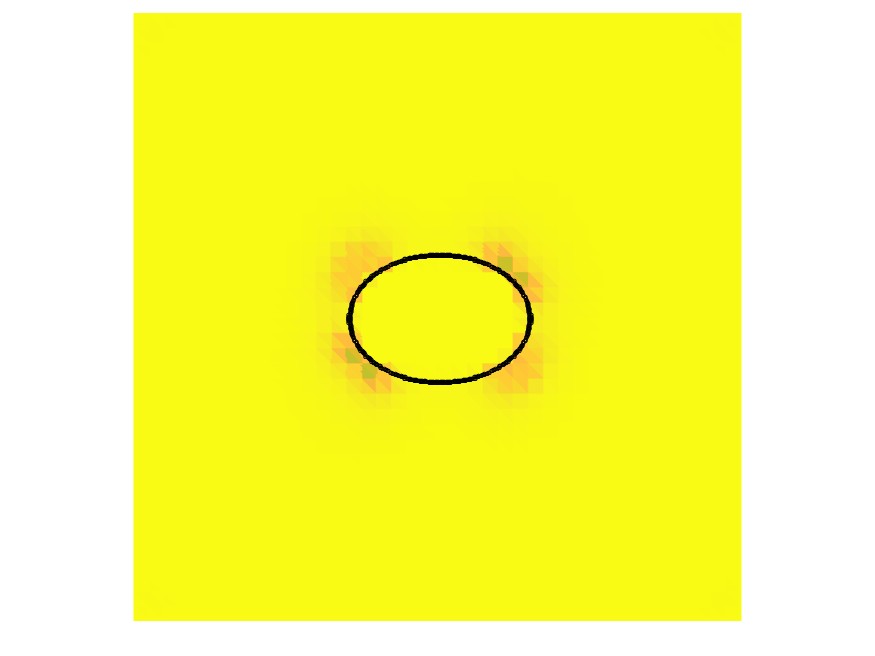}
&\includegraphics[height = \fht, trim = {3.5cm 0cm 3.5cm 0cm}, clip]{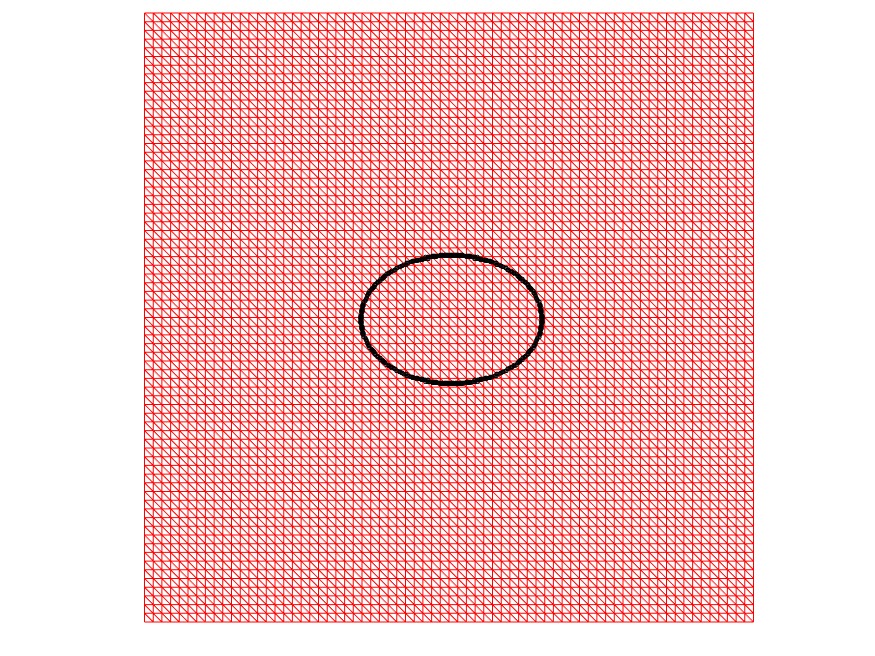}
&\includegraphics[height = \fht, trim = {3.5cm 0cm 3.5cm 0cm}, clip]{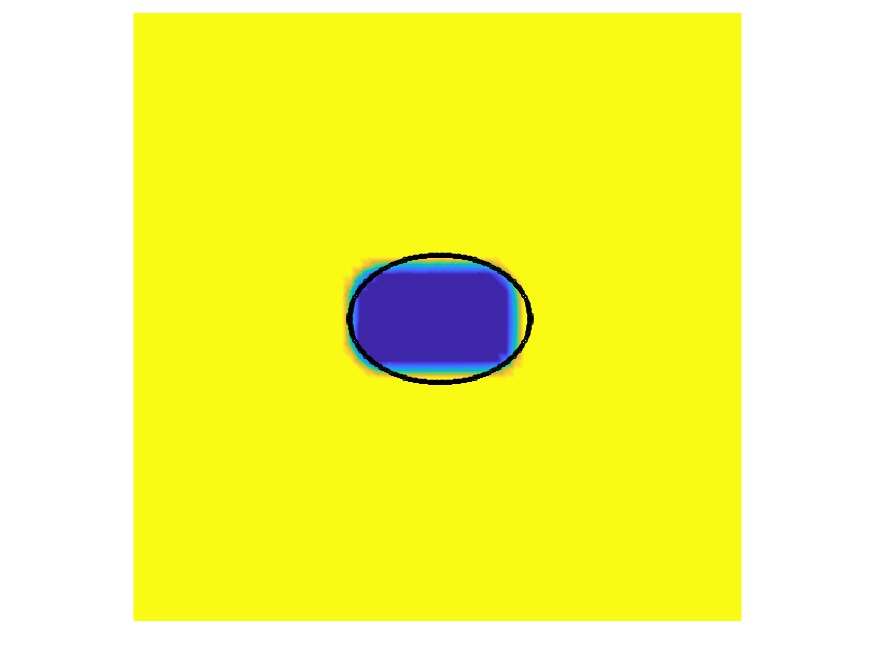}
&\includegraphics[height = \fht, trim = {3.5cm 0cm 3.5cm 0cm}, clip]{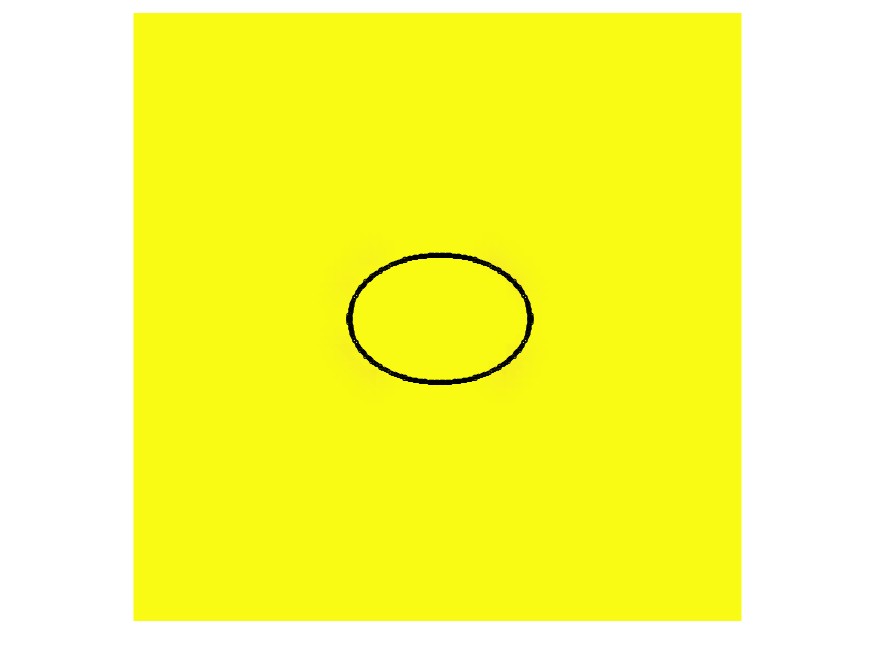}
\end{tabular}
\caption{
The results by uniform mesh refinements for the noisy data $y^\delta$.
From left to the right are figures of mesh, recovered inclusion and indicator function.
The number of nodes at each step is 676, 1681, 3136 and 5041.
}
\label{fig:ellipnoiseuniform}
\end{figure}

In Figs. \ref{fig:ellipnoiseadaptive} and \ref{fig:ellipnoiseuniform}, we present the numerical results for recovering an ellipse inclusion from the noisy data $y^\delta$. The results indicate that on the coarse mesh $\mathcal{T}_0$, the reconstructed shape is not accurate near the interface $\partial\omega$ due to the large discretization errors, and the reconstruction accuracy improves steadily as the mesh refinement proceeds, especially around the interface. The magnitude of the error indicator $\widetilde\eta_k$ always decreases, indicating the convergence of the adaptive algorithm, which agrees with Theorem \ref{thm:conv}.  In the end, the ellipse shape is well resolved by both adaptive and uniform mesh refinements. In sum, these observations resemble closely those for the circular inclusion in Figs. \ref{fig:circlenoiseadaptive} and \ref{fig:circlenoiseuniform}.

\begin{figure}[hbt!]
\centering
\setlength{\tabcolsep}{0pt}
\begin{tabular}{cccccc}
\includegraphics[height = \fht, trim = {3.5cm 0cm 3.5cm 0cm}, clip]{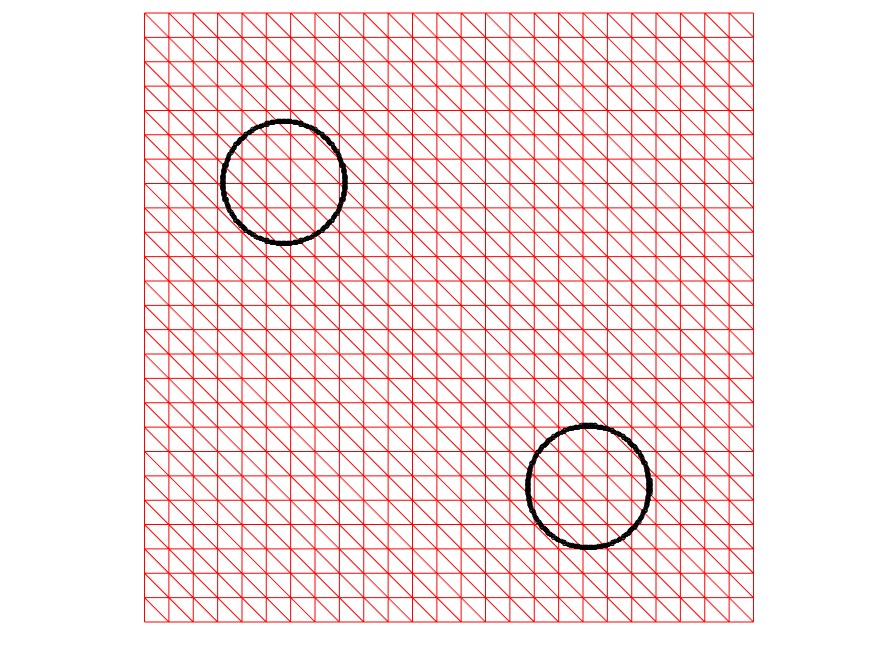}
&\includegraphics[height = \fht, trim = {3.5cm 0cm 3.5cm 0cm}, clip]{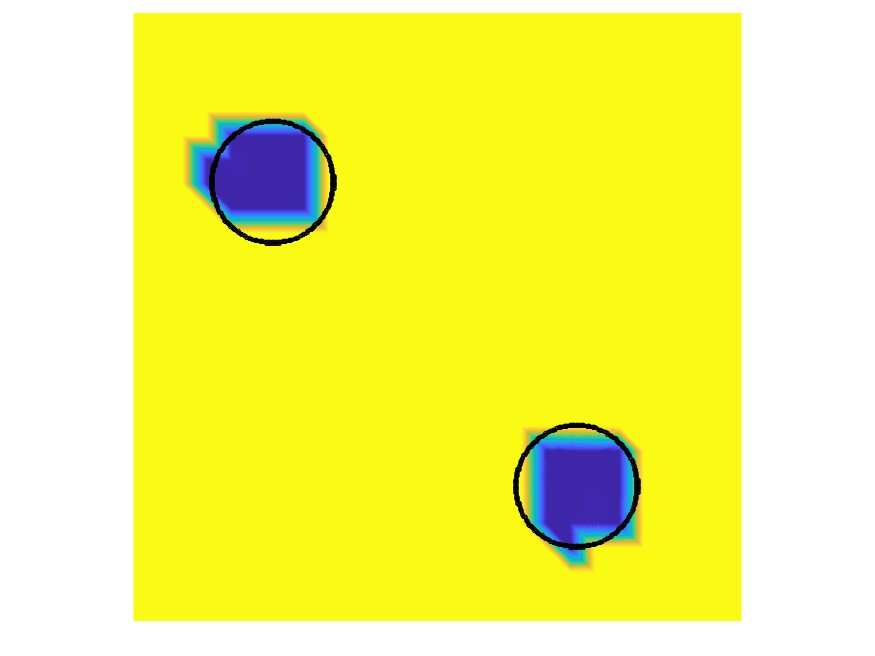}
&\includegraphics[height = \fht, trim = {3.5cm 0cm 3.5cm 0cm}, clip]{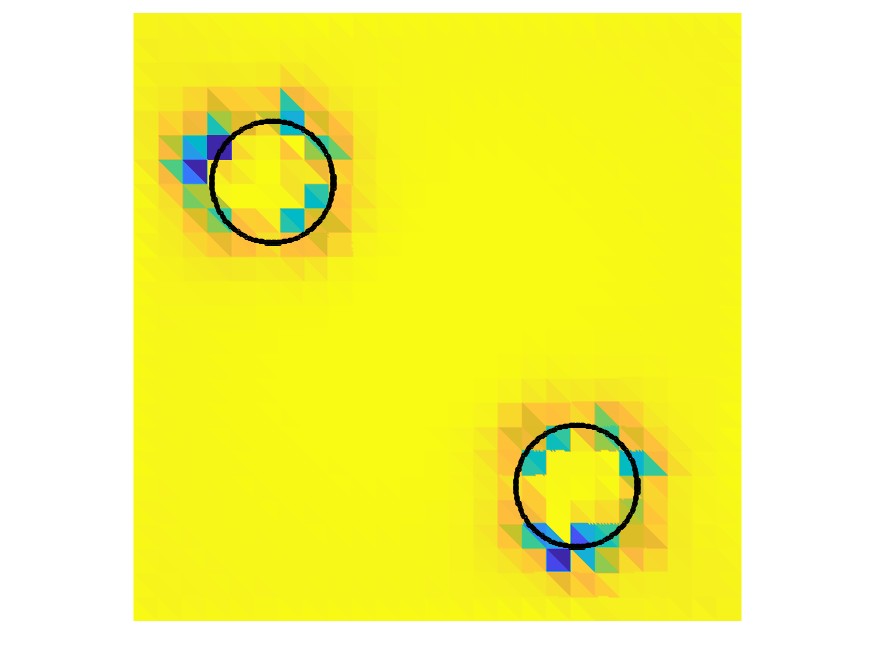}
&\includegraphics[height = \fht, trim = {3.5cm 0cm 3.5cm 0cm}, clip]{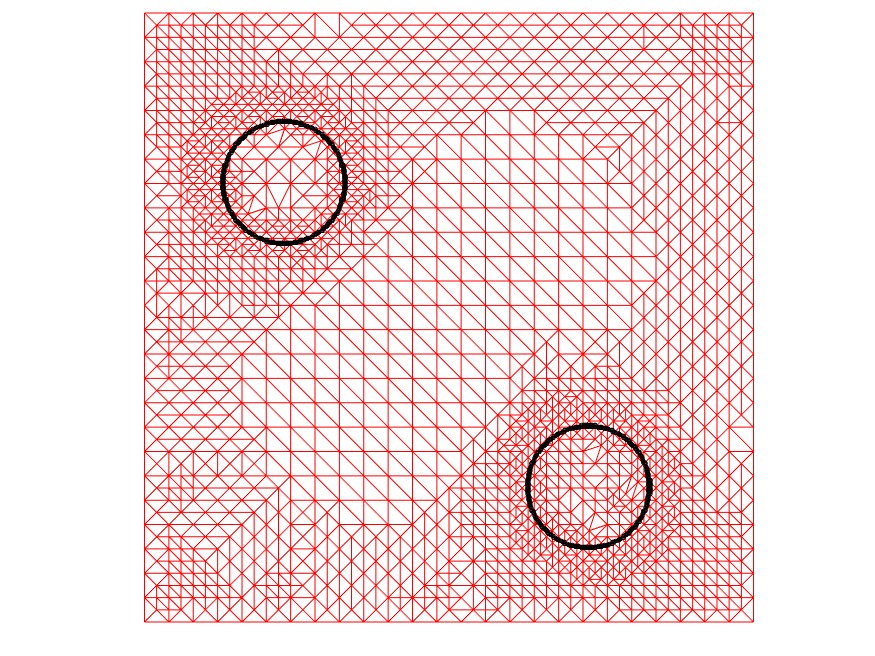}
&\includegraphics[height = \fht, trim = {3.5cm 0cm 3.5cm 0cm}, clip]{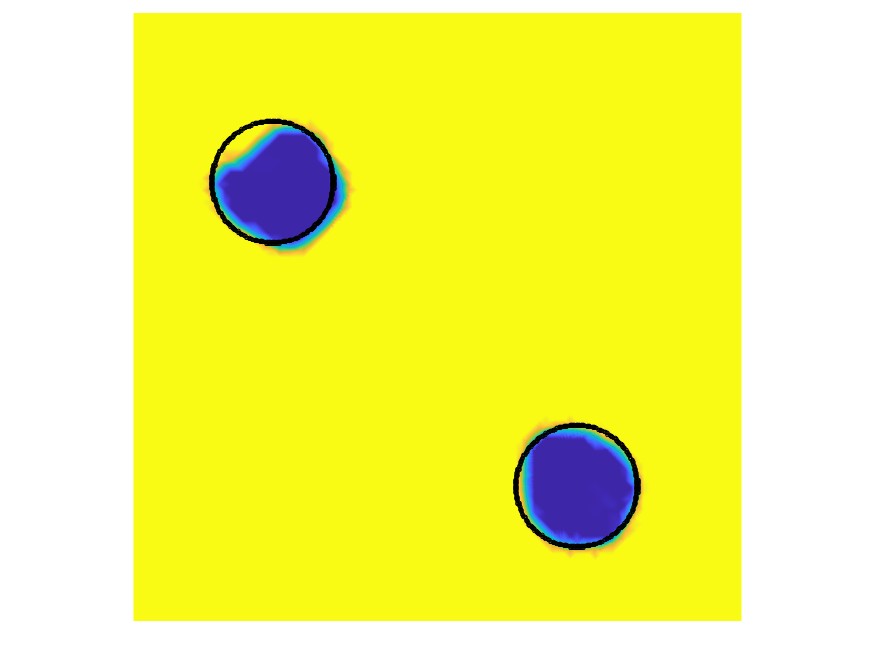}
&\includegraphics[height = \fht, trim = {3.5cm 0cm 3.5cm 0cm}, clip]{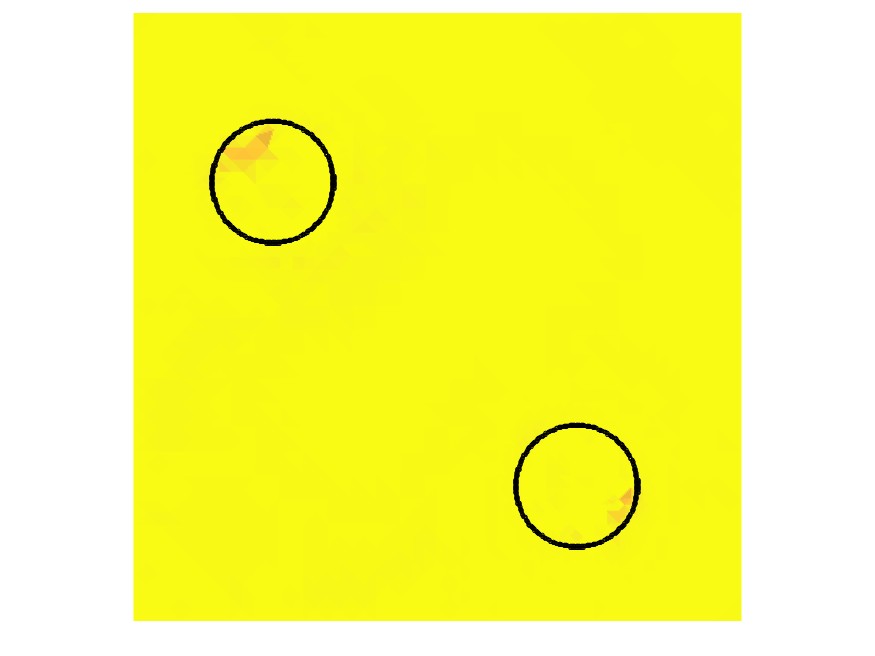}\\
\includegraphics[height = \fht, trim = {3.5cm 0cm 3.5cm 0cm}, clip]{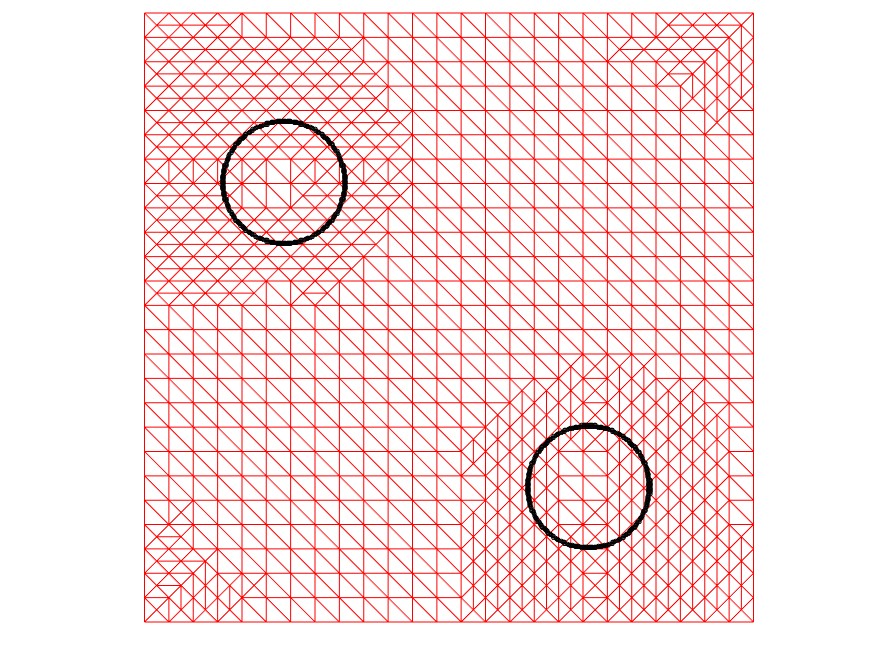}
&\includegraphics[height = \fht, trim = {3.5cm 0cm 3.5cm 0cm}, clip]{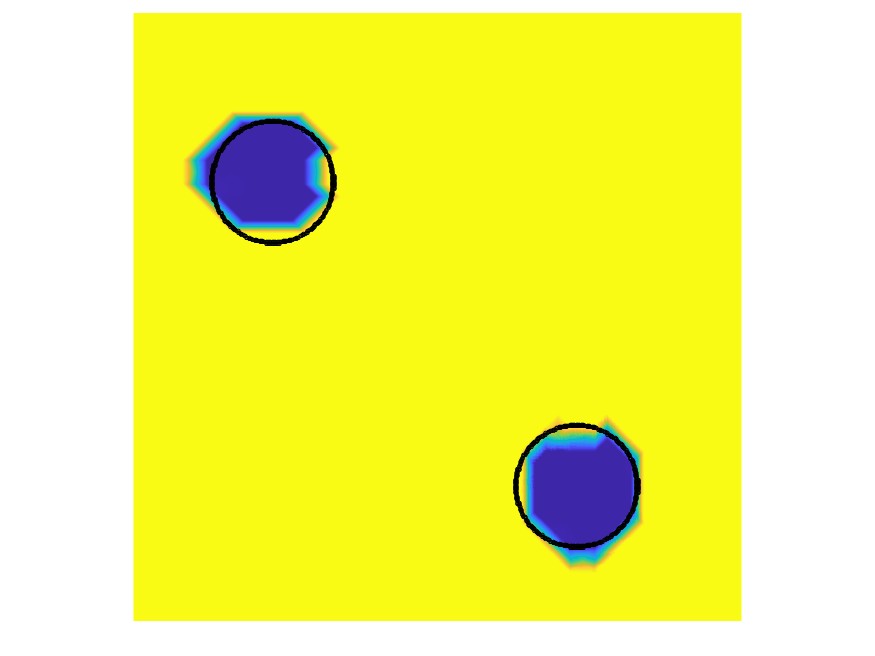}
&\includegraphics[height = \fht, trim = {3.5cm 0cm 3.5cm 0cm}, clip]{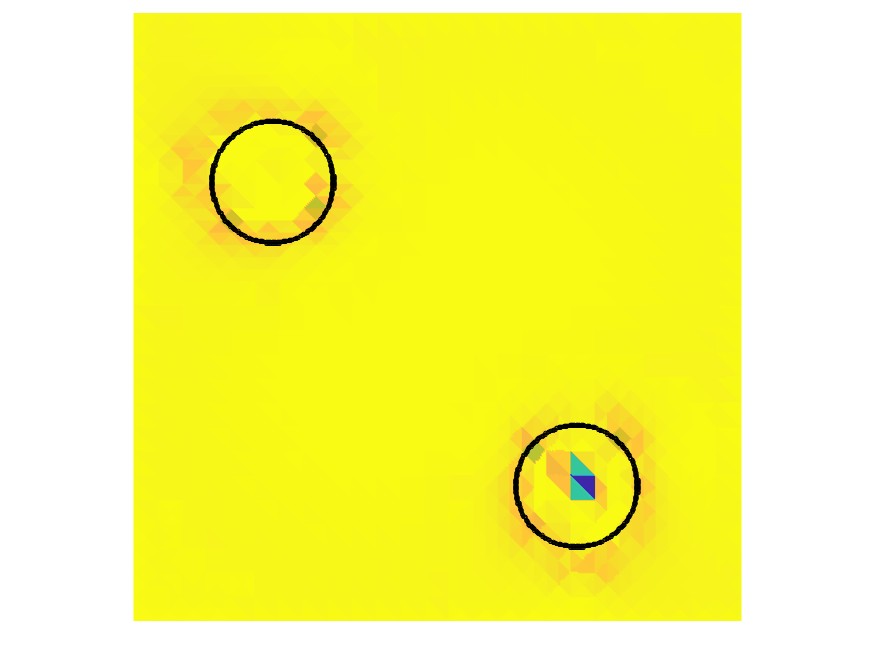}
&\includegraphics[height = \fht, trim = {3.5cm 0cm 3.5cm 0cm}, clip]{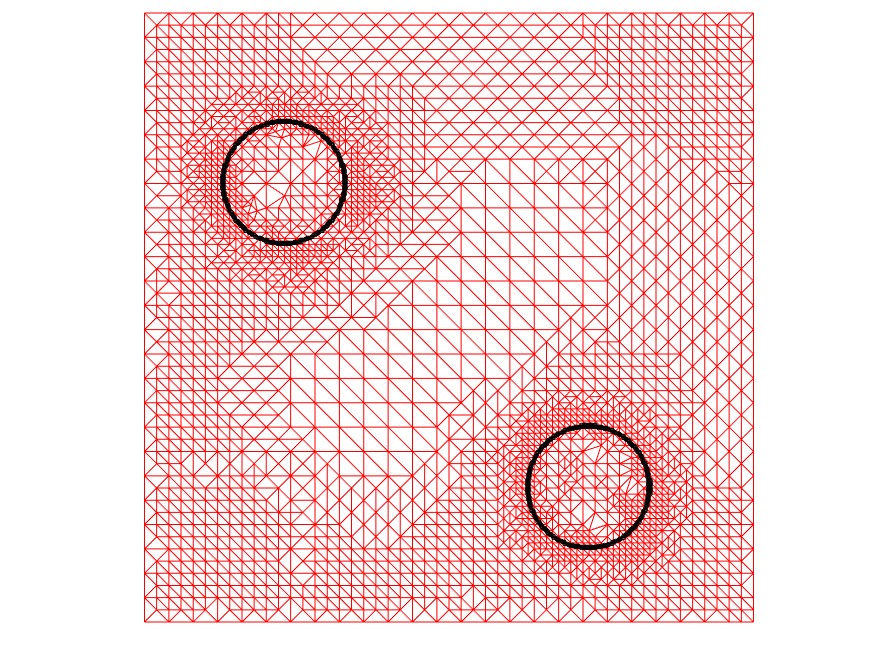}
&\includegraphics[height = \fht, trim = {3.5cm 0cm 3.5cm 0cm}, clip]{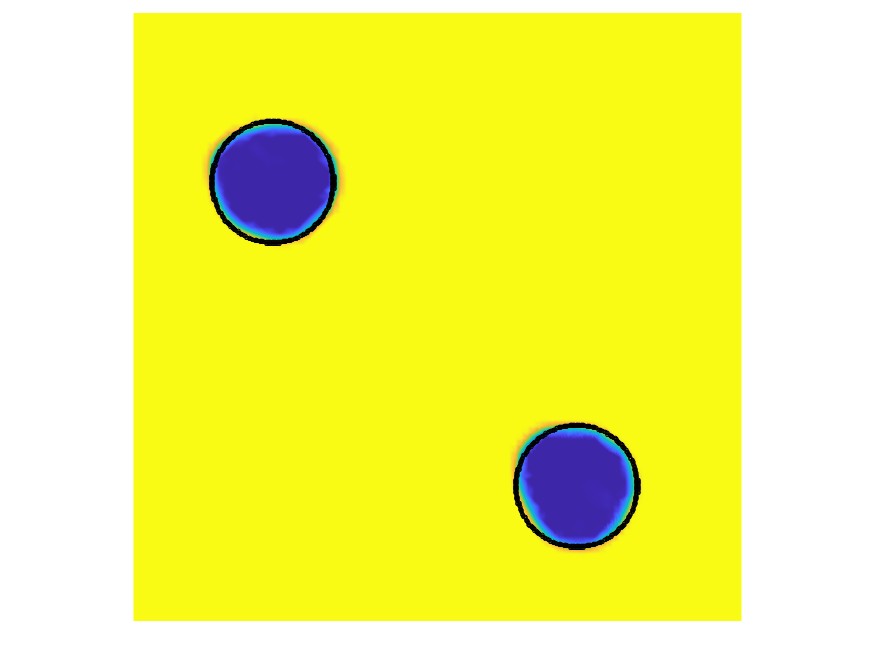}
&\includegraphics[height = \fht, trim = {3.5cm 0cm 3.5cm 0cm}, clip]{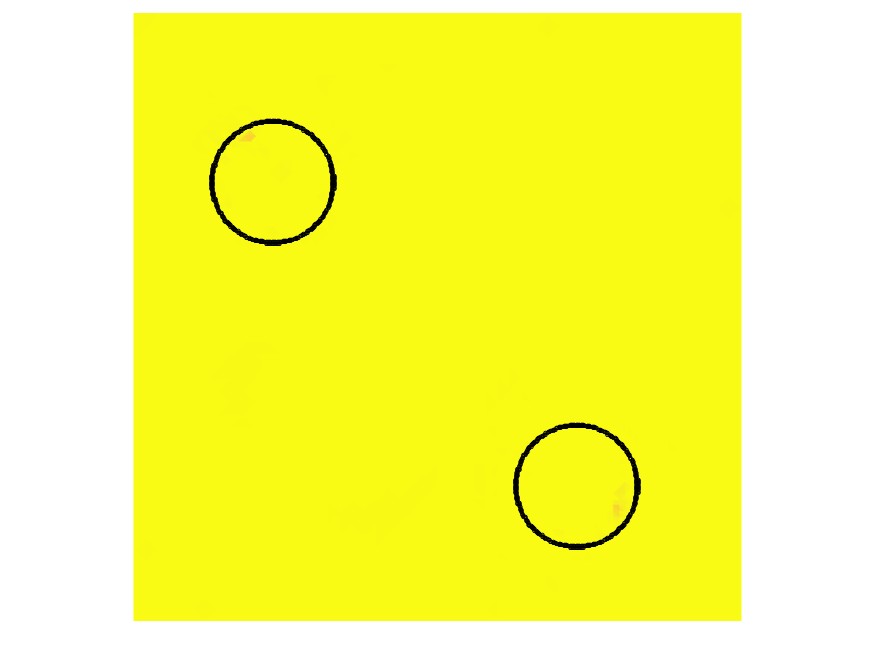}\\
\includegraphics[height = \fht, trim = {3.5cm 0cm 3.5cm 0cm}, clip]{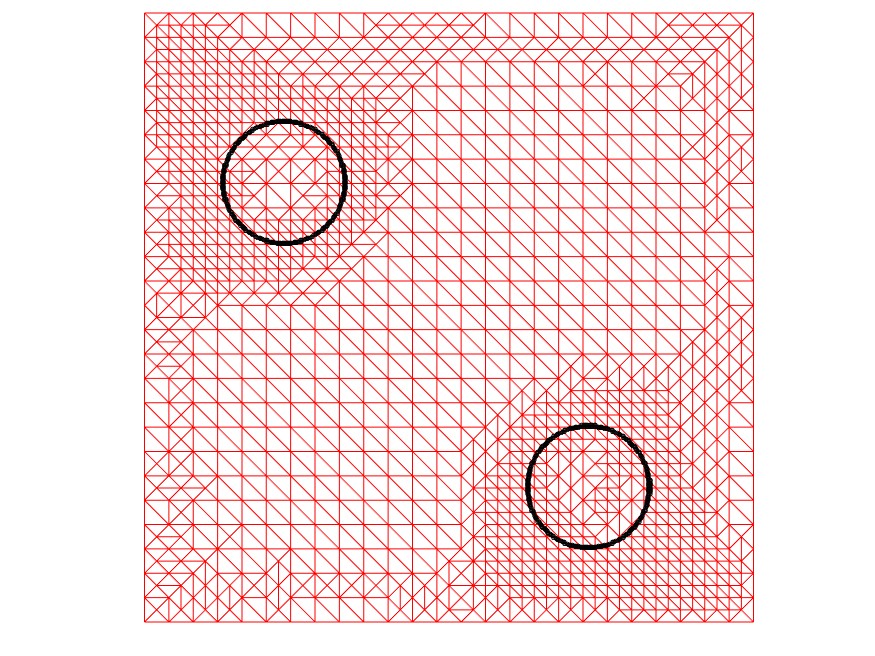}
&\includegraphics[height = \fht, trim = {3.5cm 0cm 3.5cm 0cm}, clip]{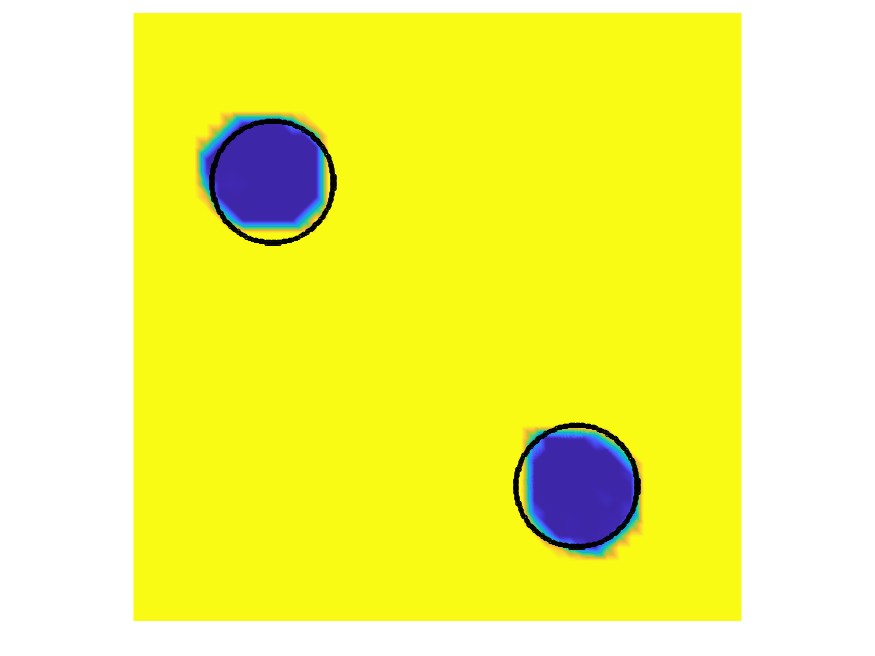}
&\includegraphics[height = \fht, trim = {3.5cm 0cm 3.5cm 0cm}, clip]{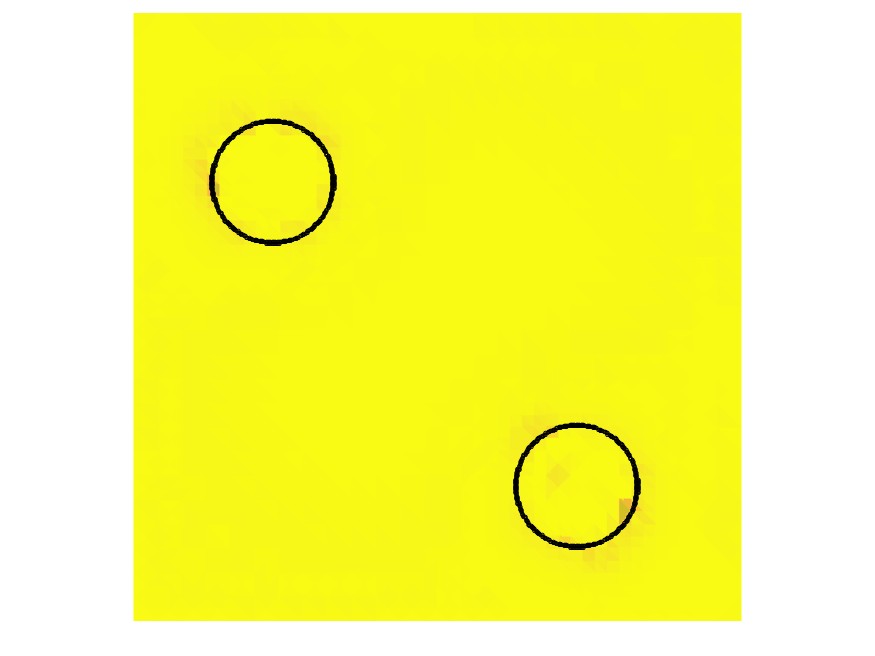}
&\includegraphics[height = \fht, trim = {3.5cm 0cm 3.5cm 0cm}, clip]{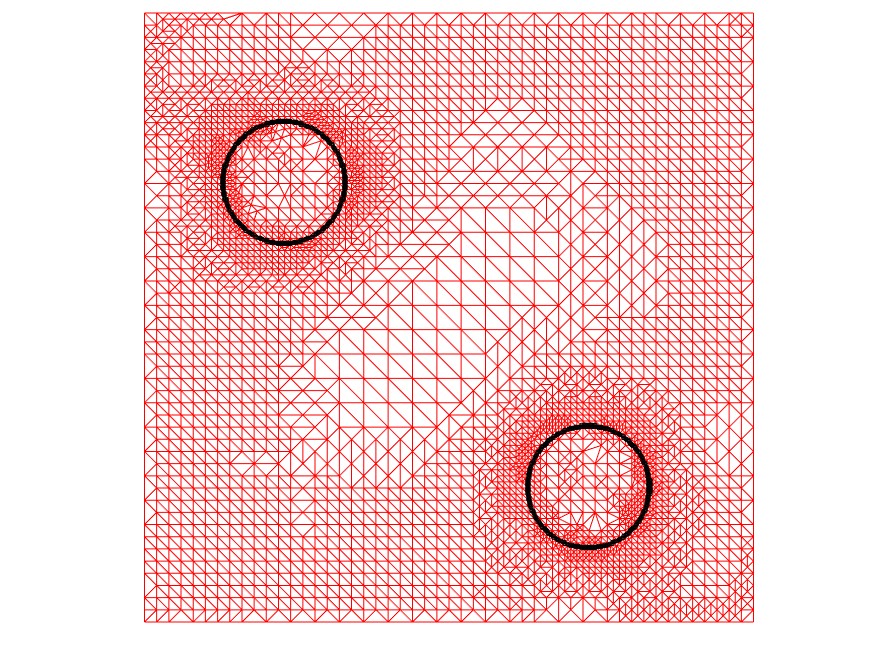}
&\includegraphics[height = \fht, trim = {3.5cm 0cm 3.5cm 0cm}, clip]{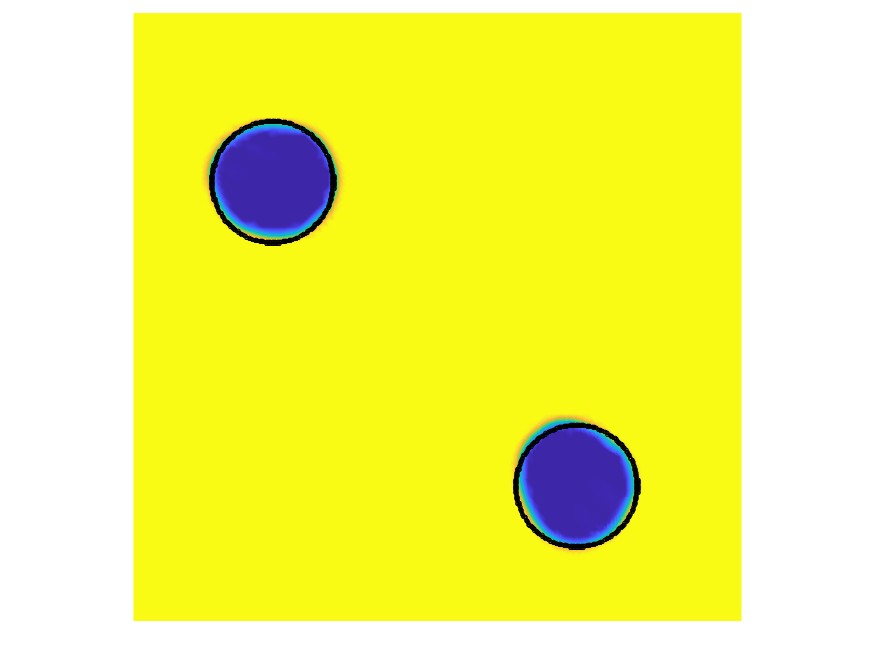}
&\includegraphics[height = \fht, trim = {3.5cm 0cm 3.5cm 0cm}, clip]{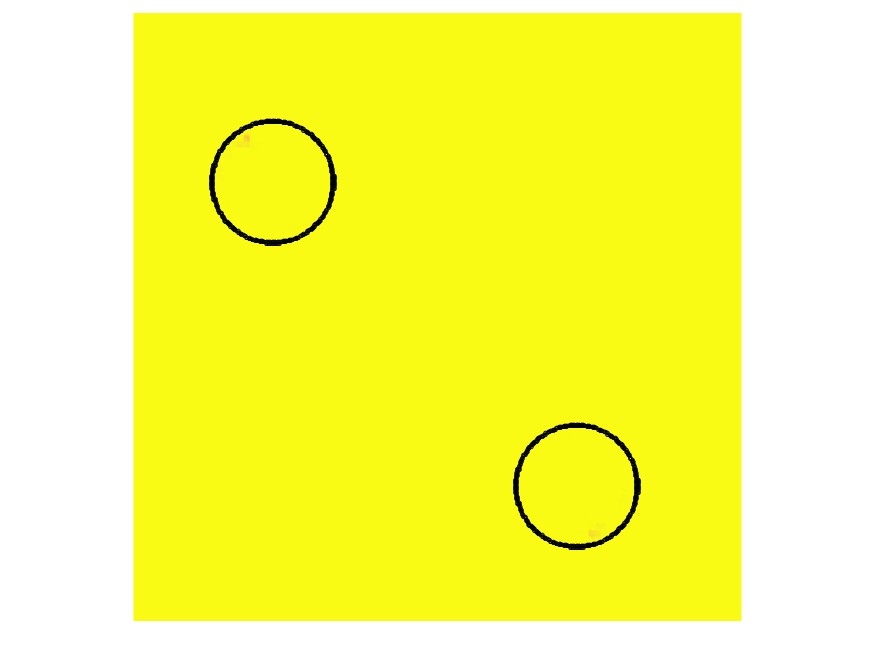}
\end{tabular}
\caption{The results by the adaptive method for the noisy data $y^\delta$.
From left to the right are the mesh, recovered inclusion and error indicator function.
The number of nodes for each step is
676, 927, 1288, 1816, 2559 and 3616.
}
\label{fig:twocirclenoiseadaptive}
\end{figure}

\begin{figure}[hbt!]
\centering
\setlength{\tabcolsep}{0pt}
\begin{tabular}{cccccc}
\includegraphics[height = \fht, trim = {3.5cm 0cm 3.5cm 0cm}, clip]{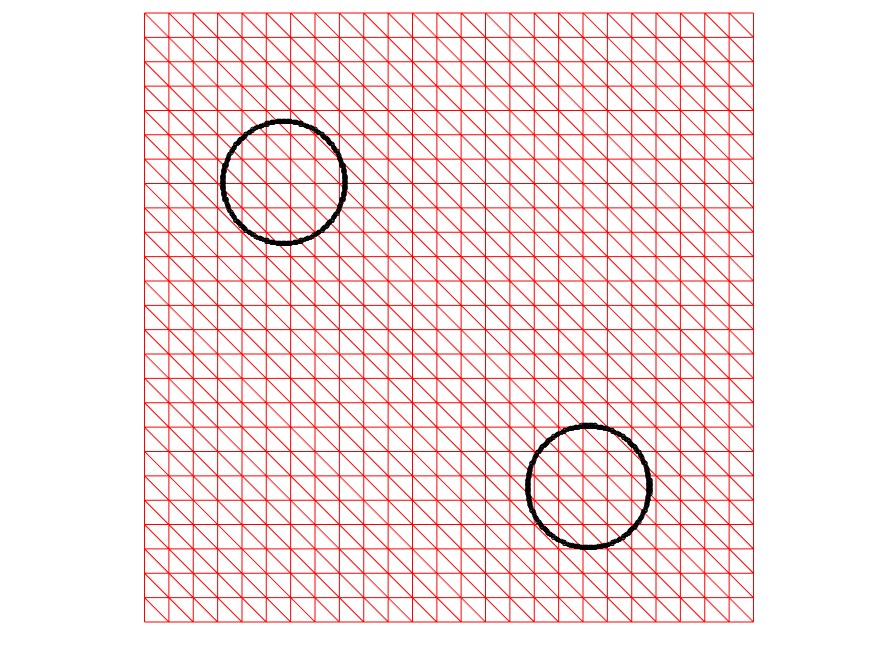}
&\includegraphics[height = \fht, trim = {3.5cm 0cm 3.5cm 0cm}, clip]{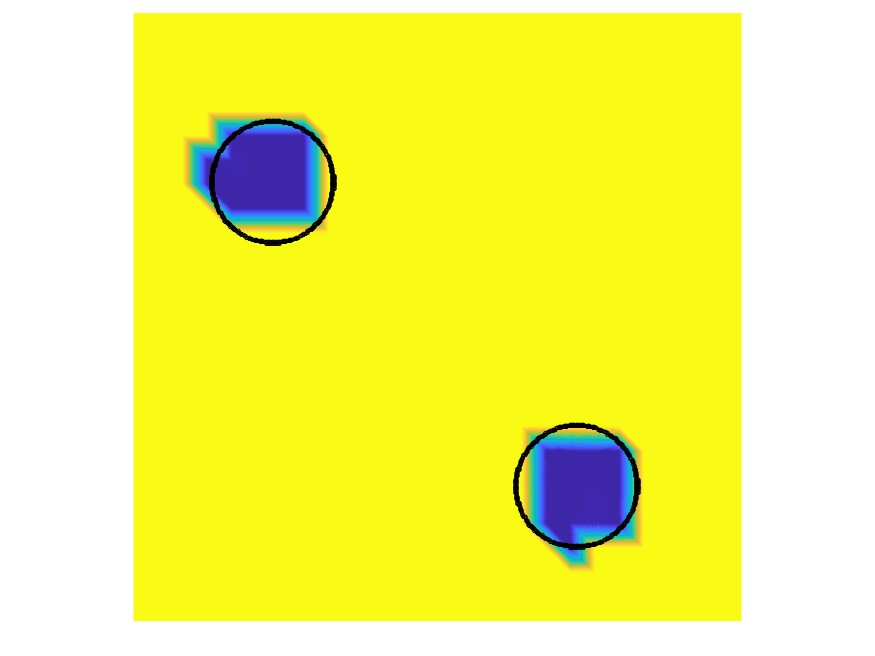}
&\includegraphics[height = \fht, trim = {3.5cm 0cm 3.5cm 0cm}, clip]{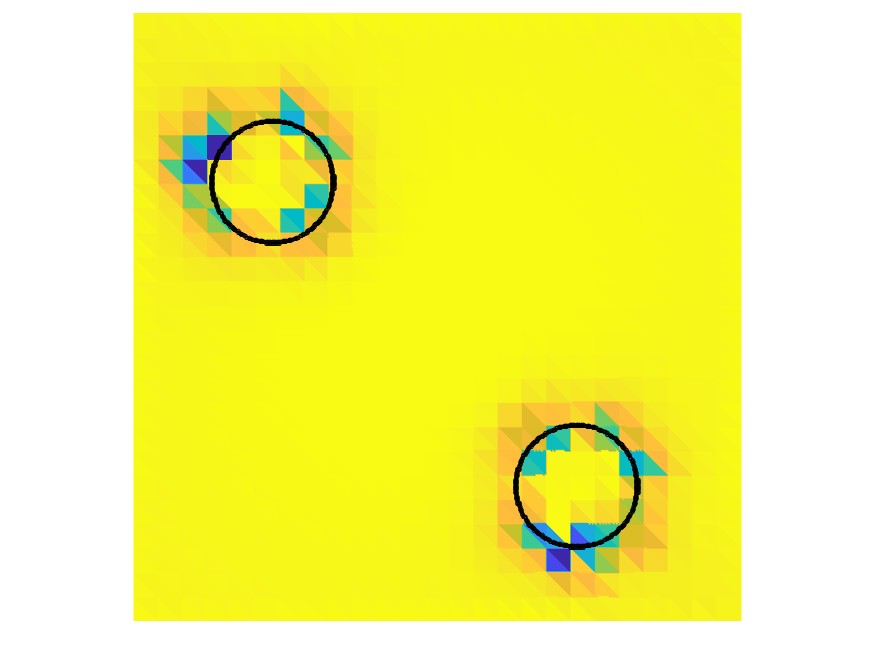}
&\includegraphics[height = \fht, trim = {3.5cm 0cm 3.5cm 0cm}, clip]{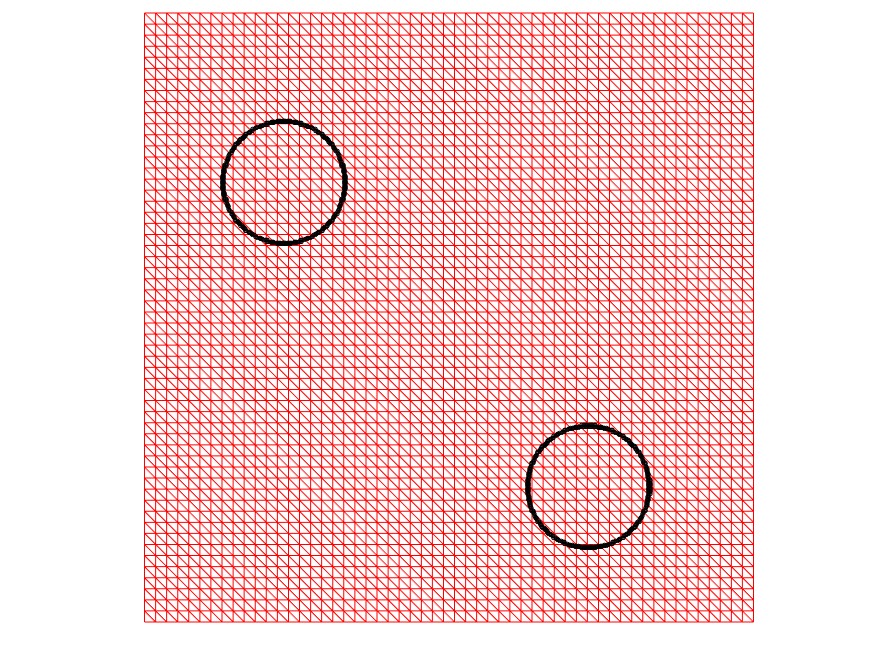}
&\includegraphics[height = \fht, trim = {3.5cm 0cm 3.5cm 0cm}, clip]{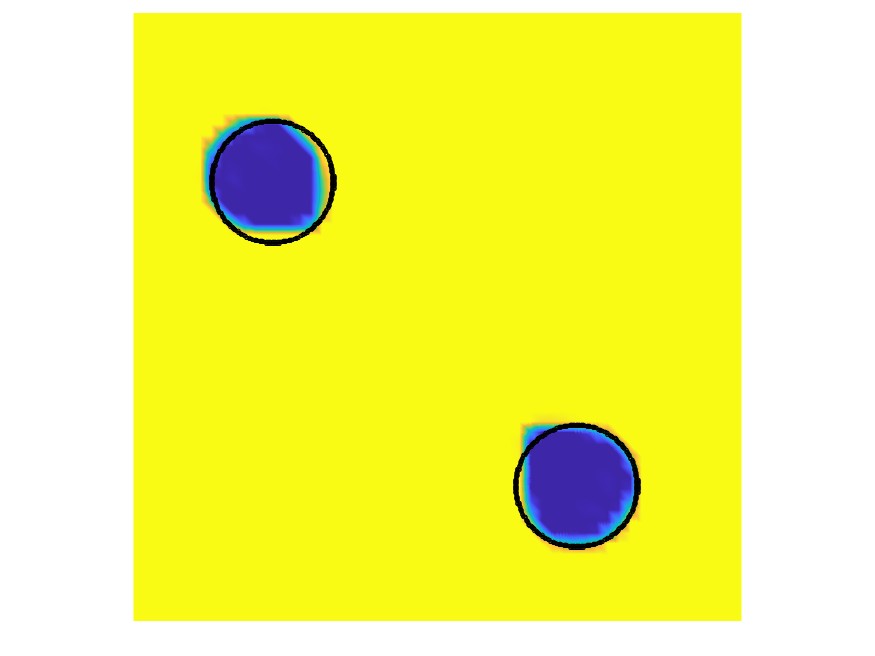}
&\includegraphics[height = \fht, trim = {3.5cm 0cm 3.5cm 0cm}, clip]{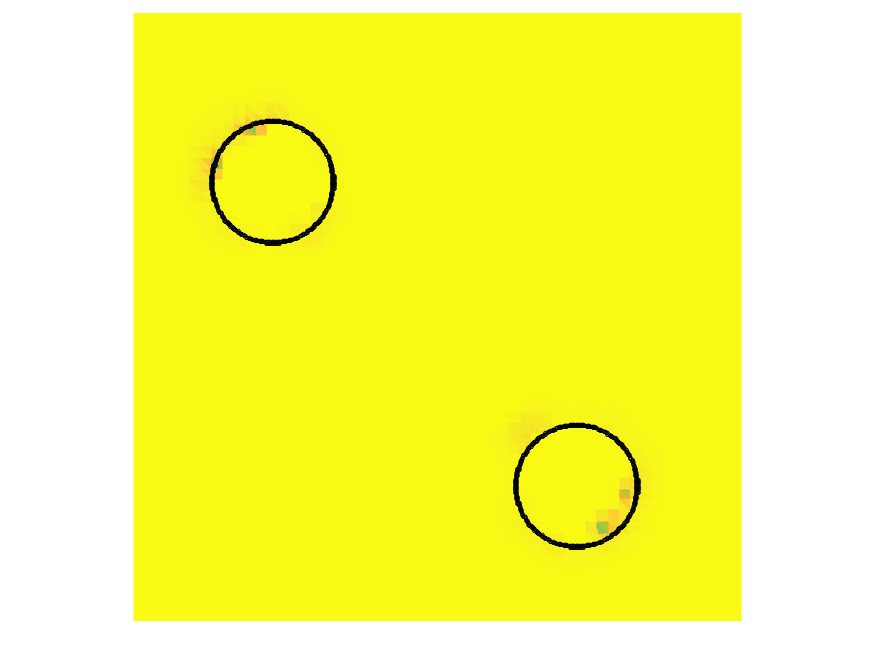}\\
\includegraphics[height = \fht, trim = {3.5cm 0cm 3.5cm 0cm}, clip]{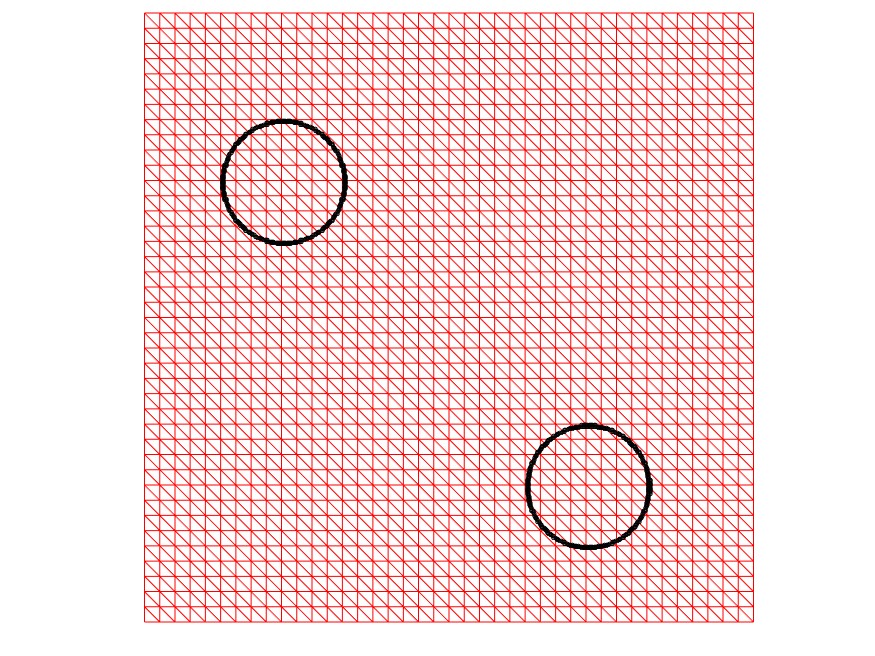}
&\includegraphics[height = \fht, trim = {3.5cm 0cm 3.5cm 0cm}, clip]{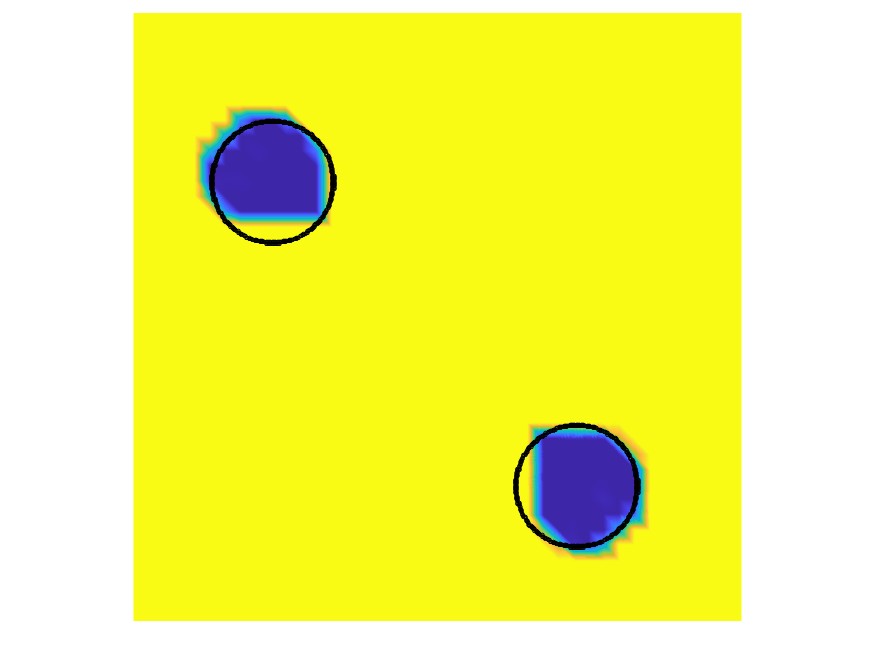}
&\includegraphics[height = \fht, trim = {3.5cm 0cm 3.5cm 0cm}, clip]{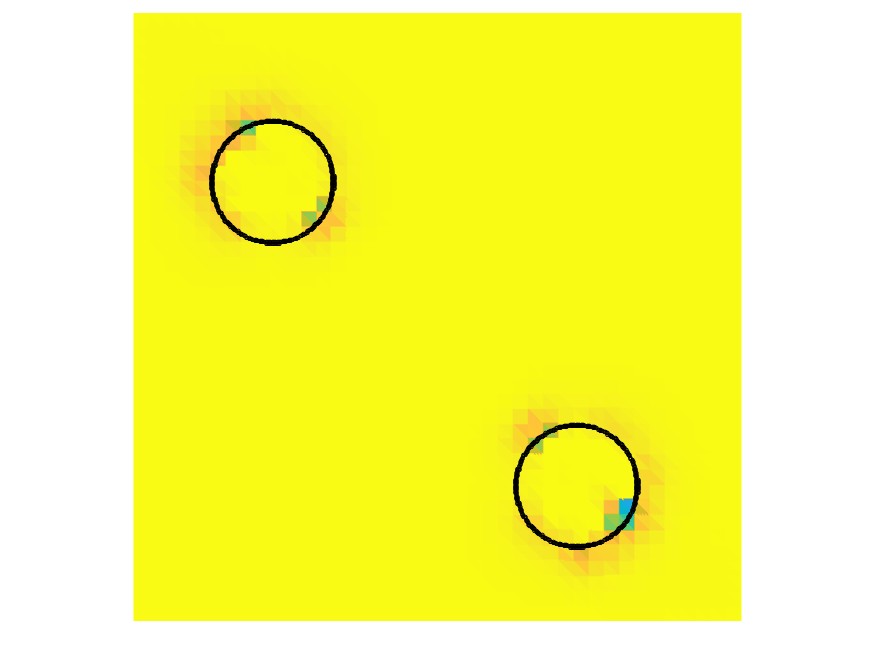}
&\includegraphics[height = \fht, trim = {3.5cm 0cm 3.5cm 0cm}, clip]{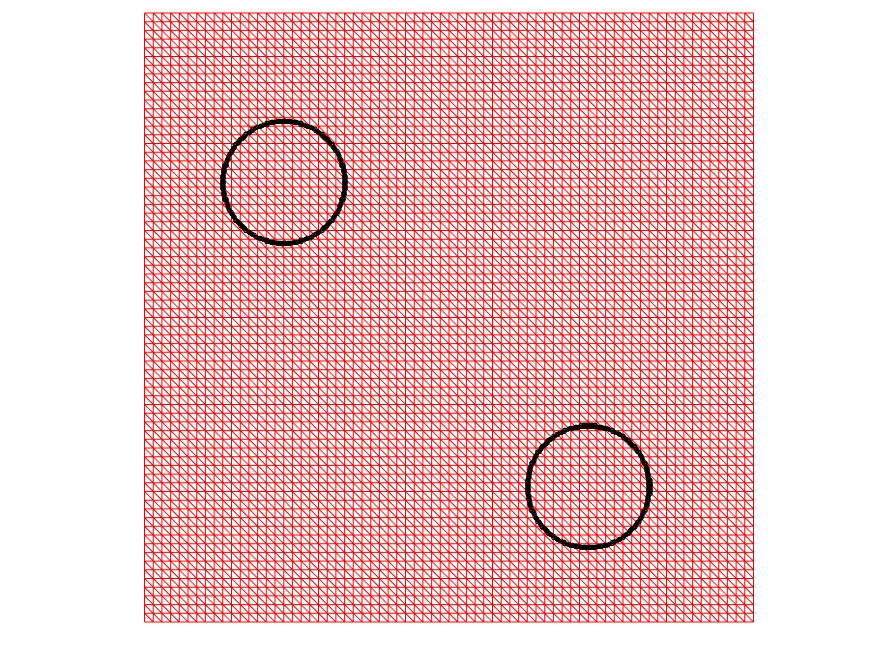}
&\includegraphics[height = \fht, trim = {3.5cm 0cm 3.5cm 0cm}, clip]{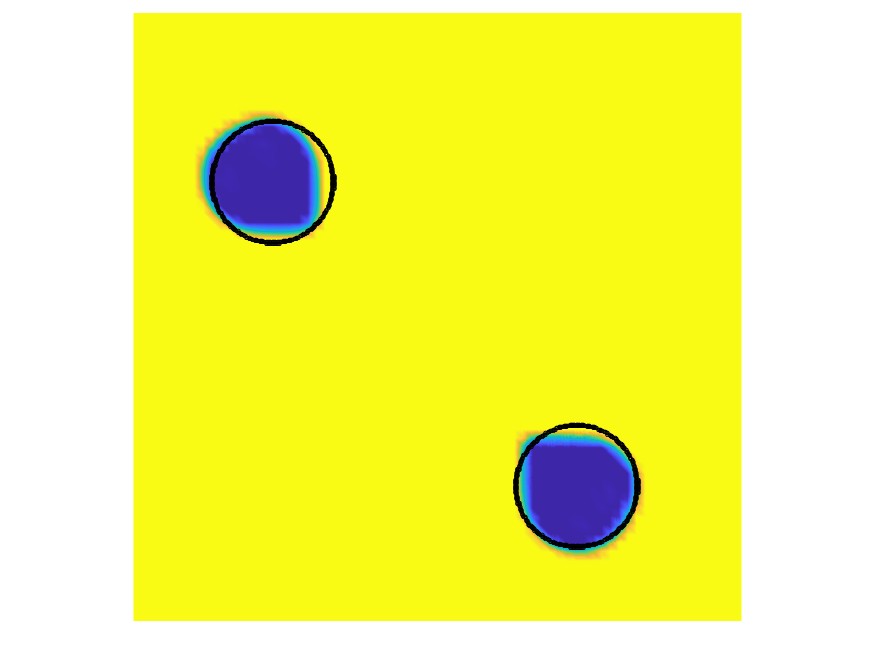}
&\includegraphics[height = \fht, trim = {3.5cm 0cm 3.5cm 0cm}, clip]{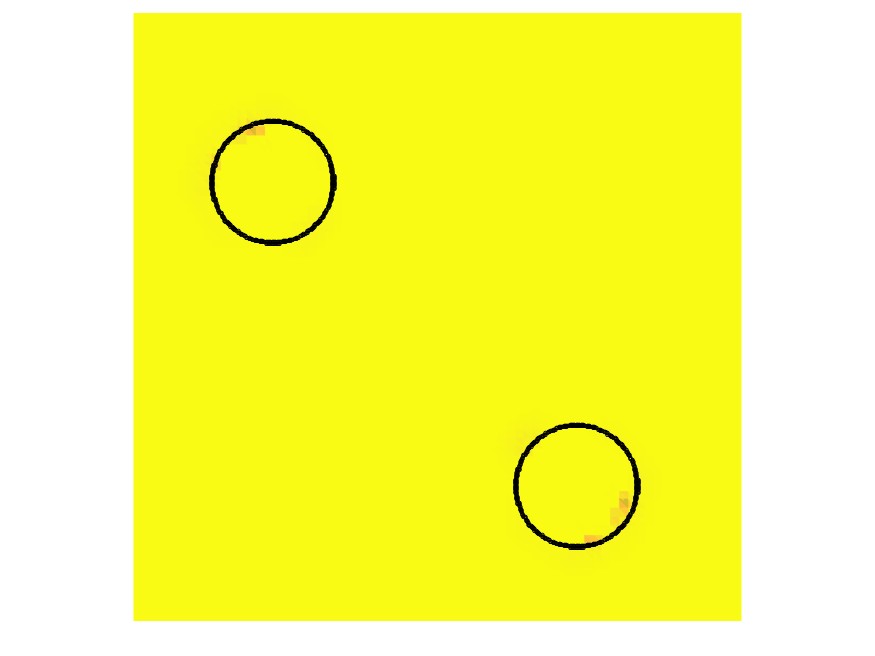}
\end{tabular}
\caption{The results by
uniform mesh refinements for the noisy data $y^\delta$.
From left to the right are the mesh, recovered inclusion and error indicator function.
The number of nodes at each step is 676, 1681, 3136 and 5041.
}
\label{fig:twocirclenoiseuniform}
\end{figure}

\begin{figure}[hbt!]
\centering
\setlength{\tabcolsep}{0pt}
\begin{tabular}{cccccc}
\includegraphics[height = \fht, trim = {3.5cm 0cm 3.5cm 0cm}, clip]{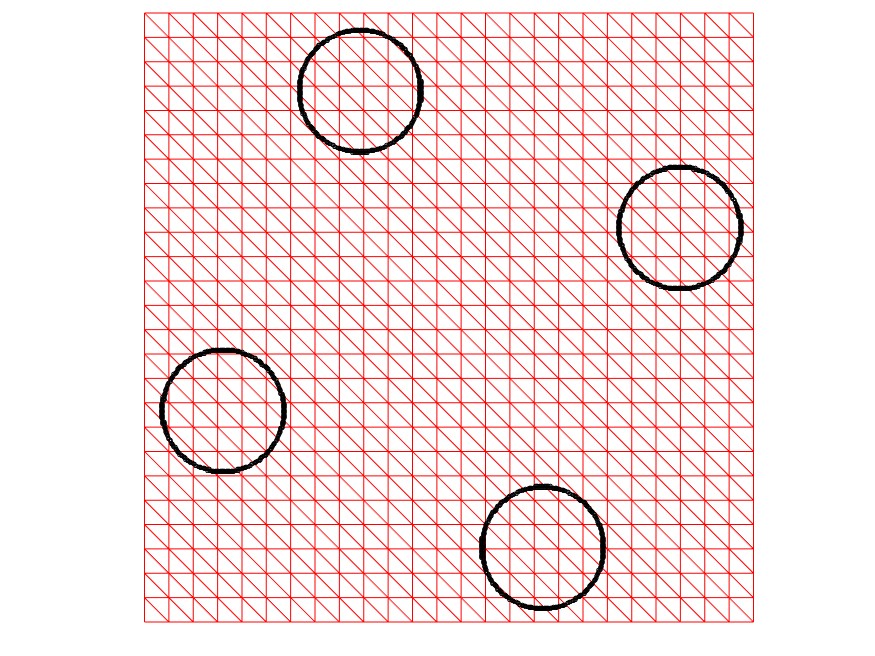}
&\includegraphics[height = \fht, trim = {3.5cm 0cm 3.5cm 0cm}, clip]{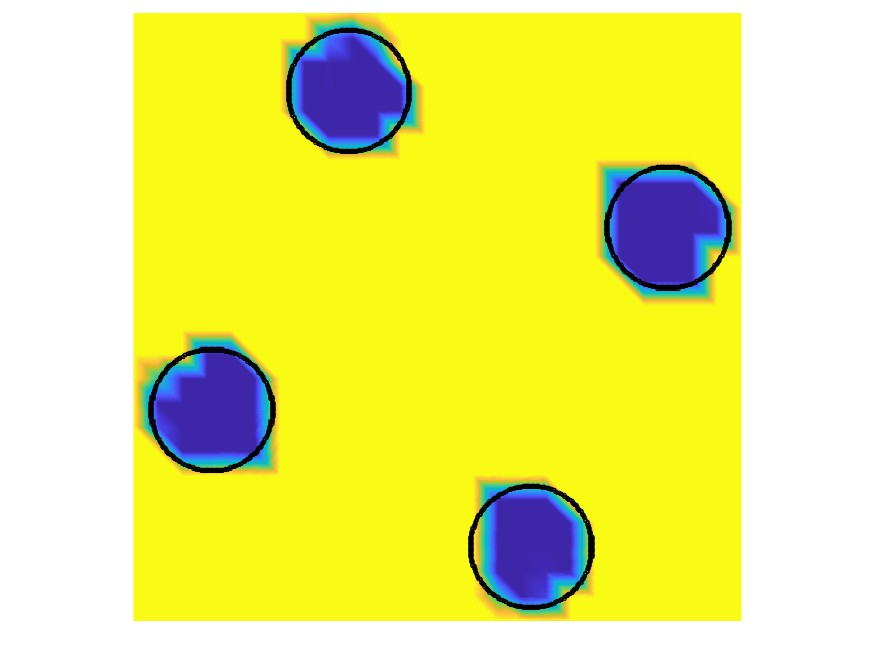}
&\includegraphics[height = \fht, trim = {3.5cm 0cm 3.5cm 0cm}, clip]{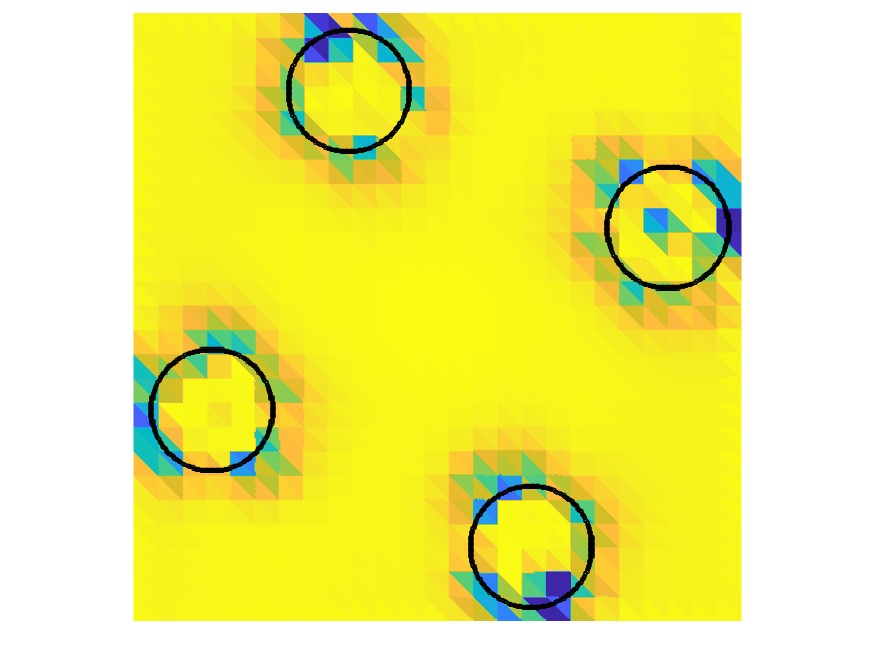}
&\includegraphics[height = \fht, trim = {3.5cm 0cm 3.5cm 0cm}, clip]{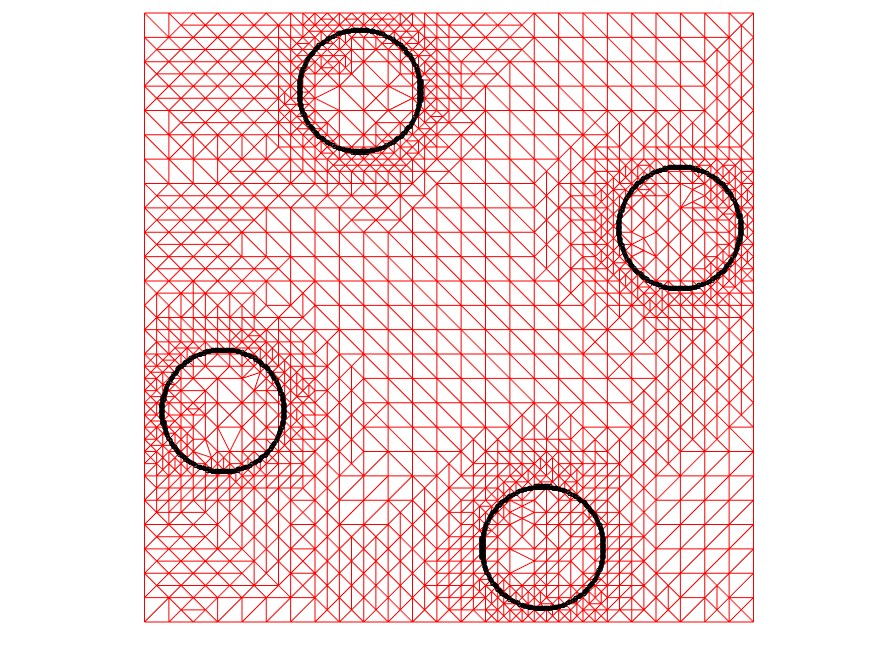}
&\includegraphics[height = \fht, trim = {3.5cm 0cm 3.5cm 0cm}, clip]{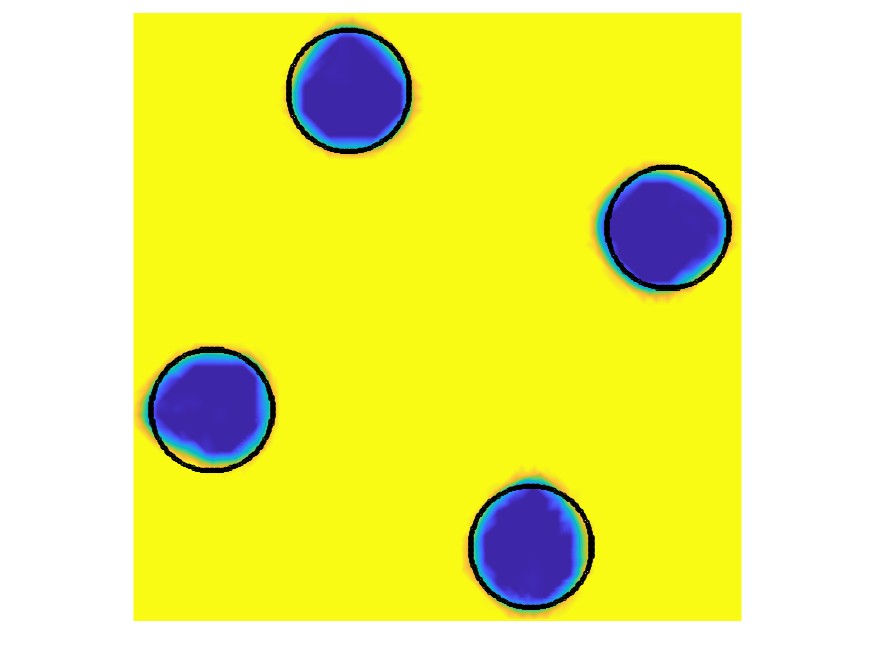}
&\includegraphics[height = \fht, trim = {3.5cm 0cm 3.5cm 0cm}, clip]{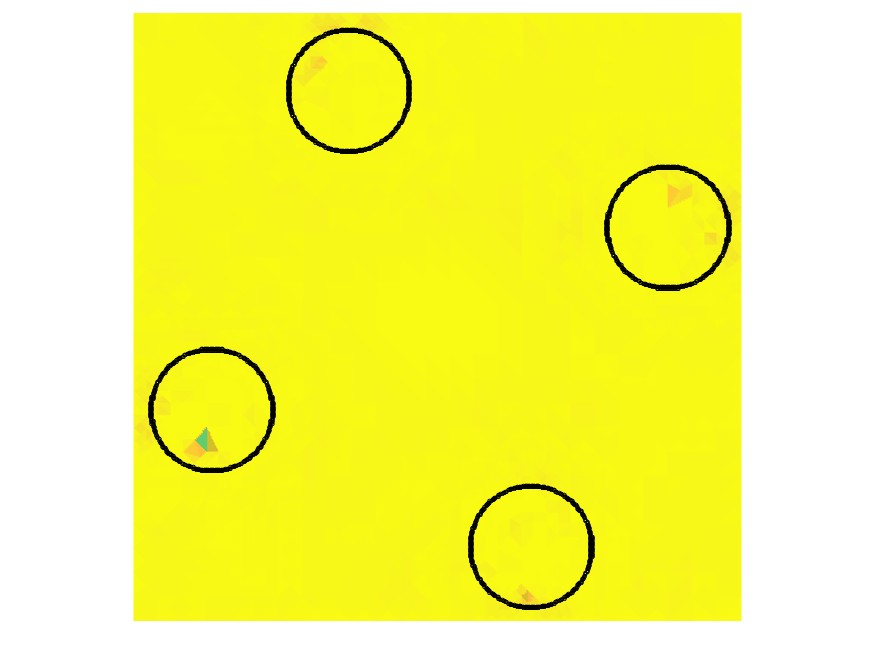}\\
\includegraphics[height = \fht, trim = {3.5cm 0cm 3.5cm 0cm}, clip]{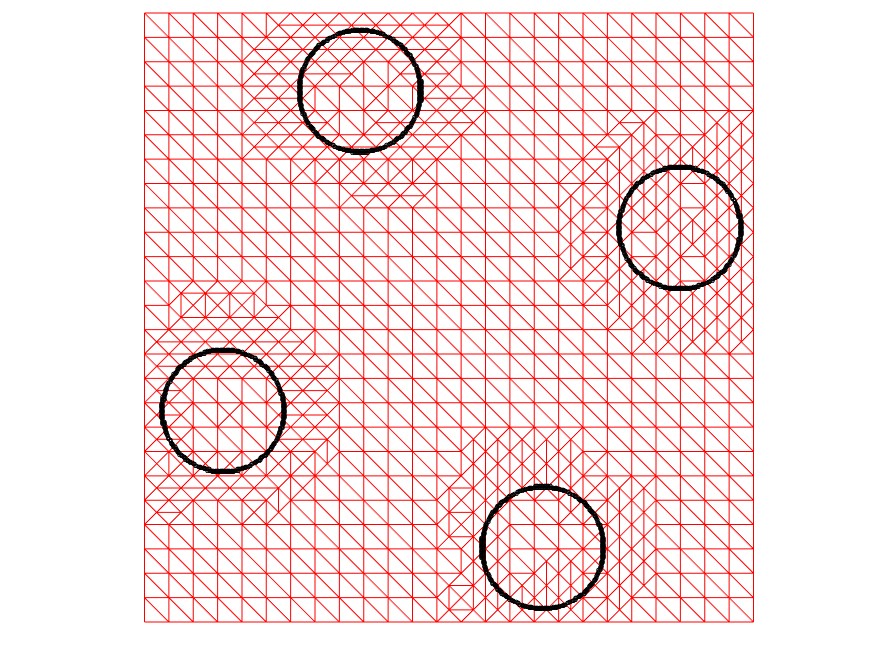}
&\includegraphics[height = \fht, trim = {3.5cm 0cm 3.5cm 0cm}, clip]{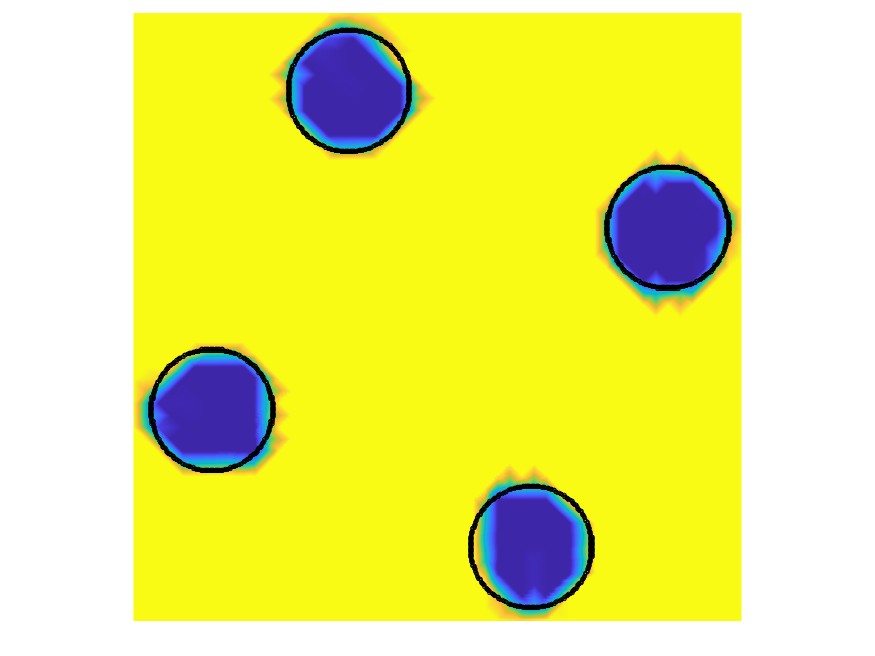}
&\includegraphics[height = \fht, trim = {3.5cm 0cm 3.5cm 0cm}, clip]{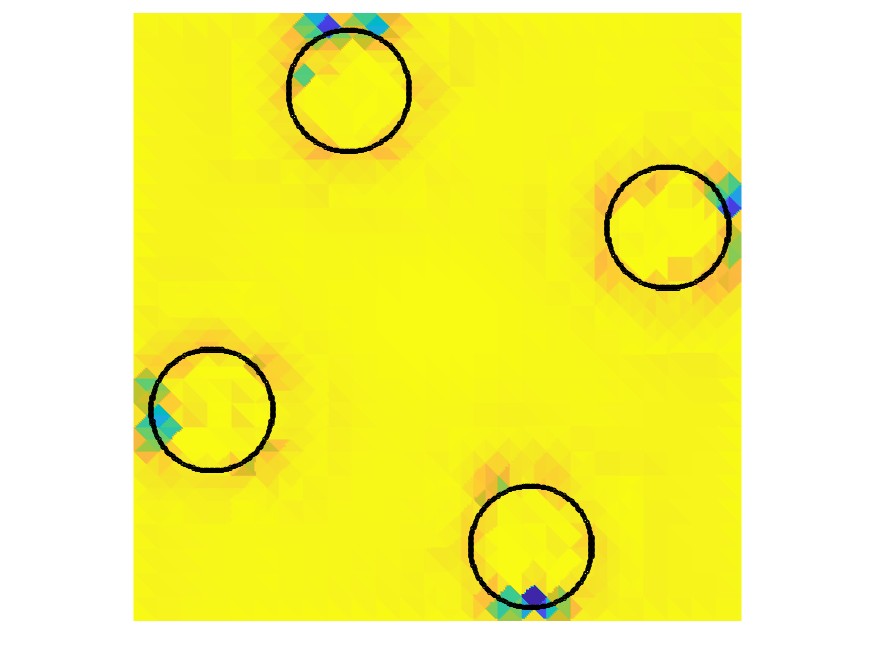}
&\includegraphics[height = \fht, trim = {3.5cm 0cm 3.5cm 0cm}, clip]{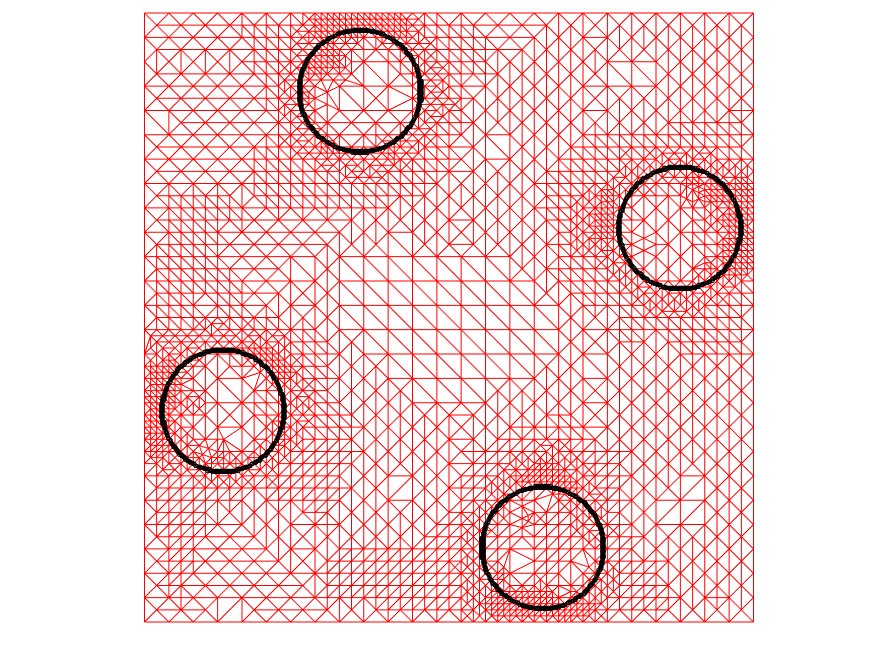}
&\includegraphics[height = \fht, trim = {3.5cm 0cm 3.5cm 0cm}, clip]{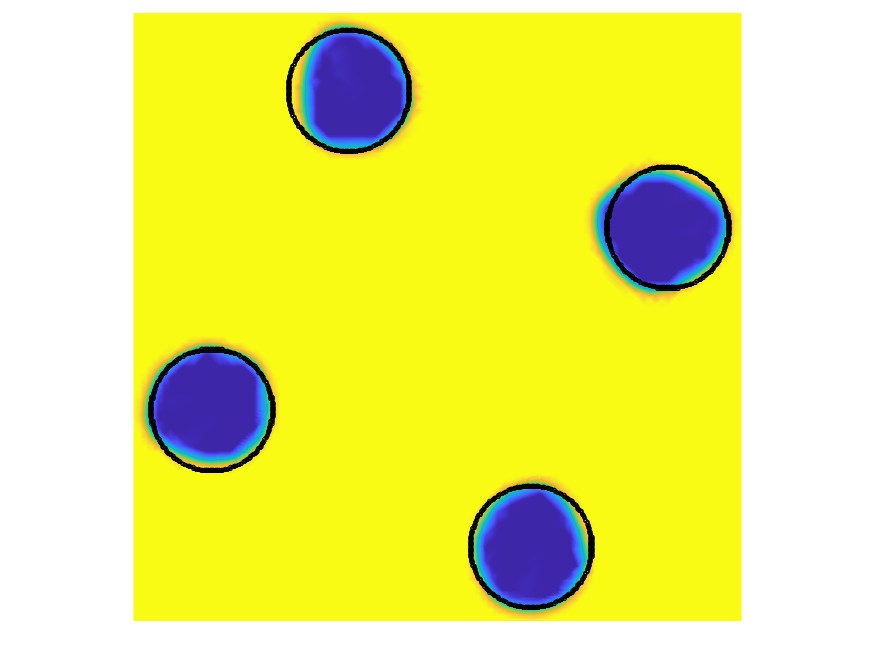}
&\includegraphics[height = \fht, trim = {3.5cm 0cm 3.5cm 0cm}, clip]{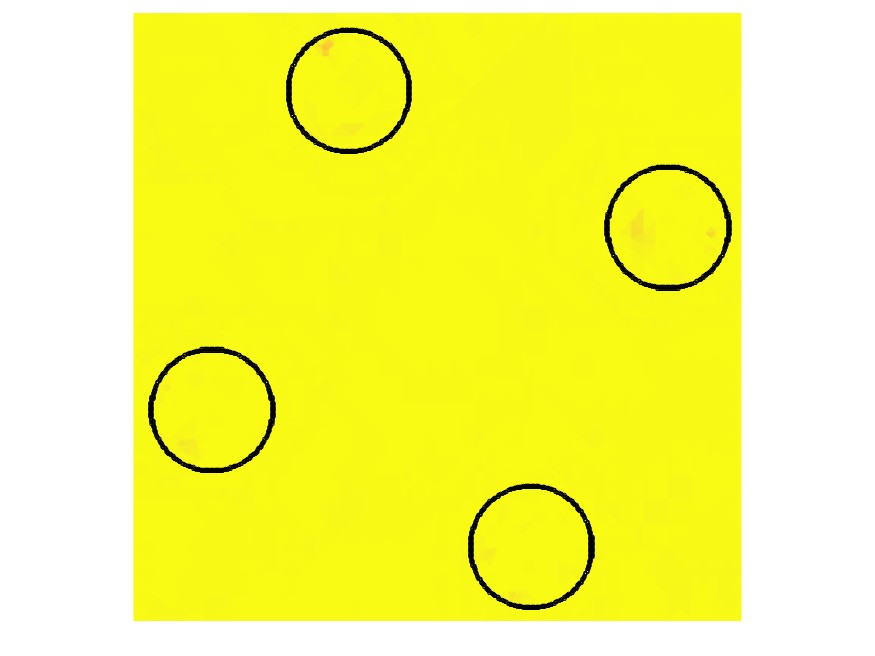}\\
\includegraphics[height = \fht, trim = {3.5cm 0cm 3.5cm 0cm}, clip]{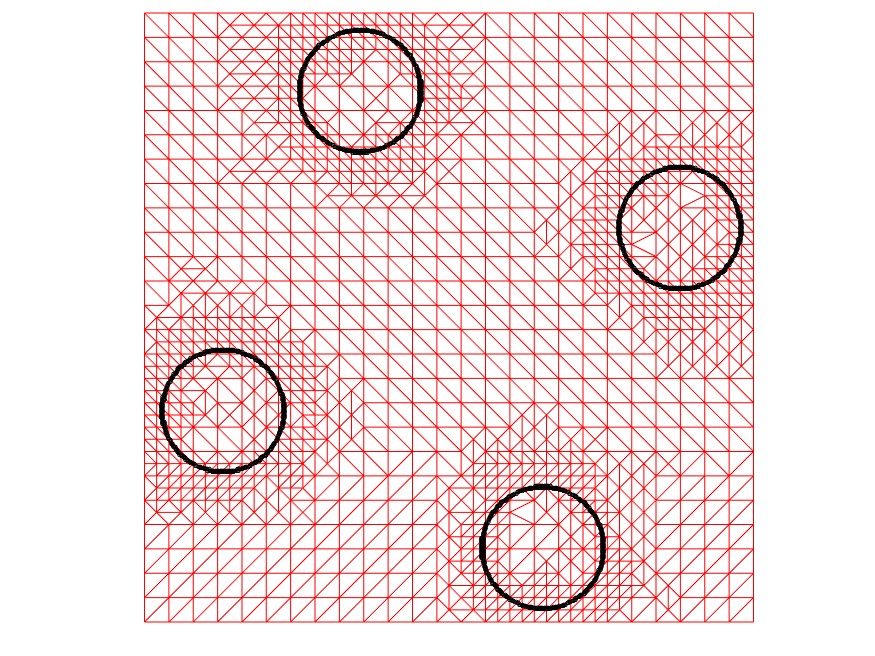}
&\includegraphics[height = \fht, trim = {3.5cm 0cm 3.5cm 0cm}, clip]{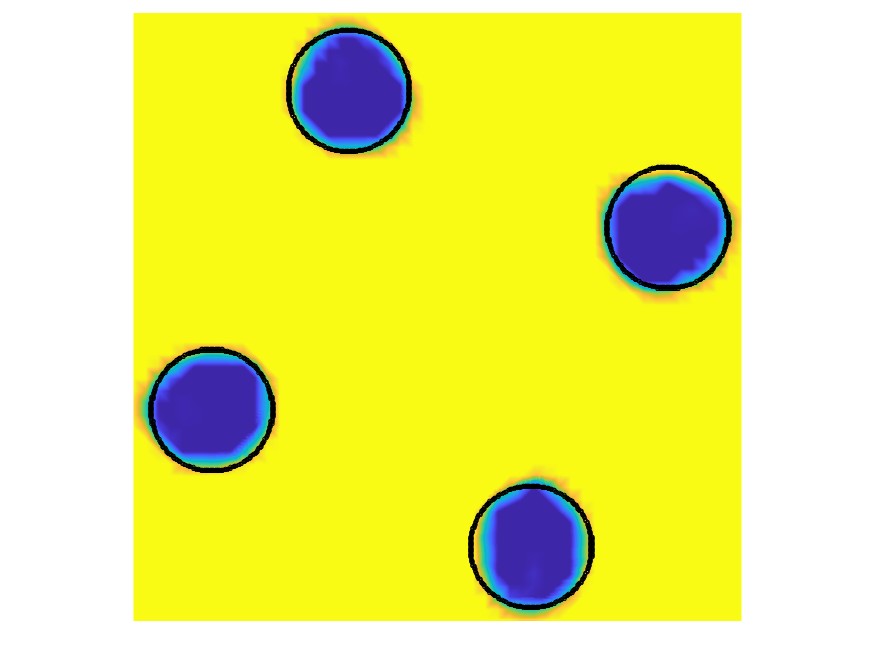}
&\includegraphics[height = \fht, trim = {3.5cm 0cm 3.5cm 0cm}, clip]{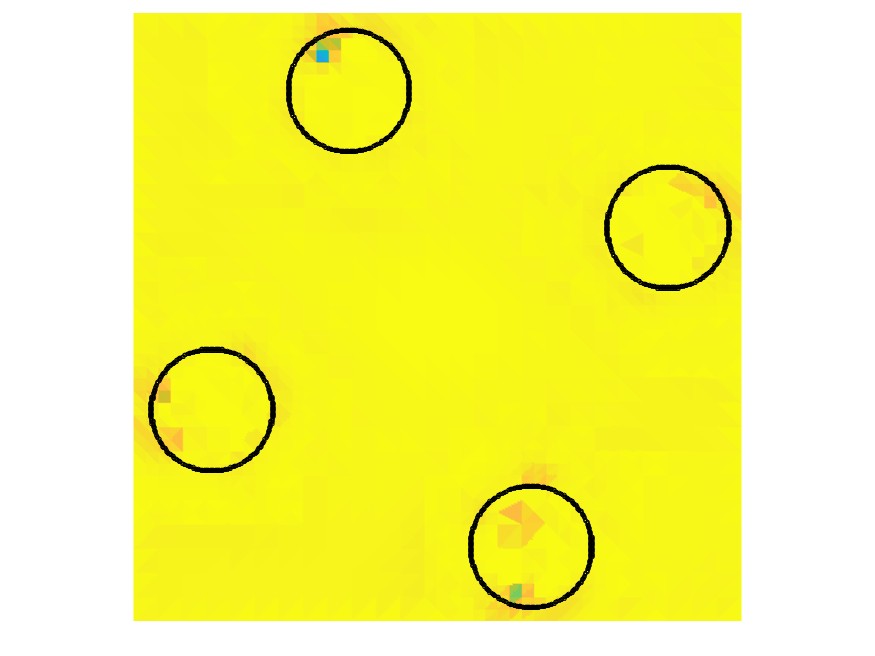}
&\includegraphics[height = \fht, trim = {3.5cm 0cm 3.5cm 0cm}, clip]{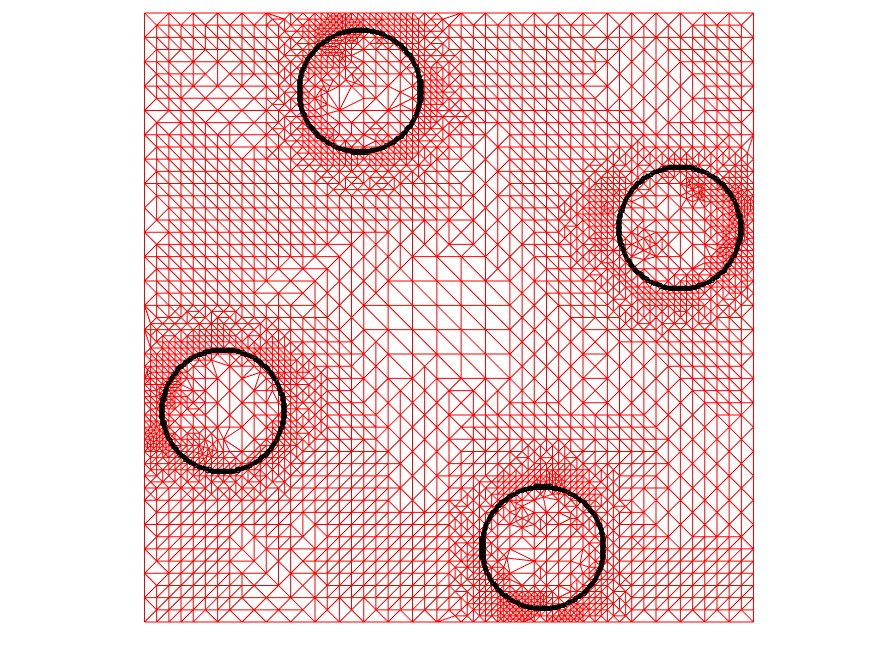}
&\includegraphics[height = \fht, trim = {3.5cm 0cm 3.5cm 0cm}, clip]{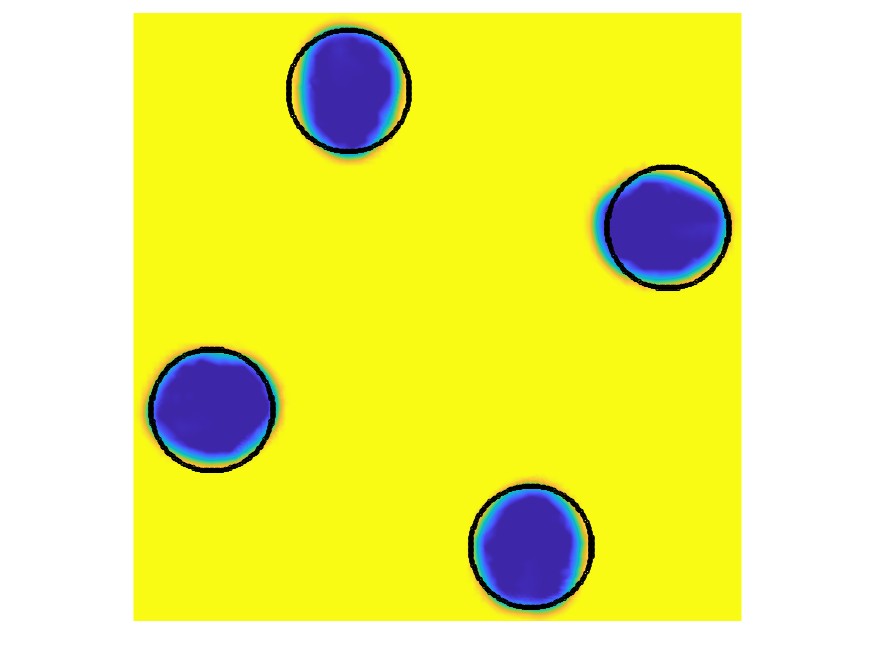}
&\includegraphics[height = \fht, trim = {3.5cm 0cm 3.5cm 0cm}, clip]{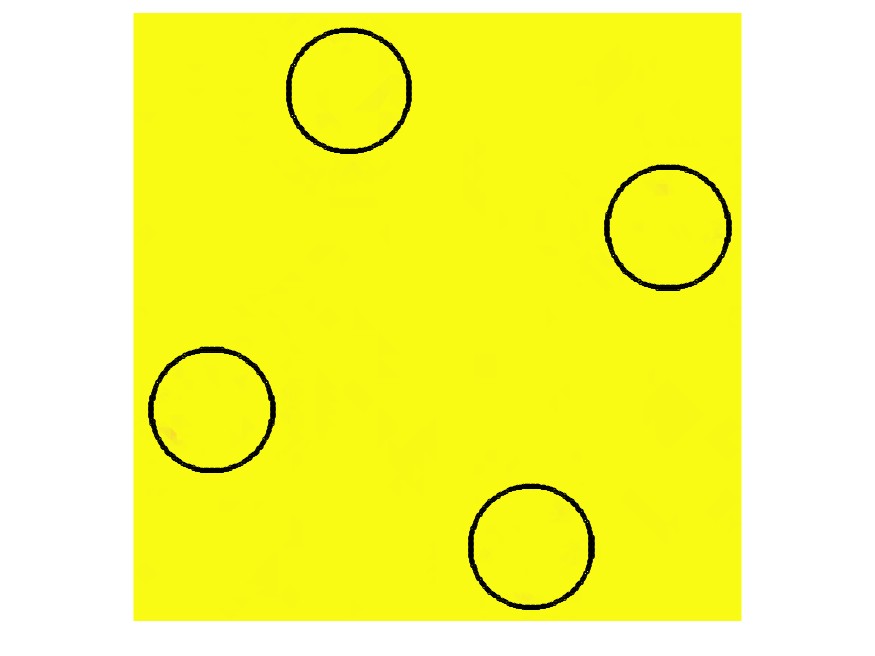}
\end{tabular}
\caption{The results by the adaptive method for the noisy data $y^\delta$.
From left to the right are the mesh, recovered inclusion and error indicator function.
The number of nodes for each step is
676, 918, 1271, 1798, 2556 and 3686.}
\label{fig:fourcirclenoiseadaptive}
\end{figure}

\begin{figure}[hbt!]
\centering
\setlength{\tabcolsep}{0pt}
\begin{tabular}{cccccc}
\includegraphics[height = \fht, trim = {3.5cm 0cm 3.5cm 0cm}, clip]{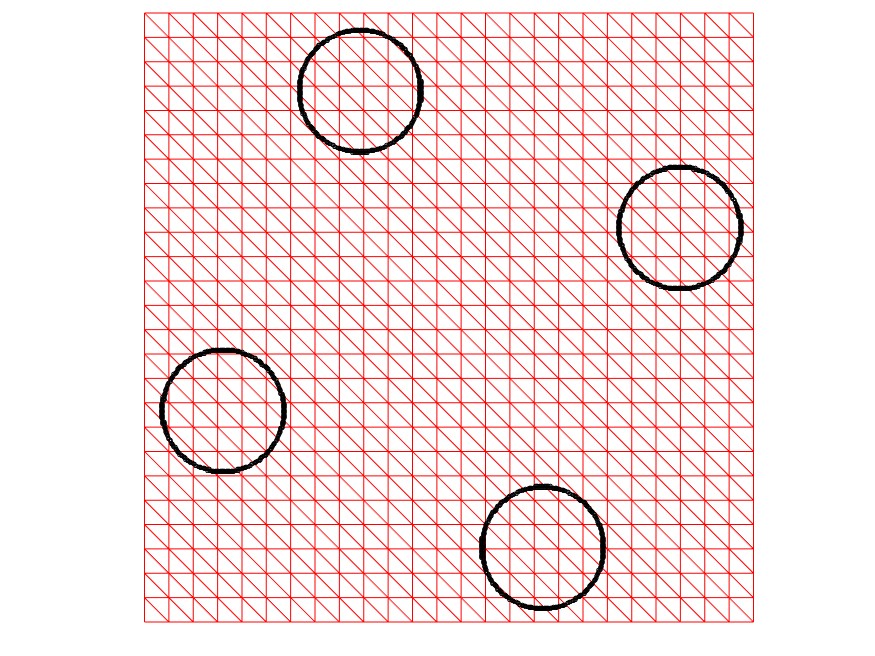}
&\includegraphics[height = \fht, trim = {3.5cm 0cm 3.5cm 0cm}, clip]{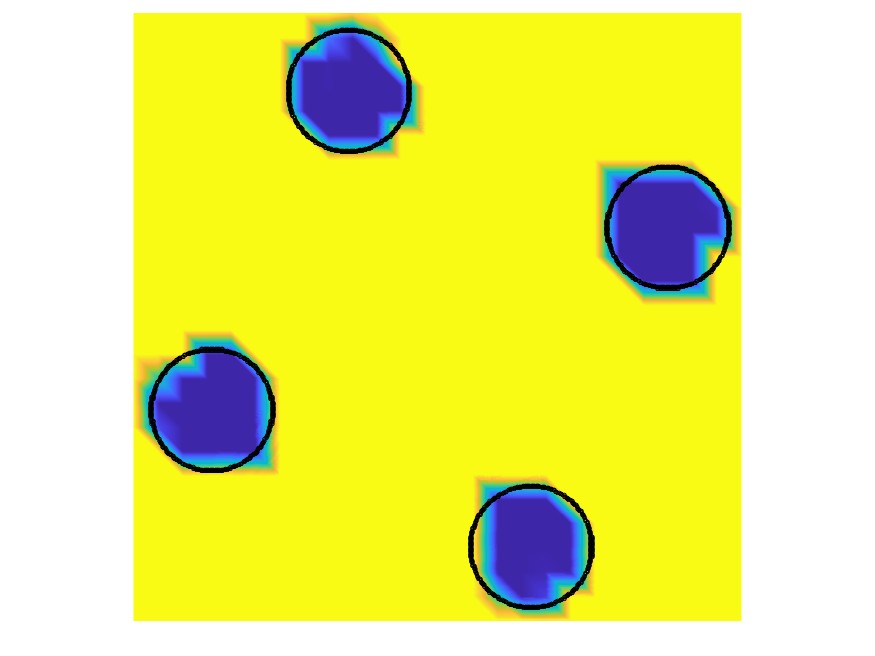}
&\includegraphics[height = \fht, trim = {3.5cm 0cm 3.5cm 0cm}, clip]{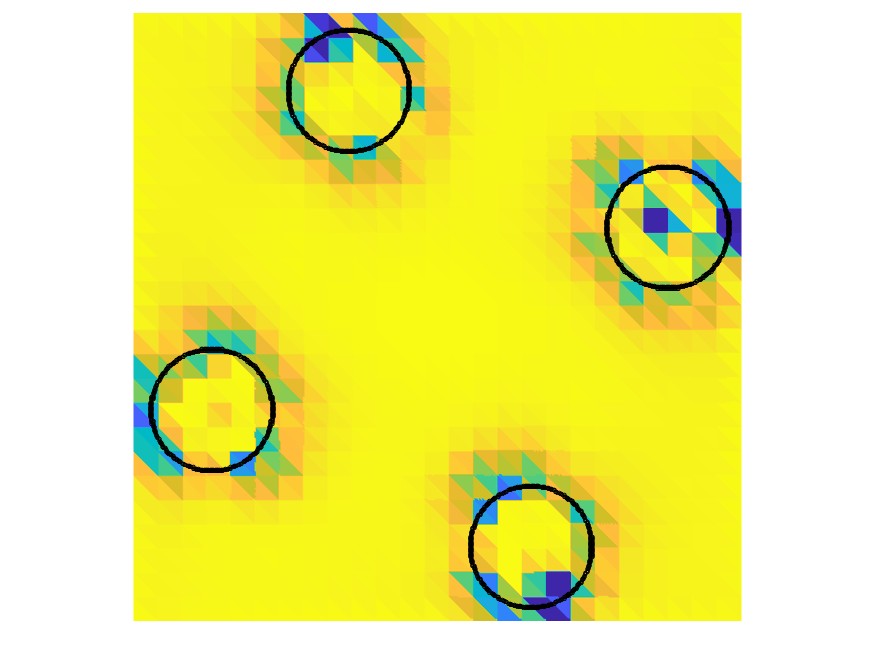}
&\includegraphics[height = \fht, trim = {3.5cm 0cm 3.5cm 0cm}, clip]{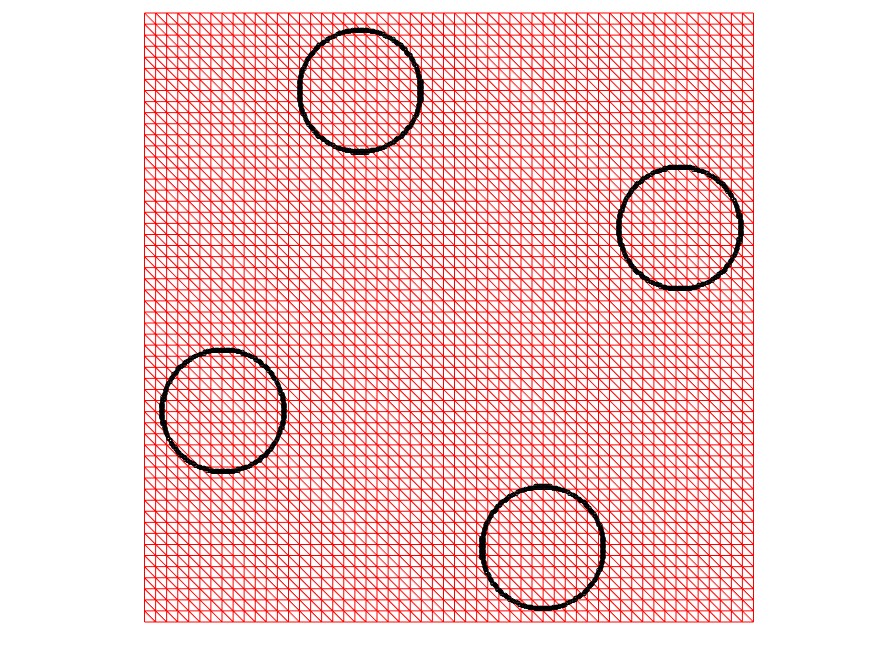}
&\includegraphics[height = \fht, trim = {3.5cm 0cm 3.5cm 0cm}, clip]{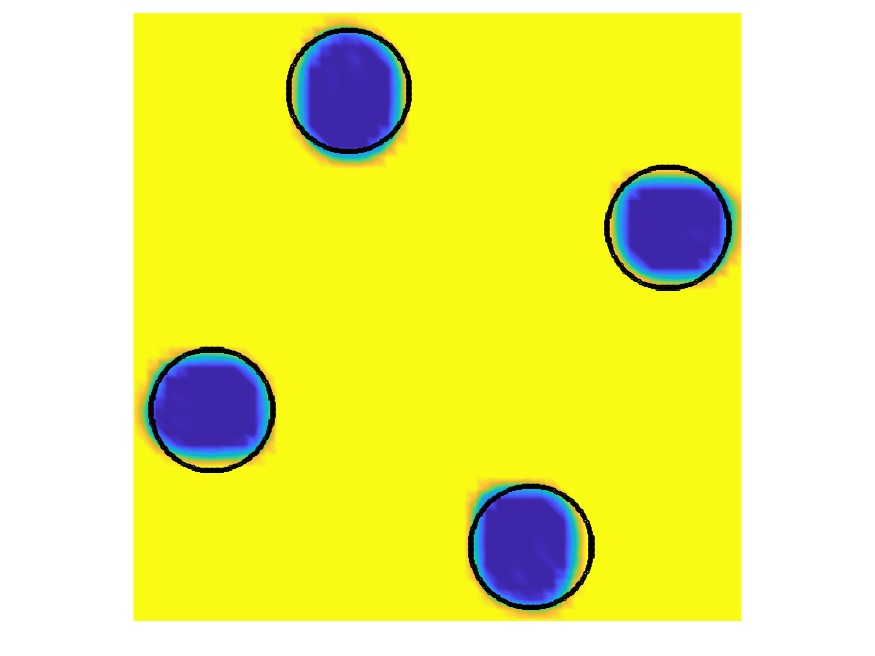}
&\includegraphics[height = \fht, trim = {3.5cm 0cm 3.5cm 0cm}, clip]{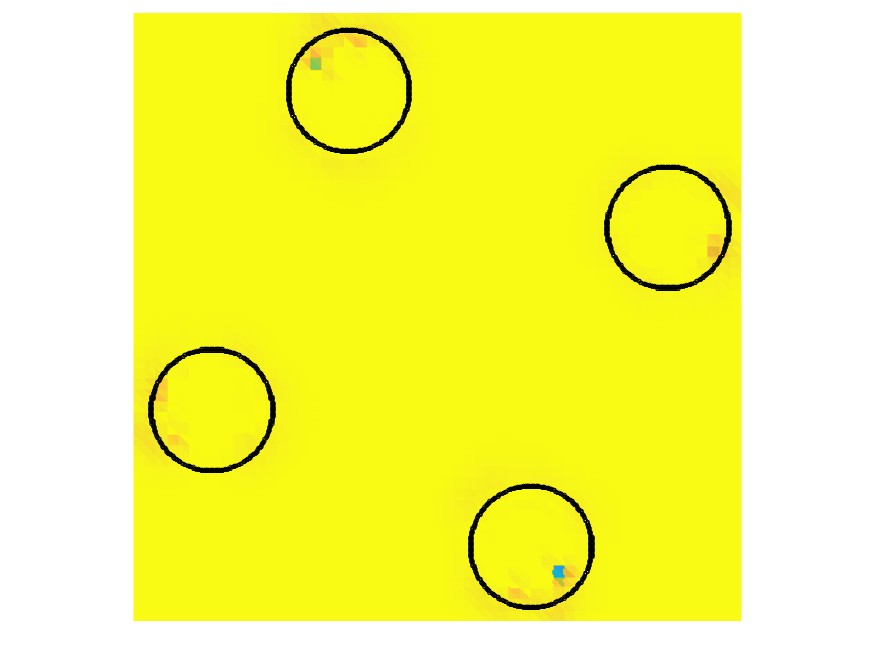}\\
\includegraphics[height = \fht, trim = {3.5cm 0cm 3.5cm 0cm}, clip]{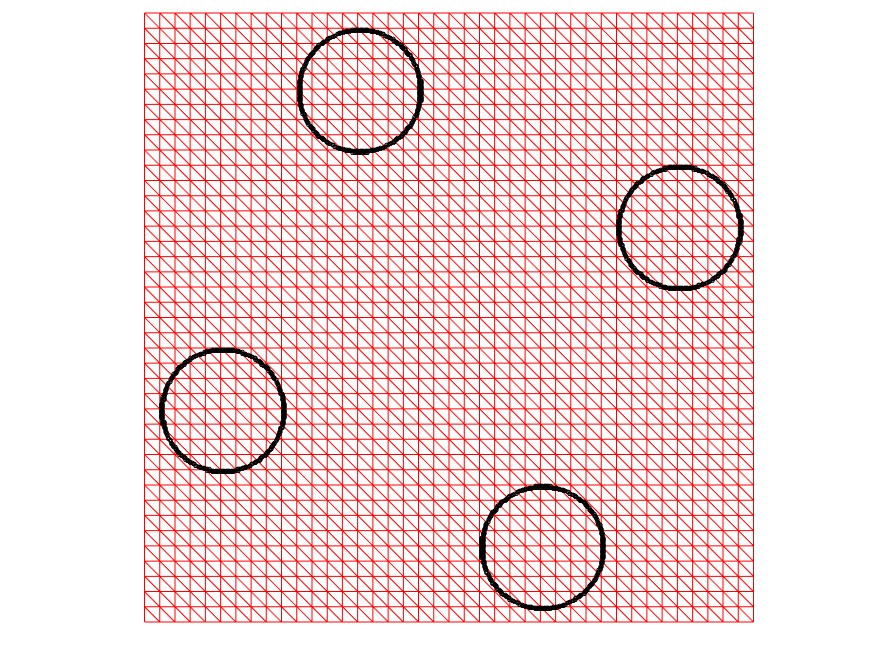}
&\includegraphics[height = \fht, trim = {3.5cm 0cm 3.5cm 0cm}, clip]{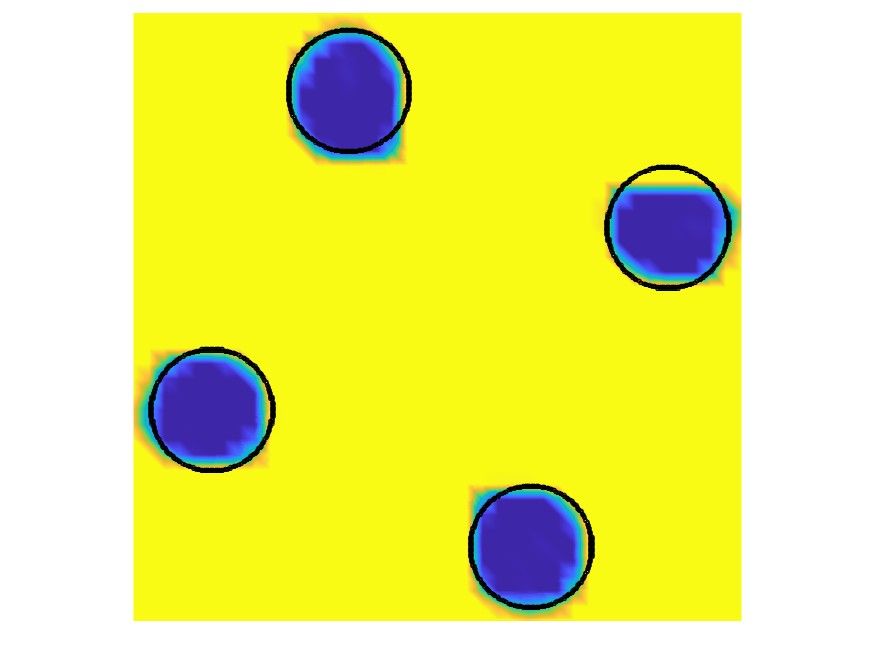}
&\includegraphics[height = \fht, trim = {3.5cm 0cm 3.5cm 0cm}, clip]{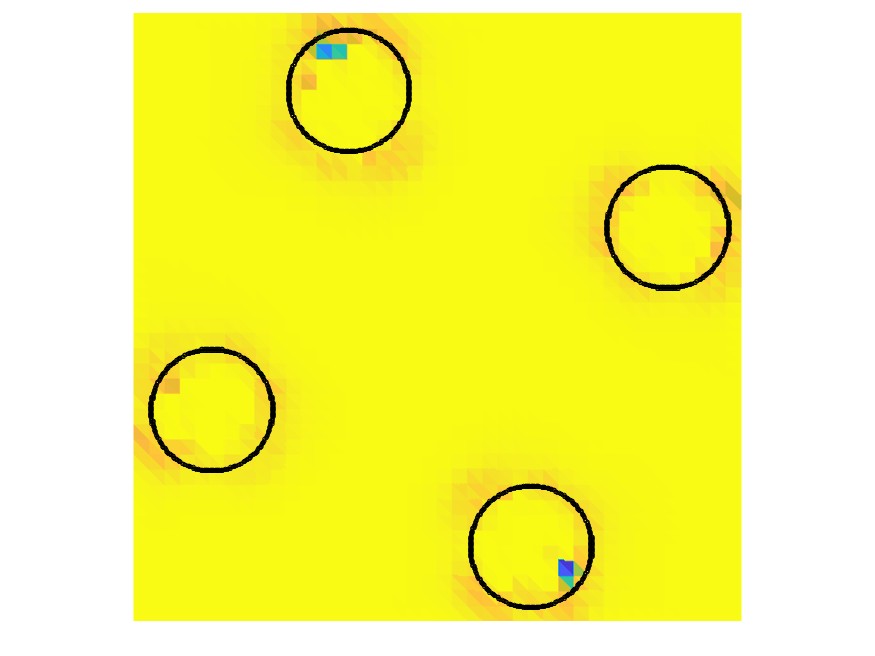}
&\includegraphics[height = \fht, trim = {3.5cm 0cm 3.5cm 0cm}, clip]{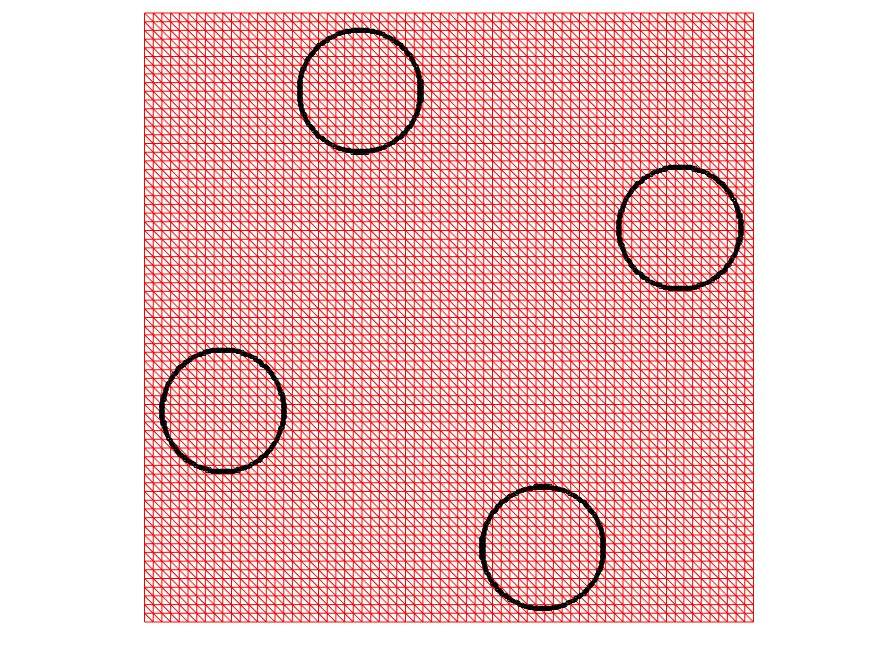}
&\includegraphics[height = \fht, trim = {3.5cm 0cm 3.5cm 0cm}, clip]{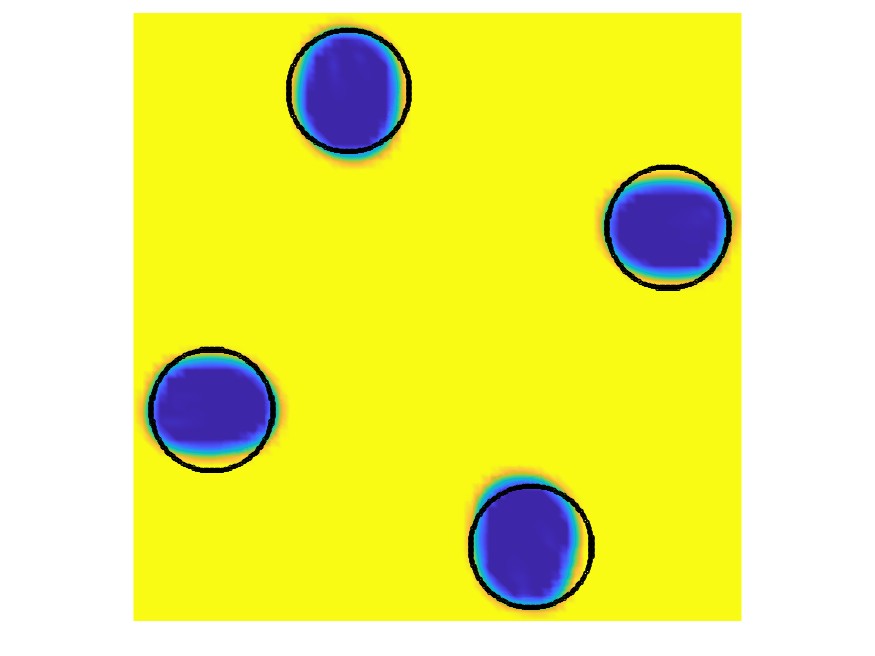}
&\includegraphics[height = \fht, trim = {3.5cm 0cm 3.5cm 0cm}, clip]
{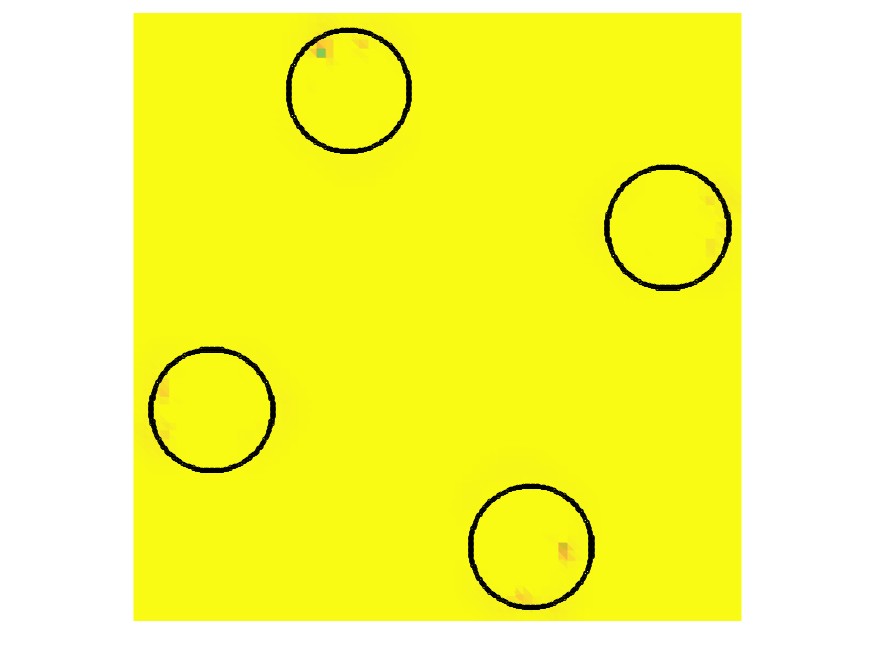}
\end{tabular}
\caption{The results by uniform mesh refinements for the noisy data $y^\delta$.
From left to the right are the mesh, recovered inclusion and error indicator function.
The number of nodes at each step is 676, 1681, 3136 and 5041.
}
\label{fig:fourcirclenoiseuniform}
\end{figure}

In Figs. \ref{fig:twocirclenoiseadaptive}--\ref{fig:fourcirclenoiseuniform}, we present the numerical results for recovering multiple circular inclusions from noisy data $y^\delta$. These two cases are more challenging (and thus require more data) due to the more complex geometric configurations. Nonetheless, when the observational data $y^\delta$ is relatively accurate, the inclusion locations and shapes can still be accurately resolved. The accuracy of the reconstructed inclusions by the adaptive algorithm and the uniform refinement strategy is largely comparable with each other at the last step. These results show the feasibility of the proposed adaptive algorithm for recovering multiple inclusions in a semilinear elliptic model.

\begin{figure}[hbt!]
\centering
\setlength{\tabcolsep}{0pt}
\begin{tabular}{cc}
\includegraphics[width = 0.45\textwidth, trim = {0.0cm 0cm 0.0cm 0cm}, clip]{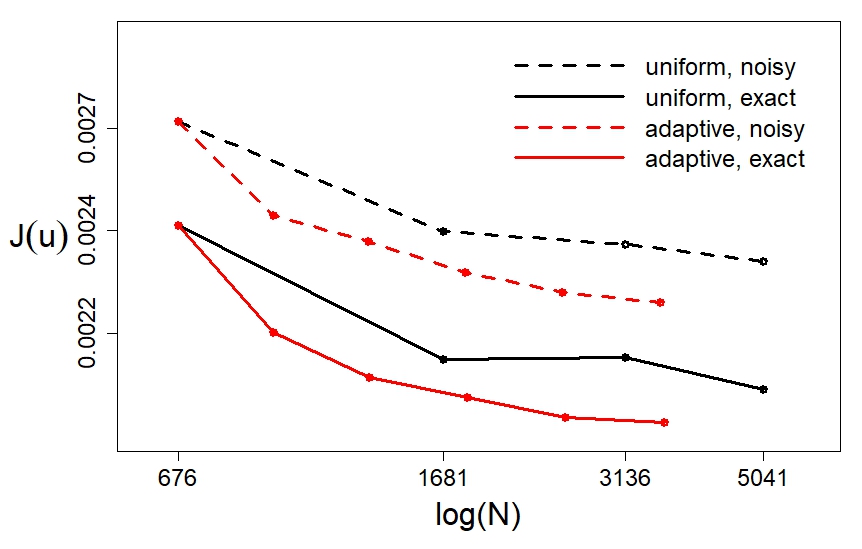}&
\includegraphics[width = 0.45\textwidth, trim = {0.0cm 0cm 0.0cm 0cm}, clip]{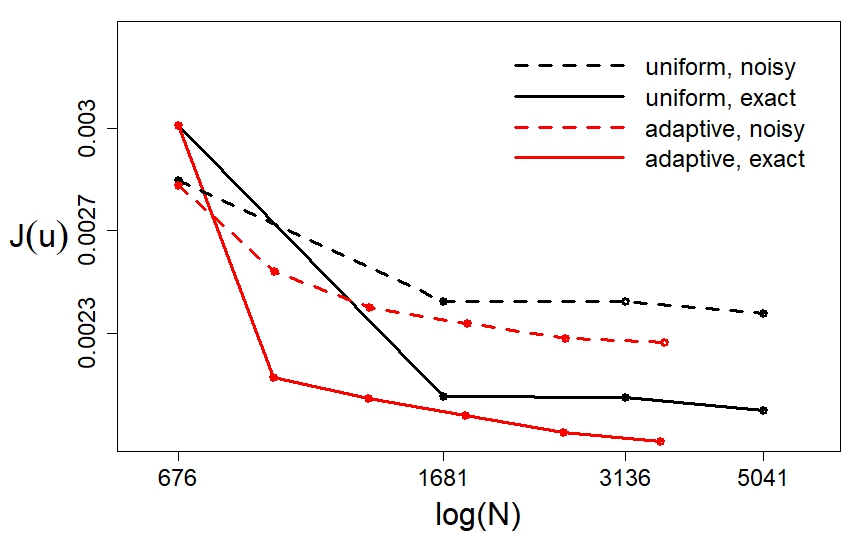}\\
(a) one circle & (b) one ellipse\\
\includegraphics[width = 0.45\textwidth, trim = {0.0cm 0cm 0.0cm 0cm}, clip]{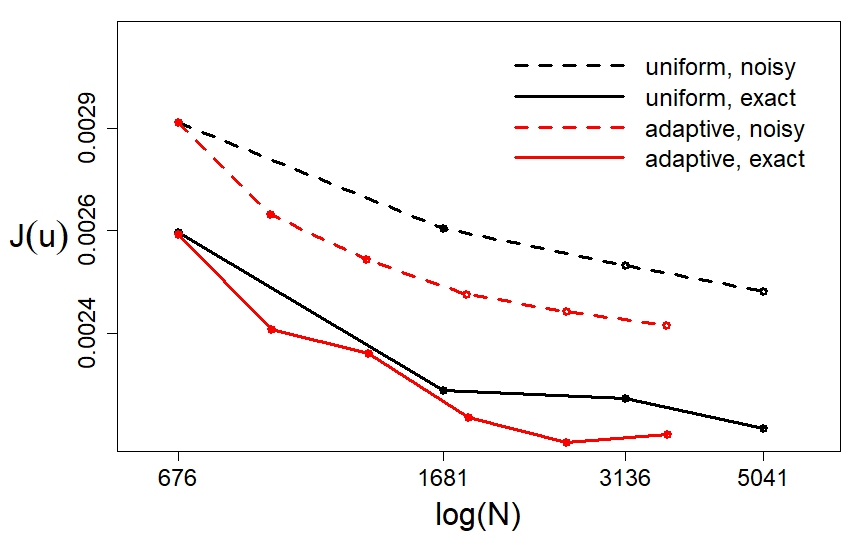}&
\includegraphics[width = 0.45\textwidth, trim = {0.0cm 0cm 0.0cm 0cm}, clip]{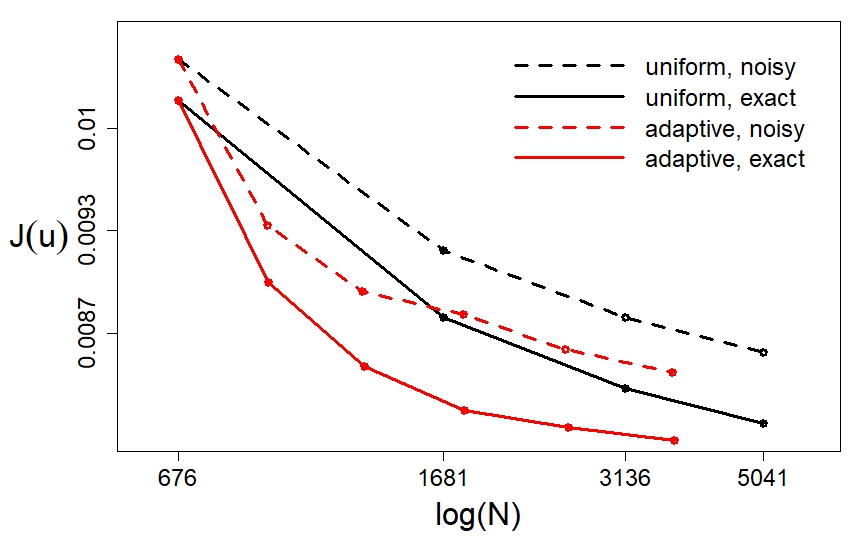}\\
(c) two circles & (d) four circles
\end{tabular}
\caption{
The convergence plots for the four cases. The $y$-axis represents the converged objective value $\mathcal{J}_{\varepsilon,k}(u_k^*)$, and the $x$-axis denotes the degrees of freedom of the mesh.
The dashed and solid lines refer to noisy and exact data, respectively.}
\label{fig:Rconv}
\end{figure}

Finally, we compare the convergence of the two mesh refinement strategies in terms of the objective $\mathcal{J}_{\varepsilon,k}(u_k^*)$ (since the two strategies give different solutions). Note that for both uniform and adaptive refinement strategies, the objective value $\mathcal{J}_{\varepsilon,k}(u_k^*)$ should converge to the minimal value of problem \eqref{min_G-L} as the mesh refinement proceeds, and thus a smaller value of the objective functional $\mathcal{J}_{\varepsilon,k}(u_k^*)$ indicates a better discretization strategy (for any fixed degree of freedom).
The convergence results for the above four cases are depicted in Fig. \ref{fig:Rconv}, which shows that each objective value $\J_{\eps,k}(u_{k}^*)$ decreases steadily as the mesh is refined either locally or globally. Moreover, in terms of the efficiency and accuracy, the adaptive algorithm has a clear advantage over the uniform refinement strategy.
For the adaptive method, the value $\J_{\eps,k}(u_{k}^*)$ decreases more rapidly as the degrees of freedom increases, indicating a faster convergence of the algorithm due to the adaptive refinements of the mesh around the inclusion interface $\partial\omega$.
The results for the noisy data $y^\delta$ also demonstrate the algorithm's robustness in the presence of measurement errors, and the overall observation is similar to that for the exact data $y^*$ except that in the presence of data noise, the converged value of the objective $J_{\eps,k}(u_k^*)$ is larger.

\begin{figure}[hbt!]
\centering
\setlength{\tabcolsep}{0pt}
\begin{tabular}{cc}
\includegraphics[width = 0.45\textwidth, trim = {0.0cm 0cm 0.0cm 0cm}, clip]{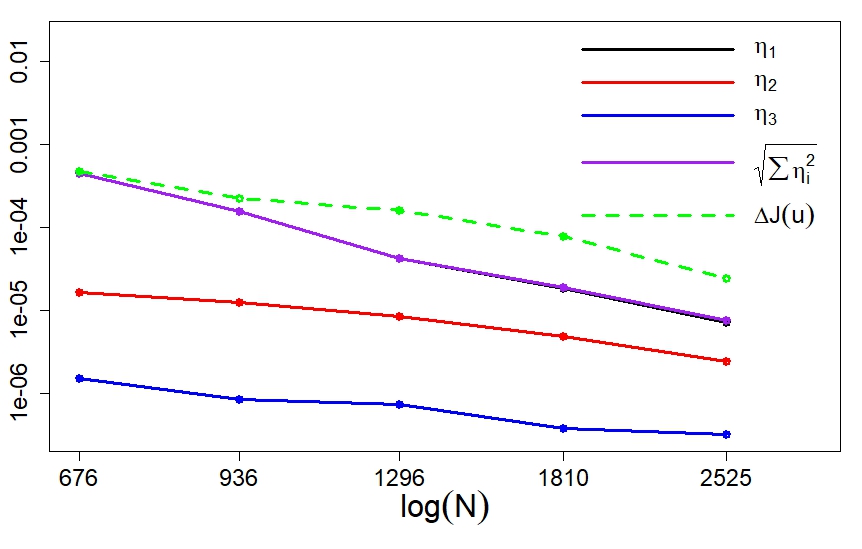}&
\includegraphics[width = 0.45\textwidth, trim = {0.0cm 0cm 0.0cm 0cm}, clip]{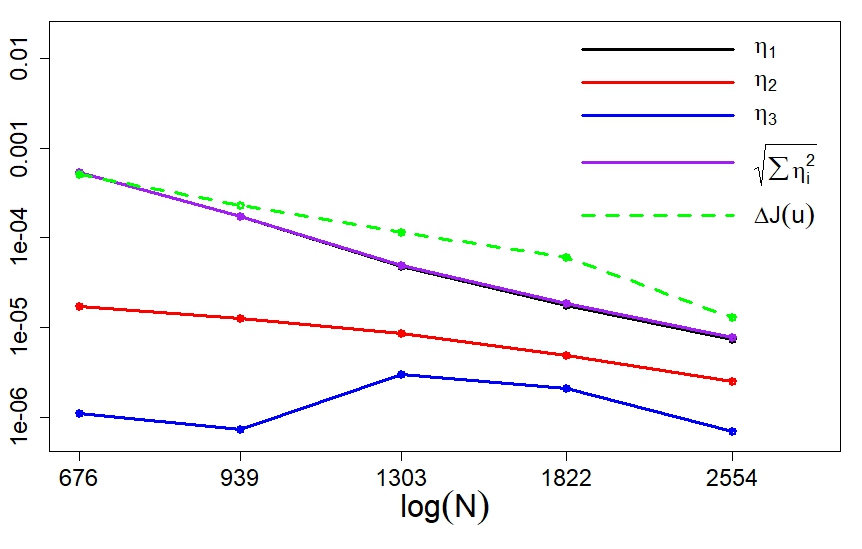}\\
(a) one circle & (b) one ellipse\\
\includegraphics[width = 0.45\textwidth, trim = {0.0cm 0cm 0.0cm 0cm}, clip]{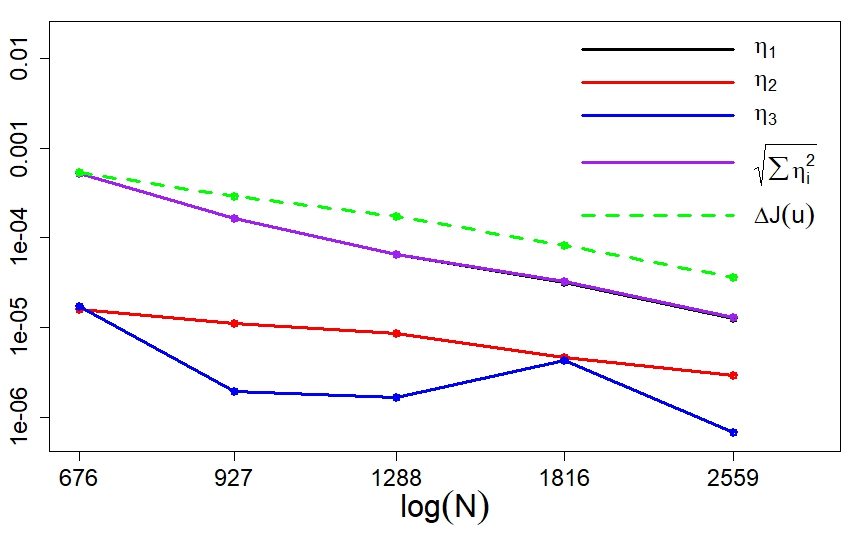}&
\includegraphics[width = 0.45\textwidth, trim = {0.0cm 0cm 0.0cm 0cm}, clip]{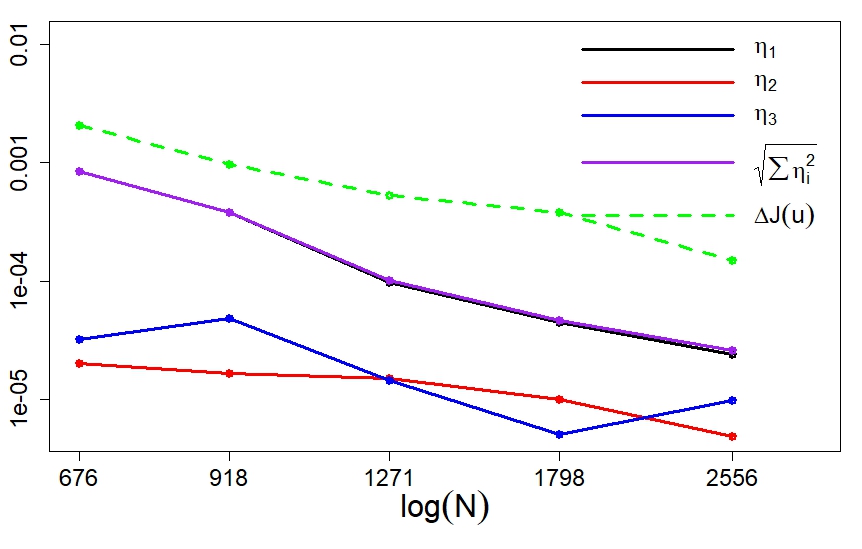}\\
(c) two circles & (d) four circles
\end{tabular}\caption{
The convergence of the error indicators $\{\eta_{k,i}(u^\ast_{k},y_{k}^\ast,T)\}_{i=1}^{3}$, $\Delta J = \mathcal{J}_{\varepsilon,k}(u_k^*) - \mathcal{J}_{\varepsilon,6}(u_6^*)$ (with $\mathcal{J}_{\varepsilon,6}(u_6^*)$ being the reference value).
The $x$-axis $N$ denotes the degrees of freedom of the mesh.}
\label{fig:Reta}
\end{figure}

To shed further insights, in Fig. \ref{fig:Reta} we present the convergence history of the error indicators $\eta_{k,i}$, $i=1,2,3$ as well as the combined estimator $\eta_k$ (computed over the entire mesh), and also the convergence history of the objective value error $\Delta J(u_k^*):= \mathcal{J}_{\varepsilon,k}(u_k^*) - \mathcal{J}_{\varepsilon,6}(u_6^*)$ (with $\mathcal{J}_{\varepsilon,6}(u_6^*)$ taken as the reference value). The plots in Fig. \ref{fig:Reta} show that all the error indicators exhibit a steady decrease as the adaptive mesh refinement proceeds. Moreover, as can be seen from the overlap of black and purple curves, in each adaptive loop the error estimator $\eta_{k,1}$ (for the state variable $y$) clearly dominates among the three components and largely drives the adaptive refinement process. The steady decreasing trend indicates the effectiveness of the adaptive algorithm in enhancing the accuracy of the numerical solution, which is reflected by more local refinements along inclusion interfaces; cf. Fig. \ref{fig:circlefreeadaptive} for an illustration in the case of one circular inclusion.
The observed convergence of the three estimators agrees well with that for the error in the objective value $\Delta J(u_k^\ast)$. However, the reliability of the estimators remains to be established theoretically.

\appendix

\section{Proof of Lemma \ref{lem:stab_est_unref}}

In this appendix, we prove the three estimates in Lemma \ref{lem:stab_est_unref} one by one.
By the triangle inequality,
\[
        \Big(\sum_{T\in\cT_{k_j}\cap\cT_{k_l}}\eta^2_{k_l,1}(u_{k_l}^\ast,y^\ast_{k_l},T) \Big)^{1/2} \leq \Big(\sum_{T\in\cT_{k_j}\cap\cT_{k_l}}\eta^2_{k_j,1}(u_{k_j}^\ast,y^\ast_{k_j},T) \Big)^{1/2} + \mathrm{I},
        \]
    with the term ${\rm I}$ given by
    \[
        \begin{aligned}
        \mathrm{I}:= &\Big(\sum_{T\in\cT_{k_l}\cap\cT_{k_j}}h_T^2\|{\nabla} \cdot (a(u_{k_l}^\ast) {\nabla} y_{k_l}^\ast - a(u_{k_j}^\ast){\nabla} y_{k_j}^\ast) + (b(u_{k_j}^\ast)(y_{k_j}^\ast)^3-b(u_{k_l}^\ast)(y_{k_l}^\ast)^3)\|^2_{L^2(T)}  \\
        &\quad + \sum_{T\in\cT_{k_l}\cap\cT_{k_j}}\sum_{F\subset\partial T}h_T\|[(a(u_{k_l}^\ast){\nabla} y_{k_l}^\ast -a(u_{k_j}^\ast){\nabla} y_{k_j}^\ast) \cdot {n}_F]\|^2_{L^2(F)}\Big)^{1/2}.
        \end{aligned}
    \]
 Next we bound the three summands of the term ${\rm I}$. First, we have the following splitting
        \begin{align*}
        {\nabla} \cdot (a(u_{k_l}^\ast) {\nabla} y_{k_l}^\ast - a(u_{k_j}^\ast) {\nabla} y_{k_j}^\ast) & = {\nabla}\cdot (a(u_{k_l}^\ast){\nabla}\overline{y}_{k_l,k_j}^\ast)
        +{\nabla}\cdot((a(u_{k_l}^\ast)-a(u_{k_j}^\ast)){\nabla}y^\ast_{k_j})\\
        &= (\sigma-1){\nabla}u_{k_l}^\ast\cdot{\nabla}\overline{y}_{k_l,k_j}^\ast+
        (\sigma-1){\nabla}\overline{u}_{k_l,k_j}^\ast\cdot{\nabla}y^\ast_{k_j}.
        \end{align*}
Since $|\sigma-1|\leq1$, by H\"{o}lder's inequality and the inverse estimate (cf. Lemma \ref{lem:inverse})
\begin{equation*} \|v_h\|_{L^p(T;\mathbb{R}^d)}\leq c h_T^{d(\frac{1}{p}-\frac12)}\|v_h\|_{L^2(T;\mathbb{R}^d)} \forall v_h\in P_0(T)^d,
\end{equation*}  we get
    \begin{align*}
     & h_T^2 \|(\sigma-1){\nabla}u_{k_l}^\ast\cdot{\nabla}\overline{y}_{k_l,k_j}^\ast\|_{L^2(T)}^2 \leq h_T^2 \|{\nabla}u_{k_l}^\ast\|_{L^4(T)}^2 \|{\nabla}\overline{y}_{k_l,k_j}^\ast\|_{L^4(T)}^2 \\
     \leq &c h_{T}^2 (h_{T}^{-1/2})^4 \|{\nabla}u_{k_l}^\ast\|_{L^2(T)}^2 \|{\nabla}\overline{y}_{k_l,k_j}^\ast\|_{L^2(T)}^2
         = c\|{\nabla}u_{k_l}^\ast\|_{L^2(T)}^2 \|{\nabla}\overline{y}_{k_l,k_j}^\ast\|_{L^2(T)}^2.
    \end{align*}
    Then summing the above inequality over the elements $T\in\cT_{k_j}\cap\cT_{k_l}$ yields
        \begin{align}
        \sum_{T\in\cT_{k_j}\cap\cT_{k_l}}  h_T^2 \|(\sigma-1){\nabla}u_{k_l}^\ast\cdot{\nabla}\overline{y}_{k_l,k_j}^\ast\|_{L^2(T)}^2
        &\leq c\|u_{k_l}^\ast\|^2_{H^1(\Omega)}\|{\nabla}\overline{y}_{k_l,k_j}^\ast\|_{L^2(\Omega)}^2.\label{pf:stab_est_unref1_01}
    \end{align}
    Similarly, we have
    \begin{align}    \sum_{T\in\cT_{k_j}\cap\cT_{k_l}} h_T^2\|(\sigma-1){\nabla}\overline{u}_{k_l,k_j}^\ast\cdot{\nabla}y^\ast_{k_l}\|^2_{L^2(T)}
        &\leq c \|y_{k_j}^\ast\|^2_{H^1(\Omega)}\|{\nabla}\overline{u}_{k_l,k_j}^\ast\|_{L^2(\Omega)}^2.\label{pf:stab_est_unref1_02}
    \end{align}
    Second, using the splitting
    \begin{equation*}
     b(u_{k_j}^\ast)(y_{k_j}^\ast)^3-b(u_{k_l}^\ast)(y_{k_l}^\ast)^3 = (b(u_{k_j}^\ast)-b(u_{k_l}^\ast))(y_{k_j}^\ast)^3+b(u_{k_l}^\ast)((y_{k_j}^\ast)^3-(y_{k_l}^\ast)^3),
    \end{equation*}
and then repeating the preceding argument with the generalized H\"{o}lder's inequality and the inverse inequality from Lemma \ref{lem:inverse} lead to
    \begin{align*}
        &h_T^2 \|(b(u_{k_j}^\ast)-b(u_{k_l}^\ast))(y_{k_j}^\ast)^3\|_{L^2(T)}^2
        = h_T^2 \|\overline{u}_{k_l,k_j}^\ast(y_{k_j}^\ast)^3\|_{L^2(T)}^2
        \leq h_T^2 \|\overline{u}_{k_l,k_j}^\ast\|^2_{L^8(T)}\|y_{k_j}^\ast\|^6_{L^8(T)}\\
        \leq& c h_T^2 (h_{T}^{-3/4})^2 (h_T^{-1/12})^6\|\overline{u}_{k_l,k_j}^\ast\|_{L^2(T)}^2\|y_{k_j}^\ast\|^6_{L^6(T)}=c\|\overline{u}_{k_l,k_j}^\ast\|_{L^2(T)}^2\|y_{k_j}^\ast\|^6_{L^6(T)}.
    \end{align*}
Similarly, the identity $a^3-b^3=(a-b)(a^2+ab+b^2)$, the generalized H\"{o}lder's inequality, and the inverse inequality in Lemma \ref{lem:inverse} imply
    \begin{align*}
            h_T^2 \|b(u_{k_l}^\ast)((y_{k_j}^\ast)^3-&(y_{k_l}^\ast)^3)\|_{L^2(T)}^2
           \leq h_T^2 \|\overline{y}_{k_l,k_j}^\ast((y_{k_l}^\ast)^2 + y_{k_j}^*y_{k_l}^*+(y_{k_j}^\ast)^2)\|_{L^2(T)}^2\\
&\leq  ch_T^2\|\overline{y}_{k_l,k_j}^\ast\|_{L^6(T)}^2(\|y_{k_l}^\ast\|_{L^6(T)}^4+  \|y_{k_j}^\ast\|_{L^6(T)}^4)\\
   &\leq ch_{T}^2(h_T^{-2/3})^2(h_T^{-1/6})^{4}\|\overline{y}_{k_l,k_j}^\ast\|_{L^2(T)}^2(\|y_{k_l}^\ast\|_{L^4(T)}^4+\|y_{k_j}^\ast\|_{L^4(T)}^4)\\
   & = c\|\overline{y}_{k_l,k_j}^\ast\|_{L^2(T)}^2(\|y_{k_l}^\ast\|_{L^4(T)}^4+\|y_{k_j}^\ast\|_{L^4(T)}^4).
    \end{align*}
    Then summing over $T\in\cT_{k_j}\cap\cT_{k_l}$ and applying Sobolev embedding theorem yield
\begin{align}
 \sum_{T\in\cT_{k_j}\cap\cT_{k_l}} h_T^2 \|(b(u_{k_j}^\ast)-b(u_{k_l}^\ast))(y_{k_j}^\ast)^3\|_{L^2(T)}^2
      &  \leq c \|y_{k_j}^\ast\|^6_{L^6(\Omega)}\|\overline{u}_{k_l,k_j}^\ast\|_{L^2(\Omega)}^2\nonumber\\
      &\leq c\|y_{k_j}^\ast\|^6_{H^1(\Omega)}\|\overline{u}_{k_l,k_j}^\ast\|_{L^2(\Omega)}^2,\label{pf:stab_est_unref1_03}\\
\sum_{T\in\cT_{k_j}\cap\cT_{k_l}} h_T^2 \|b(u_{k_l}^\ast)((y_{k_j}^\ast)^3-(y_{k_l}^\ast)^3)\|_{L^2(T)}^2
       &\leq c (\|y_{k_l}^\ast\|_{H^1(\Omega)}^4+\|y_{k_j}^\ast\|_{H^1(            \Omega)}^4)\|\overline{y}_{k_l,k_j}^\ast\|_{L^2(\Omega)}^2.\label{pf:stab_est_unref1_04}
    \end{align}
Last, for the boundary term, we use the splitting
\[
     [(a(u_{k_l}^\ast) \nabla y_{k_l}^\ast - a(u_{k_j}^\ast) \nabla y_{k_j}^\ast) \cdot {n}_F] = [a(u_{k_l}^\ast) {\nabla} \overline{y}_{k_l,k_j}^\ast \cdot {n}_F] + [ ( a(u_{k_l}^\ast) - a(u_{k_j}^\ast) ){\nabla} y_{k_j}^\ast \cdot{n}_F].
\]
The estimate $|a(u_{k_l}^\ast)| \leq 1$ a.e. in $\Omega$, the scaled trace theorem $$\|v_h\|_{L^p(F;\mathbb{R}^d)}\leq c h_T^{-\frac1p+d(\frac1p-\frac1r)}\|v_h\|_{L^r(T;\mathbb{R}^d)},\quad \forall v_h \in P_1(T)^d$$
and the inverse estimate (cf. Lemma \ref{lem:inverse}) imply
\begin{align*}
  h_T \| [a(u_{k_l}^\ast) {\nabla} \overline{y}_{k_l,k_j}^\ast \cdot {n}_F]\|^2_{L^2(F)}&\leq
            c \|{\nabla}\overline{y}_{k_l,k_j}^\ast \|^2_{L^2(\omega(F))},\\
  h_T \|[ ( a(u_{k_l}^\ast) - a(u_{k_j}^\ast) ){\nabla} y_{k_j}^\ast \cdot {n}_F]\|_{L^2(F)}^2 
            & \leq c \|\overline{u}_{k_l,k_j}^\ast {\nabla} y_{k_j}^\ast\|^2_{L^2(\omega(F))}.
\end{align*}
By combining the last two estimates and the triangle inequality, we obtain
\begin{align}
&\sum_{T\in\cT_{k_j}\cap\cT_{k_l}} \sum_{F\subset\p T} h_T \|[a(u_{k_l}^\ast) {\nabla} \overline{y}_{k_l,k_j}^\ast \cdot {n}_F]\|^2_{L^2(F)}           \leq c \|{\nabla}\overline{y}_{k_l,k_j}^\ast \|^2_{L^2(\Omega)},\label{pf:stab_est_unref1_05}\\
& \sum_{T\in\cT_{k_j}\cap\cT_{k_l}} \sum_{F\subset\p T} h_T \|[ ( a(u_{k_l}^\ast) - a(u_{k_j}^\ast) ){\nabla} y_{k_j}^\ast \cdot {n}_F]\|_{L^2(F)}^2
\leq c ( \|{\nabla}\overline{y}_{k_j}^\ast\|^2_{L^2(\Omega)} + \|\overline{u}_{k_l,k_j}^\ast {\nabla} y_{\infty}^\ast\|^2_{L^2(\Omega)}).\label{pf:stab_est_unref1_06}
\end{align}
By Theorem \ref{thm:conv_medmin}, the sequences $\{u_{k_j}^\ast\}_{j\geq0}$ and $\{y_{k_j}^\ast\}_{j\geq0}$ are both uniformly bounded
in $H^1(\Omega)$. Then collecting estimates \eqref{pf:stab_est_unref1_01}--\eqref{pf:stab_est_unref1_06} leads
to the estimate \eqref{stab_est_unref1}. Likewise, for the estimate \eqref{stab_est_unref2}, we have
\[
   \Big(\sum_{T\in\cT_{k_j}\cap\cT_{k_l}}\eta^2_{k_l,2}(u_{k_l}^\ast,y^\ast_{k_l},p_{k_l}^\ast,T) \Big)^{1/2} \leq \Big(\sum_{T\in\cT_{k_j}\cap\cT_{k_l}}\eta^2_{k_j,2}(u_{k_j}^\ast,y^\ast_{k_j},p_{k_j}^\ast,T) \Big)^{1/2} + \mathrm{II},
\]
with the term ${\rm II}$ given by
\begin{align*}
        \mathrm{II}:= &\Big(\sum_{T\in\cT_{k_l}\cap\cT_{k_j}}h_T^2\|{\nabla} \cdot ( a(u_{k_l}^\ast){\nabla} p_{k_l}^\ast - a(u_{k_j}^\ast){\nabla} p_{k_j}^\ast ) - 3b(u_{k_l}^\ast) (y_{k_l}^\ast)^2 p_{k_l}^\ast + 3b(u_{k_j}^\ast) (y_{k_j}^\ast)^2 p_{k_j}^\ast \|^2_{L^2(T)} \\
        &\quad + \sum_{T\in\cT_{k_l}\cap\cT_{k_j}}\sum_{F\subset\partial T}h_T\left\|J_{F,2}(u_{k_l}^\ast, y_{k_l}^\ast, p_{k_l}^\ast) - J_{F,2}(u_{k_j}^\ast, y_{k_j}^\ast, p_{k_j}^\ast)\right\|_{L^2(F)}^2\Big)^{1/2}.
        \end{align*}
Next we bound the three summands of ${\rm II}$. First, repeating the argument
    for the estimate \eqref{pf:stab_est_unref1_01} yields
    \begin{align*}
         &\sum_{T\in\cT_{k_l}\cap\cT_{k_j}}h_T^2\|{\nabla} \cdot( a(u_{k_l}^\ast){\nabla} p_{k_l}^\ast - a(u_{k_j}^\ast){\nabla} p_{k_j}^\ast)\|_{L^2(T)}^2\nonumber\\
         \leq & c(\|u_{k_l}^\ast\|^2_{H^1(\Omega)}\|{\nabla}\overline{p}_{k_l,k_j}^\ast \|_{L^2(\Omega)}^2+\|p_{k_j}^\ast\|^2_{H^1(\Omega)}\|{\nabla}\overline{u}_{k_l,k_j}^\ast\|_{L^2(\Omega)}^2).
    \end{align*}
    Second, in the splitting
   \begin{align*}
       &b(u_{k_j}^\ast) (y_{k_j}^\ast)^2 p_{k_j}^\ast - 3b(u_{k_l}^\ast) (y_{k_l}^\ast)^2 p_{k_l}^\ast\\
         =&  \overline{u}_{k_j,k_l}^\ast (y_{k_j}^\ast)^2p_{k_j}^\ast +  b(u_{k_l}^\ast)( (y_{k_j}^\ast)^2 - (y_{k_l}^\ast)^2) p_{k_j}^\ast + b(u_{k_l}^\ast)(y_{k_l}^\ast)^2\overline{p}_{k_j,k_l}^\ast ,
    \end{align*}
  using the argument for the estimates \eqref{pf:stab_est_unref1_03} and \eqref{pf:stab_est_unref1_04} yields
  \begin{align*}
        &\quad h_T^2  \|3 \overline{u}_{k_j,k_l}^\ast (y_{k_j}^\ast)^2p_{k_j}^\ast\|_{L^2(T)}^2
        \leq 9 h_T^2 \|\overline{u}_{k_l,k_j}^\ast\|_{L^8(T)}^2\|y_{k_j}^\ast\|_{L^8(T)}^4\|p_{k_j}^\ast\|_{L^8(T)}^2\nonumber\\
        & \leq c h_T^2 (h_{T}^{-3/4})^2 (h_T^{-1/12})^{4} (h_T^{-1/12})^2 \|\overline{u}_{k_l,k_j}^\ast \|_{L^2(T)}^2\|y_{k_j}^\ast\|_{L^6(T)}^4\|p_{k_j}^\ast\|_{L^6(T)}^2\\
        & \leq c\|\overline{u}_{k_l,k_j}^\ast \|_{L^2(T)}^2\|y_{k_j}^\ast\|_{L^6(\Omega)}^4\|p_{k_j}^\ast\|_{L^6(\Omega)}^2.
\end{align*}
Likewise, we derive
\begin{align*}
         h_T^2  \|3 b(u_{k_l}^\ast)( (y_{k_j}^\ast)^2 - (y_{k_l}^\ast)^2) p_{k_j}^\ast\|_{L^2(T)}^2
        \leq& 9h_T^2 \|\overline{y}_{k_j,k_l}^\ast \|_{L^6(T)}^2 \|  y_{k_j}^\ast + y_{k_l}^\ast \|_{L^6(T)}^2 \| p_{k_j}^\ast\|_{L^6(T)}^2\\
       \leq &c  h_T^2 (h_T^{-2/3})^2(h_T^{-1/6})^4 \|\overline{y}_{k_l,k_j}^\ast \|_{L^2(T)}^2 \|  y_{k_j}^\ast + y_{k_l}^\ast \|_{L^4(\Omega)}^2 \| p_{k_j}^\ast\|_{L^4(\Omega)}^2 \\
       =& c\|\overline{y}_{k_l,k_j}^\ast \|_{L^2(T)}^2 \|  y_{k_j}^\ast + y_{k_l}^\ast \|_{L^4(\Omega)}^2 \| p_{k_j}^\ast\|_{L^4(\Omega)}^2,\\
    h_T^2\|3 b(u_{k_l}^\ast)(y_{k_l}^\ast)^2\overline{p}_{k_j,k_l}^\ast \|_{L^2(T)}^2 \leq& 9 h_T^2 \|y_{k_l}^\ast\|_{L^6(T)}^4 \|\overline{p}_{k_l,k_j}^\ast\|_{L^6(T)}^2
    \leq c \|y_{k_l}^\ast\|_{L^4(T)}^4 \|\overline{p}_{k_l,k_j}^\ast\|_{L^2(T)}^2\\
    \leq& c \|y_{k_l}^\ast\|_{L^4(\Omega)}^4 \|\overline{p}_{k_l,k_j}^\ast\|_{L^2(T)}^2. 
\end{align*}
Consequently, the following three estimates hold
    \begin{align}
         \sum_{T\in\cT_{k_l}\cap\cT_{k_j}}h_T^2  \|3 \overline{ u}_{k_j,k_l}^\ast  (y_{k_j}^\ast)^2p_{k_j}^\ast\|_{L^2(T)}^2
        &\leq c \|y_{k_j}^\ast\|_{H^1(\Omega)}^4\|p_{k_j}^\ast\|_{H^1(\Omega)}^2 \|\overline{u}_{k_l,k_j}^\ast \|_{L^2(\Omega)}^2,\label{pf:stab_est_unref1_08}\\
      \sum_{T\in\cT_{k_l}\cap\cT_{k_j}}h_T^2  \|3 b(u_{k_l}^\ast)( (y_{k_j}^\ast)^2 - (y_{k_l}^\ast)^2) p_{k_j}^\ast\|_{L^2(T)}^2 &\leq c
      \|  y_{k_j}^\ast + y_{k_l}^\ast \|_{H^1(\Omega)}^2 \| p_{k_j}^\ast\|_{H^1(\Omega)}^2 \|\overline{ y}_{k_l,k_j}^\ast   \|_{L^2(\Omega)}^2,\nonumber\\
       \sum_{T\in\cT_{k_l}\cap\cT_{k_j}}h_T^2 \|3 b(u_{k_l}^\ast)(y_{k_l}^\ast)^2 \overline{p}_{k_j,k_l}^\ast \|_{L^2(T)}^2 &\leq c \| y_{k_l}^\ast\|_{H^1(\Omega)}^4 \| \overline{p}_{k_l,k_j}^\ast \|_{L^2(\Omega)}^2.\nonumber
\end{align}
Third, the argument for the estimates \eqref{pf:stab_est_unref1_05} and \eqref{pf:stab_est_unref1_06} gives
    \begin{align*}
      &\quad \sum_{T\in\cT_{k_l}\cap\cT_{k_j}}\sum_{F\subset\partial T}h_T\|J_{F,2}(u_{k_l}^\ast, y_{k_l}^\ast, p_{k_l}^\ast)  - J_{F,2}(u_{k_j}^\ast, y_{k_j}^\ast, p_{k_j}^\ast)\|_{L^2(F)}^2  \nonumber\\
      &\leq c\Big( \|{\nabla}\overline{p}_{k_l,k_j}^\ast \|^2_{L^2(\Omega)} +  \|{\nabla}\overline{p}_{k_j}^\ast\|^2_{L^2(\Omega)}
       + \|\overline{u}_{k_l,k_j}^\ast {\nabla} p_{\infty}^\ast\|^2_{L^2(\Omega)} + \|\overline{y}_{k_l,k_j}^\ast\|_{L^2(\Omega)}^2\Big).
    \end{align*}
    Since the sequences $\{u_{k_j}^\ast\}_{j\geq0}$, $\{y_{k_j}^\ast\}_{j\geq0}$ and $\{p_{k_j}^\ast\}_{j\geq0}$ are all uniformly bounded in $H^1(\Omega)$, the
    estimate \eqref{stab_est_unref2} follows from the last estimates. Last, for the estimate \eqref{stab_est_unref3}, the triangle inequality  gives
    \[
        \Big(\sum_{T\in\cT_{k_j}\cap\cT_{k_l}}\eta^2_{k_l,3}(u_{k_l}^\ast,y^\ast_{k_l},p_{k_l}^\ast,T) \Big)^{1/2} \leq \Big(\sum_{T\in\cT_{k_j}\cap\cT_{k_l}}\eta^2_{k_j,3}(u_{k_j}^\ast,y^\ast_{k_j},p_{k_j}^\ast,T) \Big)^{1/2} + \mathrm{III},
    \]
    with the term ${\rm III}$ given by
        \begin{align*}
        \mathrm{III}:= &\Big(\sum_{T\in\cT_{k_l}\cap\cT_{k_j}}h_T^2\| ( 1 - \sigma )({\nabla} y_{k_l}^\ast \cdot {\nabla} p_{k_l}^\ast - {\nabla} y_{k_j}^\ast \cdot {\nabla} p_{k_j}^\ast) +  (y_{k_l}^\ast)^3 p_{k_l}^\ast - (y_{k_j}^\ast)^3 p_{k_j}^\ast + 2\alpha\varepsilon^{-1}(u^\ast_{k_j} - u^\ast_{k_l})\|^2_{L^2(T)} \\
        &\quad + \sum_{T\in\cT_{k_l}\cap\cT_{k_j}}\sum_{F\subset\partial T}h_T\|2\alpha\varepsilon[{\nabla}(u_{k_l}^\ast - u_{k_j}^\ast) \cdot{n}_F]\|^2_{L^2(F)}\Big)^{1/2}.
        \end{align*}
    It suffices to bound the two summands of the term ${\rm III}$.
    By proceeding like before, we can derive
    \begin{align}
      &\quad \sum_{T\in\cT_{k_l}\cap\cT_{k_j}}h_T^2\| ( 1 - \sigma )({\nabla} y_{k_l}^\ast \cdot{\nabla} p_{k_l}^\ast - {\nabla} y_{k_j}^\ast \cdot{\nabla} p_{k_j}^\ast)\|_{L^2(T)}^2 \nonumber\\
      &\leq c(\|y_{k_l}^\ast\|^2_{H^1(\Omega)}\|{\nabla}\overline{p}_{k_l,k_j}^\ast\|_{L^2(\Omega)}^2
      +\|p_{k_j}^\ast\|^2_{H^1(\Omega)}\|{\nabla}\overline{y}_{k_l,k_j}^\ast\|_{L^2(\Omega)}^2),\label{pf:stab_est_unref1_09}\\
       & \sum_{T\in\cT_{k_l}\cap\cT_{k_j}}h_T^2\left\|2\alpha\varepsilon^{-1}\overline{u}^\ast_{k_l,k_j} \right\|_{L^2(T)}^2 \leq c \|\overline{u}^\ast_{k_l,k_j} \|_{L^2(\Omega)}^2,\nonumber\\
        & \sum_{T\in\cT_{k_l}\cap\cT_{k_j}}\sum_{F\subset\partial T}h_T\left\|2\alpha\varepsilon[{\nabla}\overline{u}_{k_l,k_j}^\ast  \cdot {n}_F]\right\|^2_{L^2(F)} \leq c \|{\nabla}\overline{u}_{k_l,k_j}^\ast\|_{L^2(\Omega)}^2.\nonumber
    \end{align}
    In the splitting $(y_{k_l}^\ast)^3 p_{k_l}^\ast - (y_{k_j}^\ast)^3 p_{k_j}^\ast = ((y_{k_l}^\ast)^3 - (y_{k_j}^\ast)^3) p_{k_l}^\ast + (y_{k_j}^\ast)^3 ( p_{k_l}^\ast - p_{k_j}^\ast) $, using the argument for the estimates \eqref{pf:stab_est_unref1_03}, \eqref{pf:stab_est_unref1_04} and \eqref{pf:stab_est_unref1_08} once again gives
    \begin{align*}
        &h_T^2 \|((y_{k_l}^\ast)^3 - (y_{k_j}^\ast)^3) p_{k_l}^\ast \|_{ L^2(T)}^2
       \leq ch_T^2\|\overline{y}_{k_l,k_j}^\ast \|_{L^8(T)}^2(\|y_{k_j}^\ast\|_{L^8(T)}^4 \|p_{k_l}^\ast\|_{L^8(T)}^2 + \|y_{k_l}^\ast\|_{L^8(T)}^4\|p_{k_l}^\ast\|_{L^8(T)}^2)\\
        \leq& c h_T^2 (h_{T}^{-3/4})^2 (h_T^{-1/12})^{4} (h_T^{-1/12})^2 \|\overline{y}_{k_l,k_j}^\ast \|_{L^2(T)}^2(\|y_{k_j}^\ast\|_{L^6(T)}^4\|p_{k_j}^\ast\|_{L^6(T)}^2 + \|y_{k_l}^\ast\|_{L^6(T)}^4\|p_{k_l}^\ast\|_{L^6(T)}^2)\\
        = & c \|\overline{y}_{k_l,k_j}^\ast \|_{L^2(T)}^2(\|y_{k_j}^\ast\|_{L^6(T)}^4\|p_{k_j}^\ast\|_{L^6(T)}^2 + \|y_{k_l}^\ast\|_{L^6(T)}^4\|p_{k_l}^\ast\|_{L^6(T)}^2),
    \end{align*}
    and hence
\begin{align*}
    \sum_{T\in\cT_{k_l}\cap\cT_{k_j}}h_T^2 \|((y_{k_l}^\ast)^3 - (y_{k_j}^\ast)^3) p_{k_l}^\ast \|_{L^2(T)}^2
        &\leq c (\|y_{k_j}^\ast\|_{H^1(\Omega)}^4\|p_{k_j}^\ast\|_{H^1(\Omega)}^2 + \|y_{k_l}^\ast\|_{H^1(\Omega)}^4\|p_{k_l}^\ast\|_{H^1(\Omega)}^2) \|\overline{y}_{k_l,k_j}^\ast \|_{L^2(\Omega)}^2,\\
        \sum_{T\in\cT_{k_l}\cap\cT_{k_j}}h_T^2 \|(y_{k_j}^\ast)^3 ( p_{k_l}^\ast - p_{k_j}^\ast) \|_{L^2(T)}^2 &\leq c \|y_{k_j}^\ast\|_{H^1(\Omega)}^6\|\overline{p}_{k_l,k_j}^\ast \|_{L^2(\Omega)}^2.
\end{align*}
Finally, the estimate \eqref{stab_est_unref3} follows from these estimates and the uniform boundedness of the sequences $\{y_{k_j}^\ast\}_{j\geq0}$ and $\{p_{k_j}^\ast\}_{j\geq0}$ in the $H^1(\Omega)$ norm.

\bibliographystyle{abbrv}
\bibliography{CE_bib}

\end{document}